\documentclass[11pt,reqno]{amsart}

\usepackage[a4paper,margin=1in]{geometry}
\usepackage{amsmath,amssymb,amsthm,mathtools,mathrsfs}
\usepackage{bm}
\usepackage{array,booktabs,tabularx,ragged2e}
\usepackage{longtable}
\usepackage{microtype}
\usepackage{enumitem}
\usepackage{needspace}
\usepackage{xcolor}
\usepackage[colorlinks=true,linkcolor=blue,citecolor=blue,urlcolor=blue]{hyperref}

\allowdisplaybreaks
\numberwithin{equation}{section}
\renewcommand{\paragraph}[1]{%
  \par\addvspace{0.6\baselineskip}\noindent\textbf{#1}\enspace\ignorespaces}
\makeatletter
\renewcommand{\subsection}{\@startsection{subsection}{2}{\z@}%
  {0.8\baselineskip\@plus0.2\baselineskip}{0.35\baselineskip}%
  {\normalfont\bfseries}}
\makeatother

\newtheorem{theorem}{Theorem}[section]
\newtheorem{proposition}{Proposition}[section]
\newtheorem{lemma}[proposition]{Lemma}
\newtheorem{corollary}[proposition]{Corollary}
\theoremstyle{remark}

\newtheorem*{remark*}{Remark}

\newcommand{\R}{\mathbb R}
\newcommand{\dd}{\,\mathrm d}
\newcommand{\dx}{\,\mathrm dx}
\newcommand{\pa}{\partial}
\newcommand{\Div}{\operatorname{div}}
\newcommand{\ip}[2]{\left\langle #1,#2\right\rangle}
\newcommand{\cE}{\mathcal E}
\newcommand{\cD}{\mathcal D}
\newcommand{\cX}{\mathcal X}

\title[Stability and nonlinear asymptotics for the NSF system]
{Global stability and nonlinear asymptotics for the three-dimensional
compressible Navier--Stokes--Fourier system with horizontal shear viscosity}

\author[Z. Feng]{Zefu Feng}
\address{Zefu Feng, School of Mathematical Sciences, Chongqing Normal University,
Chongqing 401331, China}
\email{zefufeng@mails.ccnu.edu.cn}

\author[G. Hong]{Guangyi Hong}
\address{Guangyi Hong, School of Mathematics, South China University of
Technology, Guangzhou 510641, China}
\email{magyhong@scut.edu.cn}

\author[Y. Wang]{Yinghui Wang}
\address{Yinghui Wang, MOE-LCSM, School of Mathematics and Statistics,
Hunan Normal University, Changsha 410081, China}
\email{yhwangmath@hunnu.edu.cn}

\author[J. Wu]{Jiahong Wu}
\address{Jiahong Wu, Department of Mathematics, University of Notre Dame,
Notre Dame, IN 46556, USA}
\email{jwu29@nd.edu}

\makeatletter
\@ifundefined{subjclassname@2020}{%
  \@namedef{subjclassname@2020}{\textup{2020} Mathematics Subject Classification}%
}{}
\makeatother

\begin{document}

\begin{abstract}
We study the three-dimensional compressible Navier--Stokes--Fourier
system with horizontal shear viscosity, a positive coefficient
$\lambda+\mu>0$ in the divergence term, and full heat conduction. Horizontal viscosity does not dissipate vertical shear derivatives, while thermal diffusion is coupled to the compressible variables. For small $H^2$ perturbations of a constant equilibrium, we
construct global strong solutions. Additional spatial localization
yields anisotropic decay through second order, convergence of the full
$H^2$ norm, and horizontal Gaussian profiles determined by an effective
velocity datum. The third-order behavior contains nonlinear
contributions: horizontal viscous heating generates a differentiated
Newton potential in temperature, while compressible transport slows
the decay of the potential velocity. We prove the corresponding bounds
with exponents $17/8$ and $2$; for each exponent, we construct a fixed
datum in the corresponding class for which every larger exponent fails.
For each fixed initial datum,
however, the corresponding time-weighted norms tend to zero as
$t\to\infty$. Smallness of the first vertical derivative of the initial velocity in the horizontal mixed norm improves the vertical decay estimates and restores the linear third-order temperature rate. Both the existence and decay results require only $H^2$ Sobolev regularity of the initial data.
\end{abstract}

\maketitle
\tableofcontents

\section{Introduction and main results}

\subsection{The model and the anisotropic dissipation}

We study global stability and large-time asymptotics for a compressible
heat-conducting gas with horizontal shear viscosity. The viscous
stress dissipates horizontal velocity gradients and divergence,
while the vertical shear remains undiffused. For small $H^2$
perturbations, we establish global strong solutions and determine
how nonlinear transport and viscous heating affect the decay rates.
Anisotropic eddy-viscosity models motivate this choice of stress;
see \cite{TemamZiane,CheminAnisotropic2000}.

Let \(x=(x_h,x_3)\in\R^2\times\R\), \(x_h=(x_1,x_2)\), and
\(\Delta_h=\pa_1^2+\pa_2^2\).  We consider the Cauchy problem
\begin{equation}\label{eq:original}
\left\{
\begin{aligned}
 &\rho_t+\Div(\rho u)=0,\\
 &\rho(u_t+u\cdot\nabla u)+\nabla(R\rho\theta)
   =\mu\Delta_hu+(\lambda+\mu)\nabla\Div u,\\
 &\rho c_v(\theta_t+u\cdot\nabla\theta)
      +R\rho\theta\,\Div u\\
 &\hspace{25mm}
   =\mu|\nabla_hu|^2+(\lambda+\mu)(\Div u)^2+\kappa\Delta\theta .
\end{aligned}\right.
\end{equation}
It is supplemented with
\begin{equation*}
 (\rho,u,\theta)|_{t=0}=(\rho_0,u_0,\theta_0),\qquad
 (\rho_0,u_0,\theta_0)(x)\longrightarrow
 (\bar\rho,0,\bar\theta)\quad\text{as }|x|\to\infty,
\end{equation*}
where the far-field condition is understood through the Sobolev perturbations
introduced below.  The coefficients are constants satisfying
\begin{equation}\label{eq:coeff}
 \bar\rho>0,\qquad \bar\theta>0,\qquad
 R>0,\quad c_v>0,\quad\mu>0,\quad\lambda+\mu>0,\quad\kappa>0.
\end{equation}
The constant state \((\bar\rho,0,\bar\theta)\) is an equilibrium.

The momentum operator and the viscous production in \eqref{eq:original} are
generated by the anisotropic eddy stress
\begin{equation*}
 \mathbb S^{\rm ed}_{ij}(u)
 :=\mu\mathbf 1_{\{j=1,2\}}\pa_j u_i
   +(\lambda+\mu)\delta_{ij}\Div u .
\end{equation*}
Indeed,
\begin{equation*}
 \Div\mathbb S^{\rm ed}(u)
 =\mu\Delta_hu+(\lambda+\mu)\nabla\Div u,
 \qquad
 \mathbb S^{\rm ed}(u):\nabla u
 =\mu|\nabla_hu|^2+(\lambda+\mu)(\Div u)^2.
\end{equation*}
The same anisotropic quadratic form therefore determines the momentum dissipation and viscous heating, consistently with the energy balance.  The natural coercivity condition is \(\lambda+\mu>0\),
and the velocity dissipation controls
\[
 \nabla_hu\quad\text{and}\quad \Div u,
\]
but it does not directly control the vertical shear derivative \(\pa_3u_h\).
This loss is the central degeneracy of the problem. Throughout the paper,
we denote the viscous production by
\begin{equation}
 Q(u):=\mu|\nabla_hu|^2+(\lambda+\mu)(\Div u)^2.
\end{equation}
The results concern the ideal-gas laws and the specified eddy stress
in \eqref{eq:original}. Full heat conduction supplies the temperature
gradients required by the $H^2$ coefficient estimates.

\subsection{Related work}

For the corresponding barotropic eddy-diffusion system, Feng, Hong, Wang and Wu
\cite[Theorems~1.1--1.2]{FengHongWangWu} prove global existence for small $H^N$ data, $N\ge3$. Their decay theorem assumes small $H^4\cap L^1$ data together with
$\partial_3^ju_0\in L^2_3L^1_h$ for $j=0,1,2$, and includes a
logarithmic loss in the $L^2$ density estimate.
For \eqref{eq:original}, Theorem~\ref{thm:global} closes the global
energy estimate in $H^2$.
The rates \eqref{eq:new-H2-main-rates}, including density decay without
a time logarithm, require only small $\mathcal J_2$ as defined in
\eqref{eq:new-data-sizes}.
The enhanced vertical estimates \eqref{eq:H2-intermediate-main}
use small $\mathcal J_{2,1}$ and require no second vertical mixed
norm of $u_0$.

For full viscosity and heat conduction, the small-data theory near
equilibrium goes back to Matsumura and Nishida
\cite{MatsumuraNishida}. Critical-space formulations were developed
by Danchin \cite{Danchin,DanchinFull} and by Chikami and Danchin
\cite{ChikamiDanchin}. The Shizuta--Kawashima mechanism
\cite{ShizutaKawashima} explains how coupling supplies dissipation
to hyperbolic variables, while the thermodynamic weak-solution
framework is discussed by Feireisl and Novotn\'y
\cite{FeireislNovotny}. Here the coupling still controls density
and the potential velocity, but leaves vertical shear without
direct diffusion.

Quantitative decay requires information beyond the energy norm.
The whole-space estimates of Deckelnick \cite{Deckelnick}, the
exterior-domain heat-conductive analysis of Kobayashi and Shibata
\cite{KobayashiShibata}, the negative Sobolev method of Guo and
Wang \cite{GuoWang}, and the critical decay estimates of Danchin
and Xu \cite{DanchinXuBarotropic,DanchinXuFull} illustrate the role
of low frequencies in fully viscous flows. In the present system,
the solenoidal and compressible components have the decay rates of two-dimensional horizontal diffusion and three-dimensional acoustic--thermal diffusion, respectively.

The shear estimates draw on tools developed for incompressible
partial dissipation. Anisotropic spaces and product estimates
appear in Chemin, Desjardins, Gallagher and Grenier
\cite{CheminAnisotropic2000}, Iftimie \cite{Iftimie}, and Chemin
and Zhang \cite{CheminZhang}. Directional decay and the enhanced
third-component behavior were studied by Ji, Wu and Yang
\cite{JiWuYang} and Xu and Zhang \cite{XuZhang}; related fractional
problems were treated by Ji, Luo and Jiang \cite{JiLuoJiang} and
Bie, Fang, Wang and Zhou \cite{BieFangWangZhou}. The product
structure is also useful in partially dissipative MHD
\cite{WuZhu}. For compressible anisotropic stresses, relevant
developments include the compactness theory of Bresch and Jabin
\cite{BreschJabin}, the classical-solution results of Wang and Guo
\cite{WangGuoAnisotropic}, Bai, Guo and Xu \cite{BaiGuoXu}, and
Xu and Du \cite{XuDu}, and the half-space vertical-viscosity limit
of Gao, Peng, Wu and Yao \cite{GaoPengWuYao}. These works concern different constitutive laws or domains.
For \eqref{eq:original}, the solenoidal velocity has only horizontal
diffusion.

At third order, \eqref{eq:r9-thermal-profile} isolates the Newton
potential generated by horizontal viscous heating.
The nonlinear decay exponents $17/8$ for temperature and $2$ for
potential velocity in \eqref{eq:r9-tau-third} and \eqref{eq:r9-w-third}
are sharp in the respective data classes by
Theorem~\ref{thm:decay-H2}(iv).

\subsection{Variables and initial localization}

We introduce the logarithmic perturbation variables
\begin{equation}
 q=\log\frac{\rho}{\bar\rho},\qquad
 \tau=\log\frac{\theta}{\bar\theta},\qquad U=(q,u,\tau),
\end{equation}
and set
\begin{equation*}
 c^2=R\bar\theta,\quad \alpha=\frac{R}{c_v},\quad
 \omega=\bar\theta c_v,\quad
 \nu_h=\frac{\mu}{\bar\rho},\quad
 \nu_c=\frac{\lambda+\mu}{\bar\rho},\quad
 \chi=\frac{\kappa}{\bar\rho c_v}.
\end{equation*}
The constant-coefficient velocity operator will be denoted by
\[
 \mathcal Lu=\nu_h\Delta_hu+\nu_c\nabla\Div u.
\]
Then \(\omega\alpha=c^2\), and \eqref{eq:original} becomes
\begin{equation}\label{eq:perturbed}
\left\{
\begin{aligned}
 q_t+\Div u&=F_1,\\
 u_t+c^2\nabla q+c^2\nabla\tau-\nu_h\Delta_hu
       -\nu_c\nabla\Div u&=F_2,\\
 \tau_t+\alpha\Div u-\chi\Delta\tau&=F_3,
\end{aligned}\right.
\end{equation}
where
\begin{align}
 F_1={}&-u\cdot\nabla q,\label{eq:F1}\\
 F_2={}&-u\cdot\nabla u
 -c^2(e^\tau-1)(\nabla q+\nabla\tau)
 +(e^{-q}-1)(\nu_h\Delta_hu+\nu_c\nabla\Div u),\label{eq:F2}\\
 F_3={}&-u\cdot\nabla\tau
 +\frac{e^{-q-\tau}}{\bar\rho\bar\theta c_v}
   \bigl\{\mu|\nabla_hu|^2+(\lambda+\mu)(\Div u)^2\bigr\}\notag\\
 &+\chi(e^{-q}-1)\Delta\tau
 +\chi e^{-q}|\nabla\tau|^2.\label{eq:F3}
\end{align}
Let \(\mathbb P=I-\nabla\Delta^{-1}\Div\) and
\(\mathbb Q=I-\mathbb P\) denote the solenoidal and potential projections.
We write
\[
 v=\mathbb Pu,\qquad w=\mathbb Qu,\qquad W=(q,w,\tau),
 \qquad S_h(t)=e^{\nu_ht\Delta_h}.
\]
Although \(-\chi\Delta\tau\) gives full temperature dissipation, projection of
the momentum equation yields the exact identity
\begin{equation}\label{eq:projected-P-intro}
 (\mathbb Pu)_t-\nu_h\Delta_h\mathbb Pu=\mathbb PF_2,
 \qquad \mathbb P\nabla(q+\tau)=0.
\end{equation}
Thus heat conduction strengthens the acoustic--thermal subsystem but cannot
create vertical diffusion for the shear mode.

All initial regularity assumptions in this paper are imposed in
$H^2(\R^3)$. The squared $L^2$ or Sobolev norm of a vector or tuple is the
sum of the squared norms of its components. Finite sums of
component norms are equivalent, with constants absorbed into $C$.
The symbols $\nabla^k$ and $\nabla_h^k$ denote tensors of all
ordered derivatives; their Fourier weights are $|\xi|^k$ and
$|\xi_h|^k$. In particular the mixed horizontal second derivative occurs
twice in its squared tensor norm. This convention fixes the
Gaussian constants stated below. We use
$\widehat f(\xi)=\int e^{-ix\cdot\xi}f(x)\dd x$; thus the
three-dimensional Plancherel factor is $(2\pi)^{-3}$.
When convenient, we write $z$ for the vertical variable $x_3$,
and $(\xi_h,\zeta)\in\R^2\times\R$ for the horizontal and vertical
Fourier variables.
We write $D=-i\nabla$, $D_h=-i\nabla_h$, and $D_3=-i\partial_3$,
and use $m(D)$ for the Fourier multiplier with symbol $m(\xi)$.
In particular, $|D_h|=(-\Delta_h)^{1/2}$ and
$\mathbf1_{\{|D_h|\ge1\}}$ has symbol $\mathbf1_{\{|\xi_h|\ge1\}}$.

The energy dissipation is
\begin{equation}\label{eq:H2-dissipation-definition}
 \cD_2(U)=\|\nabla q\|_{H^1}^2+\|\nabla_hu\|_{H^2}^2
                    +\|\Div u\|_{H^2}^2+\|\nabla\tau\|_{H^2}^2.
\end{equation}
It contains third derivatives in the dissipative directions,
integrated in time. These are not additional assumptions on the
initial data, and $\partial_3^3u_h$ is not controlled by this
dissipation.

For $1\le p,q<\infty$, mixed norms are defined by
\[
 \|f\|_{L^p_hL^q_3}
 =\left(\int_{\R^2}
       \left(\int_\R|f(x_h,x_3)|^q\dd x_3\right)^{p/q}
       \dd x_h\right)^{1/p}.
\]
The reversed order defines $L^p_3L^q_h$; an infinite exponent denotes
the corresponding essential supremum. Unless specified otherwise,
$\|\cdot\|_p$ is the $L^p(\R^3)$ norm. Unless stated otherwise, $C$ and $c_*$ denote positive constants, which may
change from line to line and depend on the fixed coefficients and
differentiation order, but not on the solution or time. The physical
parameter $c$ is fixed by $c^2=R\bar\theta$.
For a measurable function $g$ on $\R^2$, we use the weak $L^2$
quasi-norm
\[
 \|g\|_{L^{2,\infty}(\R^2)}
 =\sup_{s>0}s\,
   \bigl|\{x_h\in\R^2:|g(x_h)|>s\}\bigr|^{1/2}.
\]
Here the measure is the two-dimensional Lebesgue measure, and
$L^{2,\infty}_h$ abbreviates $L^{2,\infty}(\R^2)$.

We use two initial-data norms with the same Sobolev regularity:
\begin{equation}\label{eq:new-data-sizes}
\mathcal J_2=\|U_0\|_{H^2}+\|U_0\|_{L^1}+\|u_0\|_{L^2_3L^1_h},
 \qquad \mathcal J_{2,1}=\mathcal J_2+\|\partial_3u_0\|_{L^2_3L^1_h}.
\end{equation}
Each smallness assumption applies to the entire corresponding norm.
We refer to the data class measured by $\mathcal J_2$ as the
$\mathcal J_2$ class; all solutions are strong $H^2$ solutions.
The assumptions and their uses are summarized below. 
\begin{center}\small
\begin{tabularx}{\linewidth}{@{}>{\RaggedRight\arraybackslash}p{.25\linewidth}>{\RaggedRight\arraybackslash}X@{}}
\toprule
Initial hypothesis & Conclusions supported by that hypothesis\\ \midrule
Small $\|U_0\|_{H^2}$ &
Global strong solution and integrated $H^2$ dissipation
(Theorem~\ref{thm:global}).\\
Small $\mathcal J_2$ &
Directional decay and $H^2$ convergence
(Theorem~\ref{thm:decay-H2}); horizontal Gaussian profiles
(Theorem~\ref{thm:r9-gaussians}).\\
Small $\mathcal J_{2,1}$ &
Enhanced decay of the first vertical derivative, mixed third-order derivatives, and temperature
(Theorem~\ref{thm:H2-intermediate-main}); vertical Gaussian profiles
(Theorem~\ref{thm:r9-gaussians}).\\
\bottomrule
\end{tabularx}
\end{center}

\subsection{Stability and basic decay}

\begin{theorem}[Global $H^2$ strong solutions]\label{thm:global}
Assume \eqref{eq:coeff}. There are constants $\delta_2,C_2>0$,
depending only on the fixed coefficients and equilibrium, such
that every $U_0\in H^2(\R^3)$ with $\|U_0\|_{H^2}\le\delta_2$
generates a unique global strong solution of \eqref{eq:original}.
It satisfies $\rho=\bar\rho e^q>0$, $\theta=\bar\theta e^\tau>0$,
\begin{align}
 U&\in C([0,\infty);H^2),\label{eq:solspace1}\\
 \nabla q&\in L^2(0,\infty;H^1),\qquad
 \nabla_hu,\ \Div u,\ \nabla\tau\in L^2(0,\infty;H^2),
                                      \\
 q_t&\in C([0,\infty);H^1),\qquad
 u_t,\tau_t\in C([0,\infty);L^2)
                    \cap L^2_{\rm loc}([0,\infty);H^1),
                                      \label{eq:solspace3}
\end{align}
and
\begin{equation}\label{eq:globalbound}
 \sup_{t\ge0}\|U(t)\|_{H^2}^2+
       \int_0^\infty\cD_2(U(t))\dd t\le C_2\|U_0\|_{H^2}^2.
\end{equation}
\end{theorem}

Additional $L^1$ and horizontal mixed-norm localization yield
quantitative decay of the global solution. The rates distinguish
the compressible variables from the solenoidal velocity.

\Needspace{8\baselineskip}
\begin{theorem}[Decay estimates and sharpness]
\label{thm:decay-H2}
Assume \eqref{eq:coeff}. There are $\delta_{d,2},C>0$ such that
the solution of Theorem~\ref{thm:global} satisfies
\textup{(i)--(iii)} for every datum with
$\mathcal J_2\le\delta_{d,2}$. The existence statements in
\textup{(iv)} hold in arbitrarily small balls of the indicated data
classes.

\smallskip\noindent\textup{(i) Decay through second order.}
For $k=0,1,2$ and $t\ge0$,
\begin{align}
 \label{eq:new-H2-main-rates}
 \|\nabla^kW(t)\|_2&\le C\mathcal J_2(1+t)^{-3/4-k/2},\notag\\
 \|\nabla_h^ku(t)\|_2&\le C\mathcal J_2(1+t)^{-(1+k)/2},
 \qquad \|\nabla_h^ku_3(t)\|_2\le C\mathcal J_2(1+t)^{-3/4-k/2}.
\end{align}
The remaining directional estimates are
\begin{align}\label{eq:H2-main-vertical-weak}
 \|\partial_3u(t)\|_2&\le C\mathcal J_2(1+t)^{-1/4},&
 \|\nabla_h\partial_3u(t)\|_2&\le C\mathcal J_2(1+t)^{-3/4},\notag\\
 \|\partial_3u_3(t)\|_2&\le C\mathcal J_2(1+t)^{-1},&
 \|\partial_3^2u_3(t)\|_2&\le C\mathcal J_2(1+t)^{-3/4}.
\end{align}
For $2\le p\le\infty$,
\begin{equation}\label{eq:new-H2-main-Lp}
 \|u(t)\|_p\le C_p\mathcal J_2(1+t)^{-7/8+3/(4p)},\qquad
 \|(W,u_3)(t)\|_p\le C_p\mathcal J_2(1+t)^{-\frac32(1-1/p)}.
\end{equation}
Furthermore,
\begin{equation}\label{eq:H2-main-mixed-and-full}
\begin{gathered}
 \|u_3(t)\|_{L^\infty_3L^2_h}\le C\mathcal J_2(1+t)^{-1},\\
 \|U(t)\|_{H^1}\le C\mathcal J_2(1+t)^{-1/4},\qquad
 \|W(t)\|_{H^2}\le C\mathcal J_2(1+t)^{-3/4}.
\end{gathered}
\end{equation}
\smallskip\noindent\textup{(ii) Fixed-data limits.}
For each fixed datum,
\begin{equation}\label{eq:r8-main-zero}
 \|U(t)\|_{H^2}\longrightarrow0,\qquad
 \|\partial_3u(t)\|_2=o(t^{-1/4}),\qquad
 \|\nabla_h\partial_3u(t)\|_2+\|\partial_3^2u_3(t)\|_2=o(t^{-3/4}).
\end{equation}

\smallskip\noindent\textup{(iii) Third-order estimates and thermal comparison.}
For $t\ge1$,
\begin{align}
 \|\nabla_h^3u(t)\|_2&\le C\mathcal J_2(1+t)^{-2},
                                      \label{eq:r9-horizontal-third}\\
 \|\nabla_h^2\partial_3u(t)\|_2&\le C\mathcal J_2(1+t)^{-5/4},
                                      \label{eq:r9-mixed-third}\\
 \|\nabla^3w(t)\|_2&\le C\mathcal J_2(1+t)^{-9/4}
                         +C\mathcal J_2^2(1+t)^{-2},\label{eq:r9-w-third}\\
 \|\nabla^3\tau(t)\|_2&\le C\mathcal J_2(1+t)^{-9/4}
                         +C\mathcal J_2^2(1+t)^{-17/8}.
                                      \label{eq:r9-tau-third}
\end{align}
For each fixed datum, the left-hand sides of
\eqref{eq:r9-mixed-third}, \eqref{eq:r9-w-third}, and
\eqref{eq:r9-tau-third} are, respectively,
$o(t^{-5/4})$, $o(t^{-2})$, and $o(t^{-17/8})$.
The thermal contribution has the more precise form
\begin{equation}\label{eq:r9-thermal-profile}
 \left\|\nabla^3\tau(t)
 -\frac{\mu}{\bar\rho\bar\theta c_v\chi}
       \nabla^3(-\Delta)^{-1}|\nabla_hu(t)|^2\right\|_2
 \le C\mathcal J_2(1+t)^{-9/4}.
\end{equation}
The differentiated Newton potential is defined by its Fourier
multiplier, with
$\|\nabla^3(-\Delta)^{-1}g\|_2=\|\nabla g\|_2$.
The physical temperature $\theta-\bar\theta$ satisfies the same
estimate with coefficient $\mu/(\bar\rho c_v\chi)$.

One also has, for $t\ge1$,
\begin{align}
 \|\nabla_hu(t)\|_\infty+\|\nabla u_3(t)\|_\infty
       &\le C\mathcal J_2(1+t)^{-11/8},\label{eq:r9-lips}\\
 \|\nabla_hu_3(t)\|_\infty+\|\nabla(w,\tau)(t)\|_\infty
       &\le C\mathcal J_2(1+t)^{-2}.\label{eq:r9-compressible-lips}
\end{align}
\smallskip\noindent\textup{(iv) Sharpness of the third-order decay rates.}
For every sufficiently small $\delta>0$, the following assertions hold.
\begin{enumerate}[label=\textup{(\alph*)},leftmargin=*]
\item There is a fixed datum with $\mathcal J_2\le\delta$ such that,
for every $\beta>17/8$,
\[
 \limsup_{t\to\infty}(1+t)^\beta\|\partial_3^3\tau(t)\|_2
 =\limsup_{t\to\infty}(1+t)^\beta
        \|\partial_3^3\tau\|_{L^2(t,t+1;L^2)}=\infty.
\]
Its linearized temperature is identically zero.
\item There is a fixed datum with $\mathcal J_{2,1}\le\delta$ such that,
for every $\beta>2$,
\[
 \limsup_{t\to\infty}(1+t)^\beta\|\partial_3^3w_3(t)\|_2
 =\limsup_{t\to\infty}(1+t)^\beta
        \|\partial_3^3w_3\|_{L^2(t,t+1;L^2)}=\infty.
\]
The linearized solution with these initial data has third-order
potential-velocity decay $O(t^{-9/4})$ in $L^2$.
\end{enumerate}
\end{theorem}

The upper bounds and fixed-data limits in parts (i)--(iii) hold
for every datum in the stated class. Part (iv) supplies fixed data
for which every strictly larger decay exponent fails. These
statements are compatible with the vanishing of the corresponding
time-weighted norms for each fixed initial datum.

The additional localization
$\partial_3u_0\in L^2_3L^1_h$ improves the vertical derivative
and temperature decay rates.

\Needspace{7\baselineskip}
\begin{theorem}[Enhanced decay estimates]
\label{thm:H2-intermediate-main}
If $\mathcal J_{2,1}$ is sufficiently small, the solution satisfies
the estimates and limits in Theorem~\ref{thm:decay-H2}(i)--(iii), and,
for every $t\ge0$,
\begin{align}\label{eq:H2-intermediate-main}
 \|u(t)\|_p&\le C_p\mathcal J_{2,1}(1+t)^{-1+1/p},
                                          &&2\le p\le\infty,\notag\\
 \|\partial_3u(t)\|_2+\|U(t)\|_{H^1}
       &\le C\mathcal J_{2,1}(1+t)^{-1/2},\notag\\
 \|\nabla_h\partial_3u(t)\|_2+\|\partial_3^2u_3(t)\|_2
       &\le C\mathcal J_{2,1}(1+t)^{-1}.
\end{align}
The estimates for the first vertical derivative have no logarithmic loss.

For $t\ge1$, one further has
\begin{equation}\label{eq:r9-strong-third}
 \|\nabla_h^2\partial_3u(t)\|_2
       \le C\mathcal J_{2,1}(1+t)^{-3/2},\qquad
 \|\nabla^3\tau(t)\|_2
       \le C\mathcal J_{2,1}(1+t)^{-9/4}.
\end{equation}
\end{theorem}
Smooth composition transfers the density and temperature estimates
to the physical perturbations.

\begin{corollary}[Physical density and temperature]\label{cor:physicaldecay}
The density and temperature estimates in
Theorems~\ref{thm:decay-H2} and \ref{thm:H2-intermediate-main}
hold for $(\rho-\bar\rho,\theta-\bar\theta)$ in place of $(q,\tau)$,
including their derivatives through order two and their $L^p$ norms.
\end{corollary}

\subsection{Linear profiles and nonlinear leading terms}\label{sec:r12-leading-terms}

The velocity is asymptotic to a horizontal heat flow, whereas higher-order derivatives of the compressible variables may have more slowly decaying nonlinear contributions.

\begin{theorem}[Gaussian profiles for horizontal and mixed derivatives]
\label{thm:r9-gaussians}
Assume $\mathcal J_2$ is sufficiently small. The effective velocity datum
\begin{equation}\label{eq:r8-physical-effective}
 u_\infty=u_0+\int_0^\infty F_2(s)\dd s
       \in H^1\cap L^2_3L^1_h
\end{equation}
is defined before Hodge projection, with $F_2$ given in \eqref{eq:F2};
the shear comparison uses $\mathbb Pu_\infty$.
Define its horizontal mean and the horizontal heat kernel by
\[
 m_\infty(x_3)=\int_{\R^2}(u_\infty)_h(x_h,x_3)\dd x_h,\qquad
 K_h(t,x_h)=\frac{1}{4\pi\nu_ht}
                    \exp\!\left(-\frac{|x_h|^2}{4\nu_ht}\right).
\]
Then $m_\infty\in L^2(\R)$ and, for $k=0,1,2,3$,
\begin{equation}\label{eq:r9-gaussian-all}
 \|\nabla_h^ku(t)-(\nabla_h^kK_h(t)m_\infty,0)\|_2
                         =o(t^{-(1+k)/2}).
\end{equation}
If $\mathcal J_{2,1}$ is sufficiently small, then $m_\infty\in H^1(\R)$
and, for $j=0,1,2$,
\begin{equation}\label{eq:r9-gaussian-vertical-all}
 \|\nabla_h^j\partial_3u(t)
       -(\nabla_h^jK_h(t)\partial_3m_\infty,0)\|_2
                         =o(t^{-(1+j)/2}).
\end{equation}
These profiles also hold in $L^2(t,t+1;L^2)$ in the same derivative ranges.
\end{theorem}

These Gaussian remainders are asymptotically negligible at the
respective decay scales. Without a moment condition, they need not
inherit the faster quantitative
rate of the comparison with the effective heat flow.

\subsection{Main difficulties and proof ideas}

\paragraph{The $H^2$ energy estimate.}
We control the vertical derivative of the third velocity component through
\[
             \partial_3u_3=\Div u-\Div_hu_h.
\]
Together with density compensation, this identity yields the
closed energy estimate in Proposition~\ref{prop:H2-closed} and its
integrated form \eqref{eq:strong-integrated-energy}.
At this regularity $\|\nabla u\|_\infty$ is not controlled by
the energy alone. The anisotropic inequalities instead combine
the $H^2$ bound with horizontal dissipation integrated in time.
We estimate the differentiated density transport in the compensated energy identity and pair the variable-viscosity terms with horizontal velocity gradients or the divergence. This closes the estimate without
$\partial_3^3u_h$. Local well-posedness and the continuity
argument then give a unique global strong solution.

\paragraph{Pointwise estimates for dissipative derivatives.}
We prove the decay estimates in the following order
\[
\begin{aligned}
 \text{second-order decay}&\ \longrightarrow\
          \text{time-integrated dissipation estimates}\\
 &\ \longrightarrow\ \text{positive-time difference quotients}
       \ \longrightarrow\ \text{spatial derivative estimates}.
\end{aligned}
\]
Positive-time averaging and time-difference quotients give
pointwise estimates for the time derivatives; the equations then
recover spatial derivatives. No derivative of $u_t(0)$ is prescribed.
Time-integrated horizontal velocity estimates control the
differentiated viscous heating in the second-order temperature
estimate. In the density energy estimate, we integrate $u\cdot\nabla\partial^2q$ by parts against the differentiated density.

\paragraph{The nonlinear thermal contribution.}
The thermal equation contains the slowly decaying horizontal
viscous heating.
The comparison \eqref{eq:r9-thermal-profile} isolates the
differentiated Newton potential of $|\nabla_hu|^2$.
A low-frequency time primitive and the damped high-frequency
equations leave a remainder with the linear decay rate.
The identity
$\|\nabla^3(-\Delta)^{-1}g\|_2=\|\nabla g\|_2$
reduces the estimate to the gradient of the heat source.
A thin-layer construction gives the potential-velocity lower bound
through a weak derivative limit; the thermal lower bound follows
from strong convergence of the thermal profile.

\paragraph{Vertical transport and heat-flow comparisons.}
Under small $\mathcal J_{2,1}$, we retain $u_3\partial_3$ in the
projected evolution to obtain the first vertical derivative estimates
without a logarithmic loss. Vertical characteristics and horizontal
Fourier estimates under vertical changes of variables give the
horizontal heat decay rates.
The nonlinear remainder relative to the linear shear flow is $r_v=v-S_h(t)\mathbb Pu_0$. The comparison $v-S_h(t)\mathbb Pu_\infty$ accounts for the time integral of the momentum source in the Gaussian profile.
Relative-energy and frequency-localized error estimates give strict
lower bounds for solutions of \eqref{eq:original}. Weak stability and
translation separation of the initial data then yield a single
datum with the asserted limsup.

\paragraph{Organization of the paper.}
Section~\ref{sec:r10-foundations} proves global existence.
Sections~\ref{sec:r10-decay}--\ref{sec:third-order-decay}
establish the basic and third-order decay, and
Sections~\ref{sec:r13-enhanced} and \ref{sec:r10-profiles}
prove the enhanced estimates and Gaussian profiles.
Appendix~\ref{app:construction-tools} contains the stability and
comparison estimates, and Appendix~\ref{sec:r10-counterexamples}
proves the nonlinear sharpness assertions. The relative-energy
identity, source estimates, and refined Green-matrix bounds are
collected in Appendices~\ref{app:r10-relative}--\ref{app:r13-green-matrix}.

\section{Global strong solutions}\label{sec:r10-foundations}

We extend the local solution by an $H^2$ energy estimate, combining
the relative energy, anisotropic product bounds, and density compensation.

\subsection{Anisotropic inequalities}

\subsubsection{Anisotropic product and mixed-norm estimates}

\begin{lemma}[Anisotropic triple products]
For smooth, rapidly decaying scalar functions \(f,g,h\),
\begin{align}
 \int_{\R^3}|fgh|\,\dd x
 &\leq C\|f\|_2^{1/2}\|\nabla_hf\|_2^{1/2}
          \|g\|_2^{1/2}\|\nabla_hg\|_2^{1/2}
          \|h\|_2^{1/2}\|\pa_3h\|_2^{1/2},\\
 \int_{\R^3}|fgh|\,\dd x
 &\leq C\|f\|_2^{1/2}\|\pa_3f\|_2^{1/2}
          \|g\|_2^{1/2}\|\nabla_hg\|_2^{1/2}
          \|h\|_2^{1/2}\|\nabla_hh\|_2^{1/2}.
\end{align}
\end{lemma}

These horizontal-gradient forms follow from the coordinatewise
anisotropic Sobolev inequalities in \cite[Lemma~1.2]{WuZhu} and
\cite[Lemma~2.2]{JiWuYang}.

\begin{lemma}[Coordinatewise anisotropic products]
Let \(i,j,k\) be pairwise distinct elements of \(\{1,2,3\}\).
Whenever the right-hand sides are finite,
\begin{align}
 \int_{\R^3}|fgh|\dx
 &\leq C
 \|f\|_2^{1/2}\|\pa_i f\|_2^{1/2}
 \|g\|_2^{1/2}\|\pa_j g\|_2^{1/2}
 \|h\|_2^{1/2}\|\pa_k h\|_2^{1/2},
 \label{eq:coordtriple}\\
 \int_{\R^3}|fgh|\dx
 &\leq C
 \|f\|_2^{1/4}\|\pa_i f\|_2^{1/4}
 \|\pa_j f\|_2^{1/4}\|\pa_i\pa_j f\|_2^{1/4}
 \|g\|_2^{1/2}\|\pa_k g\|_2^{1/2}\|h\|_2,
 \label{eq:coordtriple2}\\
 \|fg\|_2
 &\leq C
 \|f\|_2^{1/4}\|\pa_i f\|_2^{1/4}
 \|\pa_j f\|_2^{1/4}\|\pa_i\pa_j f\|_2^{1/4}
 \|g\|_2^{1/2}\|\pa_k g\|_2^{1/2}.
 \label{eq:coordproduct}
\end{align}
The estimates remain valid after a simultaneous permutation of coordinate
directions and factors.
\end{lemma}

For \eqref{eq:coordtriple} and \eqref{eq:coordtriple2}, see
\cite[Lemma~1.2]{WuZhu} and \cite[Lemma~2.2]{JiWuYang}.
The product estimate \eqref{eq:coordproduct} follows from
\eqref{eq:coordtriple2} by $L^2$ duality.

\begin{lemma}[Mixed horizontal--vertical norms]\label{lem:plane-trace}
For \(f,g\in H^1(\R^3)\),
\begin{align}
 \|f\|_{L_h^4L_3^2}^2
 &\leq C\|f\|_2\|\nabla_hf\|_2,\\
 \|g\|_{L_h^2L_3^\infty}^2
 &\leq C\|g\|_2\|\pa_3g\|_2,\label{eq:mixedv}\\
 \int_{\R^3}|f|^2|g|\dx
 &\leq C\|f\|_2\|\nabla_hf\|_2
       \|g\|_2^{1/2}\|\pa_3g\|_2^{1/2}.
\end{align}
\end{lemma}

The two-dimensional Gagliardo--Nirenberg inequality, vertical Sobolev
embedding and Minkowski's inequality give these estimates; see
\cite[Section~2]{JiWuYang,WuZhu}.
We also record the plane trace and reversed mixed-norm estimates:
\[
 \|a\|_{L^\infty_3L^4_h}\le C\|\nabla a\|_2,\qquad
 \|b\|_{L^2_3L^4_h}\le C\|b\|_2^{1/2}\|\nabla_hb\|_2^{1/2},
\]
for $a\in H^1$ and $b,\nabla_hb\in L^2$.
For smooth compactly supported $a$,
$\sup_{x_3}\int_h|a|^4\le4\int|a|^3|\partial_3a|
\le4\|a\|_6^3\|\partial_3a\|_2\le C\|\nabla a\|_2^4$.
The second estimate follows from horizontal Gagliardo--Nirenberg and
Cauchy--Schwarz in $x_3$; density gives the stated spaces.

\subsubsection{Differential products and structural estimates}

\begin{lemma}[Differential products and compositions]
Let \(|\beta|=2\), and
\(f,g\in H^2(\R^3)\cap L^\infty(\R^3)\), with the right-hand sides below
finite. Then
\begin{align}
 \|\pa^\beta(fg)\|_2
 &\leq C\bigl(\|f\|_\infty\|g\|_{H^2}
          +\|g\|_\infty\|f\|_{H^2}\bigr),\\
 \|[\pa^\beta,f]g\|_2
 &\leq C\bigl(\|\nabla f\|_\infty\|g\|_{H^1}
          +\|f\|_{H^2}\|g\|_\infty\bigr).
\end{align}
If \(G\in C^3(I)\), \(G(0)=0\), and the range of
\(z\in H^2\cap L^\infty\) is contained in a compact subset of \(I\), then
\begin{align}
 \|G(z)\|_{H^2}&\leq C_{G,2}(\|z\|_\infty)\|z\|_{H^2},
 \label{eq:Mosercomp}\\
 \|\nabla G(z)\|_{H^1}&\leq
 C_{G,2}(\|z\|_\infty)\|\nabla z\|_{H^1}.
 \label{eq:Mosercompgrad}
\end{align}
Moreover,
\begin{equation}\label{eq:H2-H1-multiplier}
 \|fg\|_{H^1}\leq C\|f\|_{H^2}\|g\|_{H^1}.
\end{equation}
\end{lemma}

The first two inequalities are standard Sobolev product and commutator
estimates, while the next two are the classical Moser composition estimates;
see, for example, \cite{KatoPonce,BahouriCheminDanchin}. The lower-order
estimate \eqref{eq:H2-H1-multiplier} is the standard Sobolev multiplication
\(H^2(\R^3)\times H^1(\R^3)\to H^1(\R^3)\); see also
\cite[Chapter~2]{BahouriCheminDanchin}.

\begin{lemma}[Control of the vertical component]
For $m=0,1,2$,
\begin{equation}\label{eq:verticalcomp}
 \|\nabla^m\nabla u_3\|_2
 \leq C\{\|\nabla^m\nabla_hu\|_2+
          \|\nabla^m\Div u\|_2\}.
\end{equation}
\end{lemma}

\begin{proof}
The horizontal derivatives of \(u_3\) are components of \(\nabla_hu\).
For the vertical derivative, commute every \(\pa^\beta\), \(|\beta|=m\),
with the exact identity
\[
 \pa_3u_3=\Div u-\pa_1u_1-\pa_2u_2.
\]
Taking \(L^2\) norms and summing proves \eqref{eq:verticalcomp}.
\end{proof}

\subsubsection{Anisotropic inequalities at \texorpdfstring{$H^2$}{H2} regularity}

\begin{lemma}[Two anisotropic estimates]
\label{lem:H2-anisotropic}
If \(f\in H^2(\R^3)\) and \(\nabla_hf\in H^2(\R^3)\), then
\begin{equation}\label{eq:H2-Lipschitz}
 \|\nabla f\|_\infty
 \leq C\big(\|f\|_{H^2}+\|\nabla_hf\|_{H^2}\big).
\end{equation}
Moreover, for smooth scalar functions, or componentwise for vector-valued
functions,
\begin{equation}\label{eq:H2-mixed-triple}
 \int_{\R^3}|abc|\,\dd x
 \leq C\|a\|_{H^1}\|b\|_2^{1/3}
       \|\nabla_hb\|_2^{2/3}\|c\|_2.
\end{equation}
\end{lemma}

\begin{proof}
With \(\langle\xi\rangle=(1+|\xi|^2)^{1/2}\), Fourier inversion
and Cauchy--Schwarz give
\[
 \|\nabla f\|_\infty
 \leq C\left(\int_{\R^3}
 \frac{|\xi|^2}{\langle\xi\rangle^4(1+|\xi_h|^2)}\,\dd\xi
 \right)^{1/2}
 \left(\int_{\R^3}\langle\xi\rangle^4(1+|\xi_h|^2)
 |\widehat f(\xi)|^2\,\dd\xi\right)^{1/2}.
\]
The first integral is finite.  Indeed, with \(r=|\xi_h|\), integration in
\(\xi_3\), including the factor \((1+r^2)^{-1}\), produces a quantity
bounded by \(C(1+r^2)^{-3/2}\); the remaining radial integrand is bounded
by \(Cr(1+r^2)^{-3/2}\).  The second factor is equivalent to the
right-hand side of \eqref{eq:H2-Lipschitz}.

For \eqref{eq:H2-mixed-triple}, the one-dimensional Sobolev
inequality and horizontal H\"older give
\[
 \|a\|_{L^3_hL^\infty_3}^2
 \leq C\bigl\|\|a(x_h,\cdot)\|_{L^2_3}\|\pa_3a(x_h,\cdot)\|_{L^2_3}\bigr\|_{L^{3/2}_h}
 \leq C\|a\|_{L^6_hL^2_3}\|\pa_3a\|_2
 \leq C\|a\|_{H^1}^2.
\]
The Hilbert-valued two-dimensional Gagliardo--Nirenberg inequality gives
\[
 \|b\|_{L^6_hL^2_3}
 \leq C\|b\|_2^{1/3}\|\nabla_hb\|_2^{2/3}.
\]
Horizontal H\"older with exponents \((3,6,2)\) and vertical H\"older with
exponents \((\infty,2,2)\) prove \eqref{eq:H2-mixed-triple}.
\end{proof}

\subsection{Local existence and uniqueness}

The local construction uses spatial regularization and compactness
for hyperbolic--parabolic systems; see \cite{KawashimaThesis,ShizutaKawashima}
and \cite[Section~3]{FengHongWangWu}. The uniform $H^2$ estimates for
the present equations are derived in Section~\ref{sec:global-proof}.

\begin{lemma}[Local existence and uniqueness]\label{lem:local}
There are $\delta_{\rm loc}>0$ and a positive time $T_{\rm loc}$,
depending only on the fixed coefficients, the equilibrium and the
prescribed $H^2$ bound, such that every $U_0\in H^2(\R^3)$ with
$\|U_0\|_{H^2}\le\delta_{\rm loc}$ generates a unique strong solution
of \eqref{eq:perturbed} on $[0,T_{\rm loc}]$.
It has the regularities \eqref{eq:solspace1}--\eqref{eq:solspace3}
on that finite interval, and $\rho=\bar\rho e^q$, $\theta=\bar\theta e^\tau$
remain uniformly positive.
A solution satisfying
$\sup_{0\le t<T_*}\|U(t)\|_{H^2}\le\delta_{\rm loc}/2$
continues beyond any finite $T_*$.
\end{lemma}

The $H^2$ energy and integrated dissipation bounds give a common
local existence time. By \eqref{eq:H2-Lipschitz},
$\nabla u\in L^1(0,T_{\rm loc};L^\infty)$, which permits spatial
regularization of the energy identity and gives strong $H^2$
continuity. The time regularity follows from the equations and
\eqref{eq:H2-H1-multiplier}. Uniqueness and finite-time $L^2$
continuous dependence follow from Lemma~\ref{lem:r10-difference}
in Appendix~\ref{app:construction-tools}; the continuation assertion
follows by restarting with the common local existence time.

\subsection{The nonlinear energy estimate}\label{sec:global-proof}

Let $U$ be a smooth solution on $[0,T]$ satisfying
$\sup_{0\le t\le T}\|U(t)\|_{H^2}\le\delta$, where $\delta>0$
is sufficiently small. We derive an energy estimate with constants
independent of $T$, and then apply it to the local solution of
Lemma~\ref{lem:local}.

\subsubsection{Zeroth-order energy estimate}

Define the relative thermal energy by
\begin{equation*}
\mathscr H(\rho,\theta\mid\bar\rho,\bar\theta)
 ={}c_v\rho\left(\theta-\bar\theta-\bar\theta
                 \log\frac{\theta}{\bar\theta}\right)
 +R\bar\theta\left(\rho\log\frac{\rho}{\bar\rho}-\rho+\bar\rho\right),
\end{equation*}
and set
\[
 \cE_{\rm th}(t)=\int_{\R^3}
 \left(\frac12\rho|u|^2+
       \mathscr H(\rho,\theta\mid\bar\rho,\bar\theta)\right)\dd x.
\]
\begin{lemma}[Zeroth-order energy estimate]\label{lem:zeroth-order-energy}
Let $(\rho,u,\theta)$ be a smooth solution of \eqref{eq:original}
on $[0,T]$ such that $\sup_{0\le t\le T}\|U(t)\|_{H^2}$ is
sufficiently small. Then
\begin{equation}\label{eq:freeenergyidentity}
 \frac{\dd}{\dd t}\cE_{\rm th}(t)
 +\bar\theta\int_{\R^3}\left(
       \frac{Q(u)}{\theta}
       +\kappa\frac{|\nabla\theta|^2}{\theta^2}\right)\dd x=0,
\end{equation}
and
\begin{equation}\label{eq:freeequiv}
 C^{-1}\|U(t)\|_2^2\le\cE_{\rm th}(t)\le C\|U(t)\|_2^2.
\end{equation}
The dissipation in \eqref{eq:freeenergyidentity} is equivalent to
$\|\nabla_hu\|_2^2+\|\Div u\|_2^2+\|\nabla\tau\|_2^2$.
Consequently, for $0\le t\le T$,
\[
 \|U(t)\|_2^2+\int_0^t
 \left(\|\nabla_hu(s)\|_2^2+\|\Div u(s)\|_2^2
                    +\|\nabla\tau(s)\|_2^2\right)\dd s
 \le C\|U_0\|_2^2.
\]
\end{lemma}
\begin{proof}
Multiply the continuity, momentum and temperature equations in
\eqref{eq:original}, respectively, by
\[
 \frac12|u|^2+\partial_\rho\mathscr H,
 \qquad u,\qquad 1-\frac{\bar\theta}{\theta},
\]
where
\[
 \partial_\rho\mathscr H
 =c_v\left(\theta-\bar\theta-\bar\theta
                     \log\frac{\theta}{\bar\theta}\right)
      +R\bar\theta\log\frac{\rho}{\bar\rho}.
\]
The continuity equation combines the time and transport terms
into the derivative of $\cE_{\rm th}$. The pressure terms cancel
after integration, since
\[
 \rho\partial_\rho\mathscr H-\mathscr H
       =R\bar\theta(\rho-\bar\rho).
\]
The viscous terms give $-\bar\theta Q(u)/\theta$, and the heat
term gives $-\kappa\bar\theta|\nabla\theta|^2/\theta^2$.
Adding and integrating by parts therefore proves
\eqref{eq:freeenergyidentity}.
These calculations can first be made with a spatial cutoff.
The relative energy fluxes are integrable in space on finite
time intervals, so the boundary terms vanish as the cutoff
radius tends to infinity. Approximation extends the identity
to the $H^2$ solution.

For small $H^2$ perturbations, $\rho$ and $\theta$ stay in a bounded
positive range. Taylor's formula gives
\[
 \int_{\R^3}\mathscr H\,\dd x
 \simeq\|\rho-\bar\rho\|_2^2+\|\theta-\bar\theta\|_2^2
 \simeq\|q\|_2^2+\|\tau\|_2^2.
\]
Together with the kinetic energy, this proves \eqref{eq:freeequiv}.
The dissipation equivalence follows from
$\nabla\theta/\theta=\nabla\tau$ and the positive upper and lower
bounds on $\theta$. Integration of \eqref{eq:freeenergyidentity}
then yields the asserted $L^2$ estimate.
\end{proof}

\subsubsection{First- and second-order energy estimates}

Set
\[
 \cE_{+,2}:=\frac12\sum_{1\leq|\beta|\leq2}
 \left(c^2\|\pa^\beta q\|_2^2+\|\pa^\beta u\|_2^2
 +\omega\|\pa^\beta\tau\|_2^2\right).
\]
For $1\le|\beta|\le2$, apply $\partial^\beta$ to
\eqref{eq:perturbed} and test the three equations against
$c^2\partial^\beta q$, $\partial^\beta u$ and
$\omega\partial^\beta\tau$, respectively. Integration by parts,
using $\omega\alpha=c^2$, yields
\begin{equation}\label{eq:r13-positive-order-identity}
\begin{aligned}
 &\frac12\frac{\dd}{\dd t}
 \left(c^2\|\partial^\beta q\|_2^2
      +\|\partial^\beta u\|_2^2
      +\omega\|\partial^\beta\tau\|_2^2\right)\\
 &\quad+\nu_h\|\nabla_h\partial^\beta u\|_2^2
        +\nu_c\|\Div\partial^\beta u\|_2^2
        +\omega\chi\|\nabla\partial^\beta\tau\|_2^2\\
 &=c^2\ip{\partial^\beta F_1}{\partial^\beta q}
   +\ip{\partial^\beta F_2}{\partial^\beta u}
   +\omega\ip{\partial^\beta F_3}{\partial^\beta\tau}.
\end{aligned}
\end{equation}
All norms below are evaluated at the same time. We use
\eqref{eq:H2-dissipation-definition} and the bounds
\[
 \|\nabla q\|_{H^1}\le\|U\|_{H^2},\qquad
 \|\nabla\tau\|_{H^1}\le
       \min\{\|U\|_{H^2},\|\nabla\tau\|_{H^2}\}.
\]
These two bounds are used when absorbing the transport and
variable-coefficient terms.

\begin{lemma}[Nonlinear energy estimate at $H^2$ regularity]
\label{prop:H2-nonlinear}
There are \(C,\delta>0\) such that, whenever \(\|U\|_{H^2}\leq\delta\),
\begin{equation}\label{eq:H2-NLenergy}
 \left|\sum_{1\leq|\beta|\leq2}\left(
 c^2\ip{\pa^\beta F_1}{\pa^\beta q}
 +\ip{\pa^\beta F_2}{\pa^\beta u}
 +\omega\ip{\pa^\beta F_3}{\pa^\beta\tau}\right)\right|
 \leq C\|U\|_{H^2}\cD_2(U).
\end{equation}
\end{lemma}

\begin{proof}
For \(|\beta|\leq2\), Leibniz's formula gives
\begin{align}
 \ip{\pa^\beta(u\cdot\nabla f)}{\pa^\beta f}
 ={}&-\frac12\int_{\R^3}\Div u\,|\pa^\beta f|^2\,\dd x
 \notag\\
 &+\sum_{0<\gamma\leq\beta}C_{\beta\gamma}
 \int_{\R^3}\pa^\gamma u\cdot
 \nabla\pa^{\beta-\gamma}f\,\pa^\beta f\,\dd x.
 \label{eq:H2-transport-identity}
\end{align}

\smallskip\noindent\emph{Step 1. Scalar transport.}
For $z=q,\tau$, take $f=z$ in \eqref{eq:H2-transport-identity}.
The principal transport term and the commutators
with one derivative on $u$ have the common bound
\[
 C\bigl(\|U\|_{H^2}+\|(\nabla_hu,\Div u)\|_{H^2}\bigr)
                                      \|\nabla z\|_{H^1}^2.
\]
For density use $\|\nabla q\|_{H^1}\le\|U\|_{H^2}$; for
temperature also use $\|\nabla\tau\|_{H^1}\le\|\nabla\tau\|_{H^2}$.
At second order, the commutators with two derivatives on $u$ have the form
\[
 \int\nabla^2u\,\nabla z\,\nabla^2z\,\dd x.
\]
Apply \eqref{eq:H2-mixed-triple} to the three factors
$\nabla z$, $\nabla^2u$ and $\nabla^2z$. The resulting bound is
\[
 C\|U\|_{H^2}^{1/3}\|(\nabla_hu,\Div u)\|_{H^2}^{2/3}
                                      \|\nabla z\|_{H^1}^2,
 \qquad z=q,\tau.
\]
Since $\|\nabla z\|_{H^1}\le\|U\|_{H^2}$ and
$\|\nabla z\|_{H^1}^2\le\cD_2(U)$, weighted Young's inequality
bounds both scalar transport contributions by
$C\|U\|_{H^2}\cD_2(U)$. This includes the
pure vertical commutator
\[
 \int\pa_3^2u_h\cdot\nabla_hz\,\pa_3^2z\,\dd x.
\]

\smallskip\noindent\emph{Step 2. Velocity convection.}
The third velocity component satisfies
\begin{equation}
 \pa_3^{j+1}u_3=\pa_3^j\Div u-\pa_3^j\Div_hu_h,\qquad
 \|\nabla^j\nabla u_3\|_2\leq C \|(\nabla_hu,\Div u)\|_{H^2},\quad0\leq j\leq2 .
\end{equation}

Take \(f=u\) in \eqref{eq:H2-transport-identity}. For \(1\leq|\beta|\leq2\),
\eqref{eq:coordtriple}, with the derivatives \(3,1,2\) assigned to
\(\Div u,\pa^\beta u,\pa^\beta u\), respectively, yields
\begin{align}
 \left|\int\Div u\,|\pa^\beta u|^2\,\dd x\right|
 &\leq C\|\Div u\|_2^{1/2}\|\pa_3\Div u\|_2^{1/2}
 \|\pa^\beta u\|_2
 \notag\\[-1mm]
 &\hspace{17mm}\times
 \|\pa_1\pa^\beta u\|_2^{1/2}
 \|\pa_2\pa^\beta u\|_2^{1/2}
 \leq C \|u\|_{H^2}\|(\nabla_hu,\Div u)\|_{H^2}^2 .
 \label{eq:H2-velocity-principal}
\end{align}
For \(k=|\beta|\), split the velocity commutator exactly as
\begin{align}
 [\pa^\beta,u\cdot\nabla]u
 ={}&[\pa^\beta,u_h\cdot\nabla_h]u
 +[\pa^\beta,u_3\pa_3]u_3\,e_3
 +[\pa^\beta,u_3\pa_3]u_h .
 \label{eq:H2-velocity-comm-split}
\end{align}
In the first two terms of \eqref{eq:H2-velocity-comm-split}, the
differentiated transport factor is controlled by $(\nabla_hu,\Div u)$,
using \eqref{eq:verticalcomp} for $\partial_3u_3$.
Apply \eqref{eq:coordtriple} with the vertical derivative on this factor
and one horizontal derivative on each of the other two factors.
Then
\begin{equation}
 \left|\int\Big\{
 [\pa^\beta,u_h\cdot\nabla_h]u
 +[\pa^\beta,u_3\pa_3]u_3\,e_3\Big\}\cdot\pa^\beta u\,\dd x\right|
 \le C\|u\|_{H^2}\|(\nabla_hu,\Div u)\|_{H^2}^2.
\end{equation}

If $\beta_h\ne0$, choose $i\in\{1,2\}$ with $\beta_i>0$.
In each commutator product, apply \eqref{eq:coordtriple} with coordinate
derivatives $1,2,3$ on $\partial^\gamma u_3$,
$\partial_3\partial^{\beta-\gamma}u_h$ and
$\partial^\beta u_h$, respectively. The vertical derivative on the
last factor contains the horizontal derivative $\partial_i$, so
\begin{equation}
 \left|\int[\pa^\beta,u_3\pa_3]u_h\cdot\pa^\beta u_h\,\dd x\right|
 \le C\|u\|_{H^2}\|(\nabla_hu,\Div u)\|_{H^2}^2.
\end{equation}
If \(\beta=(0,0,k)\), then
\begin{align}
 &\int[\pa_3^k,u_3\pa_3]u_h\cdot\pa_3^ku_h\,\dd x
 =\sum_{\ell=1}^k\binom{k}{\ell}
 \int\pa_3^{\ell-1}(\Div u-\Div_hu_h)\,
 \pa_3^{k+1-\ell}u_h\cdot\pa_3^ku_h\,\dd x ,\\
 &\left|\int[\pa_3^k,u_3\pa_3]u_h\cdot\pa_3^ku_h\,\dd x\right|
 \notag\\
 &\quad\leq C\sum_{\ell=1}^k
 \|\pa_3^{\ell-1}(\Div u-\Div_hu_h)\|_2^{1/2}
 \|\pa_3^\ell(\Div u-\Div_hu_h)\|_2^{1/2}
 \notag\\[-1mm]
 &\hspace{17mm}\times
 \|\pa_3^{k+1-\ell}u_h\|_2^{1/2}
 \|\pa_1\pa_3^{k+1-\ell}u_h\|_2^{1/2}
 \|\pa_3^ku_h\|_2^{1/2}
 \|\pa_2\pa_3^ku_h\|_2^{1/2}
 \notag\\
 &\quad\leq C\|u\|_{H^2}\|(\nabla_hu,\Div u)\|_{H^2}^2 .
 \label{eq:H2-velocity-comm-vertical}
\end{align}
Equations \eqref{eq:H2-transport-identity} and
\eqref{eq:H2-velocity-principal}--\eqref{eq:H2-velocity-comm-vertical} give
\begin{equation}
 \left|\sum_{1\leq|\beta|\leq2}
 \ip{\pa^\beta(u\cdot\nabla u)}{\pa^\beta u}\right|
 \leq C\|u\|_{H^2}\|(\nabla_hu,\Div u)\|_{H^2}^2 .
\end{equation}
Since $\|u\|_{H^2}\le\|U\|_{H^2}$, this contribution is bounded by
$C\|U\|_{H^2}\cD_2(U)$.

\smallskip\noindent\emph{Step 3. Temperature-dependent pressure.}
For $|\beta|=2$, integrate the principal term
$(e^\tau-1)\nabla\pa^\beta(q+\tau)$ by parts once.
Its contribution is bounded by
\[
 C\|U\|_{H^2}
 (\|\nabla q\|_{H^1}+\|\nabla\tau\|_{H^2})
                         \|(\nabla_hu,\Div u)\|_{H^2}.
\]
The remaining commutators contain one or two derivatives of
the temperature coefficient and satisfy
\[
\begin{aligned}
 &\left|\int\nabla(e^\tau-1)\,\nabla^2(q+\tau)\,\nabla^2u\dx\right|
 \le\|\nabla(e^\tau-1)\|_\infty
                \|\nabla^2(q+\tau)\|_2\|\nabla^2u\|_2\\
 &\quad\le C\|U\|_{H^2}\|\nabla\tau\|_{H^2}
                 (\|\nabla q\|_{H^1}+\|\nabla\tau\|_{H^2}),\\[2mm]
 &\left|\int\nabla^2(e^\tau-1)\,\nabla(q+\tau)\,\nabla^2u\dx\right|
 \le\|\nabla^2(e^\tau-1)\|_6
                  \|\nabla(q+\tau)\|_3\|\nabla^2u\|_2\\
 &\quad\le C\|U\|_{H^2}\|\nabla\tau\|_{H^2}
                 (\|\nabla q\|_{H^1}+\|\nabla\tau\|_{H^2}).
\end{aligned}
\]
The coefficient bounds used here are
\[
 \|\nabla(e^\tau-1)\|_\infty+\|\nabla^2(e^\tau-1)\|_6
       \le C\|\nabla\tau\|_{H^2},
 \qquad \|\nabla(q+\tau)\|_3
       \le C(\|\nabla q\|_{H^1}+\|\nabla\tau\|_{H^2}).
\]
With one derivative, the same integration by parts uses
$\|\nabla(e^\tau-1)\|_6\|\nabla(q+\tau)\|_3\|\nabla u\|_2$.
Young's inequality bounds the full pressure contribution by
$C\|U\|_{H^2}\cD_2(U)$.

\smallskip\noindent\emph{Step 4. Variable viscosity and heat conductivity.}
Write
\[
 \mathbb G_{ij}(u)=\nu_h\mathbf1_{\{j=1,2\}}\pa_ju_i
 +\nu_c\delta_{ij}\Div u,\qquad \mathcal Lu_i=\pa_j\mathbb G_{ij}(u).
\]
In the \(|\beta|=2\) pairing, integrate only the
principal \((e^{-q}-1)\mathcal L\pa^\beta u\) term. Its quadratic
form is bounded by $C\|U\|_{H^2}\|(\nabla_hu,\Div u)\|_{H^2}^2$.
The term containing $\nabla(e^{-q}-1)$ and the
\(|\gamma|=1\) commutators are bounded by
\begin{equation}\label{eq:H2-variable-viscosity-top}
\begin{aligned}
 &C\|\nabla(e^{-q}-1)\|_{H^1}\|\nabla^2u\|_2^{1/3}
     \|\nabla_h\nabla^2u\|_2^{2/3}\|(\nabla_hu,\Div u)\|_{H^2}\\
 &\qquad\le C\|\nabla q\|_{H^1}\|U\|_{H^2}^{1/3}
                            \|(\nabla_hu,\Div u)\|_{H^2}^{5/3}.
\end{aligned}
\end{equation}
The factors in \eqref{eq:H2-variable-viscosity-top} satisfy
\begin{equation*}
 \begin{aligned}
 \|\nabla(e^{-q}-1)\|_{H^1}&\le C\|\nabla q\|_{H^1},
 \qquad \|\nabla^2u\|_2\le\|U\|_{H^2},\\
 \|\nabla_h\nabla^2u\|_2+\|\nabla^2\Div u\|_2
                  &\le C\|(\nabla_hu,\Div u)\|_{H^2}.
\end{aligned}
\end{equation*}
For \(|\beta|=2\), the third-order velocity derivatives generated by
\(\mathbb G(\pa^\beta u)\) belong to
\begin{equation*}
 \{\pa_h\pa^\beta u:\ |\beta|=2\}
 \cup
 \{\pa^\beta\Div u:\ |\beta|=2\},
\end{equation*}
and \(\pa_3^{j+1}u_3=\pa_3^j\Div u-\pa_3^j\Div_hu_h\) for \(j=0,1,2\).
For \(|\gamma|=2\), apply \eqref{eq:H2-mixed-triple} with
\(a=\mathcal Lu\), \(b=\nabla^2u\), and \(c=\nabla^2(e^{-q}-1)\).  Since
\(\|\mathcal Lu\|_{H^1}\leq C\|(\nabla_hu,\Div u)\|_{H^2}\) and \(\|\nabla^2(e^{-q}-1)\|_2\leq C\|\nabla q\|_{H^1}\),
\begin{equation*}
 \int_{\R^3}|\mathcal Lu|\,|\nabla^2u|\,|\nabla^2(e^{-q}-1)|\dd x
 \le C\|\nabla q\|_{H^1}\|U\|_{H^2}^{1/3}
                            \|(\nabla_hu,\Div u)\|_{H^2}^{5/3}.
\end{equation*}
Since $\|\nabla q\|_{H^1}\le\|U\|_{H^2}$,
\[
 \|\nabla q\|_{H^1}\|U\|_{H^2}^{1/3}
              \|(\nabla_hu,\Div u)\|_{H^2}^{5/3}
 \le \|U\|_{H^2}\|\nabla q\|_{H^1}^{1/3}
              \|(\nabla_hu,\Div u)\|_{H^2}^{5/3}.
\]
Weighted Young's inequality therefore gives
\begin{equation*}
 \left|\sum_{1\leq|\beta|\leq2}
 \ip{\pa^\beta[(e^{-q}-1)\mathcal Lu]}{\pa^\beta u}\right|
 \leq C\|U\|_{H^2}(\|\nabla q\|_{H^1}^2+\|(\nabla_hu,\Div u)\|_{H^2}^2).
\end{equation*}
For the variable heat coefficient, the exact expansion is
\[
 \begin{aligned}
 \ip{\pa^\beta((e^{-q}-1)\Delta\tau)}{\pa^\beta\tau}
 ={}&-\int (e^{-q}-1)|\nabla\pa^\beta\tau|^2\,\dd x
     -\int\nabla (e^{-q}-1)\cdot\nabla\pa^\beta\tau\,\pa^\beta\tau\,\dd x\\
 &+\sum_{0<\gamma\leq\beta}C_{\beta\gamma}
   \int\pa^\gamma (e^{-q}-1)\,\Delta\pa^{\beta-\gamma}\tau\,
            \pa^\beta\tau\,\dd x,\qquad 1\leq|\beta|\leq2 .
 \end{aligned}
\]
For the coefficient-gradient term and $|\gamma|=1$, use
\eqref{eq:H2-mixed-triple} with \(a=\nabla (e^{-q}-1)\), \(b=\pa^\beta\tau\),
and $c=\nabla\partial^\beta\tau$ for the coefficient-gradient
term and $c=\Delta\partial^{\beta-\gamma}\tau$ for $|\gamma|=1$.
For \(|\gamma|=2\), use \(a=\Delta\tau\), \(b=\pa^\beta\tau\),
and \(c=\pa^\gamma (e^{-q}-1)\).  Since
\[
 \|\nabla (e^{-q}-1)\|_{H^1}\leq C\|\nabla q\|_{H^1},\quad
 \|\pa^\beta\tau\|_2\leq \|U\|_{H^2},\quad
 \|\nabla_h\pa^\beta\tau\|_2+\|\Delta\tau\|_{H^1}\leq C\|\nabla\tau\|_{H^2},
\]
these terms are bounded by \(C\|\nabla q\|_{H^1}\|U\|_{H^2}^{1/3}\|\nabla\tau\|_{H^2}^{5/3}\).
Applying the same weighted Young inequality with
$\|\nabla\tau\|_{H^2}$ in place of
$\|(\nabla_hu,\Div u)\|_{H^2}$ gives
\begin{equation*}
 \begin{aligned}
 &\left|\sum_{1\le|\beta|\le2}
    \ip{\pa^\beta[(e^{-q}-1)\Delta\tau]}{\pa^\beta\tau}\right|
 \le C\|U\|_{H^2}\|\nabla\tau\|_{H^2}^2
       +C\|\nabla q\|_{H^1}\|U\|_{H^2}^{1/3}\|\nabla\tau\|_{H^2}^{5/3}\\
 &\quad\le C\|U\|_{H^2}
                    (\|\nabla q\|_{H^1}^2+\|\nabla\tau\|_{H^2}^2).
\end{aligned}
\end{equation*}

\smallskip\noindent\emph{Step 5. Viscous heating and the quadratic temperature source.}
The $H^2$ algebra estimate, with each exponential written as
a constant plus an $H^2$ function, gives
\[
 \begin{aligned}
 \|e^{-q-\tau}Q(u)\|_{H^2}
    &\le C(1+\|U\|_{H^2})\|(\nabla_hu,\Div u)\|_{H^2}^2,\\
 \|e^{-q}|\nabla\tau|^2\|_{H^2}
    &\le C(1+\|U\|_{H^2})\|\nabla\tau\|_{H^2}^2.
 \end{aligned}
\]
Pairing with $\partial^\beta\tau$, $1\le|\beta|\le2$, bounds
the production terms by $C\|U\|_{H^2}\cD_2(U)$.
This proves \eqref{eq:H2-NLenergy}.
\end{proof}

\subsubsection{Density compensation and the closed energy estimate}

The positive-order energy does not contain density dissipation.
The following cross term recovers $\|\nabla q\|_{H^1}^2$ without
requiring a third density derivative.

\begin{lemma}[Density compensation]

For a smooth solution under the smallness assumption of
Lemma~\ref{prop:H2-nonlinear}, set
\[
 \cX_2:=\sum_{|\beta|\le1}
              \ip{\pa^\beta u}{\nabla\pa^\beta q}.
\]
The sources satisfy
\begin{equation}\label{eq:H2-NLsource}
 \|F_1\|_{H^1}+\|F_2\|_{H^1}
 \le C\|U\|_{H^2}\cD_2(U)^{1/2}.
\end{equation}
The cross term satisfies
\begin{equation}\label{eq:r13-density-compensation}
\cX_2'+\frac{c^2}{2}\|\nabla q\|_{H^1}^2
 \le C\bigl(\|(\nabla_hu,\Div u)\|_{H^2}^2
                          +\|\nabla\tau\|_{H^2}^2\bigr)
 +C\|U\|_{H^2}\cD_2(U).
\end{equation}
\end{lemma}
\begin{proof}
\smallskip\noindent\emph{Step 1. Source estimates in $H^1$.}
We first prove \eqref{eq:H2-NLsource}. For velocity convection,
\begin{equation}\label{eq:H2-velocity-convection}
 \|u\cdot\nabla u\|_{H^1}
 \le C\|U\|_{H^2}\|(\nabla_hu,\Div u)\|_{H^2}.
\end{equation}
To prove \eqref{eq:H2-velocity-convection}, use
$H^2\cdot H^1\subset H^1$ and
$\pa_3u_3=\Div u-\Div_hu_h$ to obtain
\[
 \|u_h\cdot\nabla_hu\|_{H^1}+\|u_3\pa_3u_3\|_{H^1}
       \le C\|U\|_{H^2}\|(\nabla_hu,\Div u)\|_{H^2}.
\]
For \(|\beta|\leq1\), \eqref{eq:coordproduct} gives
\[
 \begin{aligned}
 \|u_3\pa_3\pa^\beta u_h\|_2
 &\leq C\bigl(\|u_3\|_2\|\pa_1u_3\|_2
             \|\pa_3u_3\|_2\|\pa_1\pa_3u_3\|_2\bigr)^{1/4}\\
 &\qquad\times
 \bigl(\|\pa_3\pa^\beta u_h\|_2
             \|\pa_2\pa_3\pa^\beta u_h\|_2\bigr)^{1/2}
 \leq C\|U\|_{H^2}^{3/4}\|(\nabla_hu,\Div u)\|_{H^2}^{5/4}.
 \end{aligned}
\]
There is also the bound
\[
 \|u_3\|_\infty\|\pa_3\pa^\beta u_h\|_2\le C\|U\|_{H^2}^2.
\]
Taking the smaller of these two bounds gives the required
source estimate, since
\[
 \begin{aligned}
 \min\bigl\{\|U\|_{H^2}^2,
       \|U\|_{H^2}^{3/4}\|(\nabla_hu,\Div u)\|_{H^2}^{5/4}\bigr\}
 \le \|U\|_{H^2}\|(\nabla_hu,\Div u)\|_{H^2}.
 \end{aligned}
\]
The differentiated coefficient term satisfies
\[
\|\pa_i u_3\,\pa_3u_h\|_2
 \le\|\pa_i u_3\|_6\|\pa_3u_h\|_3
 \le C\|\nabla u_3\|_{H^1}\|u\|_{H^2}
 \le C\|U\|_{H^2}\|(\nabla_hu,\Div u)\|_{H^2}.
\]

The other terms in $F_1,F_2$ are estimated by
\eqref{eq:H2-H1-multiplier}:
\begin{align*}
 \|u\cdot\nabla q\|_{H^1}
 &\le C\|U\|_{H^2}\|\nabla q\|_{H^1},\\
 \|(e^\tau-1)\nabla(q+\tau)\|_{H^1}
 +\|(e^{-q}-1)\mathcal Lu\|_{H^1}
 &\le C\|U\|_{H^2}\cD_2(U)^{1/2}.
\end{align*}
Here the $H^2$ norms of $e^\tau-1$ and $e^{-q}-1$ are bounded by
$C\|U\|_{H^2}$, while
$\|\mathcal Lu\|_{H^1}\le C\|(\nabla_hu,\Div u)\|_{H^2}$.
Together with \eqref{eq:H2-velocity-convection}, these bounds prove
\eqref{eq:H2-NLsource}.

\smallskip\noindent\emph{Step 2. The density--velocity cross term.}
For $|\beta|\le1$, differentiation of the pairing and the
first two equations in \eqref{eq:perturbed} give
\begin{align}
 \frac{\dd}{\dd t}\ip{\pa^\beta u}{\nabla\pa^\beta q}
 +c^2\|\nabla\pa^\beta q\|_2^2
={}&-c^2\ip{\nabla\pa^\beta\tau}{\nabla\pa^\beta q}\notag\\
 &+\ip{\pa^\beta(\nu_h\Delta_hu+\nu_c\nabla\Div u)}
                 {\nabla\pa^\beta q}\notag\\
 &+\|\Div\pa^\beta u\|_2^2
 +\ip{\pa^\beta F_2}{\nabla\pa^\beta q}
 -\ip{\Div\pa^\beta u}{\pa^\beta F_1}.
 \label{eq:crosscalculation}
\end{align}
The last term follows from
\(\ip{\pa^\beta u}{\nabla\pa^\beta F_1}
=-\ip{\Div\pa^\beta u}{\pa^\beta F_1}\);
after summing over $|\beta|\le1$, only $H^1$ norms of the
sources occur.

Consequently, the nonlinear terms in \eqref{eq:crosscalculation} satisfy
\begin{align*}
 &\left|\sum_{|\beta|\le1}\left(
 \ip{\pa^\beta F_2}{\nabla\pa^\beta q}
 -\ip{\Div\pa^\beta u}{\pa^\beta F_1}\right)\right|\\
 &\quad\le C\bigl(\|F_2\|_{H^1}\|\nabla q\|_{H^1}
                 +\|F_1\|_{H^1}\|\Div u\|_{H^1}\bigr)
 \le C\|U\|_{H^2}\cD_2(U).
\end{align*}

The temperature and viscosity terms in \eqref{eq:crosscalculation}
are bounded by
\[
 \frac{c^2}{2}\|\nabla q\|_{H^1}^2
 +C\bigl(\|(\nabla_hu,\Div u)\|_{H^2}^2
                          +\|\nabla\tau\|_{H^2}^2\bigr).
\]
Summing \eqref{eq:crosscalculation} over $|\beta|\le1$ proves
\eqref{eq:r13-density-compensation}.
\end{proof}

\begin{proposition}[Closed \(H^2\) energy]
\label{prop:H2-closed}
There are \(\eta_2,c_2,\delta_2'>0\) and an energy
\[
 \mathfrak E_2:=\cE_{\rm th}+\cE_{+,2}+\eta_2\cX_2
 \simeq\|U\|_{H^2}^2
\]
such that, whenever \(\|U\|_{H^2}\leq\delta_2'\),
\begin{equation}\label{eq:H2-closed}
 \frac{\dd}{\dd t}\mathfrak E_2+c_2\cD_2(U)\leq0.
\end{equation}
The estimate also holds for the local $H^2$ solution of
Lemma~\ref{lem:local}, in the integrated form
\begin{equation}\label{eq:strong-integrated-energy}
 \mathfrak E_2(U(t))+c_2\int_s^t\cD_2(U(r))\dd r
       \le\mathfrak E_2(U(s)),\qquad 0\le s\le t\le T,
\end{equation}
provided $\sup_{[0,T]}\|U\|_{H^2}\le\delta_2'$.
\end{proposition}

\begin{proof}
By Cauchy--Schwarz,
\begin{equation}\label{eq:H2-cross-equivalence}
 |\cX_2|\le C\|u\|_{H^1}\|\nabla q\|_{H^1}
             \le C\|U\|_{H^2}^2.
\end{equation}
First choose $\eta_2>0$ sufficiently small. Then
\eqref{eq:H2-cross-equivalence} gives
$\mathfrak E_2\simeq\|U\|_{H^2}^2$.
Sum \eqref{eq:r13-positive-order-identity}, add the zeroth-order
estimate of Lemma~\ref{lem:zeroth-order-energy} and $\eta_2$ times
\eqref{eq:r13-density-compensation}, and use \eqref{eq:H2-NLenergy}.
The viscous and thermal terms from the compensation estimate are
absorbed for this choice of $\eta_2$, giving
\[
 \mathfrak E_2'+c_*\cD_2(U)
       \le C\|U\|_{H^2}\cD_2(U)
\]
for fixed $c_*,C>0$.
Next choose $0<\delta_2'\le\delta$ so that
\begin{equation}
 C\delta_2'\le c_*/2.
\end{equation}
The nonlinear term is then absorbed, proving \eqref{eq:H2-closed}
with $c_2=c_*/2$.
For the local $H^2$ solution, spatial regularization gives
\eqref{eq:strong-integrated-energy}. The nonlinear pairings are
integrable by \eqref{eq:H2-NLenergy} and \eqref{eq:H2-NLsource},
and the transport commutator is controlled by
\eqref{eq:H2-Lipschitz}. Strong $H^2$ continuity permits passage
at the endpoints $s,t$.
\end{proof}

\subsection{Completion of the proof of Theorem~\ref{thm:global}}

\begin{proof}[Proof of Theorem~\ref{thm:global}]
Let $U$ be the local solution of Lemma~\ref{lem:local}, and let
$[0,T_*)$ be its maximal existence interval.
Choose $\delta_*\le\min(\delta_2',\delta_{\rm loc})$ so that
Proposition~\ref{prop:H2-closed} applies, and fix
$c_E,C_E>0$ with
$c_E\|U\|_{H^2}^2\le\mathfrak E_2(U)\le C_E\|U\|_{H^2}^2$.
Choose $\delta_2\le\frac12\delta_*\sqrt{c_E/C_E}$.

On every interval on which $\|U\|_{H^2}\le\delta_*$,
apply \eqref{eq:strong-integrated-energy} with $s=0$ to obtain
\begin{equation}\label{eq:H2-global-exact}
 c_E\|U(t)\|_{H^2}^2+c_2\int_0^t\cD_2(U(s))\dd s
       \le C_E\|U_0\|_{H^2}^2.
\end{equation}
Thus $\|U(t)\|_{H^2}\le\delta_*/2$.
The continuity $U\in C([0,T_*);H^2)$ keeps the solution in this
smaller ball by a first-exit argument. The logarithmic variables
also give uniform positive upper and lower bounds for $\rho,\theta$.

If $T_*<\infty$, the continuation assertion of
Lemma~\ref{lem:local} extends the solution beyond $T_*$,
a contradiction. Hence $T_*=\infty$.
Letting $t\to\infty$ in \eqref{eq:H2-global-exact} proves
\eqref{eq:globalbound}; local uniqueness and regularity give
\eqref{eq:solspace1}--\eqref{eq:solspace3} on the full interval.
\end{proof}

\section{Basic anisotropic decay through second order}\label{sec:r10-decay}

Throughout Sections~\ref{sec:r10-decay}--\ref{sec:third-order-decay},
$U$ is the solution of Theorem~\ref{thm:global} with sufficiently small
$\mathcal J_2$. We first establish decay through second order, then
full $H^2$ convergence and third-order estimates.
Appendix~\ref{sec:source-estimate-index} records the source bounds
and their uses.

\subsection{Linear and Duhamel estimates}

The Hodge decomposition separates the solenoidal horizontal heat
flow from the compressible acoustic--thermal system. We establish
semigroup bounds for both components and then estimate the
compressible Duhamel integrals with the source regularity needed
in Section~\ref{sec:H2-decay}.

\subsubsection{Hodge decomposition and the acoustic--thermal block}

Throughout the linear analysis, \(\varrho=|\xi|\) denotes the frequency radius.
The constants \(C,c_*>0\) depend only on the fixed coefficients and may
vary between estimates; in particular, \(c_*\) is distinct from the
physical coefficient \(c^2=R\bar\theta\).

With the Hodge projectors \(\mathbb P,\mathbb Q\) used in
\eqref{eq:projected-P-intro}, put \(\Lambda=(-\Delta)^{1/2}\) and
\begin{equation*}
 d:=\Lambda^{-1}\Div u.
\end{equation*}
The homogeneous linearization of \eqref{eq:perturbed} is
\begin{equation*}
 \left\{
 \begin{aligned}
 q_t+\Div u&=0,\\
 u_t+c^2\nabla(q+\tau)-\nu_h\Delta_hu-\nu_c\nabla\Div u&=0,\\
 \tau_t+\alpha\Div u-\chi\Delta\tau&=0.
 \end{aligned}\right.
\end{equation*}
Application of \(\mathbb P\) gives the exact identity
\begin{equation*}
 (\mathbb Pu)_t-\nu_h\Delta_h\mathbb Pu=0.
\end{equation*}
For \(\varrho=|\xi|\) and
\begin{equation*}
 A(\xi):=\nu_h|\xi_h|^2+\nu_c\varrho^2,
\end{equation*}
the remaining Fourier variables solve
\begin{equation}\label{eq:fourierblock}
 \left\{
 \begin{aligned}
 \widehat q_t+\varrho\widehat d&=0,\\
 \widehat d_t-c^2\varrho(\widehat q+\widehat\tau)+A(\xi)\widehat d&=0,\\
 \widehat\tau_t+\alpha \varrho\widehat d+\chi \varrho^2\widehat\tau&=0.
 \end{aligned}\right.
\end{equation}

\subsubsection{Compensated Fourier energy}

\begin{lemma}[Compensated Fourier energy]\label{lem:fourier-energy}
There is a fixed sufficiently small \(\eta>0\), depending only on the
coefficients, such that the energy
\begin{equation*}
 \mathscr H_\eta(\xi):=
 c^2|\widehat q|^2+|\widehat d|^2+\omega|\widehat\tau|^2
 -2\eta\frac{\varrho}{1+\varrho^2}
      \Re(\widehat q\,\overline{\widehat d})
\end{equation*}
satisfies, for constants \(c_0,C_0>0\),
\begin{align}
 C_0^{-1}|(\widehat q,\widehat d,\widehat\tau)|^2
 &\leq\mathscr H_\eta(\xi)
 \leq C_0|(\widehat q,\widehat d,\widehat\tau)|^2,
 \label{eq:fourier-equivalence}\\
 \frac{\dd}{\dd t}\mathscr H_\eta
 +c_0\left\{
 A(\xi)|\widehat d|^2+\chi \varrho^2|\widehat\tau|^2
 +\frac{\varrho^2}{1+\varrho^2}|\widehat q|^2\right\}&\leq0.
 \label{eq:fourier-dissipation}
\end{align}
Consequently,
\begin{equation}\label{eq:blockmultiplier}
 |(\widehat q,\widehat d,\widehat\tau)(t,\xi)|
 \leq C_0\exp\left(-c_0\frac{\varrho^2}{1+\varrho^2}t\right)
 |(\widehat q,\widehat d,\widehat\tau)(0,\xi)|.
\end{equation}
\end{lemma}

\begin{proof}
The density dissipation is obtained by adding a compensating cross
term to the Fourier energy of \eqref{eq:fourierblock}.
Multiplication of the three equations in \eqref{eq:fourierblock} by
\(c^2\overline{\widehat q}\), \(\overline{\widehat d}\), and
\(\omega\overline{\widehat\tau}\), respectively, followed by taking real
parts and using \(\omega\alpha=c^2\), gives
\begin{equation}\label{eq:basic-fourier-exact}
 \frac12\frac{\dd}{\dd t}
 \{c^2|\widehat q|^2+|\widehat d|^2+
       \omega|\widehat\tau|^2\}
 +A(\xi)|\widehat d|^2+\omega\chi \varrho^2|\widehat\tau|^2=0.
\end{equation}
Furthermore,
\begin{equation}\label{eq:cross-fourier-exact}
 \frac{\dd}{\dd t}\Re(\widehat q\,\overline{\widehat d})
 =-\varrho|\widehat d|^2+c^2\varrho|\widehat q|^2
 +c^2\varrho\Re(\widehat q\,\overline{\widehat\tau})
 -A(\xi)\Re(\widehat q\,\overline{\widehat d}).
\end{equation}
To estimate the cross terms in \eqref{eq:cross-fourier-exact}, use
\begin{equation}\label{eq:Acomparison}
 \nu_c \varrho^2\leq A(\xi)\leq(\nu_h+\nu_c)\varrho^2,
 \qquad 0\leq \frac{\varrho}{1+\varrho^2}\leq\frac12,
\end{equation}
Young's inequality yields
\begin{align}
 c^2\frac{\varrho}{1+\varrho^2}\varrho|\widehat q\widehat\tau|
 &\leq\frac{c^2}{4}\frac{\varrho^2}{1+\varrho^2}|\widehat q|^2
 +C \varrho^2|\widehat\tau|^2,
 \label{eq:qtaucross}\\
 A(\xi)\frac{\varrho}{1+\varrho^2}|\widehat q\widehat d|
 &\leq\frac{c^2}{4}\frac{\varrho^2}{1+\varrho^2}|\widehat q|^2
 +C A(\xi)|\widehat d|^2,
 \label{eq:qdcross}\\
 \frac{\varrho}{1+\varrho^2}\varrho|\widehat d|^2
 &\leq C A(\xi)|\widehat d|^2.
 \label{eq:dbad}
\end{align}
Indeed, for \eqref{eq:qdcross},
\begin{equation*}
 \frac{A(\xi)^2\frac{\varrho^2}{(1+\varrho^2)^2}}{\varrho^2/(1+\varrho^2)}
 =\frac{A(\xi)^2}{1+\varrho^2}
 \leq(\nu_h+\nu_c)A(\xi).
\end{equation*}
Subtract \(\eta \frac{\varrho}{1+\varrho^2}\) times
\eqref{eq:cross-fourier-exact} from \eqref{eq:basic-fourier-exact}.
The bounds \eqref{eq:qtaucross}--\eqref{eq:dbad} allow the cross terms
to be absorbed into the dissipation if we choose
\begin{equation*}
 0<\eta\leq\min\left\{
 \frac{c^2}{8C},\frac{\omega\chi}{4C},
 \frac{1}{8C},\frac14\min(c^2,1,\omega)
 \right\}.
\end{equation*}
This proves \eqref{eq:fourier-equivalence} and
\eqref{eq:fourier-dissipation}.  By \eqref{eq:Acomparison}, the dissipation
in \eqref{eq:fourier-dissipation} is bounded below by
\begin{equation*}
 c_*\frac{\varrho^2}{1+\varrho^2}
 |(\widehat q,\widehat d,\widehat\tau)|^2.
\end{equation*}
Gronwall's inequality and \eqref{eq:fourier-equivalence} give
\eqref{eq:blockmultiplier}.
\end{proof}

\subsubsection{Basic anisotropic semigroup estimates}
For initial data \(V=(v_1,v_2,v_3)\), with scalar first and third
components and a vector-valued second component, let \(S_c(t)V\) denote the solution
\((q,\mathbb Qu,\tau)\) of the compressible linear subsystem whose initial
Fourier vector in \eqref{eq:fourierblock} is
\begin{equation}\label{eq:physical-to-block}
 \left(\widehat v_1,
       |\xi|^{-1}i\xi\cdot\widehat v_2,
       \widehat v_3\right).
\end{equation}
At \(\xi=0\) the middle entry is set equal to zero.  Put
\(S_h(t)=e^{\nu_ht\Delta_h}\). The projected shear evolution is $S_h(t)\mathbb P$.
Let $P_L^h$ be a smooth cutoff in horizontal frequency,
equal to one near $\xi_h=0$ and with compact frequency support,
and put $P_H^h=I-P_L^h$.
Constants may depend on this fixed cutoff.

\begin{lemma}[Linear anisotropic estimates]\label{lem:linear-semigroup}
For nonnegative integers \(k,m,n\) with \(k\leq2\) and \(m+n\leq2\), and physical
triples \(V=(v_1,v_2,v_3)\) and \(G=(G_1,G_2,G_3)\),
\begin{align}
 \|\nabla^kS_c(t)V\|_2
 &\leq C(1+t)^{-3/4-k/2}\|V\|_1
 +Ce^{-c_*t}\|\nabla^kV\|_2,
 \label{eq:Scestimate}\\
 \|\nabla_h^m\pa_3^nS_h(t)\mathbb Pf\|_2
 &\leq C(1+t)^{-1/2-m/2}
       \|\pa_3^nf\|_{L^2_3L^1_h}
 +Ce^{-c_*t}\|\nabla_h^m\pa_3^nf\|_2.
 \label{eq:Shestimate}
\end{align}
For every integer $m\ge0$ and $f\in L^2\cap L^2_3L^1_h$, the fixed
low-horizontal-frequency cutoff also satisfies
\begin{equation*}\tag{3.13a}\label{eq:Sh-low-all-orders}
 \|\nabla_h^mP_L^hS_h(t)\mathbb Pf\|_2
       \le C_m(1+t)^{-(m+1)/2}\|f\|_{L^2_3L^1_h}.
\end{equation*}
No derivative of $f$ is required in this localized estimate.
The time-integrated forms of \eqref{eq:Scestimate} and
\eqref{eq:Shestimate} are
\begin{align}
 \left\|\nabla^k\int_0^tS_c(t-s)G(s)\,\dd s\right\|_2
 &\leq C\int_0^t(1+t-s)^{-3/4-k/2}\|G(s)\|_1\,\dd s
 \notag\\
 &\quad+C\int_0^te^{-c_*(t-s)}\|\nabla^kG(s)\|_2\,\dd s,
 \label{eq:ScDuhamel}\\
 \left\|\nabla_h^m\pa_3^n\int_0^tS_h(t-s)\mathbb PG_2(s)\,\dd s\right\|_2
 &\leq C\int_0^t(1+t-s)^{-1/2-m/2}
 \|\pa_3^nG_2(s)\|_{L^2_3L^1_h}\,\dd s
 \notag\\
 &\quad+C\int_0^te^{-c_*(t-s)}
 \|\nabla_h^m\pa_3^nG_2(s)\|_2\,\dd s.
 \label{eq:ShDuhamel}
\end{align}
Moreover, for $0\le a<t$, $m=0,1$, and every scalar- or vector-valued
source $G\in L^2(a,t;H^m)$, one has
\begin{equation}\label{eq:r8-horizontal-maximal}
 \left\|\nabla_h\int_a^t S_h(t-s)G(s)\dd s\right\|_{H^m}
       \le (2\nu_h)^{-1/2}\|G\|_{L^2(a,t;H^m)}.
\end{equation}

For a vector-valued source $G$ and an integer $k\ge1$, we also have
\begin{equation}\label{eq:new-maximal-horizontal}
\left\|\nabla_h^kP_H^h\int_0^tS_h(t-s)\mathbb PG(s)\dd s\right\|_2^2
 \le C\int_0^t e^{-c_*(t-s)}
                    \|\nabla_h^{k-1}G(s)\|_2^2\dd s,
\end{equation}
provided $\nabla_h^{k-1}G\in L^2(0,t;L^2)$.
\end{lemma}

\begin{proof}
\smallskip\noindent\emph{Step 1. Compressible and solenoidal semigroup bounds.}
By \eqref{eq:physical-to-block},
\begin{equation}\label{eq:riesz-fourier-pointwise}
 \left||\xi|^{-1}i\xi\cdot\widehat v_2(\xi)\right|
 \leq |\widehat v_2(\xi)|,
\end{equation}
so no \(L^1\)-boundedness of a Riesz transform is being used.  By
\eqref{eq:blockmultiplier}, on \(\{\varrho\leq1\}\) the multiplier of
\(\nabla^kS_c(t)\) is bounded by \(C\varrho^ke^{-c_*\varrho^2t}\), while on
\(\{\varrho\geq1\}\) it is bounded by \(Ce^{-c_*t}\varrho^k\).  Hence Plancherel's
theorem, \eqref{eq:riesz-fourier-pointwise}, and
\(\|\widehat v_j\|_\infty\leq\|v_j\|_1\) give
\begin{align*}
 \|\nabla^kS_c(t)V\|_2^2
 &\leq C\|V\|_1^2\int_{\varrho\leq1}\varrho^{2k}e^{-2c_*\varrho^2t}\,\dd\xi
 +Ce^{-2c_*t}\|\nabla^kV\|_2^2\\
 &\leq C(1+t)^{-3/2-k}\|V\|_1^2
 +Ce^{-2c_*t}\|\nabla^kV\|_2^2,
\end{align*}
which proves \eqref{eq:Scestimate}.  Replacing \(v_2\) by the physical
source \(G_2(s)\) gives
\begin{equation}
 \left||\xi|^{-1}i\xi\cdot\widehat G_2(s,\xi)\right|
 \leq |\widehat G_2(s,\xi)|
 \leq\|G_2(s)\|_1,
\end{equation}
which is the low-frequency estimate used in \eqref{eq:ScDuhamel}.

For \eqref{eq:Shestimate}, let \(\mathcal F_h\) be the horizontal
Fourier transform. For every \(\xi_h\), \(\mathbb P(\xi_h,D_3)\) is
an \(L^2_{x_3}\)-multiplier with
operator norm at most one.  Moreover,
\begin{equation*}
 \sup_{\xi_h}\|\mathcal F_hf(\xi_h,\cdot)\|_{L^2_{x_3}}
 \leq\left\|\int_{\R^2}|f(x_h,\cdot)|\,\dd x_h\right\|_{L^2_{x_3}}
 =\|f\|_{L^2_3L^1_h}.
\end{equation*}
Therefore
\begin{align*}
 \|\nabla_h^mS_h(t)\mathbb Pf\|_2^2
 &\leq C\int_{|\xi_h|\leq1}|\xi_h|^{2m}
 e^{-2\nu_h|\xi_h|^2t}\,\dd\xi_h\,
 \|f\|_{L^2_3L^1_h}^2\\
 &\quad+Ce^{-2c_*t}\|\nabla_h^mf\|_2^2
 \leq C(1+t)^{-1-m}\|f\|_{L^2_3L^1_h}^2
 +Ce^{-2c_*t}\|\nabla_h^mf\|_2^2.
\end{align*}
Applying this horizontal Fourier estimate to \(\pa_3^nf\) proves
\eqref{eq:Shestimate}, since \(\pa_3\) commutes with
\(S_h(t)\mathbb P\).
The same contraction and mixed-norm bound prove
\eqref{eq:Sh-low-all-orders} on the fixed support $|\xi_h|\le R$, since
$\int_{|\xi_h|\le R}|\xi_h|^{2m}e^{-2\nu_ht|\xi_h|^2}\dd\xi_h
\le C_m(1+t)^{-m-1}$ for every integer $m\ge0$.

\smallskip\noindent\emph{Step 2. Duhamel bounds and horizontal smoothing.}
Minkowski's inequality in time, followed by
\eqref{eq:Scestimate} and \eqref{eq:Shestimate}, proves
\eqref{eq:ScDuhamel} and \eqref{eq:ShDuhamel}, respectively.

To prove \eqref{eq:r8-horizontal-maximal}, use Fourier Cauchy--Schwarz:
\[
 \left|\int_a^t |\xi_h|e^{-\nu_h(t-s)|\xi_h|^2}
                    \widehat G(s,\xi)\dd s\right|^2
 \le \frac{1}{2\nu_h}\int_a^t|\widehat G(s,\xi)|^2\dd s.
\]
Multiplication by $(1+|\xi|^2)^m$ and integration in $\xi$
prove \eqref{eq:r8-horizontal-maximal}.

It remains to prove the high-horizontal-frequency estimate
\eqref{eq:new-maximal-horizontal}. Fourier Cauchy--Schwarz and
$\|\mathbb P(\xi)\|\le1$ give
\[
\left|\int_0^t|\xi_h|^k e^{-\nu_h|\xi_h|^2(t-s)}
                      \mathbb P(\xi)\widehat G(s,\xi)\dd s\right|^2
 \le\nu_h^{-1}\int_0^t|\xi_h|^{2k-2}
          e^{-\nu_h|\xi_h|^2(t-s)}|\widehat G(s,\xi)|^2\dd s.
\]
Here we used
$\int_0^t|\xi_h|^2e^{-\nu_h|\xi_h|^2(t-s)}\dd s\le\nu_h^{-1}$.
On the support of $P_H^h$, the exponential is bounded by
$e^{-c_*(t-s)}$ for a fixed $c_*>0$.
Tonelli's theorem and Plancherel now prove
\eqref{eq:new-maximal-horizontal}.
\end{proof}

\subsubsection{Third-component and mixed-norm heat estimates}

\begin{lemma}[Third-component and mixed-norm heat estimates]
\label{lem:r13-third-row-heat}
For $a>0$ and an $\R^3$-valued source $G$, the following bounds
hold whenever their right-hand sides are finite:
\begin{align}\label{eq:new-third-kernels-all}
 \|[S_h(a)\mathbb PG]_3\|_\infty&\le Ca^{-3/2}\|G\|_1,
 &\|[S_h(a)\mathbb PG]_3\|_\infty&\le Ca^{-3/4}\|G\|_2,\notag\\
 \|[S_h(a)\mathbb PG]_3\|_{L^\infty_3L^2_h}&\le Ca^{-1}\|G\|_1,
 &\|[S_h(a)\mathbb PG]_3\|_{L^\infty_3L^2_h}&\le Ca^{-1/4}\|G\|_2.
\end{align}
Each additional horizontal derivative increases the time decay
exponent by $1/2$. Moreover,
\begin{align}\label{eq:new-mixed-kernel-all}
 \|\nabla_hS_h(a)\mathbb PG\|_\infty
 &\le Ca^{-3/2}\{\|G\|_{L^2_3L^1_h}\|\partial_3G\|_{L^2_3L^1_h}\}^{1/2},\notag\\
 \|\nabla_hS_h(a)\mathbb PG\|_\infty
 &\le C\|G\|_{H^1}\begin{cases}a^{-3/4},&0<a\le1,\\a^{-1},&a\ge1,\end{cases}\notag\\
 \|\nabla_hS_h(a)\mathbb PG\|_{L^\infty_3L^2_h}
 &\le Ca^{-1}\{\|G\|_{L^2_3L^1_h}\|\partial_3G\|_{L^2_3L^1_h}\}^{1/2},\notag\\
 \|\nabla_hS_h(a)\mathbb PG\|_{L^\infty_3L^2_h}&\le Ca^{-1/2}\|G\|_{H^1}.
\end{align}
\end{lemma}
\begin{proof}
Fix $a>0$ and a horizontal frequency $\xi_h$, and write
$r=|\xi_h|>0$. The inverse Fourier transforms in $\xi_3$ of
the third row of $\mathbb P(\xi)$ are
\[
 K_{3j}(\xi_h,z)=-\frac{i\xi_j}{2}\operatorname{sgn}(z)e^{-r|z|}
 \quad(j=1,2),\qquad K_{33}(\xi_h,z)=\frac r2e^{-r|z|}.
\]
Here $z=x_3-y_3$ is the vertical convolution variable.
The horizontal heat factor is $e^{-\nu_har^2}$.
The bounds $\|K_{3j}(\xi_h,\cdot)\|_\infty\le Cr$ and
$\|K_{3j}(\xi_h,\cdot)\|_2\le Cr^{1/2}$ imply
\eqref{eq:new-third-kernels-all}.

For the $L^1$ sources, convolution and Minkowski's inequality
give the first and third Gaussian factors below. For the $L^2$
sources, Cauchy--Schwarz gives the second and fourth factors.
Plancherel's theorem in the horizontal variables gives both mixed-norm estimates.
In the order of the four bounds in \eqref{eq:new-third-kernels-all},
these factors are
\[
\begin{aligned}
 &\int_{\R^2}r e^{-\nu_har^2}\,\dd\xi_h=Ca^{-3/2},
 &&\left(\int_{\R^2}r e^{-2\nu_har^2}\,\dd\xi_h\right)^{1/2}
       =Ca^{-3/4},\\
 &\left(\int_{\R^2}r^2e^{-2\nu_har^2}\,\dd\xi_h\right)^{1/2}
       =Ca^{-1},
 &&\sup_{r>0}r^{1/2}e^{-\nu_har^2}\le Ca^{-1/4}.
\end{aligned}
\]
For the last factor, first apply Cauchy--Schwarz to the vertical
convolution at each $\xi_h$, using $\|K_{3j}\|_2\le Cr^{1/2}$,
and then apply horizontal Plancherel.
An additional horizontal derivative inserts one factor $r$,
so each time exponent in \eqref{eq:new-third-kernels-all}
increases by $1/2$.

To prove \eqref{eq:new-mixed-kernel-all}, we use vertical Sobolev
embedding and horizontal Fourier estimates.

For each fixed $\xi_h$, horizontal Fourier transformation gives
$\|\mathcal F_hG(\xi_h,\cdot)\|_{L^2_3}
\le\|G\|_{L^2_3L^1_h}$ and the same bound for
$\partial_3G$.
The vertical multiplier
$\mathcal P_{\xi_h}=\mathcal F_3^{-1}\mathbb P(\xi_h,\xi_3)\mathcal F_3$
is an $L^2_3$ contraction and commutes with $\partial_3$.
The one-dimensional Sobolev inequality therefore gives
\[
 \sup_z|\mathcal P_{\xi_h}\mathcal F_hG(\xi_h,z)|
 \le C\|\mathcal P_{\xi_h}\mathcal F_hG\|_{L^2_3}^{1/2}
       \|\partial_3\mathcal P_{\xi_h}\mathcal F_hG\|_{L^2_3}^{1/2}
 \le C\bigl(\|G\|_{L^2_3L^1_h}\|\partial_3G\|_{L^2_3L^1_h}\bigr)^{1/2}.
\]
Horizontal inversion with
$\int_{\R^2}r e^{-\nu_har^2}\dd\xi_h=Ca^{-3/2}$ gives the
first inequality in \eqref{eq:new-mixed-kernel-all}.
Horizontal Plancherel with
$\int_{\R^2}r^2e^{-2\nu_har^2}\dd\xi_h=Ca^{-2}$ gives the third.
For the second bound, Fourier Cauchy--Schwarz gives
\[
 \|\nabla_hS_h(a)\mathbb PG\|_\infty^2
 \le C\|G\|_{H^1}^2
       \int_0^\infty r^3e^{-2\nu_har^2}(1+r^2)^{-1/2}\dd r.
\]
The integral is bounded by $Ca^{-3/2}$ for $0<a\le1$ and by
$Ca^{-2}$ for $a\ge1$. This proves the second inequality in
\eqref{eq:new-mixed-kernel-all}.
For the fourth inequality in \eqref{eq:new-mixed-kernel-all}, \eqref{eq:mixedv} and Minkowski's inequality imply
$\|F\|_{L^\infty_3L^2_h}\le C\|F\|_2^{1/2}\|\partial_3F\|_2^{1/2}$.
Since $\sup_{r>0}r e^{-\nu_har^2}\le Ca^{-1/2}$, Plancherel gives
\[
\begin{aligned}
 \|\nabla_hS_h(a)\mathbb PG\|_{L^\infty_3L^2_h}
 &\le C\|\nabla_hS_h(a)\mathbb PG\|_2^{1/2}\|\partial_3\nabla_hS_h(a)\mathbb PG\|_2^{1/2}\\
 &\le Ca^{-1/2}\bigl(\|G\|_2\|\partial_3G\|_2\bigr)^{1/2}
       \le Ca^{-1/2}\|G\|_{H^1}.
\end{aligned}
\]
This proves the fourth bound in \eqref{eq:new-mixed-kernel-all}.
The projection estimates used here are $L^2$ estimates;
all mixed source norms are taken on $G$.
The bounds extend by approximation whenever their right-hand
sides are finite.

\end{proof}

\subsubsection{Estimates for the Duhamel integrals}
Choose \(\chi_\ell\in C_c^\infty(\R^3)\) such that
\(\chi_\ell=1\) on \(\{|\xi|\leq R_0\}\) and
\(\chi_\ell=0\) on \(\{|\xi|\geq2R_0\}\), where \(R_0\geq2\) is fixed by the spectral decomposition in
Appendix~\ref{app:r13-green-matrix}, and put
\begin{equation*}
 P_L=\chi_\ell(D),\qquad P_H=I-P_L.
\end{equation*}
Here $P_L$ and $P_H$ denote the low- and high-frequency cutoffs
with respect to $|\xi|$.

\begin{lemma}[Estimates for the Duhamel integrals]
\label{lem:refined-compressible-green}
Let \(F=(F_1,F_2,F_3)\), and let \((q_F,\mathbb Q u_F,\tau_F)\)
be the zero-initial-data solution of the compressible linear system with
right-hand side \(F\).  Then
\begin{align}
 &\|\nabla^2(q_F,\tau_F)(t)\|_2
   +\|\nabla\Div\mathbb Q u_F(t)\|_2
 \notag
 \leq C\int_0^t(1+t-r)^{-7/4}\|F(r)\|_1\,\dd r
 \notag\\
 &\qquad+C\int_0^t e^{-c_*(t-r)}\{1+(t-r)^{-1/2}\}
 \big\{\|(F_1,F_3)(r)\|_{H^2}+\|F_2(r)\|_{H^1}\big\}\,\dd r.
 \label{eq:refined-second-Duhamel}
\end{align}
Moreover, for \(j=0,1\),
\begin{align}
 \|\nabla^j(q_F,\mathbb Q u_F,\tau_F)(t)\|_2
 &\leq C\int_0^t(1+t-r)^{-3/4-j/2}\|F(r)\|_1\,\dd r
 \notag\\
 &\quad+C\int_0^te^{-c_*(t-r)}\|F(r)\|_{H^j}\,\dd r.
 \label{eq:refined-low-order-Duhamel}
\end{align}
\end{lemma}

\begin{proof}
At low frequencies, \eqref{eq:blockmultiplier} gives the Gaussian
bounds in \eqref{eq:refined-second-Duhamel} and
\eqref{eq:refined-low-order-Duhamel}.
At high frequencies, the slow density mode and the two fast modes
have the componentwise Green-matrix bounds
\eqref{eq:component-Green-matrix}, proved in
Appendix~\ref{app:r13-green-matrix}.
The momentum column gains a factor $\varrho^{-1}$ in the scalar
rows; the fast velocity row is controlled by the parabolic estimate
\eqref{eq:parabolic-half-gain}. Thus the second-order bound requires
$F_2$ only in $H^1$, while the scalar sources are estimated in $H^2$.
The reconstruction in that appendix proves
\eqref{eq:refined-second-Duhamel};
\eqref{eq:refined-low-order-Duhamel} follows from
\eqref{eq:ScDuhamel} with $k=0,1$.
\end{proof}

\begin{remark*}[Time convolution and unit-interval estimates]
Let $h\ge0$ be locally integrable, and fix $c_*>0$ and
$p,\beta\ge0$. If a constant $A\ge0$, independent of
$a$ and $t$, satisfies
\[
                  \int_a^{a+1}h(s)\dd s
                         \le A(1+a)^{-\beta}\qquad(a\ge0),
\]
then
\begin{equation}\label{eq:new-window-convolution}
 \int_0^te^{-c_*(t-s)}(1+s)^{-p}h(s)\dd s
       \le C_{c_*,p,\beta}A(1+t)^{-p-\beta},\qquad t\ge0.
\end{equation}
In particular, if $\int_0^\infty h(s)\dd s\le A$, then
\[
 \int_0^te^{-c_*(t-s)}(1+s)^{-p}h(s)\dd s
                       \le C_{c_*,p}A(1+t)^{-p},\qquad t\ge0.
\]
For a nonnegative function $h$ satisfying
$\int_a^{a+1}h(s)^2\dd s\le A^2(1+a)^{-\beta}$ for all $a\ge0$,
with $\beta>2$, one also has
\begin{equation}\label{eq:r13-unit-interval-tail}
 \int_t^\infty h(s)\dd s
       \le C_\beta A(1+t)^{1-\beta/2},\qquad t\ge0.
\end{equation}

For $t\ge2$, split the integral at $t/2$. Summation over unit
intervals bounds the contribution from $[0,t/2]$ by
$CA(1+t)e^{-c_*t/2}$. On $[t/2,t]$, each interval of length
at most one has integral
at most $CA(1+t)^{-\beta}$. The weight is at most
$C(1+t)^{-p}$, and the exponential factors sum to a geometric
series. For $t\le2$, at most two unit intervals suffice.
This proves \eqref{eq:new-window-convolution}. The total-integral
hypothesis implies $\int_a^{a+1}h(s)\dd s\le A$, so the second
estimate follows by taking $\beta=0$.

For \eqref{eq:r13-unit-interval-tail}, Cauchy--Schwarz on each
unit interval gives
\[
 \int_t^\infty h(s)\dd s
 \le A\sum_{k=0}^\infty(1+t+k)^{-\beta/2}
 \le C_\beta A(1+t)^{1-\beta/2},
\]
since $\beta/2>1$.
\end{remark*}

\subsection{Decay estimates up to second order}\label{sec:H2-decay}

We prove \eqref{eq:new-H2-main-rates}--\eqref{eq:H2-main-mixed-and-full}
by successive improvements of the decay at orders zero, one, and two.
Throughout this subsection, $v=\mathbb Pu$, $w=\mathbb Qu$, and
$W=(q,w,\tau)$. The global estimate in Theorem~\ref{thm:global}
provides the energy and dissipation bounds
\begin{equation}\label{eq:new-base-energy}
 \sup_{t\ge0}\|U(t)\|_{H^2}\le C\mathcal J_2,
 \qquad \int_0^\infty\cD_2(U(t))\,\dd t\le C\mathcal J_2^2.
\end{equation}
The third derivatives in $\cD_2(U)$ occur only in the dissipation:
they do not include $\partial_3^3u_h$.

Let $P_L,P_H$ be the fixed spatial-frequency cutoffs introduced before
Lemma~\ref{lem:refined-compressible-green}, with radius $R_0$
chosen sufficiently large.
The horizontal cutoffs $P_L^h,P_H^h$ are those used in
Lemma~\ref{lem:linear-semigroup}. We write
$f_L=P_Lf$ and $f_H=P_Hf$.

\subsubsection{Preliminary decay estimates}

The low-order energy inequality below supplies the dissipation tail
needed for the preliminary decay bootstrap. Its proof uses only
\eqref{eq:new-base-energy}.

\paragraph{A low-order energy and its dissipation tail.}
Define the density dissipation and the low-order dissipation by
\begin{equation}\label{eq:new-d0}
 d_q^2=\ip q{(I-(1-\Delta)^{-1})q},\qquad
 \cD_0(U)=d_q^2+\|(\nabla_hu,\Div u)\|_2^2+\|\nabla\tau\|_2^2.
\end{equation}
The Fourier weight in $d_q^2$ is $|\xi|^2/(1+|\xi|^2)$.

\begin{lemma}[Low-order energy and dissipation]\label{lem:new-low-energy}
There is an energy $E_0\simeq\|U\|_2^2$ such that
\begin{equation}\label{eq:new-low-energy}
 E_0'+c_*\cD_0(U)\le C \cD_2(U)E_0.
\end{equation}
If $\|U(s)\|_2\le M(1+s)^{-1/2}$ on $[0,t]$, then
\begin{equation}\label{eq:new-low-tail}
\begin{aligned}
 \int_a^t\cD_0(U(s))\dd s&\le CM^2(1+a)^{-1},
                                      &&0\le a\le t,\\
 \int_0^t(1+s)^r\cD_0(U(s))\dd s&\le C_rM^2,
                                      &&0\le r<1,
\end{aligned}
\end{equation}
provided $M\ge\mathcal J_2$.
\end{lemma}
\begin{proof}
Spatial regularization justifies the energy calculations by
Theorem~\ref{thm:global} and Lemma~\ref{lem:H2-anisotropic}.
For sufficiently small fixed $\eta>0$, take
\[
 E_0=\cE_{\rm th}+\eta\ip u{\nabla(1-\Delta)^{-1}q},
\]
where $\cE_{\rm th}$ is the energy in
Lemma~\ref{lem:zeroth-order-energy}. Direct differentiation of the cross term gives
\begin{equation}\label{eq:new-cross0}
\begin{aligned}
 &\frac{\dd}{\dd t}\ip u{\nabla(1-\Delta)^{-1}q}+c^2d_q^2
 =-c^2\ip{\nabla\tau}{\nabla(1-\Delta)^{-1}q}
         +\ip{\mathcal Lu}{\nabla(1-\Delta)^{-1}q}\\
 &\qquad+\ip{\Div u}{(1-\Delta)^{-1}\Div u}
         +\ip{F_2}{\nabla(1-\Delta)^{-1}q}
         -\ip{(1-\Delta)^{-1}\Div u}{F_1}.
\end{aligned}
\end{equation}
The linear terms are bounded by $c^2d_q^2/4+C(\|(\nabla_hu,\Div u)\|_2^2+
\|\nabla\tau\|_2^2)$. We now estimate the nonlinear terms $F_2$ and
$F_1$ in \eqref{eq:new-cross0}, using their definitions
\eqref{eq:F2} and \eqref{eq:F1}. The velocity convection satisfies
\[
 \|u\cdot\nabla u\|_2
 \le\|u_h\|_\infty\|\nabla_hu\|_2
       +\|u_3\|_6\|\partial_3u\|_3\le C\|U\|_{H^2}\cD_0(U)^{1/2},
\]
where the bound
$\|\nabla u_3\|_2\le C\|(\nabla_hu,\Div u)\|_2$ is
\eqref{eq:verticalcomp} with $m=0$.
Integrating the variable-viscosity term in \eqref{eq:F2} by parts
once gives $C\|U\|_{H^2}\cD_0(U)$;
its coefficient remainder uses
$\|\nabla (e^{-q}-1)\|_3\|(\nabla_hu,\Div u)\|_2\|\nabla (1-\Delta)^{-1}q\|_6\le C\|U\|_{H^2}\cD_0(U)$.

For the pressure contribution containing $(e^\tau-1)\nabla q$, use
\[
 |\ip{(e^\tau-1)q}{(I-(1-\Delta)^{-1})q}|
 \le C\|(e^\tau-1)\|_\infty d_q^2+
 C\|\nabla (e^\tau-1)\|_\infty\|q\|_2d_q.
\]
Indeed, the kernel $K(x)=e^{-|x|}/(4\pi|x|)$ of $(1-\Delta)^{-1}$ is
nonnegative, even, has mass one and a finite second moment, and
\[
\begin{aligned}
 &\ip{(e^\tau-1)q}{(I-(1-\Delta)^{-1})q}
 =\frac12\iint K(x-y)
       \bigl[(e^{\tau(x)}-1)q(x)-(e^{\tau(y)}-1)q(y)\bigr]\\
 &\hspace{45mm}\times(q(x)-q(y))\dd x\dd y.
\end{aligned}
\]
Splitting the first difference and applying Cauchy--Schwarz gives
the bound for the pairing of $(e^\tau-1)q$ with
$(I-(1-\Delta)^{-1})q$. The contribution containing
$(e^\tau-1)\nabla\tau$ in \eqref{eq:F2} is treated by the same
integration by parts.
Young's inequality and $\|\nabla (e^\tau-1)\|_\infty^2\le C\cD_2(U)$
then give a small part of $d_q^2$ plus $C\cD_2(U)\|q\|_2^2$.
Finally write $F_1=-\Div(uq)+q\Div u$. The flux term is at most
\[
 C\|q\|_\infty\|u\|_2\cD_0(U)^{1/2}
 \le\eta \cD_0(U)+C_\eta\|q\|_\infty^2\|u\|_2^2,
 \qquad
 \|q\|_\infty^2\le C\|\nabla q\|_2\|\nabla^2q\|_2\le C\cD_2(U).
\]
The term $q\Div u$ in $F_1$ contributes at most
$C\|U\|_{H^2}\cD_0(U)$. Choose the compensation parameter and
$\|U\|_{H^2}$ small, and combine \eqref{eq:new-cross0} with
\eqref{eq:freeenergyidentity}. This proves \eqref{eq:new-low-energy};
\eqref{eq:freeequiv} gives $E_0\simeq\|U\|_2^2$ for small $\eta$.

For the dissipation tail, the assumed decay gives
$E_0(s)\le CM^2(1+s)^{-1}$. Integrating
\eqref{eq:new-low-energy} on $[a,t]$ and using
\eqref{eq:new-base-energy} yields
\[
 \int_a^t\cD_0(U(s))\dd s
 \le CE_0(a)+C\int_a^t\cD_2(U(s))E_0(s)\dd s
 \le CM^2(1+a)^{-1}.
\]
To prove the weighted bound, multiply \eqref{eq:new-low-energy}
by $(1+s)^r$ and integrate on $[0,t]$:
\[
\begin{aligned}
 \int_0^t(1+s)^r\cD_0(U(s))\dd s
 &\le CE_0(0)+C_r\int_0^t(1+s)^{r-1}E_0(s)\dd s\\
 &\quad+C\int_0^t(1+s)^r\cD_2(U(s))E_0(s)\dd s
 \le C_rM^2.
\end{aligned}
\]
Here $r<1$ ensures that $\int_0^\infty(1+s)^{r-2}\dd s$ is finite
and $(1+s)^{r-1}\le1$; the last integral also uses
\eqref{eq:new-base-energy}. This proves both bounds in
\eqref{eq:new-low-tail}.
\end{proof}

\paragraph{Preliminary decay.}
The dissipation bounds \eqref{eq:new-low-tail} and
Lemma~\ref{lem:linear-semigroup} yield the following preliminary rates.
\begin{lemma}[Preliminary decay]\label{lem:new-seed}
For sufficiently small $\mathcal J_2$,
\begin{equation}\label{eq:new-seed}
 \|u(t)\|_2\le C\mathcal J_2 (1+t)^{-1/2},\qquad
 \|(W,u_3)(t)\|_2\le C\mathcal J_2 (1+t)^{-5/8}.
\end{equation}
Moreover $\int_t^\infty \cD_0(U)\le C\mathcal J_2^2(1+t)^{-1}$ and
$\int_0^\infty (1+s)^r \cD_0(U)\dd s\le C_r\mathcal J_2^2$ for $r<1$.
\end{lemma}
\begin{proof}
We estimate $W$, $v$, and $v_3$ in that order and recover $u=v+w$
and $u_3=v_3+w_3$.
Set $\beta=5/8$, $r=\vartheta=31/32$, and on a finite interval let
\[
 M_0(t)=\mathcal J_2+\sup_{s\le t}
       \bigl((1+s)^{1/2}\|u(s)\|_2+(1+s)^\beta\|(W,u_3)(s)\|_2\bigr).
\]
In the interpolation and source bounds below, norms without a time
argument are evaluated at $s\in[0,t]$. The quantity
$M_0(t)\ge\mathcal J_2$ is fixed.
Assume temporarily $M_0(t)\le1$. Interpolation between the $L^2$
bounds in the definition of $M_0(t)$ and \eqref{eq:new-base-energy}
yields
\begin{equation}\label{eq:new-seed-ledger}
\begin{aligned}
 \|u\|_\infty+\|\nabla u\|_3&\le CM_0(t)(1+s)^{-1/8},\\
 \|\partial_3u\|_2&\le CM_0(t)(1+s)^{-1/4},\\
 \|u\|_{L^\infty_3L^2_h}&\le CM_0(t)(1+s)^{-3/8},\\
 \|u_3\|_{L^\infty_3L^2_h}&\le CM_0(t)^{1/2}(1+s)^{-\beta/2}\cD_0(U)^{1/4},\\
 \|g\|_\infty+\|\nabla g\|_3+\|\nabla g\|_{L^\infty_3L^2_h}
       &\le CM_0(t)(1+s)^{-\beta/4},\\
 \|g\|_{L^\infty_3L^2_h}&\le CM_0(t)(1+s)^{-3\beta/4},
\end{aligned}
\end{equation}
where $g=q,\tau,(e^{-q}-1),(e^\tau-1)$. The mixed-norm bounds use
$\|g\|_{L^\infty_3L^2_h}\le C\|g\|_2^{1/2}\|\partial_3g\|_2^{1/2}$,
which follows from \eqref{eq:mixedv} and Minkowski's inequality;
for $g=u_3$, use \eqref{eq:verticalcomp} and \eqref{eq:new-d0}.
For each fixed $t$, apply Lemma~\ref{lem:new-low-energy} with
$M=M_0(t)$ to obtain the weighted and tail integral estimates
\eqref{eq:new-low-tail} on $[0,t]$. In particular,
\[
 \int_0^t(1+s)^r \cD_0(U(s))\dd s\le CM_0(t)^2,
 \qquad \int_0^t\cD_2(U(s))\dd s\le C\mathcal J_2^2.
\]
When a time-integrated estimate is used, extend
$\cD_0(U)^{1/2}$ by zero
outside $[0,t]$; the finite-interval tail bound then gives
the required integral bound without information beyond $t$.
For $t\le2$, \eqref{eq:new-seed} follows from
\eqref{eq:new-base-energy}. We therefore take $t\ge2$ in the decay
estimates that follow.

\paragraph{Compressible low frequencies.} We estimate $P_LW$ using
\eqref{eq:ScDuhamel}. Decompose the sources
\eqref{eq:F1}--\eqref{eq:F3} as
$F=F^{(0)}+\sum_j\partial_jF^{(j)}$, with
\begin{equation}\label{eq:new-basic-flux}
 F^{(j)}=\left(-u_jq,
 -u_ju+\nu_h\mathbf1_{j\le2}(e^{-q}-1)\partial_ju
       +\nu_c(e^{-q}-1)(\Div u)e_j,
 -u_j\tau+\chi (e^{-q}-1)\partial_j\tau\right).
\end{equation}
The remainder $F^{(0)}$ contains coefficient-gradient products,
transport and pressure remainders, and the quadratic source terms
in \eqref{eq:F3}. H\"older's inequality and
\eqref{eq:new-seed-ledger} give
\begin{align}\label{eq:new-seed-low-sources}
 \|F^{(0)}\|_1&\le CM_0(t) (1+s)^{-\beta/2}\cD_0(U)^{1/2}
                  +CM_0(t) (1+s)^{-\beta}\cD_2(U)^{1/2}+C\cD_0(U),\notag\\
 \sum_j\|F^{(j)}\|_1&\le CM_0(t)^2(1+s)^{-1}+CM_0(t) (1+s)^{-\beta}\cD_0(U)^{1/2}.
\end{align}
The pressure term uses $\|\tau\|_2\|\nabla(q,\tau)\|_2
\le CM_0(t) (1+s)^{-\beta}\cD_2(U)^{1/2}$.
In \eqref{eq:ScDuhamel}, the low-frequency kernels for $F^{(0)}$ and
$\partial_jF^{(j)}$ decay as $(1+t-s)^{-3/4}$ and
$(1+t-s)^{-5/4}$, respectively.
We estimate separately the terms containing $\cD_0(U)^{1/2}$ and
$\cD_2(U)^{1/2}$ in \eqref{eq:new-seed-low-sources}. Weighted Cauchy--Schwarz and
\eqref{eq:new-low-tail} give
\[
\begin{aligned}
 &M_0(t)\int_0^t(1+t-s)^{-3/4}(1+s)^{-\beta/2}\cD_0(U(s))^{1/2}\dd s\\
 &\quad\le CM_0(t)^2
       \left(\int_0^t(1+t-s)^{-3/2}
                          (1+s)^{-\beta-r}\dd s\right)^{1/2}
 \le CM_0(t)^2(1+t)^{-3/4},
\end{aligned}
\]
because $\beta+r=51/32>3/2$. For the term containing
$\cD_2(U)^{1/2}$ in \eqref{eq:new-seed-low-sources}, use the total
dissipation bound \eqref{eq:new-base-energy}:
\[
\begin{aligned}
 &M_0(t)\int_0^t(1+t-s)^{-3/4}(1+s)^{-\beta}
                                      \cD_2(U(s))^{1/2}\dd s\\
 &\quad\le C\mathcal J_2 M_0(t)
       \left(\int_0^t(1+t-s)^{-3/2}
                              (1+s)^{-2\beta}\dd s\right)^{1/2}\\
 &\quad\le C\mathcal J_2 M_0(t)(1+t)^{-\beta}
       \le CM_0(t)^2(1+t)^{-\beta}.
\end{aligned}
\]
Here $2\beta=5/4$ and $M_0(t)\ge\mathcal J_2$.
Splitting at $t/2$ proves the displayed convolution bounds:
the two time weights are integrable, and the smaller power
determines the decay. The remaining source and flux terms obey
\[
\begin{aligned}
 \int_0^t(1+t-s)^{-3/4}\cD_0(U(s))\dd s
       &\le CM_0(t)^2(1+t)^{-3/4},\\
 \int_0^t(1+t-s)^{-5/4}
               \sum_j\|F^{(j)}(s)\|_1\dd s
       &\le CM_0(t)^2
            \bigl[(1+t)^{-1}+(1+t)^{-\beta-r/2}\bigr].
\end{aligned}
\]
For the $\cD_0(U)$ convolution, use $r>3/4$ in
\eqref{eq:new-low-tail}. For the flux convolution, integrate the
$M_0(t)^2(1+s)^{-1}$ term of \eqref{eq:new-seed-low-sources}
directly and apply weighted Cauchy--Schwarz to the term containing $\cD_0(U)^{1/2}$.
Combining these source convolutions with the initial-data bound
\eqref{eq:Scestimate} gives
$\|P_LW(t)\|_2\le C\mathcal J_2(1+t)^{-3/4}
                  +CM_0(t)^2(1+t)^{-\beta}$.

\paragraph{High frequencies and density transport.} We next estimate
$P_HW$. To avoid differentiating the density source, keep the
transport term from \eqref{eq:F1} on the left:
\begin{equation}\label{eq:new-density-high}
 (q_H)_t+\Div w_H+u\cdot\nabla q_H=C_q,
 \qquad C_q=P_L(u\cdot\nabla q)-u\cdot\nabla q_L.
\end{equation}
The fixed cutoff and \eqref{eq:new-seed-ledger} give
$\|C_q\|_2\le CM_0(t)^2(1+s)^{-\beta-1/8}+CM_0(t) (1+s)^{-\beta}\cD_0(U)^{1/2}$.
The unprojected sources satisfy
\begin{equation}\label{eq:new-seed-negative}
 \sum_{j=1}^3\|F_j\|_{H^{-1}}
       \le CM_0(t)^2(1+s)^{-5/8}+CM_0(t) (1+s)^{-1/8}\cD_0(U)^{1/2}.
\end{equation}
To prove \eqref{eq:new-seed-negative}, write the convection terms in
\eqref{eq:F1}--\eqref{eq:F3} as divergences plus multiples of $\Div u$.
The velocity flux obeys
$\|u\otimes u\|_2\le C\|u\|_2^{5/4}\|\nabla^2u\|_2^{3/4}$.
For pressure use $(e^\tau-1)\nabla q=\nabla((e^\tau-1)q)-q\nabla (e^\tau-1)$; for diffusion
integrate by parts once, obtaining
$C(\|(e^{-q}-1)\|_\infty+\|\nabla (e^{-q}-1)\|_3)\cD_0(U)^{1/2}$ in $H^{-1}$.
The quadratic sources use $\|(\nabla_hu,\Div u)\|_3\|(\nabla_hu,\Div u)\|_2$ and
$\|\nabla\tau\|_3\|\nabla\tau\|_2$. Multiplication by the
smooth coefficients is bounded on $H^{-1}$ by duality from
$H^2\cdot H^1\subset H^1$.

Apply the compensation of Lemma~\ref{lem:fourier-energy} to the
high-frequency system, with density transport given by
\eqref{eq:new-density-high}. The resulting energy
$E_H\simeq\|W_H\|_2^2$ satisfies
\[
 E_H'+c_*(E_H+\|\nabla(w_H,\tau_H)\|_2^2)
 \le C\cD_2(U)E_H+C\left(\left[\sum_{j=1}^3\|F_j\|_{H^{-1}}\right]^2+\|C_q\|_2^2\right).
\]
In the density test, the transport contributes
$\frac12\int\Div u\,q_H^2$; use
$\|\Div u\|_\infty\le\eta+C_\eta \cD_2(U)$ and absorb $\eta$
with the spectral gap. The other sources pair in $H^{-1}$--$H^1$.
The integrating factor is bounded by \eqref{eq:new-base-energy}.
Applying \eqref{eq:new-window-convolution} with
\eqref{eq:new-low-tail} gives
\[
 \int_0^te^{-c_*(t-s)}(1+s)^{-p}\cD_0(U(s))\dd s\le CM_0(t)^2(1+t)^{-p-1}.
\]
Substitute \eqref{eq:new-seed-negative} and the estimate for the
cutoff commutator $C_q$, defined in \eqref{eq:new-density-high},
into the energy inequality. Exponential convolution then gives
\[
 E_H(t)\le C\mathcal J_2^2e^{-c_*t}
                +CM_0(t)^4(1+t)^{-5/4}.
\]
Since $E_H\simeq\|P_HW\|_2^2$, this is the required
$(1+t)^{-\beta}$ bound for the high-frequency part of $W$.

\paragraph{Decay of the solenoidal velocity.}
To obtain the $(1+t)^{-1/2}$ decay of $u=v+w$, we estimate
$v=\mathbb Pu$ by Duhamel's formula:
\[
 v(t)=S_h(t)\mathbb Pu_0+
       \int_0^tS_h(t-s)\mathbb PF_2(s)\dd s.
\]
For $t\ge2$, we split the time integral at $t-1$, using
mixed-norm heat estimates away from $s=t$ and $L^2$ estimates
near it. By \eqref{eq:Shestimate}, the homogeneous term $S_h(t)\mathbb Pu_0$
is bounded by $C\mathcal J_2(1+t)^{-1/2}$, using the initial mixed
norm $\|u_0\|_{L^2_3L^1_h}$. For the Duhamel integral, the exact
source decomposition is
\begin{equation}\label{eq:new-projected-source}
 \mathbb PF_2=\mathbb PN+\mathbb P\Div_h\bigl(\nu_h(e^{-q}-1)\nabla_hu\bigr),
 \quad N=-u\cdot\nabla u-c^2\Pi_p+\Pi_v,
\end{equation}
where
\begin{equation}\label{eq:new-N-definitions}
\Pi_p=-(q+\tau)e^\tau\nabla\tau,
 \qquad \Pi_v=-\nu_h\nabla_h(e^{-q}-1)\cdot\nabla_hu
                 -\nu_c(\Div u)\nabla(e^{-q}-1).
\end{equation}
The gradient removed from pressure is
$\nabla[(e^\tau-1)(q+\tau)]$; the removed bulk-viscosity term is
$\nu_c\nabla((e^{-q}-1)\Div u)$. In particular no $L^1$ boundedness of
$\mathbb P$ is being used. The mixed bounds are
\[
\begin{aligned}
 \|N\|_{L^2_3L^1_h}
 &\le CM_0(t)(1+s)^{-\beta/4}\cD_0(U)^{1/2}
       +CM_0(t)^{3/2}(1+s)^{-\beta/2-1/4}\cD_0(U)^{1/4},\\
 \|\nu_h(e^{-q}-1)\nabla_hu\|_{L^2_3L^1_h}
 &\le CM_0(t)(1+s)^{-3\beta/4}\cD_0(U)^{1/2}.
\end{aligned}
\]
The bounds for $N$ and the horizontal viscosity flux follow from
\eqref{eq:new-seed-ledger}, \eqref{eq:new-N-definitions}, and
$\|ab\|_{L^2_3L^1_h}\le\|a\|_{L^\infty_3L^2_h}\|b\|_2$.
\par\smallskip\noindent\emph{The interval $[0,t-1]$.}\enspace
Here $t-s\ge1$. The horizontal Fourier calculation in the proof of
\eqref{eq:Shestimate} gives heat-kernel powers $1/2$ for $N$ and
$1$ for the horizontal flux in \eqref{eq:new-projected-source}.
For the term $M_0(t)(1+s)^{-\beta/4}\cD_0(U)^{1/2}$ in the bound for $N$, we obtain
\[
\begin{aligned}
 &M_0(t)\int_0^{t-1}(t-s)^{-1/2}(1+s)^{-\beta/4}\cD_0(U(s))^{1/2}\dd s\\
 &\quad\le CM_0(t)^2
       \left(\int_0^{t-1}(t-s)^{-1}
                              (1+s)^{-\beta/2-r}\dd s\right)^{1/2}\\
 &\quad\le CM_0(t)^2
       \bigl[(1+t)^{-1/2}
          +(1+t)^{-\beta/4-r/2}\sqrt{\log(e+t)}\bigr].
\end{aligned}
\]
The exterior factor uses the weighted dissipation bound
\eqref{eq:new-low-tail}. For the term
$M_0(t)^{3/2}(1+s)^{-\beta/2-1/4}\cD_0(U)^{1/4}$ in the bound for $N$,
H\"older with exponents $4$ and $4/3$ gives
\[
\begin{aligned}
 &M_0(t)^{3/2}\int_0^{t-1}(t-s)^{-1/2}
                (1+s)^{-\beta/2-1/4}\cD_0(U(s))^{1/4}\dd s\\
 &\quad\le CM_0(t)^2
       \left(\int_0^{t-1}(t-s)^{-2/3}
              (1+s)^{-\frac43(\beta/2+1/4+r/4)}\dd s\right)^{3/4}\\
 &\quad\le CM_0(t)^2
       \bigl[(1+t)^{-1/2}+(1+t)^{-\beta/2-r/4}\bigr].
\end{aligned}
\]
Here the exterior amplitude is
$M_0(t)^{3/2}(M_0(t)^2)^{1/4}=M_0(t)^2$.
The time weights inside the two integrals have respective
powers $41/32$ and $103/96$, both greater than one.
On $s\le t/2$ keep the heat-kernel decay; on
$s\ge t/2$ integrate $(t-s)^{-1}$ or $(t-s)^{-2/3}$.
Since
\[
 \frac{\beta}{4}+\frac r2=\frac{41}{64}>\frac12,
 \qquad
 \frac{\beta}{2}+\frac r4=\frac{71}{128}>\frac12,
\]
both complete contributions are bounded by
$CM_0(t)^2(1+t)^{-1/2}$.
For the horizontal viscosity flux in \eqref{eq:new-projected-source},
weighted Cauchy--Schwarz and \eqref{eq:new-low-tail} give
\[
\begin{aligned}
 &M_0(t)\int_0^{t-1}(t-s)^{-1}
                       (1+s)^{-3\beta/4}\cD_0(U(s))^{1/2}\dd s\\
 &\quad\le CM_0(t)^2
       \left(\int_0^{t-1}(t-s)^{-2}
                              (1+s)^{-3\beta/2-r}\dd s\right)^{1/2}
 \le CM_0(t)^2(1+t)^{-61/64}.
\end{aligned}
\]
The last estimate uses $3\beta/2+r=61/32$.

\par\smallskip\noindent\emph{The interval $[t-1,t]$.}\enspace
The mixed-norm kernels are singular at $s=t$. We instead use
$L^2$ contraction for convection and pressure, together with
$\|u\cdot\nabla u\|_2\le CM_0(t) (1+s)^{-1/8}\cD_0(U)^{1/2}$.
By \eqref{eq:new-low-tail}, their integrals on this interval are
bounded by $CM_0(t)^2(1+t)^{-\beta}$.
For $\Pi_v$ in \eqref{eq:new-N-definitions} and the horizontal
viscosity flux in \eqref{eq:new-projected-source}, combine the
mixed-norm heat estimate from the proof of
Lemma~\ref{lem:linear-semigroup} with $L^2$ contraction:
\[
\begin{aligned}
 \|S_h(a)\mathbb Pg\|_2
 &\le C\min\{a^{-1/2}\|g\|_{L^2_3L^1_h},\|g\|_2\},\\
 \|S_h(a)\mathbb P\Div_h[\nu_h(e^{-q}-1)\nabla_hu]\|_2
 &\le C\min\bigl\{a^{-1/2}\|\nu_h(e^{-q}-1)\nabla_hu\|_2,\\
 &\hspace{30mm}\|\nu_h(e^{-q}-1)\nabla_hu\|_{H^1}\bigr\}.
\end{aligned}
\]
Here $\|\Pi_v\|_{L^2_3L^1_h}+\|\bigl(\nu_h(e^{-q}-1)\nabla_hu\bigr)\|_2\le CM_0(t) (1+s)^{-\beta/4}\cD_0(U)^{1/2}$
and $\|\Pi_v\|_2+\|\bigl(\nu_h(e^{-q}-1)\nabla_hu\bigr)\|_{H^1}\le CM_0(t)^2(1+s)^{-\beta/4}$.
Interpolating these bounds with parameter $\vartheta<1$ gives
an integrable time kernel:
\[
 M_0(t)^{2-\vartheta}(1+t)^{-\beta/4}
 \int_{t-1}^t(t-s)^{-\vartheta/2}\cD_0(U(s))^{\vartheta/2}\dd s
 \le CM_0(t)^2(1+t)^{-\beta/4-\vartheta/2}.
\]
The H\"older dual kernel has exponent $\vartheta/(2-\vartheta)<1$.
Since $\beta/4+\vartheta/2-\beta=1/64>0$, this contribution is
also bounded by $CM_0(t)^2(1+t)^{-\beta}$. This stronger bound
also applies to the third-component estimate in the next paragraph.
Combining both time intervals with the homogeneous term gives
$\|v(t)\|_2\le C\mathcal J_2 (1+t)^{-1/2}+CM_0(t)^2(1+t)^{-1/2}$.
The compressible estimates for $P_LW$ and the energy $E_H$ give
$\|w(t)\|_2\le(C\mathcal J_2+CM_0(t)^2)(1+t)^{-\beta}$.
Since $u=v+w$ and $\beta>1/2$, it follows that
\[
 \|u(t)\|_2\le(C\mathcal J_2+CM_0(t)^2)(1+t)^{-1/2}.
\]

\paragraph{Decay of the third velocity component.}
We estimate $v_3$ at the rate $(1+t)^{-\beta}$, $\beta=5/8$,
and then recover $u_3=v_3+w_3$.
The third row of the projection satisfies
\begin{equation}\label{eq:new-third-row}
 \sum_j|\mathbb P_{3j}(\xi)|^2=\frac{|\xi_h|^2}{|\xi|^2},
 \quad \int_\R\frac{|\xi_h|^2}{|\xi_h|^2+\xi_3^2}\dd\xi_3
       =\pi|\xi_h|,
 \quad \partial_3[\mathbb PG]_3=-\Div_h[\mathbb PG]_h
\end{equation}
For $a>0$, Plancherel and $|\widehat G|\le C\|G\|_1$ give
\[
\|[S_h(a)\mathbb PG]_3\|_2^2
 \le C\|G\|_1^2\int_{\R^3}
       e^{-2\nu_ha|\xi_h|^2}\frac{|\xi_h|^2}{|\xi|^2}\dd\xi
 \le Ca^{-3/2}\|G\|_1^2.
\]
The last inequality uses the vertical integral in
\eqref{eq:new-third-row}, followed by horizontal scaling. Thus
$\|[S_h(a)\mathbb PG]_3\|_2\le Ca^{-3/4}\|G\|_1$.
In particular, the homogeneous term is bounded by
$C\mathcal J_2(1+t)^{-3/4}$, since $u_0\in L^1$.
We split the Duhamel integral at $t-1$ as in the estimate for $v$.

\par\smallskip\noindent\emph{The interval $[0,t-1]$.}\enspace
Here $a=t-s\ge1$. The unprojected momentum source \eqref{eq:F2} satisfies
\[
 \|F_2\|_1\le CM_0(t) (1+s)^{-1/2}\cD_0(U)^{1/2}+CM_0(t)^2(1+s)^{-\beta-1/4}
                         +CM_0(t) (1+s)^{-\beta}\cD_2(U)^{1/2}.
\]
The three terms arise from horizontal convection, vertical
convection, and pressure together with variable viscosity,
respectively, using \eqref{eq:new-seed-ledger}.
Apply the third-row kernel obtained from \eqref{eq:new-third-row}
to each contribution. For horizontal convection, weighted
Cauchy--Schwarz and \eqref{eq:new-low-tail} yield
\[
\begin{aligned}
 &M_0(t)\int_0^{t-1}(t-s)^{-3/4}(1+s)^{-1/2}\cD_0(U(s))^{1/2}\dd s\\
 &\quad\le CM_0(t)^2
       \left(\int_0^{t-1}(t-s)^{-3/2}
                        (1+s)^{-1-r}\dd s\right)^{1/2}
 \le CM_0(t)^2(1+t)^{-3/4},
\end{aligned}
\]
since $1+r=63/32>3/2$. For vertical convection,
\[
 M_0(t)^2\int_0^{t-1}(t-s)^{-3/4}(1+s)^{-\beta-1/4}\dd s
                  \le CM_0(t)^2(1+t)^{-\beta}.
\]
Indeed, splitting at $t/2$ gives the power
$3/4+(\beta+1/4)-1=\beta$, with $3/4<1$ and
$\beta+1/4=7/8<1$.
For the pressure and variable-viscosity contribution, use the total
dissipation bound \eqref{eq:new-base-energy}:
\[
\begin{aligned}
 &M_0(t)\int_0^{t-1}(t-s)^{-3/4}(1+s)^{-\beta}
                         \cD_2(U(s))^{1/2}\dd s\\
 &\quad\le C\mathcal J_2M_0(t)
       \left(\int_0^{t-1}(t-s)^{-3/2}
                         (1+s)^{-2\beta}\dd s\right)^{1/2}
 \le CM_0(t)^2(1+t)^{-\beta}.
\end{aligned}
\]
Here $1<2\beta=5/4<3/2$ determines the convolution decay,
and $\mathcal J_2\le M_0(t)$ gives the last amplitude.
Since $3/4>\beta$, the entire contribution from $[0,t-1]$
is bounded by $CM_0(t)^2(1+t)^{-\beta}$.

\par\smallskip\noindent\emph{The interval $[t-1,t]$.}\enspace
The $L^1$--$L^2$ kernel cannot be paired directly with the
$L^2$ time bound \eqref{eq:new-low-tail} for $\cD_0(U)^{1/2}$ near $s=t$.
For convection and pressure, use $L^2$ contraction as in the
$[t-1,t]$ estimate for $v$. For viscosity, use the interpolation
between the mixed-norm and $L^2$ bounds for the two terms in
\eqref{eq:new-projected-source}. The convection and pressure
integrals are bounded by
$CM_0(t)^2(1+t)^{-\beta}$. The viscosity contribution has exponent
$\beta/4+\vartheta/2=41/64>\beta$. Applying these vector estimates
to the third component gives
\[
 \left\|\left[\int_{t-1}^t
      S_h(t-s)\mathbb PF_2(s)\dd s\right]_3\right\|_2
                  \le CM_0(t)^2(1+t)^{-\beta}.
\]
Combining both intervals with the homogeneous term proves
\[
 \|v_3(t)\|_2\le C\mathcal J_2(1+t)^{-3/4}
                         +CM_0(t)^2(1+t)^{-\beta}.
\]
The bound for $P_LW$ from \eqref{eq:ScDuhamel} and
\eqref{eq:new-seed-low-sources}, together with the bound for
$E_H\simeq\|P_HW\|_2^2$, gives
$\|w_3(t)\|_2\le(C\mathcal J_2+CM_0(t)^2)(1+t)^{-\beta}$.
Since $u_3=v_3+w_3$ and $3/4>\beta$, it follows that
\[
 \|u_3(t)\|_2\le(C\mathcal J_2+CM_0(t)^2)(1+t)^{-\beta}.
\]
Combining the bounds for $W$, $u$, and $u_3$ gives
\begin{equation}\label{eq:r12-seed-closure}
                 M_0(t)\le C_0\mathcal J_2+C_1M_0(t)^2.
\end{equation}
The constants $C_0,C_1$ are independent of $t$. The function $M_0$ is continuous, and its initial value is
at most $C_0\mathcal J_2$ after enlarging $C_0$.
Choose $\mathcal J_2$ so that
$2C_0\mathcal J_2<\min\{1,(2C_1)^{-1}\}$.
At a first point where $M_0=2C_0\mathcal J_2$, the right side
of \eqref{eq:r12-seed-closure} would be strictly smaller
than $2C_0\mathcal J_2$. Thus the bootstrap closes on every
finite interval. In particular, $M_0(t)\le2C_0\mathcal J_2$;
its definition gives both decay estimates in \eqref{eq:new-seed}.
The zero datum gives the zero solution.
Finally apply \eqref{eq:new-low-tail} with $M=C\mathcal J_2$ and
let the finite endpoint tend to infinity. This proves the
dissipation tail and weighted integral stated in the lemma.
\end{proof}

\subsubsection{Zero- and first-order decay}

The next lemma also supplies \eqref{eq:new-first-start}, the initial
gradient bounds for Lemma~\ref{lem:new-first}.
\begin{lemma}[Zeroth-order decay and preliminary gradient bounds]\label{lem:new-zero}
The estimates in \eqref{eq:new-seed} improve to
\begin{equation}\label{eq:new-zero}
 \|(W,u_3)(t)\|_2\le C\mathcal J_2 (1+t)^{-3/4}.
\end{equation}
The proof also gives the preliminary gradient bounds
\begin{equation}\label{eq:new-first-start}
 \|\nabla(q,w)\|_2\le C\mathcal J_2 (1+t)^{-5/8},\qquad
 \|(\nabla_hu,\Div u)\|_2+\|\nabla\tau\|_2\le C\mathcal J_2 (1+t)^{-1/2}.
\end{equation}
\end{lemma}
\begin{proof}
For $0\le t\le2$, \eqref{eq:new-base-energy} gives
\eqref{eq:new-zero} and \eqref{eq:new-first-start}; hence assume $t\ge2$. We first estimate $(q,w,u_3)$ in Steps~1--3, then derive preliminary gradient bounds in Step~4 and use them to estimate the temperature in Step~5.
Until Step~5, use only the preliminary bound \eqref{eq:new-seed}
for $\tau$ and the $H^{-1}$ estimate for $\tau_H$; the thermal
$L^2$ estimate requires \eqref{eq:new-first-start} from Step~4.

\smallskip\noindent\emph{Step 1. Low-frequency estimates.}
For the rest of this subsection, write
$F=F^{(0)}+\sum_{j=1}^3\partial_jF^{(j)}$ with the pressure
gradient included in the flux. Thus, for $j=1,2,3$, the
fluxes in \eqref{eq:new-basic-flux} are replaced by
\begin{equation}\label{eq:new-decay-flux}
\begin{aligned}
 F^{(j)}_1&=-u_jq,\\
 F^{(j)}_2&=-u_ju-c^2(e^\tau-1)(q+\tau)e_j
 +\nu_h\mathbf1_{j\le2}(e^{-q}-1)\partial_ju
              +\nu_c(e^{-q}-1)(\Div u)e_j,\\
 F^{(j)}_3&=-u_j\tau+\chi(e^{-q}-1)\partial_j\tau.
\end{aligned}
\end{equation}
The corresponding remainder is
\begin{equation}\label{eq:new-flux-remainders}
\begin{aligned}
 F^{(0)}_1&=q\Div u,\\
 F^{(0)}_2&=u\Div u+c^2(q+\tau)\nabla(e^\tau-1)\\
 &\quad-\nu_h\nabla_h(e^{-q}-1)\cdot\nabla_hu
              -\nu_c(\Div u)\nabla(e^{-q}-1),\\
 F^{(0)}_3&=\tau\Div u-\chi\nabla(e^{-q}-1)\cdot\nabla\tau
 +\frac{e^{-q-\tau}}{\bar\rho\bar\theta c_v}Q(u)
              +\chi e^{-q}|\nabla\tau|^2.
\end{aligned}
\end{equation}
The bounds \eqref{eq:new-seed} and \eqref{eq:new-seed-ledger} give
\[
 \begin{aligned}
 \|F^{(0)}\|_1&\le C\mathcal J_2(1+s)^{-5/16}\cD_0(U)^{1/2}+C\cD_0(U),\\
 \sum_j\|F^{(j)}\|_1&\le C\mathcal J_2^2(1+s)^{-1}
       +C\mathcal J_2(1+s)^{-5/8}\cD_0(U)^{1/2}.
 \end{aligned}
\]
The low-frequency multiplier bound \eqref{eq:blockmultiplier}
and Duhamel's formula give
\[
\begin{aligned}
 \|P_LW(t)\|_2
 &\le C\mathcal J_2(1+t)^{-3/4}
   +C\int_0^t(1+t-s)^{-3/4}\|F^{(0)}(s)\|_1\dd s\\
 &\quad+C\int_0^t(1+t-s)^{-5/4}
                         \sum_j\|F^{(j)}(s)\|_1\dd s.
\end{aligned}
\]
The extra derivative in each flux term gives the kernel exponent
$5/4$. Take $r=31/32$ as in Lemma~\ref{lem:new-seed}; thus $r<1$ is
admissible in \eqref{eq:new-low-tail}, and $5/8+r>3/2$. Lemma~\ref{lem:new-seed} gives
$\int_0^\infty(1+s)^r \cD_0(U(s))\dd s\le C_r\mathcal J_2^2$.
Thus weighted Cauchy--Schwarz yields
\[
\begin{aligned}
 &\mathcal J_2\int_0^t(1+t-s)^{-3/4}(1+s)^{-5/16}\cD_0(U(s))^{1/2}\dd s\\
 &\quad\le\mathcal J_2
       \left(\int_0^t(1+s)^r \cD_0(U(s))\dd s\right)^{1/2}
       \left(\int_0^t(1+t-s)^{-3/2}
                         (1+s)^{-5/8-r}\dd s\right)^{1/2}\\
 &\quad\le C_r\mathcal J_2^2(1+t)^{-3/4}.
\end{aligned}
\]
The last bound follows by splitting the time integral at $t/2$,
since $5/8+r>3/2>1$.
For the $\cD_0(U)$ contribution to $F^{(0)}$, use the total and
tail dissipation bounds in Lemma~\ref{lem:new-seed}. For the
two flux contributions, use direct time convolution and
weighted Cauchy--Schwarz, respectively. This gives
\[
\begin{aligned}
 \int_0^t(1+t-s)^{-3/4}\cD_0(U(s))\dd s
   &\le C\mathcal J_2^2(1+t)^{-3/4},\\
 \mathcal J_2^2\int_0^t(1+t-s)^{-5/4}(1+s)^{-1}\dd s
   &\le C\mathcal J_2^2(1+t)^{-1},\\
 \mathcal J_2\int_0^t(1+t-s)^{-5/4}(1+s)^{-5/8}\cD_0(U(s))^{1/2}\dd s
   &\le C_r\mathcal J_2^2(1+t)^{-5/8-r/2}.
\end{aligned}
\]
Since $5/8+r/2>3/4$, all source and flux contributions are bounded
by $C\mathcal J_2^2(1+t)^{-3/4}$. Adding the homogeneous term
and using $\mathcal J_2\le1$ proves
\begin{equation}\label{eq:new-zero-low}
 \|P_LW(t)\|_2\le C\mathcal J_2 (1+t)^{-3/4}.
\end{equation}
The fixed Fourier support of $P_LW$ also implies from
\eqref{eq:new-zero-low} that each fixed number of its derivatives
is bounded by $C\mathcal J_2(1+t)^{-3/4}$ in
$L^2,L^3,L^\infty$.

\smallskip\noindent\emph{Step 2. High-frequency density and potential velocity.}
We estimate $q_H,w_H$ in $L^2$, but initially measure $\tau_H$
in $H^{-1}$. Testing the heat equation against
$(1-\Delta)^{-1}\tau_H$ then requires only an $H^{-2}$ estimate
of $F_3$. In particular, the quadratic source terms can be
bounded in $L^1$. This energy will give pointwise decay for
$q_H,w_H$ while requiring only the preliminary bound for $\tau$.

\smallskip\noindent\emph{Energy and linear couplings.}
Keep density transport on the left, as in
\eqref{eq:new-density-high}. The high-frequency equations are
\[
\left\{
\begin{aligned}
 (q_H)_s+\Div w_H+u\cdot\nabla q_H&=C_q,\\
 (w_H)_s+c^2\nabla(q_H+\tau_H)-\mathcal Lw_H
                                      &=P_H\mathbb QF_2,\\
 (\tau_H)_s+\alpha\Div w_H-\chi\Delta\tau_H&=P_HF_3.
\end{aligned}
\right.
\]
Choose $A_\tau$ large, then the energy cutoff radius $R_0$ large,
and finally $\eta>0$ small, to meet the finitely many absorption
conditions here and in the first- and second-order energies below.
Fix these choices throughout; low-frequency constants may depend on $R_0$,
and the data smallness is chosen afterwards. The small Green-function
cutoff used later is separate from this energy cutoff. Set
\[
 \mathscr H=\tfrac12(c^2\|q_H\|_2^2+\|w_H\|_2^2
          +A_\tau\ip{\tau_H}{(1-\Delta)^{-1}\tau_H})
          +\eta\ip{w_H}{\nabla (1-\Delta)^{-1}q_H}.
\]
Testing these equations with $c^2q_H$, $w_H$ and
$A_\tau(1-\Delta)^{-1}\tau_H$, respectively, cancels the
linear density--velocity coupling. The principal dissipative
terms are
\[
 \nu_h\|\nabla_hw_H\|_2^2+\nu_c\|\Div w_H\|_2^2
       +A_\tau\chi\|\nabla(1-\Delta)^{-1/2}\tau_H\|_2^2.
\]
Since $w_H=\mathbb Qw_H$, $\|\Div w_H\|_2=\|\nabla w_H\|_2$.
On the high-frequency support $|\xi|\ge R_0$, the thermal
dissipation is equivalent to $A_\tau\|\tau_H\|_2^2$.
The two uncancelled thermal couplings satisfy
\[
\begin{aligned}
 c^2|\ip{\tau_H}{\Div w_H}|
   &\le\frac{\nu_c}{8}\|\nabla w_H\|_2^2+C\|\tau_H\|_2^2,\\
 A_\tau\alpha|\ip{(1-\Delta)^{-1}\tau_H}{\Div w_H}|
   &\le\frac{\nu_c}{8}\|\nabla w_H\|_2^2
          +CA_\tau^2(1+R_0^2)^{-2}\|\tau_H\|_2^2.
\end{aligned}
\]
Choose $A_\tau$ large enough to absorb the first thermal remainder,
then $R_0$ large enough to absorb the second.

The cross term supplies the missing density dissipation.
Using $(q_H)_s+\Div w_H=P_HF_1$, its derivative is
\[
\begin{aligned}
 &\frac{\dd}{\dd s}\ip{w_H}{\nabla(1-\Delta)^{-1}q_H}
       +c^2\|\nabla(1-\Delta)^{-1/2}q_H\|_2^2\\
 &\quad=\ip{\mathcal Lw_H}{\nabla(1-\Delta)^{-1}q_H}
        -c^2\ip{\nabla\tau_H}{\nabla(1-\Delta)^{-1}q_H}
        +\ip{\Div w_H}{(1-\Delta)^{-1}\Div w_H}\\
 &\qquad+\ip{P_H\mathbb QF_2}{\nabla(1-\Delta)^{-1}q_H}
        -\ip{P_HF_1}{(1-\Delta)^{-1}\Div w_H}.
\end{aligned}
\]
The density term on the left is equivalent to $\|q_H\|_2^2$.
Integrating the viscous pairing by parts bounds it by
$C\|\nabla w_H\|_2\|q_H\|_2$; the other linear pairings
are bounded by $C\|\tau_H\|_2\|q_H\|_2+C\|w_H\|_2^2$.
After multiplication by a sufficiently small $\eta$, these terms
are absorbed by the preceding dissipation and a part of the
density term. The same small choice of $\eta$ gives
\[
 \mathscr H\simeq\|q_H\|_2^2+\|w_H\|_2^2
                         +\|\tau_H\|_{H^{-1}}^2.
\]
Thus its dissipation controls
$\mathscr H+\|\nabla w_H\|_2^2+\|\tau_H\|_2^2$.

\smallskip\noindent\emph{Coefficient decomposition and source bounds.}
Split the density coefficient as
\begin{equation}\label{eq:new-coefficient-split}
 (e^{-q}-1)=b_{\ell}+b_{\mathrm r},
 \quad b_{\ell}=e^{-q_L}-1,\quad b_{\mathrm r}=e^{-q}-e^{-q_L}.
\end{equation}
The mean-value formula
$b_{\mathrm r}=-q_H\int_0^1e^{-q_L-\sigma q_H}\dd\sigma$
gives $\|b_{\mathrm r}\|_2\le C\|q_H\|_2$. The low-frequency estimate
\eqref{eq:new-zero-low}, together with the chain rule, gives
\[
 \|b_{\ell}\|_{W^{2,\infty}}+\|\nabla b_{\ell}\|_3
                       \le C\mathcal J_2(1+s)^{-3/4}.
\]
Here $b_{\ell}$ is defined from the low-frequency density $q_L$;
no Fourier support property is assumed for $b_{\ell}$ or $b_{\mathrm r}$.
We estimate $F_1,F_2,F_3$, using the decomposition $q=q_L+q_H$.
The boundedness of $P_H$ and $\mathbb Q$ on $H^{-1}$ and $H^{-2}$
then gives the corresponding bounds for $P_HF_1$,
$P_H\mathbb QF_2$, and $P_HF_3$.

The required source bounds are
\begingroup
\allowdisplaybreaks[0]
\begin{align}\label{eq:new-weakened-sources}
 \|F_1\|_{H^{-1}}^2&\le C\mathcal J_2^2\|q_H\|_2^2
                           +C\mathcal J_2^4(1+s)^{-5/2},\notag\\
 \|F_2\|_{H^{-1}}^2&\le C(\mathcal J_2^2+\cD_2(U))\|q_H\|_2^2
        +C\mathcal J_2^2(1+s)^{-3/4}\cD_0(U)+C\mathcal J_2^4(1+s)^{-25/16},\notag\\
 \|F_3\|_{H^{-2}}^2&\le C\cD_2(U)\|q_H\|_2^2
        +C\mathcal J_2^2(1+s)^{-1}\cD_0(U)+C\mathcal J_2^4(1+s)^{-11/4}+C\cD_0(U)^2.
\end{align}
\endgroup

For $F_1$, write $-u\cdot\nabla q=-u\cdot\nabla q_H
-u\cdot\nabla q_L$. Integration by parts against an $H^1$
test function and \eqref{eq:H2-H1-multiplier} yield
\[
 \|u\cdot\nabla q_H\|_{H^{-1}}
       \le C\|u\|_{H^2}\|q_H\|_2,\qquad
 \|u\cdot\nabla q_L\|_2
       \le\|u\|_2\|\nabla q_L\|_\infty
       \le C\mathcal J_2^2(1+s)^{-5/4}.
\]
For the convection term in $F_2$, defined in \eqref{eq:F2},
write $-u\cdot\nabla u=-\Div(u\otimes u)+u\Div u$.
The anisotropic $L^4$ inequality, the velocity bounds in
\eqref{eq:new-seed} and \eqref{eq:new-seed-ledger}, and the
definition of $\cD_0(U)^{1/2}$ in \eqref{eq:new-d0} give
\[
 \|u\otimes u\|_2^2=\|u\|_4^4
       \le C\|u\|_2\|\partial_3u\|_2\|\nabla_hu\|_2^2
       \le C\mathcal J_2^2(1+s)^{-3/4}\cD_0(U).
\]
For $u\Div u$, interpolation between \eqref{eq:new-seed} and
\eqref{eq:new-base-energy} gives
$\|u\|_3\le C\|u\|_2^{3/4}\|\nabla^2u\|_2^{1/4}
\le C\mathcal J_2(1+s)^{-3/8}$.
Thus $L^{6/5}\hookrightarrow H^{-1}$ yields
\[
 \|u\Div u\|_{H^{-1}}
 \le C\|u\|_3\|\Div u\|_2
 \le C\mathcal J_2(1+s)^{-3/8}\cD_0(U)^{1/2}.
\]
For pressure, split $\nabla q$ into its low and high parts.
The high part has $H^{-1}$ norm at most
$C\|e^\tau-1\|_{H^2}\|q_H\|_2\le C\mathcal J_2\|q_H\|_2$,
and the low part is bounded by
$\|e^\tau-1\|_2\|\nabla q_L\|_\infty
\le C\mathcal J_2^2(1+s)^{-11/8}$.
For the thermal part use
$(e^\tau-1)\nabla\tau=\nabla(e^\tau-1-\tau)$ and
\[
 \|(e^\tau-1)\nabla\tau\|_{H^{-1}}
   \le C\|\tau\|_4^2
   \le C\|\tau\|_2^{5/4}\|\tau\|_{H^2}^{3/4}
   \le C\mathcal J_2^2(1+s)^{-25/32}.
\]

For the variable viscosity, duality with $H^1\hookrightarrow L^6$
and $\|\mathcal Lu\|_3\le C\mathcal D_2(U)^{1/2}$ give
\[
 \|b_{\mathrm r}\mathcal Lu\|_{H^{-1}}
   \le C\|b_{\mathrm r}\|_2\|\mathcal Lu\|_3
   \le C\mathcal D_2(U)^{1/2}\|q_H\|_2.
\]
For $b_{\ell}\mathcal Lu$, integrate by parts once. The resulting flux
and coefficient-gradient terms satisfy
\[
 \|b_{\ell}\mathcal Lu\|_{H^{-1}}
   \le C\bigl(\|b_{\ell}\|_\infty+\|\nabla b_{\ell}\|_3\bigr)\cD_0(U)^{1/2}
   \le C\mathcal J_2(1+s)^{-3/4}\cD_0(U)^{1/2}.
\]
For $F_3$, scalar transport is bounded directly in $L^1$:
$\|u\cdot\nabla\tau\|_{H^{-2}}
\le C\|u\|_2\|\nabla\tau\|_2
\le C\mathcal J_2(1+s)^{-1/2}\cD_0(U)^{1/2}$.
Since $H^2\hookrightarrow L^\infty$, duality gives
$L^1\hookrightarrow H^{-2}$ and hence
\[
 \|b_{\mathrm r}\Delta\tau\|_{H^{-2}}
   \le C\|b_{\mathrm r}\Delta\tau\|_1
   \le C\|b_{\mathrm r}\|_2\|\Delta\tau\|_2
   \le C\mathcal D_2(U)^{1/2}\|q_H\|_2.
\]
For $b_{\ell}=e^{-q_L}-1$, use
$b_{\ell}\Delta\tau=\Delta(b_{\ell}\tau)-2\Div(\tau\nabla b_{\ell})
+\tau\Delta b_{\ell}$ to obtain
\[
 \|b_{\ell}\Delta\tau\|_{H^{-2}}
   \le C\|b_{\ell}\|_{W^{2,\infty}}\|\tau\|_2
   \le C\mathcal J_2^2(1+s)^{-11/8}.
\]
The quadratic source terms $e^{-q-\tau}Q(u)$ and
$e^{-q}|\nabla\tau|^2$ in \eqref{eq:F3} have $L^1$ norm at
most $C\cD_0(U)$. Their coefficients are bounded by
\eqref{eq:new-base-energy}, and \eqref{eq:new-d0} controls
the quadratic factors $|\nabla_hu|^2$, $(\Div u)^2$ and
$|\nabla\tau|^2$.
The preceding product estimates prove
\eqref{eq:new-weakened-sources} after squaring and summation.

\smallskip\noindent\emph{The nonlinear energy inequality.}
The density transport test contributes
$\frac{c^2}{2}\int\Div u\,q_H^2\dd x+c^2\ip{C_q}{q_H}$.
To estimate $C_q$, write
\[
 C_q=\Div P_L(uq)-P_L(q\Div u)-u\cdot\nabla q_L.
\]
The fixed operators $P_L,\nabla P_L$ map $L^1$ into $L^2$.
Thus \eqref{eq:new-seed} and \eqref{eq:new-zero-low} imply
\[
\begin{aligned}
 \|C_q\|_2
  &\le C\|u\|_2\|q\|_2+C\|q\|_2\|\Div u\|_2
                          +\|u\|_2\|\nabla q_L\|_\infty\\
  &\le C\mathcal J_2^2(1+s)^{-9/8}
                         +C\mathcal J_2(1+s)^{-5/8}\cD_0(U)^{1/2}.
\end{aligned}
\]
Using $\|\Div u\|_\infty\le C\mathcal D_2(U)^{1/2}$ and
Young's inequality, the transport contribution is bounded by
a small part of $\|q_H\|_2^2$ plus
$C\mathcal D_2(U)\|q_H\|_2^2+C\|C_q\|_2^2$.

The velocity source pairs in $H^{-1}$--$H^1$. For the heat source,
the test function satisfies
$\|(1-\Delta)^{-1}\tau_H\|_{H^2}=\|\tau_H\|_2$,
so its pairing is at most
$\|F_3\|_{H^{-2}}\|\tau_H\|_2$.
In the cross term, the test functions satisfy
$\|\nabla(1-\Delta)^{-1}q_H\|_{H^1}\le C\|q_H\|_2$
and $\|(1-\Delta)^{-1}\Div w_H\|_{H^1}\le C\|w_H\|_2$.
Young's inequality therefore bounds all source pairings by
a small part of the dissipation plus
\[
 C\bigl(\|F_1\|_{H^{-1}}^2+\|F_2\|_{H^{-1}}^2
                   +\|F_3\|_{H^{-2}}^2+\|C_q\|_2^2\bigr).
\]
Inserting \eqref{eq:new-weakened-sources}, absorb the terms
$C\mathcal J_2^2\|q_H\|_2^2$ by smallness of $\mathcal J_2$.
The terms with $\mathcal D_2(U)\|q_H\|_2^2$ are bounded by
$C\mathcal D_2(U)\mathscr H$.
Interpolation gives $\cD_0(U)^{1/2}\le C\mathcal J_2(1+s)^{-1/4}$, and hence
\[
 \cD_0(U)^2\le C\mathcal J_2^2(1+s)^{-1/2}\cD_0(U),\qquad
 \|C_q\|_2^2\le C\mathcal J_2^4(1+s)^{-9/4}
                         +C\mathcal J_2^2(1+s)^{-5/4}\cD_0(U).
\]
The terms proportional to $\mathcal J_2^4$ in
\eqref{eq:new-weakened-sources} and the bound for
$\|C_q\|_2^2$ decay at least as fast as $(1+s)^{-25/16}$.
The terms containing $\cD_0(U)$, including the bound for $\cD_0(U)^2$,
are at most $C\mathcal J_2^2(1+s)^{-1/2}\cD_0(U)$.
The resulting inequality for $\mathscr H$ is
\[
 \mathscr H'+c_*(\mathscr H+\|\nabla w_H\|_2^2+\|\tau_H\|_2^2)
 \le C\cD_2(U)\mathscr H+C\mathcal J_2^4(1+s)^{-25/16}
                         +C\mathcal J_2^2(1+s)^{-1/2}\cD_0(U).
\]

\smallskip\noindent\emph{Time decay.}
The integrating factor is bounded by \eqref{eq:new-base-energy}.
Lemma~\ref{lem:new-seed} gives
$\int_a^{a+1}\cD_0(U(s))\dd s\le C\mathcal J_2^2(1+a)^{-1}$.
Thus \eqref{eq:new-window-convolution}, with $p=1/2$ and
$\beta=1$, bounds the weighted dissipation forcing by
$C\mathcal J_2^4(1+t)^{-3/2}$.
Since $25/16>3/2$ and $\mathscr H(0)\le C\mathcal J_2^2$,
\[
 \mathscr H(t)\le C\mathcal J_2^2e^{-c_*t}
                    +C\mathcal J_2^4(1+t)^{-3/2}.
\]
The energy equivalence now gives
\[
 \|(q_H,w_H)(t)\|_2
   \le C\mathcal J_2e^{-c_*t/2}
                    +C\mathcal J_2^2(1+t)^{-3/4}.
\]
Adding the low-frequency bound \eqref{eq:new-zero-low} gives
\[
 \|(q,w)(t)\|_2\le C\mathcal J_2(1+t)^{-3/4}.
\]
Here the energy controls $\tau_H$ only in $H^{-1}$.
Its pointwise $L^2$ decay will follow in Step~5 from the
preliminary gradient estimates of Step~4.

\smallskip\noindent\emph{Step 3. Third velocity component.}
Since $u_3=v_3+w_3$ and $w_3$ has the required decay by Step~2,
it remains to improve the decay of $v_3$. Duhamel's formula is
\[
 v_3(t)=[S_h(t)\mathbb Pu_0]_3
              +\int_0^t[S_h(t-s)\mathbb PF_2(s)]_3\dd s.
\]
The homogeneous term is bounded by
$C\mathcal J_2(1+t)^{-3/4}$, as proved using
\eqref{eq:new-third-row} in Lemma~\ref{lem:new-seed}.
We split the source integral at $t-1$.

\par\smallskip\noindent\emph{The interval $[0,t-1]$.}\enspace
Here $t-s\ge1$. The identities in \eqref{eq:new-third-row} and
the horizontal heat estimates give, for the appropriate vector
or tensor source $G$ and $a>0$,
\[
\begin{aligned}
 \|[S_h(a)\mathbb P\Div_hG]_3\|_2
       &\le Ca^{-5/4}\|G\|_1,\\
 \|[S_h(a)\mathbb P\partial_3G]_3\|_2
       &\le Ca^{-1}\|G\|_{L^2_3L^1_h},
 \qquad \|[S_h(a)\mathbb PG]_3\|_2
       \le Ca^{-3/4}\|G\|_1.
\end{aligned}
\]
The first estimate adds one horizontal derivative to the
third-component $L^1$--$L^2$ estimate. The second uses
$\partial_3[\mathbb PG]_3=-\Div_h[\mathbb PG]_h$, so the
vertical derivative is estimated by a horizontal derivative
of the solenoidal field. For this mixed-norm estimate,
Plancherel in the vertical variable and the horizontal
$L^1$ Fourier bound give
$\|\widehat G(\xi_h,\cdot)\|_{L^2_{\xi_3}}
\le C\|G\|_{L^2_3L^1_h}$, uniformly in $\xi_h$. Consequently,
\[
\begin{aligned}
 \|[S_h(a)\mathbb P\partial_3G]_3\|_2^2
 &\le C\int_{\R^2}e^{-2\nu_ha|\xi_h|^2}|\xi_h|^2
              \|\widehat G(\xi_h,\cdot)\|_{L^2_{\xi_3}}^2\dd\xi_h\\
 &\le C\|G\|_{L^2_3L^1_h}^2
             \int_{\R^2}e^{-2\nu_ha|\xi_h|^2}|\xi_h|^2\dd\xi_h
 \le Ca^{-2}\|G\|_{L^2_3L^1_h}^2.
\end{aligned}
\]
The projection is thus handled inside the Fourier estimate;
the argument does not require its boundedness on $L^1$.

Decompose convection as
\[
 -u\cdot\nabla u
     =-\Div_h(u\otimes u_h)-\partial_3(u_3u)+u\Div u.
\]
The horizontal flux, vertical flux, and remainder use the
$L^1$, mixed-norm, and $L^1$ semigroup bounds, respectively.
Their source norms satisfy
\[
\begin{aligned}
 \|u\otimes u_h\|_1
       &\le\|u\|_2^2\le C\mathcal J_2^2(1+s)^{-1},\\
 \|u_3u\|_{L^2_3L^1_h}
       &\le\|u_3\|_{L^\infty_3L^2_h}\|u\|_2
        \le C\mathcal J_2^{3/2}(1+s)^{-13/16}\cD_0(U)^{1/4},\\
 \|u\Div u\|_1
       &\le\|u\|_2\|\Div u\|_2
        \le C\mathcal J_2(1+s)^{-1/2}\cD_0(U)^{1/2}.
\end{aligned}
\]
To estimate the vertical flux $u_3u$, use the preliminary
bound for $u_3$ in \eqref{eq:new-seed} and the derivative
identity \eqref{eq:verticalcomp}:
\[
 \|u_3\|_{L^\infty_3L^2_h}
    \le C\|u_3\|_2^{1/2}\|\partial_3u_3\|_2^{1/2}
    \le C\mathcal J_2^{1/2}(1+s)^{-5/16}\cD_0(U)^{1/4}.
\]
Multiplication by $\|u\|_2\le C\mathcal J_2(1+s)^{-1/2}$
gives decay exponent $13/16$.
For the time convolutions, take $r=31/32$ and use
$\int_0^\infty(1+s)^r \cD_0(U(s))\dd s\le C\mathcal J_2^2$
from Lemma~\ref{lem:new-seed}.

For the horizontal flux, splitting at $t/2$ gives
\[
\mathcal J_2^2\int_0^{t-1}(t-s)^{-5/4}(1+s)^{-1}\dd s
 \le C\mathcal J_2^2
       \bigl((1+t)^{-5/4}\log(e+t)+(1+t)^{-1}\bigr)
        \le C\mathcal J_2^2(1+t)^{-1}.
\]
For the vertical flux, H\"older's inequality with exponents
$4$ and $4/3$ gives
\[
\begin{aligned}
 &\mathcal J_2^{3/2}\int_0^{t-1}(t-s)^{-1}
                         (1+s)^{-13/16}\cD_0(U(s))^{1/4}\dd s\\
 &\quad\le\mathcal J_2^{3/2}
       \left(\int_0^{t-1}(1+s)^r \cD_0(U(s))\dd s\right)^{1/4}
       \left(\int_0^{t-1}(t-s)^{-4/3}
                    (1+s)^{-13/12-r/3}\dd s\right)^{3/4}\\
 &\quad\le C\mathcal J_2^2(1+t)^{-1}.
\end{aligned}
\]
Indeed $13/12+r/3=45/32>4/3$; splitting at $t/2$ bounds
the last convolution integral by $C(1+t)^{-4/3}$.

For pressure, write
\[
 -c^2(e^\tau-1)\nabla(q+\tau)
   =-c^2\nabla\bigl((e^\tau-1)(q+\tau)\bigr)
                    +c^2(q+\tau)e^\tau\nabla\tau.
\]
The first term vanishes after applying $\mathbb P$.
For the remaining pressure term, Step~2 gives the $3/4$ decay
of $q$, while Lemma~\ref{lem:new-seed} gives the preliminary
$5/8$ decay of $\tau$. Hence
$\|(q+\tau)e^\tau\nabla\tau\|_1
\le C\|(q,\tau)\|_2\|\nabla\tau\|_2
\le C\mathcal J_2(1+s)^{-5/8}\cD_0(U)^{1/2}$.
The convection remainder $u\Div u$ and this pressure remainder
use the same weighted Cauchy--Schwarz estimate, with
$\gamma=1/2$ and $\gamma=5/8$, respectively:
\[
\begin{aligned}
 &\mathcal J_2\int_0^{t-1}(t-s)^{-3/4}(1+s)^{-\gamma}\cD_0(U(s))^{1/2}\dd s\\
 &\quad\le\mathcal J_2
       \left(\int_0^{t-1}(1+s)^r \cD_0(U(s))\dd s\right)^{1/2}
       \left(\int_0^{t-1}(t-s)^{-3/2}(1+s)^{-2\gamma-r}\dd s\right)^{1/2}\\
 &\quad\le C\mathcal J_2^2(1+t)^{-3/4}.
\end{aligned}
\]
For both values, $2\gamma+r\ge1+r=63/32>3/2$,
so splitting at $t/2$ bounds the last convolution integral
by $C(1+t)^{-3/2}$.

Finally, the variable-viscosity source satisfies
\[
 \|(e^{-q}-1)\mathcal Lu\|_1
   \le C\|q\|_2\|\mathcal Lu\|_2
   \le C\mathcal J_2(1+s)^{-3/4}\mathcal D_2(U)^{1/2}.
\]
Use the total dissipation estimate \eqref{eq:new-base-energy}
in Cauchy--Schwarz to obtain
\[
\begin{aligned}
 &\mathcal J_2\int_0^{t-1}(t-s)^{-3/4}(1+s)^{-3/4}
                              \mathcal D_2(U(s))^{1/2}\dd s\\
 &\quad\le\mathcal J_2
       \left(\int_0^{t-1}\mathcal D_2(U(s))\dd s\right)^{1/2}
       \left(\int_0^{t-1}(t-s)^{-3/2}(1+s)^{-3/2}\dd s\right)^{1/2}
 \le C\mathcal J_2^2(1+t)^{-3/4}.
\end{aligned}
\]
Adding the estimates of the horizontal flux $u\otimes u_h$,
the vertical flux $u_3u$, the convection remainder $u\Div u$,
the projected pressure remainder and the variable-viscosity
source gives
\[
 \left\|\left[\int_0^{t-1}S_h(t-s)\mathbb PF_2(s)\dd s\right]_3\right\|_2
                        \le C\mathcal J_2^2(1+t)^{-3/4}.
\]

\par\smallskip\noindent\emph{The interval $[t-1,t]$.}\enspace
On the last unit interval, we estimate the complete source
$F_2$ in $H^{-1}$. This avoids the singular time kernels used
for the separate source terms on $[0,t-1]$.
For $0\le a<t$ with $t-a\le1$, Fourier Cauchy--Schwarz gives
\begin{equation}\label{eq:new-third-short}
 \left\|\left[\int_a^t e^{\nu_h(t-s)\Delta_h}\mathbb PG(s)\dd s\right]_3\right\|_2^2
 \le C\int_a^t\|G(s)\|_{H^{-1}}^2\dd s.
\end{equation}
The kernel integral is
$\int_0^1e^{-\nu_h|\xi_h|^2s}(1+|\xi|^2)|\xi_h|^2/|\xi|^2\dd s\le C$.
Indeed,
\[
 (1+|\xi|^2)\frac{|\xi_h|^2}{|\xi|^2}\le1+|\xi_h|^2,
 \qquad
 \int_0^1e^{-\nu_h|\xi_h|^2s}\dd s
          \le\min\{1,(\nu_h|\xi_h|^2)^{-1}\},
\]
and the product of these bounds is uniform in $\xi$.
To apply \eqref{eq:new-third-short} with $G=F_2$ and $a=t-1$,
we estimate the time integral of $\|F_2\|_{H^{-1}}^2$.
Since $t\ge2$, $1+s\simeq1+t$ on $[t-1,t]$.

The high-frequency density estimate in Step~2, the total
dissipation bound \eqref{eq:new-base-energy}, and the
dissipation tail in Lemma~\ref{lem:new-seed} give, respectively,
\[
\begin{aligned}
 \sup_{t-1\le s\le t}\|q_H(s)\|_2^2
       &\le C\mathcal J_2^2(1+t)^{-3/2},
 \qquad \int_{t-1}^t\mathcal D_2(U(s))\dd s
       \le\int_0^\infty\mathcal D_2(U(s))\dd s
         \le C\mathcal J_2^2,\\
 \int_{t-1}^t \cD_0(U(s))\dd s
       &\le\int_{t-1}^\infty \cD_0(U(s))\dd s
         \le C\mathcal J_2^2(1+t)^{-1}.
\end{aligned}
\]
For the first term in the $F_2$ bound of
\eqref{eq:new-weakened-sources}, take $\|q_H(s)\|_2^2$ outside
the integral by its supremum:
\[
\begin{aligned}
 &\int_{t-1}^t\bigl(\mathcal J_2^2+\mathcal D_2(U(s))\bigr)
                            \|q_H(s)\|_2^2\dd s\\
 &\quad\le C\mathcal J_2^2(1+t)^{-3/2}
       \left(\mathcal J_2^2+\int_{t-1}^t\mathcal D_2(U(s))\dd s\right)
 \le C\mathcal J_2^4(1+t)^{-3/2}.
\end{aligned}
\]
The weighted $\cD_0(U)$ term uses the dissipation tail, while the
explicit time weight is integrated over an interval of length one:
\[
\begin{aligned}
 \mathcal J_2^2\int_{t-1}^t(1+s)^{-3/4}\cD_0(U(s))\dd s
   &\le C\mathcal J_2^2(1+t)^{-3/4}\int_{t-1}^t \cD_0(U(s))\dd s
 \le C\mathcal J_2^4(1+t)^{-7/4},\\
 \mathcal J_2^4\int_{t-1}^t(1+s)^{-25/16}\dd s
   &\le C\mathcal J_2^4(1+t)^{-25/16}.
\end{aligned}
\]
Since $7/4>3/2$ and $25/16>3/2$, summing the three contributions
in \eqref{eq:new-weakened-sources} gives
\[
 \int_{t-1}^t\|F_2(s)\|_{H^{-1}}^2\dd s
                       \le C\mathcal J_2^4(1+t)^{-3/2}.
\]
Using \eqref{eq:new-third-short} and taking the square root now yields
\[
 \left\|\left[\int_{t-1}^t
     e^{\nu_h(t-s)\Delta_h}\mathbb PF_2(s)\dd s\right]_3\right\|_2
                       \le C\mathcal J_2^2(1+t)^{-3/4}.
\]
The estimates on $[0,t-1]$ and $[t-1,t]$, together with
the homogeneous term, prove the $(1+t)^{-3/4}$ bound for
$v_3$. Adding the bound for $w_3$ established in Step~2 and
using $u_3=v_3+w_3$ proves the $u_3$ assertion of
\eqref{eq:new-zero}:
\[
 \|u_3(t)\|_2\le C\mathcal J_2(1+t)^{-3/4}.
\]

\smallskip\noindent\emph{Step 4. Preliminary gradient estimates.}
To estimate $\tau_H$ in Step~5, we need $H^{-1}$ bounds
for the quadratic source terms in \eqref{eq:F3}. These bounds
depend on the pointwise $L^2$ norms of
$(\nabla_hu,\Div u,\nabla\tau)$.
We first prove \eqref{eq:new-first-start}, using
\eqref{eq:new-seed}, \eqref{eq:new-zero-low}, the estimates
of $(q,w)$ in Step~2 and the bound for $u_3$ in Step~3.
A compensated high-frequency energy estimates $\nabla(q,w)$;
a second energy argument converts the time averages of
$(\nabla_hu,\Div u,\nabla\tau)$ into pointwise bounds.

\smallskip\noindent\emph{The gradient bound for $q,w$.}
Use the energy
\[
 K_1=\tfrac12(c^2\|q_H\|_{H^1}^2+\|w_H\|_{H^1}^2
             +A_\tau\|\tau_H\|_2^2)+\eta\ip{w_H}{\nabla q_H}
\]
with the parameters fixed in Step~2.
Then
\[
 K_1\simeq\|q_H\|_{H^1}^2+\|w_H\|_{H^1}^2+\|\tau_H\|_2^2.
\]
Test the three high-frequency equations with
$c^2(1-\Delta)q_H$, $(1-\Delta)w_H$ and $A_\tau\tau_H$,
and differentiate the compensating cross term.
As in Step~2, this supplies density dissipation, while the
viscous and thermal terms control
$\|w_H\|_{H^2}^2+\|\tau_H\|_{H^1}^2$.
The uncancelled thermal couplings satisfy
\[
\begin{aligned}
 c^2|\ip{\nabla\tau_H}{(1-\Delta)w_H}|
 &\le\varepsilon\|w_H\|_{H^2}^2
                             +C_\varepsilon\|\tau_H\|_{H^1}^2,\\
 A_\tau\alpha|\ip{\Div w_H}{\tau_H}|
 &\le\varepsilon\|w_H\|_{H^2}^2
                 +C_\varepsilon A_\tau^2R_0^{-4}\|\tau_H\|_{H^1}^2.
\end{aligned}
\]
The fixed choices of $A_\tau$ and $R_0$ absorb both
thermal remainders. The small cross term supplies the
dissipation of $\|q_H\|_{H^1}^2$.
The temperature has one fewer derivative in this energy;
its source is paired in $H^{-1}$--$H^1$.
For density transport, integration by parts gives
\[
 |\ip{u\cdot\nabla q_H}{q_H}_{H^1}|
       \le C\|\nabla u\|_\infty\|q_H\|_{H^1}^2.
\]
Only first derivatives of $q_H$ remain in this bound.
In the cross term, $F_1$ pairs with $\Div w_H$ in $L^2$;
the bound
$\|F_1\|_2\le C\mathcal J_2\|q_H\|_{H^1}
+C\mathcal J_2^2(1+s)^{-5/4}$ follows by splitting
$q=q_L+q_H$. Its first contribution is absorbed by smallness.

For $C_q=\Div P_L(uq)-P_L(q\Div u)-u\cdot\nabla q_L$,
the fixed cutoff bounds give
\[
\begin{aligned}
 \|C_q\|_{H^1}
 &\le C\|q\|_2(\|u\|_2+\cD_0(U)^{1/2})
                  +C\|u\|_{H^1}\|\nabla q_L\|_{W^{1,\infty}}\\
 &\le C\mathcal J_2^2(1+s)^{-1}
                  +C\mathcal J_2(1+s)^{-3/4}\cD_0(U)^{1/2}.
\end{aligned}
\]
Here $\|u\|_{H^1}\le C\mathcal J_2(1+s)^{-1/4}$ follows
by interpolation, and the low-frequency density bound is
\eqref{eq:new-zero-low}.
For the variable viscosity, the decomposition
\eqref{eq:new-coefficient-split} gives
$\nabla b_{\mathrm r}=-e^{-q}\nabla q_H-b_{\mathrm r}\nabla q_L$.
The mean-value bound $\|b_{\mathrm r}\|_2\le C\|q_H\|_2$ from
\eqref{eq:new-coefficient-split} and the bounds for
$q_L$ in \eqref{eq:new-zero-low} then yield
$\|b_{\mathrm r}\|_{H^1}+\|b_{\mathrm r}\|_6\le C\|q_H\|_{H^1}$ and
\[
\begin{aligned}
 \|b_{\mathrm r}\mathcal Lu\|_2&\le C\cD_2(U)^{1/2}\|q_H\|_{H^1},\\
 \|b_{\ell}\mathcal Lu\|_2&\le C\mathcal J_2(1+s)^{-3/4}\cD_2(U)^{1/2}.
\end{aligned}
\]
The first line uses $L^6\cdot L^3\subset L^2$ and
$\|\mathcal Lu\|_3\le C\cD_2(U)^{1/2}$.
The second uses $\|b_{\ell}\|_\infty\le C\mathcal J_2(1+s)^{-3/4}$.

For the quadratic source terms, $L^{6/5}\hookrightarrow H^{-1}$ and
interpolation give
\[
 \|(\nabla_hu,\Div u)\|_{12/5}^2
 \le\|(\nabla_hu,\Div u)\|_2^{3/2}
                  \|(\nabla_hu,\Div u)\|_6^{1/2}.
\]
The $L^6$ factor is bounded by $C\mathcal J_2$, so
\[
 \|Q(u)\|_{H^{-1}}
       \le C\mathcal J_2^{1/2}\|(\nabla_hu,\Div u)\|_2^{3/2}.
\]
Applying $L^{6/5}\hookrightarrow H^{-1}$ and the
$L^{12/5}$ interpolation inequality to $\nabla\tau$ also
gives $\||\nabla\tau|^2\|_{H^{-1}}
\le C\mathcal J_2^{1/2}\|\nabla\tau\|_2^{3/2}$.
By \eqref{eq:new-seed-ledger}, $\cD_0(U)^{1/2}\le C\mathcal J_2(1+s)^{-1/4}$,
so the squares of both quadratic source estimates are at most
$C\mathcal J_2\cD_0(U)^{3/2}\le C\mathcal J_2^2(1+s)^{-1/4}\cD_0(U)$.
The smooth exponential coefficients preserve the $H^{-1}$
bounds by \eqref{eq:H2-H1-multiplier} and duality.

The squared $L^2$ norm of velocity convection also contributes
the weight $(1+s)^{-1/4}\cD_0(U)$. Indeed,
\eqref{eq:verticalcomp} and \eqref{eq:new-seed-ledger} give
\[
 \|u\cdot\nabla u\|_2
 \le\|u_h\|_\infty\|\nabla_hu\|_2
                   +\|u_3\|_6\|\partial_3u\|_3
 \le C\mathcal J_2(1+s)^{-1/8}\cD_0(U)^{1/2}.
\]
Scalar transport is bounded in $H^{-1}$ by
$\|u\|_\infty\|\nabla\tau\|_2$.
In pressure, the factor $e^\tau-1$ multiplying
$\nabla q_H$ is small and its contribution is absorbed;
the terms with $\nabla q_L$ and $\nabla\tau$ use
\eqref{eq:new-zero-low} and \eqref{eq:new-seed-ledger}.
For thermal diffusion, the $b_{\mathrm r}$ term is bounded by
$C\cD_2(U)^{1/2}\|q_H\|_{H^1}$ in $H^{-1}$,
and the $b_{\ell}$ term is integrated by parts once.
Finally, Lemma~\ref{lem:H2-anisotropic} gives
$\|\nabla u\|_\infty\le C\mathcal J_2+C\cD_2(U)^{1/2}$.
Young's inequality bounds the transport commutators by
a small part of the dissipation plus $C\cD_2(U)K_1$.
Combining these bounds yields
\[
 \begin{aligned}
 K_1'+c_*(K_1+\|w_H\|_{H^2}^2+\|\tau_H\|_{H^1}^2)
 &\le C\cD_2(U)K_1+C\mathcal J_2^4(1+s)^{-2}\\
 &\quad+C\mathcal J_2^2(1+s)^{-1/4}\cD_0(U)
       +C\mathcal J_2^2(1+s)^{-3/2}\cD_2(U).
 \end{aligned}
\]
The three source terms on the right come, respectively,
from the explicit cutoff and low-frequency products, the
transport and quadratic source bounds, and $b_{\ell}\mathcal Lu$.

The integrating factor is bounded by \eqref{eq:new-base-energy}.
Lemma~\ref{lem:new-seed} supplies the unit-interval bound
$\int_a^{a+1}\cD_0(U(s))\dd s\le C\mathcal J_2^2(1+a)^{-1}$,
while \eqref{eq:new-base-energy} supplies the total mass of
$\cD_2(U)$. Applying \eqref{eq:new-window-convolution} with
$(p,\beta)=(1/4,1)$ and $(3/2,0)$, respectively, gives decay
powers $5/4$ and $3/2$ for the weighted dissipation terms.
The explicit forcing has power $2$. Hence
\[
 K_1(t)\le C\mathcal J_2^2e^{-c_*t}
             +C\mathcal J_2^4(1+t)^{-5/4}
       \le C\mathcal J_2^2(1+t)^{-5/4}.
\]
Taking square roots controls $\nabla(q_H,w_H)$;
adding the differentiated low-frequency bound
\eqref{eq:new-zero-low} gives
$\|\nabla(q,w)(t)\|_2\le C\mathcal J_2(1+t)^{-5/8}$.

\smallskip\noindent\emph{The remaining gradient bounds.}
To estimate $(\nabla_hu,\Div u,\nabla\tau)$, use
\begin{equation}
E_{\rm str}=\tfrac12\bigl(c^2\ip q{\mathcal Mq}
                   -\ip u{\mathcal Lu}+\omega\ip\tau{\mathcal M\tau}\bigr),
 \qquad \mathcal M=-\nu_h\Delta_h-\nu_c\Delta.
\end{equation}
Test the three equations in \eqref{eq:perturbed} against
$c^2\mathcal Mq$, $-\mathcal Lu$ and $\omega\mathcal M\tau$.
The identities
$-\mathcal L\nabla=\nabla\mathcal M$,
$-\Div\mathcal L=\mathcal M\Div$ and $\omega\alpha=c^2$
cancel the density--velocity and temperature--velocity couplings.
Integration by parts also gives
\[
 -\ip u{\mathcal Lu}
      =\nu_h\|\nabla_hu\|_2^2+\nu_c\|\Div u\|_2^2,
 \qquad
 \ip q{\mathcal Mq}
      =\nu_h\|\nabla_hq\|_2^2+\nu_c\|\nabla q\|_2^2,
\]
and the analogous identity for $\tau$.
The linear dissipative terms are
$\|\mathcal Lu\|_2^2+\omega\chi\ip{-\Delta\tau}{\mathcal M\tau}$;
the second controls $\|\nabla^2\tau\|_2^2$.

For scalar transport, integration by parts bounds the
pairing by $C\|\nabla u\|_\infty$ times the scalar gradient energy.
The convection estimate in the $K_1$ argument of this step
implies
$\|u\cdot\nabla u\|_2\le C\mathcal J_2\|(\nabla_hu,\Div u)\|_2$.
Pairing the convection term with $\mathcal Lu$ and applying Young's
inequality controls it by $E_{\rm str}$ and its dissipation.
The pressure coefficients are
bounded by $C\mathcal J_2$ in $L^\infty$.
The velocity diffusion remainder pairs exactly as
$-\int(e^{-q}-1)|\mathcal Lu|^2\dd x$, so smallness absorbs it
without differentiating the coefficient.
The thermal diffusion remainder is bounded by
$C\mathcal J_2\|\nabla^2\tau\|_2^2$ and is absorbed as well.
For the quadratic source terms, use
\[
\begin{aligned}
 \|Q(u)\|_2^2+\||\nabla\tau|^2\|_2^2
 &\le C\bigl(\|(\nabla_hu,\Div u)\|_\infty^2
                         +\|\nabla\tau\|_\infty^2\bigr)\\
 &\quad\times\bigl(\|(\nabla_hu,\Div u)\|_2^2+\|\nabla\tau\|_2^2\bigr)
 \le C\cD_2(U)E_{\rm str}.
\end{aligned}
\]
The smooth coefficients in \eqref{eq:F3} are bounded in
$L^\infty$. Combining the pairings gives
\begin{equation}\label{eq:new-structural-differential}
\begin{gathered}
 E_{\rm str}'+c_*(\|\mathcal Lu\|_2^2+\|\nabla^2\tau\|_2^2)
       \le C(1+\cD_2(U))E_{\rm str},
 \\ E_{\rm str}\simeq\|\nabla q\|_2^2+\|(\nabla_hu,\Div u)\|_2^2+\|\nabla\tau\|_2^2.
\end{gathered}
\end{equation}

To obtain pointwise decay from \eqref{eq:new-structural-differential},
combine the bound for $\nabla q$ obtained from $K_1$ with the
dissipation tail in Lemma~\ref{lem:new-seed}:
\[
\begin{aligned}
 \int_{t-1}^tE_{\rm str}(s)\dd s
 &\le C\int_{t-1}^t\|\nabla q(s)\|_2^2\dd s
                              +C\int_{t-1}^t\cD_0(U(s))\dd s\\
 &\le C\mathcal J_2^2(1+t)^{-5/4}
                              +C\mathcal J_2^2(1+t)^{-1}
 \le C\mathcal J_2^2(1+t)^{-1}.
\end{aligned}
\]
For every $a\in[t-1,t]$, Gronwall's inequality gives
\[
 E_{\rm str}(t)
 \le E_{\rm str}(a)\exp\left(C(t-a)+C\int_a^t\cD_2(U(s))\dd s\right)
 \le CE_{\rm str}(a),
\]
where \eqref{eq:new-base-energy} bounds the exponential
uniformly. Integrating this inequality with respect to $a$
over $[t-1,t]$ converts the average into the pointwise bound
$E_{\rm str}(t)\le C\mathcal J_2^2(1+t)^{-1}$.
Combining the pointwise bound for $E_{\rm str}$ with the
$(1+t)^{-5/8}$ estimate of $\nabla(q,w)$ obtained from $K_1$
proves \eqref{eq:new-first-start}.

Integrating \eqref{eq:new-structural-differential} over
$[t,t+1]$ and discarding $E_{\rm str}(t+1)\ge0$ also gives
\[
\begin{aligned}
 \int_t^{t+1}\bigl(\|\mathcal Lu(s)\|_2^2
                         +\|\nabla^2\tau(s)\|_2^2\bigr)\dd s
 &\le CE_{\rm str}(t)
        +C\int_t^{t+1}(1+\cD_2(U(s)))E_{\rm str}(s)\dd s\\
 &\le C\mathcal J_2^2(1+t)^{-1}.
\end{aligned}
\]

\smallskip\noindent\emph{Step 5. High-frequency temperature and completion.}
The low-frequency bound \eqref{eq:new-zero-low} already
controls $\tau_L$. To complete the $W$ assertion of
\eqref{eq:new-zero}, we now use
\eqref{eq:new-first-start} to estimate $\tau_H$ in $L^2$.
Testing the high-frequency heat equation with $\tau_H$ gives
\[
 \frac12\frac{\dd}{\dd s}\|\tau_H\|_2^2
      +\chi\|\nabla\tau_H\|_2^2
   =-\alpha\ip{\Div w_H}{\tau_H}+\ip{P_HF_3}{\tau_H}.
\]
Since $\tau_H$ has only high frequencies,
$\|\tau_H\|_{H^1}\le C\|\nabla\tau_H\|_2$.
Thus each $H^{-1}$ source pairing is bounded by a small
part of the thermal dissipation plus the square of its
$H^{-1}$ norm. We estimate the terms of $F_3$ separately.

For the linear coupling, integration by parts and Step~2 give
\[
 \alpha|\ip{\Div w_H}{\tau_H}|
 \le\varepsilon\|\nabla\tau_H\|_2^2+C_\varepsilon\|w_H\|_2^2
 \le\varepsilon\|\nabla\tau_H\|_2^2
                      +C_\varepsilon\mathcal J_2^2(1+s)^{-3/2}.
\]
For transport, use
$-u\cdot\nabla\tau=-\Div(u\tau)+\tau\Div u$.
The preliminary bounds in Lemma~\ref{lem:new-seed} give
$\|u\|_\infty\|\tau\|_2
\le C\mathcal J_2^2(1+s)^{-1/8-5/8}$.
Moreover, \eqref{eq:new-first-start} and interpolation give
$\|\tau\|_3\le C\|\tau\|_2^{1/2}\|\nabla\tau\|_2^{1/2}
\le C\mathcal J_2(1+s)^{-9/16}$.
Consequently
\[
\begin{aligned}
 \|u\cdot\nabla\tau\|_{H^{-1}}
 &\le C\|u\tau\|_2+C\|\tau\Div u\|_{6/5}
 \le C\mathcal J_2^2(1+s)^{-3/4}
                        +C\mathcal J_2(1+s)^{-9/16}\cD_0(U)^{1/2},\\
 \|u\cdot\nabla\tau\|_{H^{-1}}^2
 &\le C\mathcal J_2^4(1+s)^{-3/2}
                        +C\mathcal J_2^2(1+s)^{-1/2}\cD_0(U).
\end{aligned}
\]
For variable thermal diffusion, retain the decomposition
\eqref{eq:new-coefficient-split}. Duality with
$H^1\hookrightarrow L^6$ gives
\[
 \|b_{\mathrm r}\Delta\tau\|_{H^{-1}}
 \le C\|b_{\mathrm r}\|_2\|\Delta\tau\|_3
 \le C\cD_2(U)^{1/2}\|q_H\|_2.
\]
Step~2 therefore bounds its squared norm by
$C\mathcal J_2^2(1+s)^{-3/2}\cD_2(U)$.
For $b_{\ell}\Delta\tau$, integrate by parts against
$g\in H^1$:
\[
 \ip{b_{\ell}\Delta\tau}{g}
      =-\int b_{\ell}\nabla\tau\cdot\nabla g\dd x
                       -\int g\nabla b_{\ell}\cdot\nabla\tau\dd x.
\]
The chain rule applied to \eqref{eq:new-zero-low} bounds
$\|b_{\ell}\|_\infty+\|\nabla b_{\ell}\|_3$. Combining this bound with
\eqref{eq:new-first-start} yields
\[
 \|b_{\ell}\Delta\tau\|_{H^{-1}}
 \le C(\|b_{\ell}\|_\infty+\|\nabla b_{\ell}\|_3)\|\nabla\tau\|_2
 \le C\mathcal J_2^2(1+s)^{-5/4}.
\]
Its squared norm is therefore bounded by
$C\mathcal J_2^4(1+s)^{-5/2}$.

For $e^{-q-\tau}Q(u)$ and $e^{-q}|\nabla\tau|^2$,
use the $H^{-1}$ bounds obtained by $L^{12/5}$ interpolation
in Step~4 and the multiplication estimate
\eqref{eq:H2-H1-multiplier}:
\[
\begin{aligned}
 &\|e^{-q-\tau}Q(u)\|_{H^{-1}}^2
                 +\|e^{-q}|\nabla\tau|^2\|_{H^{-1}}^2
 \le C\mathcal J_2
       \bigl(\|(\nabla_hu,\Div u)\|_2^3+\|\nabla\tau\|_2^3\bigr)\\
 &\quad\le C\mathcal J_2^2(1+s)^{-1/2}
       \bigl(\|(\nabla_hu,\Div u)\|_2^2+\|\nabla\tau\|_2^2\bigr)
 \le C\mathcal J_2^2(1+s)^{-1/2}\cD_0(U).
\end{aligned}
\]
The passage from cubic to quadratic gradient norms uses
\eqref{eq:new-first-start}.

Absorb the small thermal dissipation terms.
The high-frequency inequality
$\|\nabla\tau_H\|_2^2\ge R_0^2\|\tau_H\|_2^2$
then supplies a positive damping term for $\|\tau_H\|_2^2$.
The estimates for transport, thermal diffusion, and the quadratic
source terms, together with $\mathcal J_2\le1$, give
\[
\begin{aligned}
 &\frac{\dd}{\dd s}\|\tau_H\|_2^2
       +c_*(\|\tau_H\|_2^2+\|\nabla\tau_H\|_2^2)\\
 &\quad\le C\mathcal J_2^2(1+s)^{-3/2}(1+\cD_2(U))
      +C\mathcal J_2^4(1+s)^{-5/2}+C\mathcal J_2^2(1+s)^{-1/2}\cD_0(U).
\end{aligned}
\]
Apply \eqref{eq:new-window-convolution} using
$\int_0^\infty\cD_2(U(s))\dd s\le C\mathcal J_2^2$ from
\eqref{eq:new-base-energy} and
$\int_a^{a+1}\cD_0(U(s))\dd s\le C\mathcal J_2^2(1+a)^{-1}$
from Lemma~\ref{lem:new-seed}.
The choices $(p,\beta)=(3/2,0)$ and $(1/2,1)$ give power
$3/2$ for both weighted dissipation terms; the explicit
forcing has powers $3/2$ and $5/2$.
Solving the scalar differential inequality and using
$\mathcal J_2\le1$ therefore gives
\[
 \|\tau_H(t)\|_2^2\le C\mathcal J_2^2e^{-c_*t}
                  +C\mathcal J_2^2(1+t)^{-3/2}
             \le C\mathcal J_2^2(1+t)^{-3/2}.
\]
Taking square roots and adding the low-frequency bound
\eqref{eq:new-zero-low} proves
$\|\tau(t)\|_2\le C\mathcal J_2(1+t)^{-3/4}$.
Together with the bounds for $q,w$ in Step~2 and $u_3$ in
Step~3, this completes \eqref{eq:new-zero}.
\end{proof}

\paragraph{First-order estimates.}

We now improve the first-order decay estimates in
\eqref{eq:new-first-start}. Lemma~\ref{lem:interpolation-dissipation}
provides the interpolation, velocity derivative, and quadratic
estimates used in the proof below.

\begin{lemma}[Anisotropic interpolation and dissipation estimates]
\label{lem:interpolation-dissipation}
Let $g,\tau\in H^2(\R^3)$ and $u\in H^2(\R^3;\R^3)$.
Write $v=\mathbb Pu$ and $w=\mathbb Qu$. The following estimates
hold, with $C$ depending only on the fixed coefficients.

\smallskip\noindent\textup{(i)} The interpolation inequalities are
\begin{equation}\label{eq:new-agmon}
\begin{aligned}
 \|g\|_\infty^2&\le C\|\nabla_hg\|_2\|\nabla^2g\|_2,
 \qquad \|g\|_\infty^2\le C\|\nabla g\|_2\|\nabla_h\nabla g\|_2,\\
 \|g\|_\infty&\le C\|g\|_2^{1/4}\|\nabla_h^2g\|_2^{1/2}
                                      \|\partial_3^2g\|_2^{1/4}.
\end{aligned}
\end{equation}
They also hold componentwise for vector-valued functions.

\smallskip\noindent\textup{(ii)} The velocity satisfies
\begin{equation}\label{eq:new-elliptic}
 \|\nabla_h^2u\|_2+\|\nabla^2w\|_2+
 \|\nabla_h\nabla u_3\|_2\le C\|\mathcal Lu\|_2,
 \qquad \|u_3\|_\infty^2\le C\|(\nabla_hu,\Div u)\|_2\|\mathcal Lu\|_2.
\end{equation}

\smallskip\noindent\textup{(iii)} For the quadratic terms
$Q(u)=\mu|\nabla_hu|^2+(\lambda+\mu)(\Div u)^2$ and
$|\nabla\tau|^2$, we have
\begin{equation}\label{eq:new-production-L2}
 \|Q(u)\|_2\le C(\|u\|_\infty+\|w\|_\infty)\|\mathcal Lu\|_2,
 \qquad \||\nabla\tau|^2\|_2\le C\|\tau\|_\infty \|\nabla^2\tau\|_2.
\end{equation}
\end{lemma}

\begin{proof}
\smallskip\noindent\textup{(i)} For every auxiliary parameter
$\varepsilon>0$, integration in spherical coordinates for the first
weight, and in $\xi_3$ first for the second, gives
\[
 \int_{\R^3}\frac{\dd\xi}{\varepsilon|\xi_h|^2+|\xi|^4}
       \le C\varepsilon^{-1/2},\qquad
 \int_{\R^3}\frac{\dd\xi}{|\xi|^2(\varepsilon+|\xi_h|^2)}
       \le C\varepsilon^{-1/2}.
\]
Fourier inversion and Cauchy--Schwarz therefore yield
\[
\|g\|_\infty^2
 \le C\varepsilon^{-1/2}
       \bigl(\varepsilon\|\nabla_hg\|_2^2+\|\nabla^2g\|_2^2\bigr),
 \qquad \|g\|_\infty^2
 \le C\varepsilon^{-1/2}
       \bigl(\varepsilon\|\nabla g\|_2^2+\|\nabla_h\nabla g\|_2^2\bigr).
\]
Choose $\varepsilon=\|\nabla^2g\|_2^2/\|\nabla_hg\|_2^2$ in the first
inequality and
$\varepsilon=\|\nabla_h\nabla g\|_2^2/\|\nabla g\|_2^2$ in the second.
These choices prove the first two inequalities in \eqref{eq:new-agmon}.

For the third inequality, scaling the horizontal and vertical
Fourier variables separately gives, for $a,b,c>0$,
\[
 \int_{\R^3}\frac{\dd\xi}{a+b|\xi_h|^4+c\xi_3^4}
       =C_0a^{-1/4}b^{-1/2}c^{-1/4}.
\]
Here $C_0<\infty$, since integration in $\xi_h$ reduces the
integral with $a=b=c=1$ to a constant multiple of
$\int_{\R}(1+\xi_3^4)^{-1/2}\,\dd\xi_3$.
Fourier inversion and Cauchy--Schwarz yield
\[
 \|g\|_\infty^2\le Ca^{-1/4}b^{-1/2}c^{-1/4}
 \bigl(a\|g\|_2^2+b\|\nabla_h^2g\|_2^2
                              +c\|\partial_3^2g\|_2^2\bigr).
\]
For $g\ne0$, choose
\[
 a=1,\qquad b=\frac{\|g\|_2^2}{\|\nabla_h^2g\|_2^2},\qquad
 c=\frac{\|g\|_2^2}{\|\partial_3^2g\|_2^2}.
\]
Taking square roots proves the third inequality in
\eqref{eq:new-agmon}. The case $g=0$ is immediate.

\smallskip\noindent\textup{(ii)} The decomposition $u=v+w$ gives
\[
 \widehat{-\mathcal Lu}(\xi)
 =\nu_h|\xi_h|^2\widehat v(\xi)
  +\bigl(\nu_h|\xi_h|^2+\nu_c|\xi|^2\bigr)\widehat w(\xi).
\]
Since $\widehat v$ and $\widehat w$ are orthogonal, Plancherel's
identity controls $\|\nabla_h^2u\|_2+\|\nabla^2w\|_2$ by
$C\|\mathcal Lu\|_2$. For the third component, use
$|\widehat v_3|\le |\xi_h||\widehat v|/|\xi|$ to obtain
\[
 \|\nabla_h\nabla v_3\|_2\le C\|\nabla_h^2v\|_2,
 \qquad
 \|\nabla_h\nabla w_3\|_2\le C\|\nabla^2w\|_2.
\]
Thus $\|\nabla_h\nabla u_3\|_2\le C\|\mathcal Lu\|_2$,
which proves the first inequality in \eqref{eq:new-elliptic}.
The second follows by applying the second inequality in
\eqref{eq:new-agmon} to $u_3$ and using \eqref{eq:verticalcomp}
to bound $\|\nabla u_3\|_2$.
These derivative bounds do not control $\partial_3^2v_h$.

\smallskip\noindent\textup{(iii)} Integration by parts in the
horizontal variables for $u$, and in all variables for $w,\tau$,
gives
\[
\begin{aligned}
 \|\nabla_hu\|_4^2
   &\le C\|u\|_\infty\|\nabla_h^2u\|_2,\\
 \|\Div u\|_4^2=\|\Div w\|_4^2
   &\le C\|w\|_\infty\|\nabla^2w\|_2,
 \qquad \|\nabla\tau\|_4^2
   \le C\|\tau\|_\infty\|\nabla^2\tau\|_2.
\end{aligned}
\]
The integration-by-parts inequalities extend from smooth,
compactly supported functions by approximation in $H^2$.
Combining them with \eqref{eq:new-elliptic} proves
\eqref{eq:new-production-L2}.
\end{proof}

\begin{lemma}[First-order decay]\label{lem:new-first}
The solution satisfies
\begin{align}\label{eq:new-first}
 \|\nabla W\|_2&\le C\mathcal J_2 (1+t)^{-5/4},&
 \|\nabla_hu\|_2&\le C\mathcal J_2 (1+t)^{-1},\notag\\
 \|\nabla_hu_3\|_2&\le C\mathcal J_2 (1+t)^{-5/4},&
 \|\partial_3u_3\|_2&\le C\mathcal J_2 (1+t)^{-1}.
\end{align}
Moreover,
\begin{equation}\label{eq:new-first-windows}
 \int_t^{t+1}\|\mathcal Lu\|_2^2\dd s\le C\mathcal J_2^2(1+t)^{-2},\qquad
 \int_t^{t+1}(\|\nabla^2w\|_2^2+\|\nabla^2\tau\|_2^2)\dd s
                         \le C\mathcal J_2^2(1+t)^{-5/2}.
\end{equation}
\end{lemma}
\begin{proof}
Starting from \eqref{eq:new-first-start}, Steps~1--3 improve
the common decay exponent for $\nabla W$ and $\nabla_hu$ through
\begin{equation}\label{eq:new-first-iteration}
 \frac12\ \longrightarrow\ \frac34\ \longrightarrow\ \frac78\ \longrightarrow\ 1.
\end{equation}
Step~4 then gives the $(1+t)^{-5/4}$ decay of $\nabla W$ and
\eqref{eq:new-first-windows}; Step~5 estimates the first
derivatives of $u_3$ using $u_3=v_3+w_3$.

\smallskip\noindent\emph{Step 1. Estimates of the nonlinear terms.}
Suppose that for some $1/2\le\sigma\le1$ we have already proved
\begin{equation}\label{eq:iteration-hypothesis}
 \|\nabla W\|_2+\|\nabla_hu\|_2
                 \le C_\sigma\mathcal J_2 (1+t)^{-\sigma}
\end{equation}
The starting case $\sigma=1/2$ follows from
\eqref{eq:new-first-start}. Since $\Div u=\Div w$, the structural
energy in \eqref{eq:new-structural-differential} is bounded by
$C_\sigma\mathcal J_2^2(1+t)^{-2\sigma}$.
Integrating \eqref{eq:new-structural-differential} over
$[t,t+1]$ therefore gives a unit-interval bound for
$(\mathcal Lu,\nabla^2\tau)$. Applying \eqref{eq:new-agmon}
to $u$ and each component of $W$, with
\eqref{eq:new-base-energy}, gives the $L^\infty$ bound.
Thus
\begin{equation}\label{eq:new-sigma-window}
 \int_t^{t+1}(\|\mathcal Lu\|_2^2+\|\nabla^2\tau\|_2^2)\dd s\le C_\sigma\mathcal J_2^2(1+t)^{-2\sigma},
 \qquad \|(u,W)\|_\infty\le C_\sigma\mathcal J_2 (1+t)^{-\sigma/2}.
\end{equation}
As in \eqref{eq:new-seed-ledger}, the $L^2$ decay
\eqref{eq:new-seed} and the uniform $H^2$ bound
\eqref{eq:new-base-energy} give, by interpolation,
\[
\begin{aligned}
 \|\partial_3u\|_2+\|u\|_{H^1}
 &\le C\|u\|_2^{1/2}\|u\|_{H^2}^{1/2}
       \le C\mathcal J_2(1+s)^{-1/4},\\
 \|u\|_{L^\infty_3L^2_h}
 &\le C\|u\|_2^{1/2}\|\partial_3u\|_2^{1/2}
       \le C\mathcal J_2(1+s)^{-3/8}.
\end{aligned}
\]
The mixed-norm bound uses the one-dimensional Sobolev
inequality in the vertical variable. We now estimate the
exponential convolutions of the velocity source $F_2$ and
temperature source $F_3$ from \eqref{eq:F2} and \eqref{eq:F3}:
\begin{equation}\label{eq:new-sigma-sources}
\begin{gathered}
 \int_0^t e^{-c_*(t-s)}\|F_2(s)\|_2^2\dd s\le C\mathcal J_2^4(1+t)^{-2\sigma-1/2},
 \\ \int_0^t e^{-c_*(t-s)}\|F_3(s)\|_2^2\dd s\le C\mathcal J_2^4(1+t)^{-3\sigma}.
\end{gathered}
\end{equation}
Horizontal convection, pressure, and temperature convection
have squared $L^2$ norms bounded by
$C_\sigma\mathcal J_2^4(1+s)^{-3\sigma}$, using the
$L^\infty$ bound in \eqref{eq:new-sigma-window}.
For vertical velocity convection, \eqref{eq:new-elliptic} gives
\[
\|u_3\partial_3u\|_2^2
 \le C\|u_3\|_\infty^2\|\partial_3u\|_2^2
 \le C_\sigma\mathcal J_2^3
                  (1+s)^{-\sigma-1/2}\|\mathcal Lu\|_2.
\]
The unit-interval bound in \eqref{eq:new-sigma-window} and
\eqref{eq:new-window-convolution} imply
\[
 \int_0^t e^{-c_*(t-s)}\|\mathcal Lu(s)\|_2^2\dd s
       \le C_\sigma\mathcal J_2^2(1+t)^{-2\sigma}.
\]
Consequently, Cauchy--Schwarz in time yields
\[
\begin{aligned}
 &\int_0^t e^{-c_*(t-s)}\|u_3\partial_3u\|_2^2\dd s\\
 &\quad\le C_\sigma\mathcal J_2^3
    \left(\int_0^t e^{-c_*(t-s)}(1+s)^{-2\sigma-1}\dd s\right)^{1/2}
    \left(\int_0^t e^{-c_*(t-s)}\|\mathcal Lu\|_2^2\dd s\right)^{1/2}\\
 &\quad\le C_\sigma\mathcal J_2^4(1+t)^{-2\sigma-1/2}.
\end{aligned}
\]
For variable diffusion and the quadratic source terms,
\eqref{eq:new-production-L2} bounds the squared norms by
\[
 C_\sigma\mathcal J_2^2(1+s)^{-\sigma}
             \bigl(\|\mathcal Lu\|_2^2+\|\nabla^2\tau\|_2^2\bigr).
\]
Their exponential convolutions are therefore bounded by
$C_\sigma\mathcal J_2^4(1+t)^{-3\sigma}$, again by
\eqref{eq:new-sigma-window} and \eqref{eq:new-window-convolution}.
Since $3\sigma\ge2\sigma+1/2$ for $\sigma\ge1/2$, vertical
convection determines the decay exponent $2\sigma+1/2$ for $F_2$;
the exponent for $F_3$ is $3\sigma$.
This proves both estimates in \eqref{eq:new-sigma-sources}.

\smallskip\noindent\emph{Step 2. Preliminary decay estimates for $\nabla W$.}
Under \eqref{eq:iteration-hypothesis}, we improve the decay
exponent for $\nabla W$ to $\min\{1,\sigma+1/4\}$ by estimating
$P_LW$ and $P_HW$ separately.

For low frequencies, use the decomposition
\eqref{eq:new-decay-flux}--\eqref{eq:new-flux-remainders}.
By \eqref{eq:iteration-hypothesis} and Lemma~\ref{lem:new-zero},
\[
 \|F^{(0)}\|_1\le C_\sigma\mathcal J_2^2(1+s)^{-\sigma-1/2},
 \qquad
 \sum_j\|F^{(j)}\|_1\le C_\sigma\mathcal J_2^2(1+s)^{-1}.
\]
After one spatial derivative, the low-frequency kernels have
powers $5/4$ for $F^{(0)}$ and $7/4$ for a flux.
Both kernels are integrable in time, and $\sigma+1/2\ge1$.
Splitting the convolutions at $t/2$ therefore gives
\begin{equation}\label{eq:new-low-first-rough}
 \|\nabla P_LW\|_2+\sum_{j=1}^2\|\nabla^jq_L\|_\infty
                                      \le C\mathcal J_2 (1+t)^{-1}.
\end{equation}
The $L^\infty$ bounds follow from the fixed frequency support
and the $L^2$ bound for $\nabla P_LW$.

For the high-frequency part $W_H$, use the compensated
$H^1$ energy
\[
 E_{1,H}=\tfrac12(c^2\|q_H\|_{H^1}^2+\|w_H\|_{H^1}^2
                 +\omega\|\tau_H\|_{H^1}^2)
                 +\eta\ip{w_H}{\nabla q_H}
           \simeq\|W_H\|_{H^1}^2,
\]
where $\eta>0$ is fixed and sufficiently small.
Keep the density transport in \eqref{eq:new-density-high}
on the left. The cutoff commutator satisfies
\[
 \|C_q\|_{H^1}\le C_\sigma\mathcal J_2^2
          \bigl((1+s)^{-\sigma-1/2}+(1+s)^{-5/4}\bigr),
\]
and the differentiated transport terms contribute at most
$(C\mathcal J_2+\eta_0+C_{\eta_0}\cD_2(U))E_{1,H}$.
In the velocity and temperature tests, integrate by parts to
place one derivative on the dissipative test function; their
source terms then require only the $L^2$ norms of $F_2,F_3$.
In the compensation term, the high-frequency density
transport $u\cdot\nabla q_H$ has $L^2$ norm at most
$C\mathcal J_2\|q_H\|_{H^1}$. Its low-frequency counterpart
$u\cdot\nabla q_L$ is bounded by
$C\mathcal J_2^2(1+s)^{-3/2}$, using \eqref{eq:new-seed}
and \eqref{eq:new-low-first-rough}.
Absorbing the small terms yields
\begin{equation}\label{eq:new-E1-differential}
\begin{aligned}
 E_{1,H}'+c_*(E_{1,H}+\|(w_H,\tau_H)\|_{H^2}^2)
 &\le C\cD_2(U)E_{1,H}\\
 &\quad+C(\|F_2\|_2^2+\|F_3\|_2^2+\|C_q\|_{H^1}^2)
       +C\mathcal J_2^4(1+s)^{-3}.
\end{aligned}
\end{equation}
The integrating factor is bounded by \eqref{eq:new-base-energy}.
Using \eqref{eq:new-sigma-sources}, the bound for $C_q$, and
$1/2\le\sigma\le1$, we obtain
\[
 E_{1,H}(t)\le C\mathcal J_2^2e^{-c_*t}
       +C_\sigma\mathcal J_2^4(1+t)^{-2\sigma-1/2}.
\]
The squared cutoff commutator $\|C_q\|_{H^1}^2$ has powers
$2\sigma+1$ and $5/2$, both at least $2\sigma+1/2$.
Taking square roots and adding \eqref{eq:new-low-first-rough} gives
\[
 \|W_H(t)\|_{H^1}\le C_\sigma\mathcal J_2(1+t)^{-\sigma-1/4},
 \qquad
 \|\nabla W(t)\|_2\le C_\sigma\mathcal J_2
                       (1+t)^{-\min\{1,\sigma+1/4\}}.
\]

\smallskip\noindent\emph{Step 3. Horizontal derivatives of $v$.}
To complete each iteration, we must also improve the bound
for $\nabla_hu$. First estimate the solenoidal part
$v=\mathbb Pu$. At high horizontal frequencies,
\eqref{eq:new-maximal-horizontal} with $k=1$ and $G=F_2$,
together with \eqref{eq:new-sigma-sources} and the homogeneous
heat estimate, gives
\[
 \|\nabla_hP_H^hv(t)\|_2^2
       \le C\mathcal J_2^2e^{-c_*t}
          +C_\sigma\mathcal J_2^4(1+t)^{-2\sigma-1/2}.
\]
For the low horizontal frequencies, use the decomposition
\eqref{eq:new-projected-source}. The spatially differentiated kernels
have powers $1$ for $N$ and $3/2$ for the horizontal flux, so
\[
\begin{aligned}
 \|\nabla_hP_L^hv(t)\|_2
 &\le C\mathcal J_2(1+t)^{-1}
     +C\int_0^t(1+t-s)^{-1}\|N(s)\|_{L^2_3L^1_h}\dd s\\
 &\quad+C\int_0^t(1+t-s)^{-3/2}
          \|\nu_h(e^{-q}-1)\nabla_hu\|_{L^2_3L^1_h}\dd s.
\end{aligned}
\]
Horizontal H\"older and vertical Sobolev give
\begin{align}\label{eq:new-sigma-mixed}
 \|u_h\cdot\nabla_hu\|_{L^2_3L^1_h}
      &\le C\mathcal J_2^2(1+s)^{-\sigma-3/8},\notag\\
 \|u_3\partial_3u\|_{L^2_3L^1_h}
      &\le C\mathcal J_2^2(1+s)^{-5/8-\sigma/2},\notag\\
 \|\Pi_v\|_{L^2_3L^1_h}
      &\le C\mathcal J_2^2(1+s)^{-3\sigma/2},\notag\\
 \|\Pi_p\|_{L^2_3L^1_h}
   +\|\nu_h(e^{-q}-1)\nabla_hu\|_{L^2_3L^1_h}
      &\le C\mathcal J_2^2(1+s)^{-3/8-3\sigma/2}.
\end{align}
Here we use
\[
 \|ab\|_{L^2_3L^1_h}
       \le\|a\|_{L^\infty_3L^2_h}\|b\|_2,
\]
together with
$\|u_3\|_{L^\infty_3L^2_h}
 \le C_\sigma\mathcal J_2(1+s)^{-3/8-\sigma/2}$ and
$\|\nabla(e^{-q}-1)\|_{L^\infty_3L^2_h}
 \le C_\sigma\mathcal J_2(1+s)^{-\sigma/2}$.
We use \eqref{eq:new-sigma-mixed} in the three cases
$\sigma=1/2,3/4,7/8$. After each estimate of $\nabla_hv$,
the identity $\nabla_hu=\nabla_hv+\nabla_hw$ and the bound
for $\nabla W$ from Step~2 verify
\eqref{eq:iteration-hypothesis} with the next exponent.

\smallskip\noindent\emph{The case $\sigma=1/2$.}
We retain $\cD_0(U)^{1/2}$ in the source estimates:
\[
\begin{aligned}
 \|N\|_{L^2_3L^1_h}
 &\le C\mathcal J_2(1+s)^{-1/4}\cD_0(U)^{1/2}
           +C\mathcal J_2^{3/2}(1+s)^{-5/8}\cD_0(U)^{1/4},\\
 \|\nu_h(e^{-q}-1)\nabla_hu\|_{L^2_3L^1_h}
 &\le C\mathcal J_2(1+s)^{-5/8}\cD_0(U)^{1/2}.
\end{aligned}
\]
On $[0,t/2]$, take $r=3/4$ in the weighted dissipation
estimate of Lemma~\ref{lem:new-seed}. Cauchy--Schwarz for
the terms containing $\cD_0(U)^{1/2}$, and H\"older with exponents $4,4/3$
for the term containing $\cD_0(U)^{1/4}$, bound this part by
$C\mathcal J_2^2(1+t)^{-1}$.
On $[t/2,t]$, the dissipation tail in
Lemma~\ref{lem:new-seed} gives the following contributions
from the $\cD_0(U)^{1/2}$ term of $N$, the $\cD_0(U)^{1/4}$ term of $N$ and the
horizontal flux, respectively:
\[
\begin{aligned}
 &\mathcal J_2\int_{t/2}^t(1+t-s)^{-1}(1+s)^{-1/4}\cD_0(U)^{1/2}\dd s
                 \le C\mathcal J_2^2(1+t)^{-3/4},\\
 &\mathcal J_2^{3/2}\int_{t/2}^t
        (1+t-s)^{-1}(1+s)^{-5/8}\cD_0(U)^{1/4}\dd s
                 \le C\mathcal J_2^2(1+t)^{-7/8},\\
 &\mathcal J_2\int_{t/2}^t
        (1+t-s)^{-3/2}(1+s)^{-5/8}\cD_0(U)^{1/2}\dd s
                 \le C\mathcal J_2^2(1+t)^{-9/8}.
\end{aligned}
\]
The kernels raised to powers $2,4/3,2$ are integrable on
$[t/2,t]$. Adding the high-horizontal-frequency estimate for
$v$ at $\sigma=1/2$ and the estimate of $\nabla W$ from
\eqref{eq:new-low-first-rough} and \eqref{eq:new-E1-differential}
gives
$\|\nabla_hv(t)\|_2+\|\nabla W(t)\|_2
 \le C\mathcal J_2(1+t)^{-3/4}$.

\smallskip\noindent\emph{The case $\sigma=3/4$.}
For $\sigma=3/4$, \eqref{eq:new-sigma-mixed} gives
$\|N\|_{L^2_3L^1_h}\le C\mathcal J_2^2(1+s)^{-1}$.
The corresponding convolution is bounded by
$C\mathcal J_2^2(1+t)^{-1}\log(e+t)$.
Since
\[
 (1+t)^{-1}\log(e+t)\le C(1+t)^{-7/8},
\]
we obtain the $(1+t)^{-7/8}$ bound for $\|\nabla_hv\|_2$.
The high-horizontal-frequency estimate for $v$ at
$\sigma=3/4$ and the estimate of $\nabla W$ obtained from
\eqref{eq:new-low-first-rough} and \eqref{eq:new-E1-differential}
also have power at least $7/8$. Thus
$\|\nabla_hv\|_2+\|\nabla W\|_2
\le C\mathcal J_2(1+t)^{-7/8}$.

\smallskip\noindent\emph{The case $\sigma=7/8$.}
For $\sigma=7/8$, \eqref{eq:new-sigma-mixed} gives
$\|N\|_{L^2_3L^1_h}\le C\mathcal J_2^2(1+s)^{-17/16}$,
where $17/16=5/8+\sigma/2>1$. Splitting at $t/2$ gives
\[
\int_0^t(1+t-s)^{-1}(1+s)^{-17/16}\dd s
 \le C\bigl((1+t)^{-1}
       +(1+t)^{-17/16}\log(e+t)\bigr)
 \le C(1+t)^{-1}.
\]
The other terms of $N$ and the horizontal flux in
\eqref{eq:new-sigma-mixed} decay faster than $(1+s)^{-17/16}$;
their convolutions are also bounded by $C(1+t)^{-1}$.
The high-horizontal-frequency estimate for $v$ at
$\sigma=7/8$ and the estimate of $\nabla W$ from
\eqref{eq:new-low-first-rough} and \eqref{eq:new-E1-differential}
therefore give $(1+t)^{-1}$ decay for
$\|\nabla_hv\|_2+\|\nabla W\|_2$.
Using $\nabla_hu=\nabla_hv+\nabla_hw$ and
$\|\nabla_hw\|_2\le\|\nabla W\|_2$ at each stage completes
the iteration in \eqref{eq:new-first-iteration}, with
\eqref{eq:iteration-hypothesis} verified at each exponent.

In particular,
\[
 \|\nabla W(t)\|_2+\|\nabla_hu(t)\|_2
                  \le C\mathcal J_2(1+t)^{-1}.
\]
This proves the horizontal velocity estimate in the first
line of \eqref{eq:new-first} and supplies the case $\sigma=1$
for the remaining argument.

\smallskip\noindent\emph{Step 4. The $(1+t)^{-5/4}$ bound for $\nabla W$.}
We now improve the estimate of $\nabla W$ from Step~3
to the first bound in \eqref{eq:new-first}, and then prove
\eqref{eq:new-first-windows}. The case $\sigma=1$ of
\eqref{eq:iteration-hypothesis} gives the following $L^2$
bound for the fluxes in \eqref{eq:new-decay-flux}:
\[
 \sum_j\|F^{(j)}\|_2
                   \le C\mathcal J_2^2(1+s)^{-11/8}.
\]
The slowest product is $u\otimes u$, because
\[
 \|u\|_4^4\le C\|u\|_2\|\partial_3u\|_2\|\nabla_hu\|_2^2
                       \le C\mathcal J_2^4(1+s)^{-11/4}.
\]
The scalar convection and pressure fluxes decay faster, using
\[
 \|(q,\tau)\|_4
    \le C\|(q,\tau)\|_2^{1/4}\|\nabla(q,\tau)\|_2^{3/4}
    \le C\mathcal J_2(1+s)^{-15/16}.
\]
The diffusion fluxes are bounded by
$C\|q\|_\infty\|(\nabla_hu,\Div u,\nabla\tau)\|_2
 \le C\mathcal J_2^2(1+s)^{-3/2}$.
For the remainder \eqref{eq:new-flux-remainders}, the
case $\sigma=1$ gives
$\|F^{(0)}(s)\|_1\le C\mathcal J_2^2(1+s)^{-3/2}$.
For the flux, use its $L^1$ bound from Step~2 on $[0,t/2]$,
and the $L^2$ bound
$\sum_j\|F^{(j)}(s)\|_2\le C\mathcal J_2^2(1+s)^{-11/8}$
on $[t/2,t]$. The spatial derivative and the divergence of the flux
together give an $L^2$ kernel bounded by $C(1+t-s)^{-1}$,
using \eqref{eq:blockmultiplier}.
The remainder and the flux contributions on $[0,t/2]$ and
$[t/2,t]$, respectively, satisfy
\[
\begin{aligned}
 &\int_0^t(1+t-s)^{-5/4}\|F^{(0)}(s)\|_1\dd s
             \le C\mathcal J_2^2(1+t)^{-5/4},\\
 &\int_0^{t/2}(1+t-s)^{-7/4}
                          \sum_j\|F^{(j)}(s)\|_1\dd s
             \le C\mathcal J_2^2(1+t)^{-7/4}\log(e+t),\\
 &\int_{t/2}^t(1+t-s)^{-1}
                          \sum_j\|F^{(j)}(s)\|_2\dd s
             \le C\mathcal J_2^2(1+t)^{-11/8}\log(e+t).
\end{aligned}
\]
All three are bounded by $C\mathcal J_2^2(1+t)^{-5/4}$.
Adding the homogeneous term gives
$\|\nabla P_LW(t)\|_2\le C\mathcal J_2(1+t)^{-5/4}$.
For high frequencies, \eqref{eq:new-sigma-sources} at
$\sigma=1$ and \eqref{eq:new-E1-differential} give
$E_{1,H}(t)\le C\mathcal J_2^2(1+t)^{-5/2}$.
Taking square roots and combining the two frequency ranges proves
the $\nabla W$ estimate in \eqref{eq:new-first}:
\[
                  \|\nabla W(t)\|_2
                        \le C\mathcal J_2(1+t)^{-5/4}.
\]

The $\mathcal Lu$ estimate in \eqref{eq:new-first-windows}
follows directly from \eqref{eq:new-sigma-window} with
$\sigma=1$. To prove the $(\nabla^2w,\nabla^2\tau)$ estimate,
first control the high frequencies. The bounds
\eqref{eq:new-sigma-sources}, evaluated at $t+1$ with
$\sigma=1$, control the source integrals over $[t,t+1]$, because
the exponential weight there is bounded below by $e^{-c_*}$.
Integrating \eqref{eq:new-E1-differential} over this interval
therefore yields
\[
\begin{aligned}
 \int_t^{t+1}\|(w_H,\tau_H)(s)\|_{H^2}^2\dd s
 &\le C E_{1,H}(t)
       +C\int_t^{t+1}\cD_2(U)E_{1,H}\dd s\\
 &\quad+C\int_t^{t+1}
     \bigl(\|F_2\|_2^2+\|F_3\|_2^2+\|C_q\|_{H^1}^2
                  +\mathcal J_2^4(1+s)^{-3}\bigr)\dd s\\
 &\le C\mathcal J_2^2(1+t)^{-5/2}.
\end{aligned}
\]
The integral of $\cD_2(U)E_{1,H}$ uses the total dissipation
bound \eqref{eq:new-base-energy} and
$E_{1,H}(s)\le C\mathcal J_2^2(1+s)^{-5/2}$, obtained from
\eqref{eq:new-E1-differential} with $\sigma=1$.
For the low frequencies, the fixed cutoff gives
\[
 \|\nabla^2(w_L,\tau_L)(s)\|_2
          \le C\|\nabla P_LW(s)\|_2
          \le C\mathcal J_2(1+s)^{-5/4}.
\]
Squaring and integrating over $[t,t+1]$, then adding the
high-frequency contribution, proves the second estimate in
\eqref{eq:new-first-windows}.

\smallskip\noindent\emph{Step 5. First derivatives of $u_3$.}
The identity
\[
 \partial_3u_3=\Div w-\Div_hu_h
\]
together with the first line of \eqref{eq:new-first},
proved in Steps~3--4, gives
\[
 \|\partial_3u_3(t)\|_2
       \le C\bigl(\|\nabla W(t)\|_2+\|\nabla_hu(t)\|_2\bigr)
       \le C\mathcal J_2(1+t)^{-1}.
\]
To estimate the horizontal derivatives, write $u_3=v_3+w_3$.
The $\nabla W$ bound in \eqref{eq:new-first}, proved in Step~4,
controls $\nabla_hw_3$ at rate $(1+t)^{-5/4}$; it remains to
estimate $\nabla_hv_3$.
For $t\ge2$, split its Duhamel integral at $t-1$.

\smallskip\noindent\emph{The interval $[0,t-1]$.}
For $a=t-s\ge1$, the third-row identities
\eqref{eq:new-third-row} give
\[
 \|\nabla_h[S_h(a)\mathbb PG]_3\|_2
       \le Ca^{-5/4}\|G\|_1,\qquad
 \|\nabla_h[S_h(a)\mathbb P\partial_3G]_3\|_2
       \le Ca^{-3/2}\|G\|_{L^2_3L^1_h}.
\]
The second bound uses
$\partial_3[\mathbb PG]_3=-\Div_h[\mathbb PG]_h$.
Decompose convection as
\[
 -u\cdot\nabla u
      =-u_h\cdot\nabla_hu-\partial_3(u_3u)+(\partial_3u_3)u.
\]
The bounds for $\nabla_hu$ from Step~3 and $\partial_3u_3$
from the beginning of this step, together with the velocity
decay in \eqref{eq:new-seed}, give
\[
 \|u_h\cdot\nabla_hu\|_1+\|(\partial_3u_3)u\|_1
                   \le C\mathcal J_2^2(1+s)^{-3/2}.
\]
For the vertical flux, vertical Sobolev gives
\[
\|u_3u\|_{L^2_3L^1_h}
 \le\|u_3\|_{L^\infty_3L^2_h}\|u\|_2
 \le C\|u_3\|_2^{1/2}\|\partial_3u_3\|_2^{1/2}\|u\|_2
       \le C\mathcal J_2^2(1+s)^{-11/8}.
\]
After removal of its gradient, pressure is represented by
$\Pi_p$ in \eqref{eq:new-N-definitions}, with
\[
 \|\Pi_p\|_1\le C\|(q,\tau)\|_2\|\nabla\tau\|_2
                    \le C\mathcal J_2^2(1+s)^{-2}.
\]
Convolve the horizontal convection $u_h\cdot\nabla_hu$,
the remainder $(\partial_3u_3)u$ and $\Pi_p$ with the
$a^{-5/4}$ kernel obtained from \eqref{eq:new-third-row}.
The vertical flux $u_3u$ uses the $a^{-3/2}$ mixed-norm
kernel. Splitting each convolution at $t/2$ gives a bound
$C\mathcal J_2^2(1+t)^{-5/4}$ for all four contributions.

For variable viscosity, use
\[
 \|(e^{-q}-1)\mathcal Lu\|_1
       \le C\mathcal J_2(1+s)^{-3/4}\|\mathcal Lu\|_2.
\]
Cauchy--Schwarz on each unit interval and the first bound in
\eqref{eq:new-first-windows} give
\[
 \int_n^{n+1}(1+s)^{-3/4}\|\mathcal Lu(s)\|_2\dd s
                     \le C\mathcal J_2(1+n)^{-7/4}.
\]
Summing the corresponding kernel contributions yields
\[
\begin{aligned}
 &\mathcal J_2\int_0^{t-1}(t-s)^{-5/4}
                      (1+s)^{-3/4}\|\mathcal Lu\|_2\dd s\\
 &\quad\le C\mathcal J_2^2
     \sum_{n=0}^{\lfloor t-1\rfloor}
                 (1+t-n)^{-5/4}(1+n)^{-7/4}
     \le C\mathcal J_2^2(1+t)^{-5/4}.
\end{aligned}
\]

\smallskip\noindent\emph{The interval $[t-1,t]$.}
Apply \eqref{eq:r8-horizontal-maximal}
with $m=0$ to $\mathbb PF_2$.
Since $\mathbb P$ is bounded on $L^2$, and
\eqref{eq:new-sigma-sources} at $\sigma=1$ controls the last
unit interval, we obtain
\[
\left\|\nabla_h\int_{t-1}^tS_h(t-s)\mathbb PF_2(s)\dd s\right\|_2^2
 \le C\int_{t-1}^t\|F_2(s)\|_2^2\dd s
 \le C\mathcal J_2^4(1+t)^{-5/2}.
\]
Taking square roots bounds this contribution by
$C\mathcal J_2^2(1+t)^{-5/4}$.
For the homogeneous term, the third-row estimate derived
from \eqref{eq:new-third-row} with $u_0\in L^1$ gives
$\|\nabla_h[S_h(t)\mathbb Pu_0]_3\|_2
 \le C\mathcal J_2(1+t)^{-5/4}$.
Adding the estimates on both time intervals proves the bound
for $\nabla_hv_3$. The $\nabla W$ estimate in Step~4 then
controls $\nabla_hw_3$, so $u_3=v_3+w_3$ gives
\[
                  \|\nabla_hu_3(t)\|_2
                        \le C\mathcal J_2(1+t)^{-5/4}.
\]
Together with the bound for $\partial_3u_3$ at the beginning
of this step, this proves the second line of
\eqref{eq:new-first}.
For $0\le t\le2$, the pointwise bounds and the unit-interval
integrals follow from the uniform $H^2$ estimate
\eqref{eq:new-base-energy}. This completes the proof.
\end{proof}

\subsubsection{Second-order decay}

\begin{lemma}[Second-order decay]\label{lem:new-joint-second}
For every $t\ge0$,
\[
 \|\nabla_h^2u(t)\|_2\le C\mathcal J_2(1+t)^{-3/2},\qquad
 \int_t^{t+1}\|\nabla_h^3u(s)\|_2^2\dd s
                       \le C\mathcal J_2^2(1+t)^{-3}.
\]
Moreover,
\begin{equation}\label{eq:new-compressible-second}
 \|\nabla^2W(t)\|_2\le C\mathcal J_2 (1+t)^{-7/4},\qquad
 \int_t^{t+1}\|\nabla^3(w,\tau)(s)\|_2^2\dd s
                              \le C\mathcal J_2^2(1+t)^{-7/2}.
\end{equation}
\end{lemma}
\begin{proof}
Steps~1--3 establish a common second-order decay rate $(1+t)^{-3/2}$. Step~4 improves the density and potential-velocity estimates, and Step~5 uses the time-integrated dissipation to improve the temperature estimate.
Set
\[
 \begin{aligned}
 A_2(s)&=\|\nabla^2q(s)\|_2+\|\nabla^2w(s)\|_2
             +\|\nabla^2\tau(s)\|_2+\|\nabla_h^2v(s)\|_2,\\
 M_2(t)&=\mathcal J_2+\sup_{0\le s\le t}(1+s)^{3/2}A_2(s),
 \end{aligned}
\]
where $A_2=A_2(s)$ when time is suppressed. Assume temporarily $M_2(t)\le1$.

\smallskip\noindent\emph{Step 1. Preliminary estimates.}
Using $u=v+w$, $\Div v=0$, and
\eqref{eq:new-elliptic}, we obtain
\begin{equation}\label{eq:new-S2-structure}
\begin{aligned}
 \|\mathcal Lu\|_2+\|\nabla_h^2u\|_2+\|\nabla_h\nabla u_3\|_2&\le C A_2,
 \\
 \|\nabla_h\partial_3u\|_2^2&\le C\mathcal J_2 A_2.
\end{aligned}
\end{equation}
For the second line of \eqref{eq:new-S2-structure}, Plancherel gives
\[
 \|\nabla_h\partial_3u\|_2^2
       \le\|\nabla_h^2u\|_2\|\partial_3^2u\|_2,
\]
and \eqref{eq:new-base-energy} bounds the last factor by
$C\mathcal J_2$. Since $\mathbb P$ is an $L^2$ contraction
and commutes with derivatives, $v=\mathbb Pu$ also satisfies
$\|\partial_3^2v\|_2\le C\mathcal J_2$. Thus
\[
 \|\nabla_h\partial_3v\|_2^2
 \le\|\nabla_h^2v\|_2\|\partial_3^2v\|_2
 \le C\mathcal J_2\|\nabla_h^2v\|_2.
\]
For $u_3$ and $(\nabla_hu,\Div u,\nabla\tau)$, we obtain
\begin{equation}\label{eq:new-G-structure}
\begin{aligned}
 \|u_3\|_\infty^2&\le C\mathcal J_2 A_2,\\
 \|\nabla(\nabla_hu,\Div u,\nabla\tau)\|_2^2
          +\|(\nabla_hu,\Div u,\nabla\tau)\|_6^2
 &\le C\mathcal J_2 A_2,\\
 \|(\nabla_hu,\Div u,\nabla\tau)\|_\infty^2
 &\le C A_2\cD_2(U)^{1/2}.
\end{aligned}
\end{equation}
The first line uses the second inequality in
\eqref{eq:new-elliptic}: its two factors are bounded by
\eqref{eq:new-base-energy} and \eqref{eq:new-S2-structure},
respectively. The second line uses both bounds in
\eqref{eq:new-S2-structure}, $\Div u=\Div w$, and the uniform
$H^2$ bound in \eqref{eq:new-base-energy}; the $L^6$ estimate then follows
from $\|g\|_6\le C\|\nabla g\|_2$.
For the last line of \eqref{eq:new-G-structure}, apply the first
inequality in \eqref{eq:new-agmon} componentwise: the first
factor is controlled by \eqref{eq:new-S2-structure}, and the
second by the dissipation \eqref{eq:H2-dissipation-definition}.

Applying \eqref{eq:new-agmon} and using
Lemmas~\ref{lem:new-seed} and \ref{lem:new-first}, together with
\eqref{eq:new-S2-structure}, we obtain
\begin{equation}\label{eq:new-second-bootstrap-ledger}
\begin{aligned}
 \|u\|_\infty&\le CM_2(t)(1+s)^{-7/8},&
 \|u_3\|_\infty&\le CM_2(t)(1+s)^{-5/4},\\
 \|W\|_\infty&\le CM_2(t)(1+s)^{-11/8},&
 \|\nabla_h\partial_3u\|_2&\le CM_2(t)(1+s)^{-3/4}.
\end{aligned}
\end{equation}
For $u$, use the third inequality in \eqref{eq:new-agmon}.
For $u_3$ and $W$, use its second and first inequalities,
respectively, with the first-order bounds \eqref{eq:new-first}
and the second-order bounds \eqref{eq:new-S2-structure}.
The last estimate in \eqref{eq:new-second-bootstrap-ledger}
follows from the second line of \eqref{eq:new-S2-structure}.

\Needspace{5\baselineskip}
\smallskip\noindent\emph{Step 2. Low-frequency estimates.}
\smallskip\noindent\emph{The compressible variables.}
Write
$F=\mathcal R-(0,\Div(u\otimes u),0)$, where
$\mathcal R_1=F_1$, $\mathcal R_3=F_3$, and
\[
 \mathcal R_2=u\Div u-c^2(e^\tau-1)\nabla(q+\tau)
                         +(e^{-q}-1)\mathcal Lu.
\]
Using the zero-order bounds \eqref{eq:new-seed} and
\eqref{eq:new-zero}, the first-order bounds \eqref{eq:new-first},
and \eqref{eq:new-second-bootstrap-ledger}, we obtain
\begin{equation}\label{eq:new-second-low-ledger}
\begin{aligned}
 \|\mathcal R\|_1&\le CM_2(t)^2(1+s)^{-7/4},&
 \|\mathcal R\|_2&\le CM_2(t)^2(1+s)^{-17/8},\\
 \|u\otimes u\|_1&\le C\mathcal J_2^2(1+s)^{-1},&
 \|u_3u\|_2&\le CM_2(t)^2(1+s)^{-7/4},\\
 \|\Div_h(u\otimes u_h)\|_2&\le CM_2(t)^2(1+s)^{-15/8}.
\end{aligned}
\end{equation}
The powers $7/4$ and $17/8$ in the remainder are attained by
scalar transport and $u\Div u$, since
\[
 \|u\|_2\|\nabla W\|_2
       \le CM_2(t)^2(1+s)^{-7/4},\qquad
 \|u\|_\infty\|\nabla W\|_2
       \le CM_2(t)^2(1+s)^{-17/8}.
\]
For coefficient diffusion, use
$\|(e^{-q}-1)\mathcal Lu\|_1
 \le C\|q\|_2\|\mathcal Lu\|_2$
and
$\|(e^{-q}-1)\mathcal Lu\|_2
 \le C\|q\|_\infty\|\mathcal Lu\|_2$.
Equations~\eqref{eq:new-zero}, \eqref{eq:new-S2-structure}, and
\eqref{eq:new-second-bootstrap-ledger} give decay
$(1+s)^{-9/4}$ in $L^1$ and $(1+s)^{-23/8}$ in $L^2$.
For $(e^{-q}-1)\Delta\tau$, use the same coefficient factors
and the bound for $\nabla^2\tau$ in the definition of $M_2(t)$.
For the quadratic sources, use \eqref{eq:new-first} in $L^1$ and
\eqref{eq:new-production-L2} in $L^2$, with
\eqref{eq:new-S2-structure} and \eqref{eq:new-second-bootstrap-ledger}.
The coefficients are bounded by \eqref{eq:new-base-energy}, so
\[
\begin{aligned}
 \|e^{-q-\tau}Q(u)\|_1&\le C\mathcal J_2^2(1+s)^{-2},
 \qquad \|e^{-q}|\nabla\tau|^2\|_1\le C\mathcal J_2^2(1+s)^{-5/2},\\
 \|e^{-q-\tau}Q(u)\|_2&\le CM_2(t)^2(1+s)^{-19/8},
 \qquad \|e^{-q}|\nabla\tau|^2\|_2\le CM_2(t)^2(1+s)^{-23/8}.
\end{aligned}
\]
Since $\mathcal J_2\le M_2(t)$, these bounds are stronger than the
corresponding remainder estimates in \eqref{eq:new-second-low-ledger}.

The Fourier proof of \eqref{eq:Scestimate} gives, for every
fixed integer $m\ge0$ and $a\ge0$,
\[
\|\nabla^mP_LS_c(a)G\|_2
       \le C_m(1+a)^{-3/4-m/2}\|G\|_1,
 \qquad \|\nabla^mP_LS_c(a)G\|_2
       \le C_m(1+a)^{-m/2}\|G\|_2.
\]
Thus three derivatives may be kept on the low-frequency kernel
when a derivative is present in the source.
Split the Duhamel integral at $s=t/2$, as in Step~4 of
Lemma~\ref{lem:new-first}. On $[0,t/2]$, estimate the remainder and the full momentum
flux in $L^1$. On $[t/2,t]$, estimate the remainder in $L^2$
and write
$\Div(u\otimes u)=\Div_h(u\otimes u_h)+\partial_3(u_3u)$.
The five contributions to $\|\nabla^2P_LW(t)\|_2$ satisfy
\[
\begin{aligned}
 &\int_0^{t/2}(1+t-s)^{-7/4}\|\mathcal R(s)\|_1\dd s
       \le CM_2(t)^2(1+t)^{-7/4},\\
 &\int_0^{t/2}(1+t-s)^{-9/4}\|u\otimes u\|_1\dd s
       \le CM_2(t)^2(1+t)^{-9/4}\log(e+t),\\
 &\int_{t/2}^t(1+t-s)^{-1}\|\mathcal R(s)\|_2\dd s
       \le CM_2(t)^2(1+t)^{-17/8}\log(e+t),\\
 &\int_{t/2}^t(1+t-s)^{-1}
                       \|\Div_h(u\otimes u_h)\|_2\dd s
       \le CM_2(t)^2(1+t)^{-15/8}\log(e+t),\\
 &\int_{t/2}^t(1+t-s)^{-3/2}\|u_3u\|_2\dd s
       \le CM_2(t)^2(1+t)^{-7/4}.
\end{aligned}
\]
On $[t/2,t]$, the kernel $(1+t-s)^{-1}$ gives a logarithm,
whereas the integrable kernel $(1+t-s)^{-3/2}$ retains the
$(1+s)^{-7/4}$ bound for $u_3u$ in
\eqref{eq:new-second-low-ledger}.
The remainder and both momentum flux contributions are
consequently bounded by $CM_2(t)^2(1+t)^{-7/4}$.
Adding $\nabla^2P_LS_c(t)W_0$, estimated by the Fourier proof
of \eqref{eq:Scestimate}, proves the first line of
\eqref{eq:new-low-second}:
\begin{equation}\label{eq:new-low-second}
\begin{aligned}
 \|\nabla^2P_LW\|_2&\le C(\mathcal J_2+M_2(t)^2)(1+t)^{-7/4},\\
 \|P_LW\|_\infty&\le C(\mathcal J_2+M_2(t)^2)(1+t)^{-3/2},\\
 \|\nabla^jP_LW\|_\infty&\le C_j(\mathcal J_2+M_2(t)^2)(1+t)^{-7/4}
 \quad(1\le j\le4).
\end{aligned}
\end{equation}
For the last two lines of \eqref{eq:new-low-second}, use
Fourier absolute integration on $[0,t/2]$ and $[t/2,t]$
with the source bounds in \eqref{eq:new-second-low-ledger}.
The corresponding kernels are
$(1+a)^{-3/2-m/2}$ for an $L^1$ source and
$(1+a)^{-3/4-m/2}$ for an $L^2$ source.
Without spatial differentiation, the remainder on $[0,t/2]$ determines
the rate $(1+t)^{-3/2}$. With one spatial derivative, the kernels
on $[t/2,t]$ are integrable, and the vertical flux determines
the rate $(1+t)^{-7/4}$.
Since $P_L$ has fixed compact frequency support, each further
derivative is bounded by a fixed multiple of the first.
This proves the last two lines of \eqref{eq:new-low-second}
and also
\[
 \|\widehat{P_LW}\|_1
       \le C(\mathcal J_2+M_2(t)^2)(1+s)^{-3/2},\qquad
 \sum_{j=1}^4\||\xi|^j\widehat{P_LW}\|_1
       \le C(\mathcal J_2+M_2(t)^2)(1+s)^{-7/4}.
\]

\smallskip\noindent\emph{Horizontal derivatives of $v$.}
Use the source decomposition \eqref{eq:new-projected-source}.
The mixed-product estimates in \eqref{eq:new-sigma-mixed}
at $\sigma=1$, the final first-order bounds
\eqref{eq:new-first}, and vertical Sobolev embedding
\eqref{eq:mixedv} give
\[
\begin{aligned}
 \|u\cdot\nabla u\|_{L^2_3L^1_h}
       &\le C\mathcal J_2^2(1+s)^{-9/8},
 \qquad \|(\Pi_p,\Pi_v)\|_{L^2_3L^1_h}
       \le C\mathcal J_2^2(1+s)^{-13/8},\\
 \|\nu_h(e^{-q}-1)\nabla_hu\|_{L^2_3L^1_h}
       &\le C\mathcal J_2^2(1+s)^{-2}.
\end{aligned}
\]
For $\Pi_v$, vertical interpolation bounds
$\|\nabla(e^{-q}-1)\|_{L^\infty_3L^2_h}$
by $C\mathcal J_2(1+s)^{-5/8}$; multiplying by
$\|(\nabla_hu,\Div u)\|_2$ gives power $13/8$.
For the horizontal diffusion flux,
$\|e^{-q}-1\|_{L^\infty_3L^2_h}$ has power $1$,
as does $\|\nabla_hu\|_2$.
For one horizontal derivative of convection, horizontal
H\"older bounds the four products by
\[
\begin{gathered}
 \|\nabla_hu\|_{L^\infty_3L^2_h}\|\nabla_hu\|_2,\qquad
 \|u\|_{L^\infty_3L^2_h}\|\nabla_h^2u\|_2,\\
 \|\nabla_hu_3\|_{L^\infty_3L^2_h}\|\partial_3u\|_2,\qquad
 \|u_3\|_{L^\infty_3L^2_h}\|\nabla_h\partial_3u\|_2.
\end{gathered}
\]
Vertical Sobolev embedding \eqref{eq:mixedv},
\eqref{eq:new-S2-structure}, and the first-order bounds
\eqref{eq:new-first} give
\[
\begin{aligned}
 \|\nabla_hu\|_{L^\infty_3L^2_h}
       &\le CM_2(t)(1+s)^{-7/8},&
 \|u\|_{L^\infty_3L^2_h}
       &\le C\mathcal J_2(1+s)^{-3/8},\\
 \|\nabla_hu_3\|_{L^\infty_3L^2_h}
       &\le CM_2(t)(1+s)^{-11/8},&
 \|u_3\|_{L^\infty_3L^2_h}
       &\le C\mathcal J_2(1+s)^{-7/8}.
\end{aligned}
\]
The four products therefore have powers
$15/8,15/8,13/8,13/8$, and
\[
 \|\nabla_h(u\cdot\nabla u)\|_{L^2_3L^1_h}
       \le CM_2(t)^2(1+s)^{-13/8}.
\]
For convection, keep two horizontal derivatives on the kernel
on $[0,t/2]$ and move one to the source on $[t/2,t]$.
The pressure and viscosity remainders keep two kernel
derivatives throughout; the horizontal diffusion flux keeps
three. Estimate~\eqref{eq:Sh-low-all-orders} gives the
required kernel powers for these fixed derivative orders.
On $[t/2,t]$, the convection kernel $(1+t-s)^{-1}$ and the
source bound $(1+s)^{-13/8}$ give
$(1+t)^{-13/8}\log(e+t)$.
The undifferentiated convection source has power $9/8$,
while the pressure and viscosity remainders have power $13/8$;
all three are integrable in time.
Together with the horizontal diffusion flux of power $2$,
they bound the nonlinear contribution by
\[
 CM_2(t)^2\bigl((1+t)^{-3/2}
        +(1+t)^{-13/8}\log(e+t)+(1+t)^{-2}\bigr).
\]
Since $(1+t)^{-13/8}\log(e+t)\le C(1+t)^{-3/2}$,
adding the homogeneous horizontal heat estimate from the
localized estimate \eqref{eq:Sh-low-all-orders} proves
\begin{equation}\label{eq:new-low-horizontal-second}
 \|\nabla_h^2P_L^hv\|_2\le C(\mathcal J_2+M_2(t)^2)(1+t)^{-3/2}.
\end{equation}

\begin{samepage}
\smallskip\noindent\emph{Step 3. High-frequency energy estimate.}
For the compressible variables and the horizontal derivatives
of $v$, use
\begin{align}\label{eq:new-joint-energy}
 E_c={}&\tfrac12\sum_{|\gamma|\le2}
 (c^2\|\partial^\gamma q_H\|_2^2+
 \|\partial^\gamma w_H\|_2^2+\omega\|\partial^\gamma\tau_H\|_2^2)
 +\eta\sum_{|\gamma|\le1}\ip{\partial^\gamma w_H}{\nabla\partial^\gamma q_H},\notag\\
 E={}&E_c+\tfrac12\|\Delta_hP_H^hv\|_2^2.
\end{align}
\end{samepage}
For sufficiently small $\eta>0$,
\[
 E\simeq\|W_H\|_{H^2}^2+\|\nabla_h^2P_H^hv\|_2^2.
\]
The compressible equations are differentiated up to order two,
and the shear equation is tested after two horizontal derivatives.
As in Step~2 of Lemma~\ref{lem:new-first}, the cross terms
in $E_c$ supply the density dissipation.
The fixed frequency cutoffs then give linear dissipation
controlling
$E+\|(w_H,\tau_H)\|_{H^3}^2+\|\nabla_h^3P_H^hv\|_2^2$.
The full second-order norms satisfy
\[
\bigl(\|\nabla^2q\|_2+\|\nabla^2w\|_2
             +\|\nabla^2\tau\|_2+\|\nabla_h^2v\|_2\bigr)^2
 \le C\bigl(E+\|\nabla^2P_LW\|_2^2
                         +\|\nabla_h^2P_L^hv\|_2^2\bigr).
\]
By \eqref{eq:new-low-second}, \eqref{eq:new-low-horizontal-second},
and the coefficient decomposition \eqref{eq:new-coefficient-split},
\[
 \|\nabla^2P_LW\|_2+\|\nabla_h^2P_L^hv\|_2
                 +\|b_{\ell}\|_{W^{1,\infty}}
       \le C(\mathcal J_2+M_2(t)^2)(1+s)^{-3/2}.
\]
We estimate the nonlinear terms contributing to the energy
$E$ in \eqref{eq:new-joint-energy}, keeping the
low-frequency terms from \eqref{eq:new-low-second} and
\eqref{eq:new-low-horizontal-second} as forcing.

\smallskip\noindent\emph{Variable diffusion.}
The coefficient decomposition \eqref{eq:new-coefficient-split}
and the composition estimate \eqref{eq:Mosercomp} give
$\|b_{\mathrm r}\|_{H^2}\le C\|q_H\|_{H^2}$.
The coefficient $b_{\ell}=e^{-q_L}-1$ is estimated in $W^{1,\infty}$.
By the dissipation definition \eqref{eq:H2-dissipation-definition},
$\|\mathcal Lu\|_{H^1}+\|\Delta\tau\|_{H^1}$
is controlled by $\cD_2(U)^{1/2}$. Hence
\begin{equation}\label{eq:new-joint-coefficient}
\begin{aligned}
 \|(e^{-q}-1)\mathcal Lu\|_{H^1}+\|(e^{-q}-1)\Delta\tau\|_{H^1}
 \le C\cD_2(U)^{1/2}
                (\|q_H\|_{H^2}+\|b_{\ell}\|_{W^{1,\infty}}).
\end{aligned}
\end{equation}

\smallskip\noindent\emph{Density transport and compensation.}
Keep $u\cdot\nabla q_H$ on the left of the density equation.
Lemma~\ref{lem:H2-anisotropic} bounds the second transport
commutator by
\[
 \|\nabla^2u\,\nabla q_H\|_2\|q_H\|_{H^2}
 \le C\mathcal J_2^{1/3}\cD_2(U)^{1/3}\|q_H\|_{H^2}^2
 \le(\delta+C_\delta\cD_2(U))E.
\]
The first transport commutator is also bounded by
$(\delta+C_\delta\cD_2(U))E$, using
Lemma~\ref{lem:H2-anisotropic}.
For the cutoff commutator $C_q$ in
\eqref{eq:new-density-high}, the low-frequency estimate
\eqref{eq:new-low-second} gives
\begin{equation}\label{eq:new-second-Cq}
 \|C_q\|_{H^2}\le C\|u\|_2\|\nabla q\|_2
 +C\|u\|_{H^2}\sum_{j=1}^3\|\nabla^jq_L\|_\infty
 \le C\mathcal J_2(\mathcal J_2+M_2(t)^2)(1+s)^{-7/4}.
\end{equation}
The compensation term requires only
$\|P_HF_1\|_{H^1}$. The $H^2$ product estimate bounds it by
$C\mathcal J_2(\|\nabla^2q\|_2+\|\nabla^2w\|_2+
\|\nabla^2\tau\|_2+\|\nabla_h^2v\|_2)$.

\smallskip\noindent\emph{Velocity convection.}
We claim
\begin{equation}\label{eq:new-projected-convection}
 \|P_H\mathbb Q(u\cdot\nabla u)\|_{H^1}
       +\|\nabla_h(u\cdot\nabla u)\|_2
 \le C\mathcal J_2(\|\nabla_h^2v\|_2+\|\nabla^2w\|_2).
\end{equation}
The projected term uses $\Div v=0$:
\begin{equation}\label{eq:new-potential-identity}
 \mathbb Q(u\cdot\nabla v)=\nabla\Delta^{-1}
                       \sum_{i,j}\partial_i u_j\partial_jv_i.
\end{equation}
On the support of $P_H$, the multiplier
$P_H\nabla\Delta^{-1}$ maps $L^2$ to $H^1$.
The horizontal products are controlled by
$\|\nabla_ha\|_4^2\le C\|a\|_\infty\|\nabla_h^2a\|_2$.
For $\partial_3u_3\partial_3v_3$, use
$\partial_3v_3=-\Div_hv_h$ and
$\partial_3u_3=\Div w-\Div_hu_h$, together with the
corresponding $L^4$ estimate for $\nabla w$.
The remaining vertical products satisfy
\[
\begin{aligned}
 \|\partial_3u_h\nabla_hv_3\|_2
 &\le\|\partial_3u_h\|_3\|\nabla_hv_3\|_6
       \le C\mathcal J_2\|\nabla_h^2v\|_2,\\
 \|\partial_3v_h\nabla_hu_3\|_2
 &\le C\mathcal J_2
                   (\|\nabla_h^2v\|_2+\|\nabla^2w\|_2).
\end{aligned}
\]
Here
$\|\nabla\nabla_hv_3\|_2\le C\|\nabla_h^2v\|_2$.
For $u\cdot\nabla w$, H\"older's inequality with exponents $3,6$
and the uniform
$H^2$ bound \eqref{eq:new-base-energy} give the $H^1$ bound
$C\mathcal J_2\|\nabla^2w\|_2$.
Combining this with the products in
\eqref{eq:new-potential-identity} proves the bound for
$P_H\mathbb Q(u\cdot\nabla u)$ in
\eqref{eq:new-projected-convection}.

For the horizontally differentiated convection, the vertical
products are $\nabla_hu_3\partial_3u$ and
$u_3\nabla_h\partial_3u$. For $\nabla_hu_3\partial_3u$,
use $L^6L^3$, the uniform $H^2$ bound
\eqref{eq:new-base-energy}, and the bound for
$\nabla_h\nabla u_3$ in \eqref{eq:new-S2-structure}.
For $u_3\nabla_h\partial_3u$,
\[
 \|u_3\nabla_h\partial_3u\|_2^2
 \le\|u_3\|_\infty^2\|\nabla_h\partial_3u\|_2^2
 \le C\mathcal J_2^2
                   (\|\nabla_h^2v\|_2+\|\nabla^2w\|_2)^2,
\]
by \eqref{eq:new-elliptic} and Plancherel.
For the horizontal products, use
$\|\nabla_ha\|_4^2\le C\|a\|_\infty\|\nabla_h^2a\|_2$
and \eqref{eq:new-base-energy}.
This proves the bound for $\nabla_h(u\cdot\nabla u)$ in
\eqref{eq:new-projected-convection}.

\smallskip\noindent\emph{Scalar transport, pressure, and quadratic sources.}
By $L^\infty\cdot L^2\subset L^2$ and
$L^3\cdot L^6\subset L^2$, the scalar transport and nonlinear
pressure terms have $H^1$ norms bounded by
$C\mathcal J_2(\|\nabla^2q\|_2+\|\nabla^2w\|_2+
\|\nabla^2\tau\|_2+\|\nabla_h^2v\|_2)$.
For the differentiated quadratic sources, \eqref{eq:new-G-structure} gives
\[
\begin{aligned}
 &\|\nabla(e^{-q-\tau}Q(u))\|_2^2
          +\|\nabla(e^{-q}|\nabla\tau|^2)\|_2^2\\
 &\quad\le C(\mathcal J_2\cD_2(U)^{1/2}+\mathcal J_2^4)
       \bigl(\|\nabla^2q\|_2+\|\nabla^2w\|_2
                    +\|\nabla^2\tau\|_2+\|\nabla_h^2v\|_2\bigr)^2.
\end{aligned}
\]
A derivative on the coefficient uses
$\|\nabla e^{-q-\tau}\|_6
 \|(\nabla_hu,\Div u,\nabla\tau)\|_6^2$.
A derivative on the quadratic factor uses
$\|(\nabla_hu,\Div u,\nabla\tau)\|_\infty
 \|\nabla(\nabla_hu,\Div u,\nabla\tau)\|_2$.
The zeroth-order quadratic sources are bounded by
\eqref{eq:new-production-L2}.

In the velocity and temperature tests, integrate by parts to
move one derivative onto the dissipative factor.
The source norms required are therefore
$\|P_H\mathbb QF_2\|_{H^1}$,
$\|P_HF_3\|_{H^1}$, and $\|\nabla_hF_2\|_2$.
Use
$C\mathcal J_2\cD_2(U)^{1/2}
 \le\delta+C_\delta\mathcal J_2^2\cD_2(U)$
and choose $\delta$ and $\mathcal J_2$ sufficiently small.
Combining \eqref{eq:new-joint-coefficient},
\eqref{eq:new-second-Cq}, \eqref{eq:new-projected-convection},
and the scalar and quadratic source estimates yields
\begin{equation}\label{eq:new-joint-differential}
\begin{aligned}
 &E'+c_*(E+\|(w_H,\tau_H)\|_{H^3}^2+\|\nabla_h^3P_H^hv\|_2^2)
 \le C\cD_2(U)E+C\|C_q\|_{H^2}^2\\
 &\quad+C(1+\cD_2(U))\bigl(\|\nabla^2P_LW\|_2^2
             +\|\nabla_h^2P_L^hv\|_2^2+\|b_{\ell}\|_{W^{1,\infty}}^2\bigr).
\end{aligned}
\end{equation}
The estimates \eqref{eq:new-low-second},
\eqref{eq:new-low-horizontal-second}, and
\eqref{eq:new-second-Cq} bound the forcing in
\eqref{eq:new-joint-differential} by
$C(\mathcal J_2+M_2(t)^2)^2(1+s)^{-3}(1+\cD_2(U))$.
The integrating factor for $C\cD_2(U)E$ is bounded by
\eqref{eq:new-base-energy}.
Exponential convolution, using the total-integral case of
\eqref{eq:new-window-convolution} for the term containing
$\cD_2(U)$, gives
\[
 E(t)\le C\mathcal J_2^2e^{-c_*t}
          +C(\mathcal J_2+M_2(t)^2)^2(1+t)^{-3}.
\]
Taking square roots of this bound for $E$ and adding
\eqref{eq:new-low-second} and \eqref{eq:new-low-horizontal-second}
yields $M_2(t)\le C_0\mathcal J_2+C_1M_2(t)^2$.
For sufficiently small $\mathcal J_2$, the continuity argument
in the proof of Lemma~\ref{lem:new-seed} gives
$M_2(t)\le2C_0\mathcal J_2$ on every finite interval. In particular,
\begin{equation}\label{eq:new-second-preliminary}
 \|\nabla^2q(t)\|_2+\|\nabla^2w(t)\|_2
       +\|\nabla^2\tau(t)\|_2+\|\nabla_h^2v(t)\|_2
       \le C\mathcal J_2(1+t)^{-3/2}.
\end{equation}

Integrate \eqref{eq:new-joint-differential} over $[t,t+1]$.
The bound for $E$ obtained from that inequality controls
the initial energy $E(t)$ and gives
$\int_t^{t+1}\cD_2(U)E\dd s
 \le C\mathcal J_2^4(1+t)^{-3}$,
using \eqref{eq:new-base-energy}.
The forcing is controlled by \eqref{eq:new-low-second},
\eqref{eq:new-low-horizontal-second}, and \eqref{eq:new-second-Cq}.
Dropping the nonnegative terminal energy gives
\[
 \int_t^{t+1}
  \bigl(\|(w_H,\tau_H)\|_{H^3}^2
                      +\|\nabla_h^3P_H^hv\|_2^2\bigr)\dd s
       \le C\mathcal J_2^2(1+t)^{-3}.
\]
For low frequencies, the fixed cutoffs bound third derivatives
by the second derivatives in \eqref{eq:new-low-second}
and \eqref{eq:new-low-horizontal-second}. Adding these terms gives
\begin{equation}\label{eq:new-joint-windows}
 \int_t^{t+1}(\|\nabla_h^3u\|_2^2+\|\nabla^3(w,\tau)\|_2^2)\dd s
                                     \le C\mathcal J_2^2(1+t)^{-3}.
\end{equation}
Combining \eqref{eq:new-second-preliminary} with
\eqref{eq:new-S2-structure} and \eqref{eq:new-joint-windows}
proves the estimates for the horizontal derivatives of $u$
in the lemma.

\smallskip\noindent\emph{Step 4. Second-order decay of $q$ and $w$.}
We now improve $\|\nabla^2(q,w)\|_2$ to
$C\mathcal J_2(1+t)^{-7/4}$.
Use the compensated energy
\[
 K_2=\tfrac12\sum_{|\gamma|\le2}
 (c^2\|\partial^\gamma q_H\|_2^2+\|\partial^\gamma w_H\|_2^2)
 +\tfrac{A_\tau}{2}\|\tau_H\|_{H^1}^2
 +\eta\sum_{|\gamma|\le1}
             \ip{\partial^\gamma w_H}{\nabla\partial^\gamma q_H}.
\]
Use the parameters fixed in Lemma~\ref{lem:new-zero}, Step~2.
Then
\[
 K_2\simeq\|q_H\|_{H^2}^2+\|w_H\|_{H^2}^2+\|\tau_H\|_{H^1}^2,
\]
and its dissipation controls
$K_2+\|w_H\|_{H^3}^2+\|\tau_H\|_{H^2}^2$.
The temperature is tested in $H^1$, so its source needs
the $L^2$ bound
$|\ip{F_3}{\tau_H}_{H^1}|
 \le C\|F_3\|_2\|\tau_H\|_{H^2}$.

\smallskip\noindent\emph{The viscous term containing $\Delta_hv$.}
For $b_{\ell}=e^{-q_L}-1$, variable viscosity contributes
$\nu_hP_H\mathbb Q(b_{\ell}\Delta_hv)$ to the equation for $w_H$.
Its $H^1$ norm is needed in the $H^2$ velocity energy estimate.
The commutator structure gives an $H^1$ bound using only
second derivatives of $v$:
\begin{equation}\label{eq:new-projection-commutator}
 \|P_H\mathbb Q(b_{\ell}\Delta_hv)\|_{H^1}
 \le C\|\widehat{\nabla b_{\ell}}\|_1\|\nabla_h^2v\|_2.
\end{equation}
Since $\mathbb Q\Delta_hv=0$, the Fourier integral contains
$\mathbb Q(\xi)-\mathbb Q(\widetilde\xi)$, with
\[
 |\mathbb Q(\xi)-\mathbb Q(\widetilde\xi)|
       \le C\frac{|\xi-\widetilde\xi|}{|\xi|+|\widetilde\xi|}.
\]
On the support of $P_H$, this denominator cancels the
$H^1$ weight $\langle\xi\rangle$.
Young's convolution inequality then proves
\eqref{eq:new-projection-commutator}.
The exponential series gives
\[
 \|\widehat{\nabla b_{\ell}}\|_1
 \le C e^{C\|\widehat q_L\|_1}
                       \|\widehat{\nabla q_L}\|_1
 \le C\mathcal J_2(1+s)^{-5/4}.
\]
The fixed Fourier support of $q_L$ and
\eqref{eq:new-base-energy} give
$\|\widehat q_L\|_1\le C\mathcal J_2$.
The same support and the first-order bound
\eqref{eq:new-first} give
$\|\widehat{\nabla q_L}\|_1
 \le C\|\nabla P_LW\|_2
 \le C\mathcal J_2(1+s)^{-5/4}$.
Combining this coefficient bound with
\eqref{eq:new-second-preliminary} gives
\[
 \|P_H\mathbb Q(b_{\ell}\Delta_hv)(s)\|_{H^1}
       \le C\mathcal J_2^2(1+s)^{-11/4}.
\]

\smallskip\noindent\emph{Convection and the density source.}
The coordinatewise product inequality \eqref{eq:coordtriple} gives
\begin{equation}\label{eq:new-vertical-L3}
 \|\partial_3u\|_3+\|\partial_3v\|_3\le C\mathcal J_2 (1+t)^{-3/8}.
\end{equation}
Indeed, apply \eqref{eq:coordtriple} to three copies of each
component of $g=\partial_3u$ or $\partial_3v$ and take cube roots:
\[
 \|g\|_3\le C\|g\|_2^{1/2}
                     \prod_{i=1}^3\|\partial_i g\|_2^{1/6}.
\]
Interpolation between \eqref{eq:new-seed} and
\eqref{eq:new-base-energy} gives
$\|\partial_3u\|_2+\|\partial_3v\|_2
 \le C\mathcal J_2(1+s)^{-1/4}$.
The two horizontal derivatives of $g$ have power $3/4$ by
\eqref{eq:new-S2-structure} and \eqref{eq:new-second-preliminary};
$\partial_3g$ is uniformly bounded by
\eqref{eq:new-base-energy}.
The product inequality therefore gives decay exponent $3/8$
in \eqref{eq:new-vertical-L3}.
In \eqref{eq:new-potential-identity}, the vertical cross
products therefore have decay exponent $15/8$.
The horizontal products decay faster.
Splitting $w=w_H+w_L$ in $u\cdot\nabla w$ and using
\eqref{eq:new-low-second}, we obtain
\[
 \|P_H\mathbb Q(u\cdot\nabla u)\|_{H^1}
       \le C\mathcal J_2\|w_H\|_{H^2}
                         +C\mathcal J_2^2(1+s)^{-7/4}.
\]
For the density compensation, split $q=q_H+q_L$:
\[
 \|P_HF_1\|_{H^1}
       \le C\mathcal J_2\|q_H\|_{H^2}
                         +C\mathcal J_2^2(1+s)^{-7/4}.
\]
The high-frequency term uses the uniform $H^2$ product bound;
the low-frequency term uses the derivatives of $q_L$
in \eqref{eq:new-low-second}.
The density transport commutators retain the bound
$(\delta+C_\delta\cD_2(U))K_2$ from the argument for
\eqref{eq:new-joint-differential}.
Substituting $M_2(t)\le C\mathcal J_2$ into
\eqref{eq:new-second-Cq} gives
$\|C_q\|_{H^2}\le C\mathcal J_2^2(1+s)^{-7/4}$.

\smallskip\noindent\emph{Diffusion, pressure, and the thermal source.}
The nonlinear pressure in $H^1$ is bounded by
\[
 C\mathcal J_2(\|q_H\|_{H^2}+\|\tau_H\|_{H^2})
                         +C\mathcal J_2^2(1+s)^{-7/4}.
\]
The coefficient remainder satisfies
$\|b_{\mathrm r}\mathcal Lu\|_{H^1}
 \le C\cD_2(U)^{1/2}\|q_H\|_{H^2}$.
For the coefficient $b_{\ell}$, write
\[
P_H\mathbb Q(b_{\ell}\mathcal Lu)
  ={}P_H\mathbb Q(b_{\ell}\mathcal Lw_H)
            +P_H\mathbb Q(b_{\ell}\mathcal Lw_L)
 +\nu_hP_H[\mathbb Q,b_{\ell}]\Delta_hv.
\]
The term $P_H\mathbb Q(b_{\ell}\mathcal Lw_H)$ is bounded in
$H^1$ by $C\mathcal J_2\|w_H\|_{H^3}$ and is absorbed
by the velocity dissipation.
For $P_H\mathbb Q(b_{\ell}\mathcal Lw_L)$, use the
low-frequency derivative bounds in \eqref{eq:new-low-second}.
The commutator $P_H[\mathbb Q,b_{\ell}]\Delta_hv$ is controlled
by \eqref{eq:new-projection-commutator} and
\eqref{eq:new-second-preliminary}, with power $11/4$.
Its contribution to the $H^2$ velocity energy is bounded by
\[
\begin{aligned}
 C\mathcal J_2^2(1+s)^{-11/4}\|w_H\|_{H^3}
 &\le \varepsilon\|w_H\|_{H^3}^2
       +C_\varepsilon\mathcal J_2^4(1+s)^{-11/2}.
\end{aligned}
\]
Absorb $\varepsilon\|w_H\|_{H^3}^2$ into the dissipation.
Since $11/2>7/2$, the remaining
$C_\varepsilon\mathcal J_2^4(1+s)^{-11/2}$ term is
included in the forcing of \eqref{eq:new-K2-differential}.

The thermal convection and diffusion in $L^2$ are bounded by
\[
 C\mathcal J_2\|\tau_H\|_{H^2}
       +C\cD_2(U)^{1/2}\|q_H\|_{H^2}
       +C\mathcal J_2^2(1+s)^{-7/4}.
\]
For the quadratic source terms, \eqref{eq:new-production-L2} and
\eqref{eq:new-second-bootstrap-ledger} give
$\|Q(u)\|_2\le C\mathcal J_2^2(1+s)^{-19/8}$ and
$\||\nabla\tau|^2\|_2
 \le C\mathcal J_2^2(1+s)^{-23/8}$.
Thus Young's inequality and smallness give
\begin{equation}\label{eq:new-K2-differential}
 K_2'+c_*(K_2+\|w_H\|_{H^3}^2+\|\tau_H\|_{H^2}^2)
       \le C\cD_2(U)K_2+C\mathcal J_2^4(1+t)^{-7/2}.
\end{equation}
The bounded integrating factor from \eqref{eq:new-base-energy}
and the exponential convolution estimate
\eqref{eq:new-window-convolution} applied to
\eqref{eq:new-K2-differential} yield
\[
 K_2(t)\le C\mathcal J_2^2e^{-c_*t}
                +C\mathcal J_2^4(1+t)^{-7/2}
          \le C\mathcal J_2^2(1+t)^{-7/2}.
\]
Integrating \eqref{eq:new-K2-differential} over $[t,t+1]$
then gives
\[
 \int_t^{t+1}
       (\|w_H\|_{H^3}^2+\|\tau_H\|_{H^2}^2)\dd s
       \le C\mathcal J_2^2(1+t)^{-7/2}.
\]
Add the low-frequency estimates \eqref{eq:new-low-second}.
The fixed cutoff controls their third derivatives by their
second derivatives, so
\begin{equation}\label{eq:new-K2-consequences}
\begin{aligned}
 \|\nabla^2(q,w)(t)\|_2
       &\le C\mathcal J_2(1+t)^{-7/4},\\
 \int_t^{t+1}
       (\|\nabla^3w\|_2^2+\|\nabla^2\tau\|_2^2)\dd s
       &\le C\mathcal J_2^2(1+t)^{-7/2}.
\end{aligned}
\end{equation}
This proves the $q,w$ estimates in
\eqref{eq:new-compressible-second}.
The unit-interval bound for $\|\nabla^2\tau\|_2^2$ will control
$\||\nabla\tau|^2\|_2$ in Step~5.

\smallskip\noindent\emph{Step 5. Second-order decay of $\tau$.}
It remains to improve the preliminary pointwise bound
$\|\nabla^2\tau\|_2\le C\mathcal J_2(1+t)^{-3/2}$
in \eqref{eq:new-second-preliminary} to $(1+t)^{-7/4}$,
and to improve its third-derivative time integral.
We first estimate the quadratic source terms in $H^1$, using
\eqref{eq:new-joint-windows} and \eqref{eq:new-K2-consequences}.

\smallskip\noindent\emph{Time-integrated gradient bounds in $L^\infty$.}
Fourier Cauchy--Schwarz with
$a+b|\xi_h|^6+c\xi_3^4$, followed by optimization over
$a,b,c>0$, gives
\begin{equation}\label{eq:new-horizontal-lipschitz-interpolation}
 \|\nabla_hg\|_\infty
 \le C\|g\|_2^{1/12}\|\nabla_h^3g\|_2^{2/3}
                         \|\partial_3^2g\|_2^{1/4}.
\end{equation}
The reciprocal-weight integral with numerator $|\xi_h|^2$
is $Ca^{-1/12}b^{-2/3}c^{-1/4}$.
After integration in $\xi_3$, its finiteness follows from
$\int_0^\infty r^3(1+r^6)^{-3/4}\dd r<\infty$.
Square \eqref{eq:new-horizontal-lipschitz-interpolation},
use \eqref{eq:new-seed} and \eqref{eq:new-base-energy},
and apply H\"older on the unit interval. By
\eqref{eq:new-joint-windows},
\begin{equation}\label{eq:new-lipschitz-window-h}
 \int_t^{t+1}\|\nabla_hu\|_\infty^2\dd s
 \le C\mathcal J_2^{2/3}(1+t)^{-1/12}
       \left(\int_t^{t+1}\|\nabla_h^3u\|_2^2\dd s\right)^{2/3}
 \le C\mathcal J_2^2(1+t)^{-25/12}.
\end{equation}
For divergence and temperature, use
$\|\nabla z\|_\infty^2
 \le C\|\nabla^2z\|_2\|\nabla^3z\|_2$.
For $w$, use the pointwise and unit-interval bounds in
\eqref{eq:new-K2-consequences}.
For $\tau$, use the preliminary bounds
\eqref{eq:new-second-preliminary} and \eqref{eq:new-joint-windows}.
Since $\Div u=\Div w$, these estimates give
\[
 \int_t^{t+1}\|\Div u\|_\infty^2\dd s
       \le C\mathcal J_2^2(1+t)^{-7/2},\qquad
 \int_t^{t+1}\|\nabla\tau\|_\infty^2\dd s
       \le C\mathcal J_2^2(1+t)^{-3}.
\]

\smallskip\noindent\emph{Viscous heating and the thermal quadratic source.}
We claim
\begin{equation}\label{eq:new-heating-window}
 \begin{aligned}
 \int_t^{t+1}\|e^{-q-\tau}Q(u)\|_{H^1}^2\dd s
       &\le C\mathcal J_2^4(1+t)^{-43/12},\\
 \int_t^{t+1}\|e^{-q}|\nabla\tau|^2\|_{H^1}^2\dd s
       &\le C\mathcal J_2^4(1+t)^{-6}.
 \end{aligned}
\end{equation}
For the horizontal viscous heating term, the product rule and
the uniform bound on $e^{-q-\tau}$ give
\[
 \|e^{-q-\tau}\nabla(\mu|\nabla_hu|^2)\|_2
       \le C\|\nabla_hu\|_\infty\|\nabla\nabla_hu\|_2.
\]
The second line of \eqref{eq:new-G-structure}, together with
\eqref{eq:new-second-preliminary}, yields
\[
 \|\nabla\nabla_hu(s)\|_2^2
       \le C\mathcal J_2^2(1+s)^{-3/2}.
\]
Taking the supremum of this factor on $[t,t+1]$ and applying
\eqref{eq:new-lipschitz-window-h}, we obtain
\[
\begin{aligned}
 &\int_t^{t+1}
       \|\nabla_hu\|_\infty^2\|\nabla\nabla_hu\|_2^2\dd s
 \le \sup_{t\le s\le t+1}\|\nabla\nabla_hu(s)\|_2^2
                    \int_t^{t+1}\|\nabla_hu\|_\infty^2\dd s\\
 &\quad\le C\mathcal J_2^2(1+t)^{-3/2}
                    \int_t^{t+1}\|\nabla_hu\|_\infty^2\dd s
 \le C\mathcal J_2^4(1+t)^{-43/12},
\end{aligned}
\]
where $43/12>7/2$.
For the term involving $\nabla(\Div u)^2$, the pointwise bound
for $w$ in \eqref{eq:new-K2-consequences} gives
$\|\nabla\Div u\|_2^2\le C\mathcal J_2^2(1+s)^{-7/2}$.
Combining it with the unit-interval bound for
$\|\Div u\|_\infty^2$ just proved gives power $7$.
A derivative on $e^{-q-\tau}$ is estimated by
\[
 \|\nabla e^{-q-\tau}\|_6\|Q(u)\|_3
       \le C\mathcal J_2^3(1+s)^{-3}.
\]
The zero-order viscous term is bounded by
\eqref{eq:new-production-L2} with power $19/8$ before squaring.

For the derivative of $|\nabla\tau|^2$, combine
\eqref{eq:new-second-preliminary} with the unit-interval bound
for $\|\nabla\tau\|_\infty^2$ obtained from
\eqref{eq:new-joint-windows}. This gives
\[
 \int_t^{t+1}
       \|\nabla\tau\|_\infty^2\|\nabla^2\tau\|_2^2\dd s
       \le C\mathcal J_2^4(1+t)^{-6}.
\]
For the zero-order term $|\nabla\tau|^2$, use
\eqref{eq:new-production-L2}, the coefficient bound
\eqref{eq:new-second-bootstrap-ledger} with
$M_2(t)\le C\mathcal J_2$, and the unit-interval bound for
$\nabla^2\tau$ in \eqref{eq:new-K2-consequences}:
\[
\int_t^{t+1}\||\nabla\tau|^2\|_2^2\dd s
 \le C\sup_{t\le s\le t+1}\|\tau(s)\|_\infty^2
                           \int_t^{t+1}\|\nabla^2\tau\|_2^2\dd s
 \le C\mathcal J_2^4(1+t)^{-25/4}.
\]
The coefficient derivative is bounded by
$\|\nabla e^{-q}\|_6\|\nabla\tau\|_6^2
 \le C\mathcal J_2^3(1+s)^{-9/2}$.
Since $25/4>6$, these estimates prove
\eqref{eq:new-heating-window}.
Thus \eqref{eq:new-heating-window} uses the preliminary temperature
bound \eqref{eq:new-second-preliminary}.

\smallskip\noindent\emph{The temperature energy inequality.}
The convection estimate uses
\eqref{eq:new-vertical-L3} and the preliminary
temperature bound \eqref{eq:new-second-preliminary}.
Splitting the diffusion coefficient as in
\eqref{eq:new-coefficient-split}, we have
\[
 \|u\cdot\nabla\tau\|_{H^1}
       \le C\mathcal J_2^2(1+s)^{-15/8},\qquad
 \|b_{\mathrm r}\Delta\tau\|_{H^1}
       \le C\mathcal J_2(1+s)^{-7/4}\cD_2(U)^{1/2},
\]
\[
 \|b_{\ell}\Delta\tau_H\|_{H^1}
       \le C\mathcal J_2\|\tau_H\|_{H^3},\qquad
 \|b_{\ell}\Delta\tau_L\|_{H^1}
       \le C\mathcal J_2^2(1+s)^{-13/4}.
\]
For $b_{\mathrm r}\Delta\tau$, use the $q$ bound in
\eqref{eq:new-K2-consequences}, which controls $\|q_H\|_{H^2}$.
The $w$ bound there also gives the linear forcing estimate
$\|\Div w_H\|_{H^1}
 \le C\mathcal J_2(1+s)^{-7/4}$.
Test the high-frequency temperature equation in $H^2$,
move one derivative to the dissipative test function, and
absorb the term $C\mathcal J_2\|\tau_H\|_{H^3}$.
This gives
\begin{align}\label{eq:new-thermal-second-energy}
 \frac12\frac{\dd}{\dd t}\|\tau_H\|_{H^2}^2+c_*(\tfrac12\|\tau_H\|_{H^2}^2+\|\tau_H\|_{H^3}^2)
 \le{}&C\mathcal J_2^2(1+t)^{-7/2}(1+\cD_2(U))+C\mathcal J_2^4(1+t)^{-15/4}\notag\\
 &+C\|e^{-q-\tau}Q(u)\|_{H^1}^2+C\|e^{-q}|\nabla\tau|^2\|_{H^1}^2.
\end{align}
Apply \eqref{eq:new-window-convolution} to the two quadratic source
terms with $p=0$, $A=C\mathcal J_2^4$, and
$\beta=43/12,6$.
For $(1+s)^{-7/2}\cD_2(U)$, use its total-integral case
with $p=7/2$ and \eqref{eq:new-base-energy}.
The remaining terms use ordinary exponential convolution.
We obtain
\[
\begin{aligned}
 \|\tau_H(t)\|_{H^2}^2
 &\le C\mathcal J_2^2e^{-c_*t}
            +C\mathcal J_2^2(1+t)^{-7/2}\\
 &\quad+C\mathcal J_2^4\bigl((1+t)^{-15/4}
                       +(1+t)^{-43/12}+(1+t)^{-6}\bigr)
 \le C\mathcal J_2^2(1+t)^{-7/2}.
\end{aligned}
\]
Finally integrate \eqref{eq:new-thermal-second-energy}
over $[t,t+1]$. Use the pointwise bound for $\|\tau_H\|_{H^2}^2$
obtained from that inequality, the total dissipation bound
\eqref{eq:new-base-energy}, and the quadratic source bounds
\eqref{eq:new-heating-window}. This gives
\[
 \int_t^{t+1}\|\tau_H(s)\|_{H^3}^2\dd s
       \le C\mathcal J_2^2(1+t)^{-7/2}.
\]
By \eqref{eq:new-low-second}, the low-frequency second
derivatives of $\tau$ have power $7/4$; the fixed cutoff
also controls $\nabla^3\tau_L$.
Adding these contributions proves both temperature estimates
in \eqref{eq:new-compressible-second}.
The $q,w$ estimates proved in Step~4 complete
\eqref{eq:new-compressible-second}.
\end{proof}

\subsubsection{The third velocity component and interpolation estimates}\label{sec:r11-second-completion}
Interpolation of \eqref{eq:new-first} and the estimates in
Lemma~\ref{lem:new-joint-second} gives the $L^\infty$ decay of
$u$ and $W$. The third row of the Hodge projection supplies the
stronger component estimates for $u_3$.

\begin{lemma}[Third velocity component and $L^\infty$ decay]
For every $t\ge0$, one has
\begin{equation}\label{eq:new-uniform-basic}
 \|u\|_\infty\le C\mathcal J_2 (1+t)^{-7/8},\qquad
 \|W\|_\infty\le C\mathcal J_2 (1+t)^{-3/2}.
\end{equation}
Furthermore,
\begin{equation}\label{eq:new-third-final}
 \begin{aligned}
 \|u_3\|_\infty&\le C\mathcal J_2 (1+t)^{-3/2},
 \qquad \|\nabla_h^2u_3\|_2\le C\mathcal J_2 (1+t)^{-7/4},\\
 \|u_3\|_{L^\infty_3L^2_h}&\le C\mathcal J_2 (1+t)^{-1}.
\end{aligned}
\end{equation}
\end{lemma}
\begin{proof}
For $0\le t\le2$, all the conclusions follow from the uniform
$H^2$ bound \eqref{eq:new-base-energy}. We therefore take $t\ge2$.

\smallskip\noindent\emph{Step 1. $L^\infty$ decay of $u$ and $W$.}
Applying the third inequality in \eqref{eq:new-agmon} to $u$, and using
\eqref{eq:new-seed}, \eqref{eq:new-base-energy},
\eqref{eq:new-S2-structure}, and
\eqref{eq:new-second-preliminary}, we obtain
\[
 \|u\|_\infty
 \le C\|u\|_2^{1/4}\|\nabla_h^2u\|_2^{1/2}
                         \|\partial_3^2u\|_2^{1/4}
 \le C\mathcal J_2(1+t)^{-7/8}.
\]
For $W$, the homogeneous Agmon inequality and
\eqref{eq:new-first}, \eqref{eq:new-compressible-second} yield
\[
 \|W\|_\infty
 \le C\|\nabla W\|_2^{1/2}\|\nabla^2W\|_2^{1/2}
 \le C\mathcal J_2(1+t)^{-3/2}.
\]
This proves \eqref{eq:new-uniform-basic}.

\smallskip\noindent\emph{Step 2. The $L^\infty$ estimate for $u_3$.}
The bound for $W$ in \eqref{eq:new-uniform-basic}
already gives the required decay of $w_3$.
Since $u_3=v_3+w_3$, it remains to estimate
\[
 v_3(t)=[S_h(t)\mathbb Pu_0]_3+
          \int_0^t[S_h(t-s)\mathbb PF_2(s)]_3\,\dd s.
\]
The kernel bounds are \eqref{eq:new-third-kernels-all} and
\eqref{eq:new-mixed-kernel-all} in Lemma~\ref{lem:r13-third-row-heat}.

\smallskip\noindent\emph{Source estimates.}
Using the convection identity and removing the pressure gradient,
we write
\begin{equation}\label{eq:new-third-flux-decomposition}
\begin{aligned}
 \mathbb PF_2
 &=\mathbb P\bigl(\mathcal R_u-\Div_h(u\otimes u_h)-\partial_3(u_3u)\bigr),\\
 \mathcal R_u&=u\Div u-c^2\Pi_p+(e^{-q}-1)\mathcal Lu,
\end{aligned}
\end{equation}
where $\Pi_p$ is defined in \eqref{eq:new-N-definitions}.
The source bounds needed below are
\begin{equation}\label{eq:new-third-source-ledger}
\begin{aligned}
 &\begin{alignedat}{2}
 \|\mathcal R_u\|_1&\le C\mathcal J_2^2(1+s)^{-7/4},\quad&
 \|\mathcal R_u\|_2&\le C\mathcal J_2^2(1+s)^{-17/8},\\
 \|u\otimes u_h\|_1&\le C\mathcal J_2^2(1+s)^{-1},&
 \|\Div_h(u\otimes u_h)\|_2&\le C\mathcal J_2^2(1+s)^{-15/8},\\
 \|u_3u\|_{L^2_3L^1_h}&\le C\mathcal J_2^2(1+s)^{-11/8},&
 \|\partial_3(u_3u)\|_{L^2_3L^1_h}&\le C\mathcal J_2^2(1+s)^{-9/8},
 \end{alignedat}\\[2pt]
 &\|u_3u\|_{H^1}\le C\mathcal J_2(1+s)^{-1/4}\|u_3\|_\infty
                            +C\mathcal J_2^2(1+s)^{-15/8}.
\end{aligned}
\end{equation}
The product bounds used in Step~2 of Lemma~\ref{lem:new-joint-second}
to prove \eqref{eq:new-second-low-ledger} control $u\Div u$,
$(e^{-q}-1)\mathcal Lu$, and the horizontal flux in
\eqref{eq:new-third-source-ledger}, with the closed bound
$M_2(t)\le C\mathcal J_2$.
For the pressure remainder $\Pi_p$ in \eqref{eq:new-N-definitions},
H\"older's inequality with \eqref{eq:new-zero},
\eqref{eq:new-first}, and \eqref{eq:new-uniform-basic}
gives the faster $L^1$ and $L^2$ powers $2$ and $11/4$.
Thus $u\Div u$ determines the powers $7/4$ and $17/8$
in the first line of \eqref{eq:new-third-source-ledger}.

The undifferentiated mixed-norm bound for $u_3u$ was proved in
Step~5 of Lemma~\ref{lem:new-first}.
For its vertical derivative and $H^1$ norm, the product rule gives
\[
\begin{aligned}
 \|\partial_3(u_3u)\|_{L^2_3L^1_h}
 &\le \|u\|_{L^\infty_3L^2_h}\|\partial_3u_3\|_2
       +\|u_3\|_{L^\infty_3L^2_h}\|\partial_3u\|_2,\\
 \|u_3u\|_{H^1}
 &\le C\|u_3\|_\infty\|u\|_{H^1}
       +C\|u\|_\infty\|\nabla u_3\|_2.
\end{aligned}
\]
Vertical Sobolev interpolation, \eqref{eq:new-seed},
\eqref{eq:new-base-energy}, \eqref{eq:new-zero}, and
\eqref{eq:new-first} give decay exponents $11/8$ and
$9/8$ in the first line.
In the second line, $\|u\|_{H^1}\le C\mathcal J_2(1+s)^{-1/4}$
by \eqref{eq:new-seed} and \eqref{eq:new-base-energy};
\eqref{eq:new-uniform-basic} and \eqref{eq:new-first}
give decay exponent $15/8$ for the remainder.
Retaining $\|u_3\|_\infty$ closes the bootstrap for the first
bound in \eqref{eq:new-third-final}.

\smallskip\noindent\emph{Duhamel estimates.}
Split the time integral at $s=t/2$.
For $\mathcal R_u$, use the $L^1$ kernel in
\eqref{eq:new-third-kernels-all} on $[0,t/2]$ and the $L^2$
kernel on $[t/2,t]$. Its contribution to $\|v_3(t)\|_\infty$ is
at most
\[
\begin{aligned}
 C\mathcal J_2^2\bigg[
 &(1+t)^{-3/2}\int_0^{t/2}(1+s)^{-7/4}\dd s\\
 &+(1+t)^{-17/8}\int_0^{t/2}a^{-3/4}\dd a\bigg]
 \le C\mathcal J_2^2\bigl((1+t)^{-3/2}+(1+t)^{-15/8}\bigr).
\end{aligned}
\]
For the horizontal flux, put its derivative on the kernel on
$[0,t/2]$, and use $\Div_h(u\otimes u_h)$ in $L^2$ on $[t/2,t]$.
The corresponding bound is
\[
 C\mathcal J_2^2\bigl((1+t)^{-2}\log(e+t)
                              +(1+t)^{-13/8}\bigr).
\]
Here the logarithm comes from
$\int_0^{t/2}(1+s)^{-1}\dd s$.

For the vertical flux, the divergence-free identity
\eqref{eq:new-third-row} gives
\[
 \partial_3[S_h(a)\mathbb P(u_3u)]_3
       =-\Div_h[S_h(a)\mathbb P(u_3u)]_h.
\]
We therefore apply \eqref{eq:new-mixed-kernel-all}.
On $[0,t/2]$, the geometric mean of the bounds for
$\|u_3u\|_{L^2_3L^1_h}$ and
$\|\partial_3(u_3u)\|_{L^2_3L^1_h}$ in
\eqref{eq:new-third-source-ledger} has power
$5/4>1$.
The first inequality in \eqref{eq:new-mixed-kernel-all}
therefore bounds this contribution by
$C\mathcal J_2^2(1+t)^{-3/2}$.
For the contribution from $[t/2,t]$, set
\[
 M_\infty(t)=\mathcal J_2+
       \sup_{0\le s\le t}(1+s)^{3/2}\|u_3(s)\|_\infty.
\]
The $H^1$ source bound in \eqref{eq:new-third-source-ledger}
is then at most
$C\mathcal J_2M_\infty(t)(1+s)^{-7/4}
 +C\mathcal J_2^2(1+s)^{-15/8}$.
The two ranges of the $H^1$ kernel satisfy
\[
 \int_0^1a^{-3/4}\dd a+\int_1^{t/2}a^{-1}\dd a
       \le C(1+\log(e+t)).
\]
The contribution of the vertical flux on $[t/2,t]$ is therefore
bounded by
\[
 C\bigl(\mathcal J_2M_\infty(t)(1+t)^{-7/4}
       +\mathcal J_2^2(1+t)^{-15/8}\bigr)(1+\log(e+t)).
\]
The first bound in \eqref{eq:new-third-kernels-all}
controls the homogeneous term by $C\|u_0\|_1t^{-3/2}$.
The bound for $W$ in \eqref{eq:new-uniform-basic}
controls $w_3$.
Since $(1+t)^{-1/4}(1+\log(e+t))\le C$, multiplication by
$(1+t)^{3/2}$ and taking the supremum give
\[
 M_\infty(t)\le C_0\mathcal J_2+
                         C_1\mathcal J_2M_\infty(t).
\]
Choosing $C_1\mathcal J_2\le1/2$ closes the bootstrap
and proves the $L^\infty$ estimate for $u_3$ in
\eqref{eq:new-third-final}.

\smallskip\noindent\emph{Step 3. The estimate for $\nabla_h^2u_3$.}
We first establish the source bound
\begin{equation}\label{eq:new-horizontal-source-window}
 \int_t^{t+1}\|\nabla_hF_2(s)\|_2^2\dd s
                             \le C\mathcal J_2^4(1+t)^{-15/4}.
\end{equation}
For vertical convection, the slowest product is
\[
 \|\nabla_hu_3\,\partial_3u\|_2
 \le\|\nabla_hu_3\|_6\|\partial_3u\|_3
 \le C\mathcal J_2^2(1+s)^{-15/8}.
\]
Sobolev embedding, \eqref{eq:new-S2-structure}, and
\eqref{eq:new-second-preliminary} give decay exponent $3/2$ for
the first factor; \eqref{eq:new-vertical-L3} gives exponent $3/8$
for the second.
The other vertical product, $u_3\nabla_h\partial_3u$, has
decay exponent $9/4$, using Step~2 and
\eqref{eq:new-S2-structure}.
For horizontal convection,
\[
 \|\nabla_hu\|_4^2+\|u\|_\infty\|\nabla_h^2u\|_2
 \le C\|u\|_\infty\|\nabla_h^2u\|_2
 \le C\mathcal J_2^2(1+s)^{-19/8}.
\]
The differentiated pressure term in \eqref{eq:F2} is bounded by
\[
\begin{aligned}
 &\|\nabla_h[(e^\tau-1)\nabla(q+\tau)]\|_2\\
 &\quad\le C\|\nabla_h\tau\|_6\|\nabla(q+\tau)\|_3
       +C\|\tau\|_\infty\|\nabla_h\nabla(q+\tau)\|_2
 \le C\mathcal J_2^2(1+s)^{-13/4},
\end{aligned}
\]
by \eqref{eq:new-first}, \eqref{eq:new-compressible-second},
and \eqref{eq:new-uniform-basic}.

For viscosity, \eqref{eq:new-uniform-basic} and
\eqref{eq:new-joint-windows} give
\[
 \int_t^{t+1}\|(e^{-q}-1)\nabla_h\mathcal Lu\|_2^2\dd s
       \le C\mathcal J_2^4(1+t)^{-6}.
\]
The derivative of the coefficient satisfies
\[
\|\nabla_h(e^{-q}-1)\,\mathcal Lu\|_2
 \le C\|\nabla_hq\|_6\|\mathcal Lu\|_2^{1/2}
                                  \|\nabla\mathcal Lu\|_2^{1/2}
 \le C\mathcal J_2^{3/2}(1+s)^{-5/2}\cD_2(U)^{1/4}.
\]
After squaring and integrating, Cauchy--Schwarz and
\eqref{eq:new-base-energy} bound this contribution by
$C\mathcal J_2^4(1+t)^{-5}$.
The convection bounds, the differentiated pressure estimate,
and these two viscosity estimates control all terms of
$\nabla_hF_2$ in \eqref{eq:F2}. After squaring and integrating,
the slowest decay exponent is $15/4$, from
$\nabla_hu_3\partial_3u$.
This proves \eqref{eq:new-horizontal-source-window}.

\smallskip\noindent\emph{High horizontal frequencies.}
Apply \eqref{eq:new-maximal-horizontal} with $k=2$.
Then \eqref{eq:new-window-convolution} and
\eqref{eq:new-horizontal-source-window} yield
\[
\|\nabla_h^2P_H^hv_3(t)\|_2^2
 \le C\mathcal J_2^2e^{-c_*t}
       +C\int_0^t e^{-c_*(t-s)}\|\nabla_hF_2(s)\|_2^2\dd s
 \le C\mathcal J_2^2e^{-c_*t}
       +C\mathcal J_2^4(1+t)^{-15/4}.
\]
Thus this part is bounded by
$C\mathcal J_2(1+t)^{-15/8}$, which is faster than the required
$(1+t)^{-7/4}$ decay.

\smallskip\noindent\emph{Low horizontal frequencies.}
For a vector field $G$, a horizontal flux tensor $B$, and $a\ge0$,
Plancherel and \eqref{eq:new-third-row} give
\[
\begin{aligned}
 \|\nabla_h^2P_L^h[S_h(a)\mathbb PG]_3\|_2
       &\le C(1+a)^{-7/4}\|G\|_1,\\
 \|\nabla_h^2P_L^h[S_h(a)\mathbb P\Div_h B]_3\|_2
       &\le C(1+a)^{-9/4}\|B\|_1,\\
 \|\nabla_h^2P_L^h\partial_3[S_h(a)\mathbb PG]_3\|_2
       &\le C(1+a)^{-2}\|G\|_{L^2_3L^1_h}.
\end{aligned}
\]
In the last line, first replace the vertical derivative by
$-\Div_h[S_h(a)\mathbb PG]_h$.
These three bounds apply, respectively, to $\mathcal R_u$,
$u\otimes u_h$, and $u_3u$ in
\eqref{eq:new-third-flux-decomposition}.
On $[0,t/2]$, the source bounds
\eqref{eq:new-third-source-ledger} give the sum
\[
 C\mathcal J_2^2\bigl((1+t)^{-7/4}
              +(1+t)^{-9/4}\log(e+t)+(1+t)^{-2}\bigr).
\]
The middle term retains the logarithm from the horizontal
flux, whose $L^1$ bound has power $1$.

On $[t/2,t]$, use $\mathcal R_u$ and
$\Div_h(u\otimes u_h)$ in $L^2$ with the kernel
$(1+a)^{-1}$.
For the vertical flux, Step~2 and \eqref{eq:new-seed} give
\[
 \|u_3u\|_2\le\|u_3\|_\infty\|u\|_2
       \le C\mathcal J_2^2(1+s)^{-2},
\]
and its $L^2$ kernel is $(1+a)^{-3/2}$.
The contributions on $[t/2,t]$ are therefore bounded by
\[
 C\mathcal J_2^2
 \bigl[((1+t)^{-17/8}+(1+t)^{-15/8})\log(e+t)
                              +(1+t)^{-2}\bigr].
\]
Both time intervals give at most
$C\mathcal J_2^2(1+t)^{-7/4}$.
The homogeneous low-frequency estimate for
$\nabla_h^2[S_h(t)\mathbb Pu_0]_3$ gives
$C\mathcal J_2(1+t)^{-7/4}$.
Adding it to the two Duhamel contributions and the
high-horizontal-frequency bound proves
$\|\nabla_h^2v_3(t)\|_2\le C\mathcal J_2(1+t)^{-7/4}$.
By \eqref{eq:new-compressible-second},
$\nabla_h^2w_3$ satisfies that bound as well.
The identity $u_3=v_3+w_3$ now proves the second assertion
in \eqref{eq:new-third-final}.

\smallskip\noindent\emph{Step 4. The $L^\infty_3L^2_h$ estimate for $u_3$.}
Using the Duhamel formula from Step~2 and
\eqref{eq:new-third-flux-decomposition}, estimate $v_3$ in
$L^\infty_3L^2_h$.
Split each time integral at $t/2$ and use the last two bounds in
\eqref{eq:new-third-kernels-all} and \eqref{eq:new-mixed-kernel-all},
together with \eqref{eq:new-third-source-ledger}.
For $\mathcal R_u$, the contributions on $[0,t/2]$ and $[t/2,t]$
are bounded by $C\mathcal J_2^2(1+t)^{-1}$ and
$C\mathcal J_2^2(1+t)^{-11/8}$, respectively.
For $-\Div_h(u\otimes u_h)$, the corresponding bounds are
$C\mathcal J_2^2(1+t)^{-3/2}\log(e+t)$ and
$C\mathcal J_2^2(1+t)^{-9/8}$.
Thus each of these two source terms contributes at most
$C\mathcal J_2^2(1+t)^{-1}$.

For $-\partial_3(u_3u)$, the divergence-free identity
\eqref{eq:new-third-row} allows us to apply
\eqref{eq:new-mixed-kernel-all} to $u_3u$.
The first line of \eqref{eq:new-third-final}, proved in Step~2,
and the last line of \eqref{eq:new-third-source-ledger} give
\[
 \|u_3u\|_{H^1}\le C\mathcal J_2^2(1+s)^{-7/4}.
\]
On $[0,t/2]$, the bounds for
$\|u_3u\|_{L^2_3L^1_h}$ and
$\|\partial_3(u_3u)\|_{L^2_3L^1_h}$ in
\eqref{eq:new-third-source-ledger} give
\[
\bigl\{\|u_3u\|_{L^2_3L^1_h}
          \|\partial_3(u_3u)\|_{L^2_3L^1_h}\bigr\}^{1/2}
 \le C\mathcal J_2^2(1+s)^{-(11/8+9/8)/2}
 =C\mathcal J_2^2(1+s)^{-5/4}.
\]
The third inequality in \eqref{eq:new-mixed-kernel-all}
therefore gives
\[
 C\mathcal J_2^2\int_0^{t/2}(t-s)^{-1}(1+s)^{-5/4}\,\dd s
       \le C\mathcal J_2^2(1+t)^{-1}.
\]
On $[t/2,t]$, the $H^1$ source bound just obtained has decay
exponent $7/4$.
The fourth inequality in \eqref{eq:new-mixed-kernel-all} then gives
\[
 C\mathcal J_2^2(1+t)^{-7/4}
          \int_0^{t/2}a^{-1/2}\dd a
       \le C\mathcal J_2^2(1+t)^{-5/4}.
\]
Hence the contribution of $-\partial_3(u_3u)$ is also bounded by
$C\mathcal J_2^2(1+t)^{-1}$.
The homogeneous term is bounded by $C\mathcal J_2(1+t)^{-1}$
using the $L^1$ mixed-norm estimate in
\eqref{eq:new-third-kernels-all}. Adding the three source
contributions proves
\[
 \|v_3(t)\|_{L^\infty_3L^2_h}
 \le C(\mathcal J_2+\mathcal J_2^2)(1+t)^{-1}
 \le C\mathcal J_2(1+t)^{-1}.
\]

For $w_3$, \eqref{eq:new-zero} and \eqref{eq:new-first} give
\[
 \|w_3\|_{L^\infty_3L^2_h}
 \le C\|w\|_2^{1/2}\|\partial_3w\|_2^{1/2}
 \le C\mathcal J_2(1+t)^{-1}.
\]
Since $u_3=v_3+w_3$, we conclude that
\[
 \|u_3(t)\|_{L^\infty_3L^2_h}
 \le\|v_3(t)\|_{L^\infty_3L^2_h}
      +\|w_3(t)\|_{L^\infty_3L^2_h}
 \le C\mathcal J_2(1+t)^{-1}.
\]
This is the last assertion of \eqref{eq:new-third-final}.
\end{proof}

\paragraph{Completion of the low-order decay estimates.}
The decay estimates \eqref{eq:new-seed},
\eqref{eq:new-zero}, \eqref{eq:new-first},
\eqref{eq:new-compressible-second},
\eqref{eq:new-second-preliminary}, and \eqref{eq:new-third-final},
together with \eqref{eq:new-S2-structure}, prove
\eqref{eq:new-H2-main-rates}.

For the weak vertical bounds, interpolate the $L^2$ estimate
for $u$ in \eqref{eq:new-seed} with the uniform
$\partial_3^2u$ bound in \eqref{eq:new-base-energy}.
This gives the estimate for $\partial_3u$ in
\eqref{eq:H2-main-vertical-weak}, while
$\|\nabla_h\partial_3u\|_2^2\le
C\|\nabla_h^2u\|_2\|\partial_3^2u\|_2$
gives its mixed second-order estimate. The identity
$\partial_3u_3=\Div w-\Div_hu_h$ and its vertical derivative
give the two remaining bounds in \eqref{eq:H2-main-vertical-weak}.
Interpolation between \eqref{eq:new-seed},
\eqref{eq:new-zero}, and the $L^\infty$ bounds
\eqref{eq:new-uniform-basic}--\eqref{eq:new-third-final}
proves \eqref{eq:new-H2-main-Lp}.

The mixed-norm estimate in \eqref{eq:H2-main-mixed-and-full}
is the last bound of \eqref{eq:new-third-final}.
For $\|U\|_{H^1}$, combine \eqref{eq:new-zero} and
\eqref{eq:new-first} with the interpolation bound for $\partial_3u$
obtained from \eqref{eq:new-seed} and \eqref{eq:new-base-energy}.
For $\|W\|_{H^2}$, sum \eqref{eq:new-zero},
\eqref{eq:new-first}, and \eqref{eq:new-compressible-second}.
This proves Theorem~\ref{thm:decay-H2}(i).

For use in Section~\ref{sec:r8-refined}, we also record
the following time-integral bounds.
By \eqref{eq:new-lipschitz-window-h},
the identity $\partial_3u_3=\Div w-\Div_hu_h$, and
\eqref{eq:new-compressible-second},
\begin{equation}\label{eq:new-lipschitz-tails}
\begin{aligned}
 \int_t^\infty(\|\nabla_hu\|_\infty+\|\nabla u_3\|_\infty)\dd s
       &\le C\mathcal J_2 (1+t)^{-1/24},\\
 \int_t^\infty(\|\nabla w\|_\infty+\|\nabla\tau\|_\infty)\dd s
       &\le C\mathcal J_2 (1+t)^{-3/4}.
\end{aligned}
\end{equation}
On $[a,a+1]$, $a\ge0$, the squared norms in the first and second
lines have integrals bounded by $C\mathcal J_2^2(1+a)^{-25/12}$
and $C\mathcal J_2^2(1+a)^{-7/2}$, respectively.
Apply \eqref{eq:r13-unit-interval-tail} to the two sums of
gradient norms with $A=C\mathcal J_2$ and
$\beta=25/12,7/2$, respectively. The resulting exponents
$-1/24$ and $-3/4$ prove \eqref{eq:new-lipschitz-tails}.
These bounds do not assert integrability of
$\|\partial_3u_h\|_\infty$.

The estimates in \eqref{eq:H2-main-mixed-and-full} control
$\|W(t)\|_{H^2}$, but $\partial_3^2v$ remains only uniformly
bounded in $L^2$. Section~\ref{sec:r8-refined} uses the
time-integrated dissipation bounds to prove decay of this
derivative for each fixed initial datum and hence the full
$H^2$ convergence in \eqref{eq:r8-main-zero}.

\section{Time-integrated dissipation and full \texorpdfstring{$H^2$}{H2} convergence}\label{sec:r8-refined}

We use the second-order decay to improve the integrated dissipation,
prove the fixed-data limits in \eqref{eq:r8-main-zero}, and estimate
the first vertical derivative of the nonlinear shear correction.

\subsection{Preliminary time-integrated estimates and nonlinear source bounds}

\begin{corollary}[Preliminary time-integrated estimates]

For sufficiently small $\mathcal J_2$, the solution satisfies
\begin{align}\label{eq:H2-final-windows}
 \int_t^{t+1}\|\mathcal Lu(s)\|_2^2\dd s&\le C\mathcal J_2^2(1+t)^{-3},\notag\\
 \int_t^{t+1}\|\nabla_h^3u(s)\|_2^2\dd s&\le C\mathcal J_2^2(1+t)^{-3},\notag\\
 \int_t^{t+1}\|\nabla^3(w,\tau)(s)\|_2^2\dd s
                                    &\le C\mathcal J_2^2(1+t)^{-7/2}.
\end{align}
Furthermore,
\begin{align}\label{eq:H2-final-lipschitz-windows}
 \int_t^{t+1}(\|\nabla_hu\|_\infty^2+\|\nabla u_3\|_\infty^2)\dd s
                                    &\le C\mathcal J_2^2(1+t)^{-25/12},\notag\\
 \int_t^{t+1}(\|\nabla w\|_\infty^2+\|\nabla\tau\|_\infty^2)\dd s
                                    &\le C\mathcal J_2^2(1+t)^{-7/2}.
\end{align}
The $L^1$ time integrals of the gradient norms in
\eqref{eq:H2-final-lipschitz-windows} satisfy
\eqref{eq:new-lipschitz-tails}. With $\cD_0(U)^{1/2}$ defined in
\eqref{eq:new-d0}, one also has
\[
 \int_t^\infty \cD_0(U(s))\dd s\le C\mathcal J_2^2(1+t)^{-1},\qquad
 \int_0^\infty(1+s)^r \cD_0(U(s))\dd s\le C_r\mathcal J_2^2,
                                      \quad 0\le r<1.
\]
\end{corollary}
\begin{proof}
The bound for $\mathcal Lu$ in \eqref{eq:H2-final-windows}
follows by integrating
$\|\mathcal Lu(s)\|_2\le C\mathcal J_2(1+s)^{-3/2}$,
which is a consequence of \eqref{eq:new-S2-structure} and
\eqref{eq:new-second-preliminary}.
The time-integral bound for $\nabla_h^3u$ is
\eqref{eq:new-joint-windows}, and the bound for
$\nabla^3(w,\tau)$ is \eqref{eq:new-compressible-second}.
This proves \eqref{eq:H2-final-windows}.

The gradient estimates in \eqref{eq:H2-final-lipschitz-windows}
follow from \eqref{eq:new-lipschitz-window-h}, the identity
$\partial_3u_3=\Div w-\Div_hu_h$, and
$\|\nabla z\|_\infty^2\le C\|\nabla^2z\|_2\|\nabla^3z\|_2$
for $z=w,\tau$, with the second- and third-order bounds in
\eqref{eq:new-compressible-second}.
This proves \eqref{eq:H2-final-lipschitz-windows}; its $L^1$
time tails are \eqref{eq:new-lipschitz-tails}.
The tail and weighted-integral bounds for $\cD_0(U)$ follow from
\eqref{eq:new-low-tail} and \eqref{eq:new-seed}, as shown in
Lemma~\ref{lem:new-seed}. Only dissipative derivatives occur
in \eqref{eq:H2-final-windows}.
\end{proof}

The decomposition
\eqref{eq:new-projected-source}--\eqref{eq:new-N-definitions}
separates the source $N$ from the horizontal viscous flux.
We estimate both terms for the dissipation bounds in
Proposition~\ref{prop:r8-windows} and the vertical correction in
Proposition~\ref{prop:r8-vertical-correction}.

\begin{lemma}[Nonlinear source estimates]\label{lem:new-profile-sources}
For sufficiently small $\mathcal J_2$, the sources satisfy
\begin{equation}\label{eq:new-profile-source-ledger}
\begin{aligned}
 \|N\|_{L^2_3L^1_h}&\le C\mathcal J_2^2(1+s)^{-5/4},\\
 \|\partial_3N\|_{L^2_3L^1_h}&\le C\mathcal J_2^2(1+s)^{-1},\\
 \|\nu_h(e^{-q}-1)\nabla_hu\|_{L^2_3L^1_h}
                  &\le C\mathcal J_2^2(1+s)^{-2},\\
 \|\partial_3[\nu_h(e^{-q}-1)\nabla_hu]\|_{L^2_3L^1_h}
                  &\le C\mathcal J_2^2(1+s)^{-7/4},\\
 \|N\|_2&\le C\mathcal J_2^2(1+s)^{-7/4},\\
 \|\nu_h(e^{-q}-1)\nabla_hu\|_2
       +\|\nabla_h[\nu_h(e^{-q}-1)\nabla_hu]\|_2
                  &\le C\mathcal J_2^2(1+s)^{-5/2}.
\end{aligned}
\end{equation}
For $F_2$ in \eqref{eq:F2} and any fixed $c_*>0$, one also has
\begin{equation}\label{eq:new-F2H1}
\begin{aligned}
 \|F_2\|_{H^1}
 &\le C\mathcal J_2^2(1+s)^{-9/8}
       +C\mathcal J_2(1+s)^{-3/2}\cD_2(U)^{1/2}\\
 &\quad+C\mathcal J_2^{3/2}(1+s)^{-5/2}\cD_2(U)^{1/4},\\
 \int_0^t e^{-c_*(t-s)}\|F_2(s)\|_{H^1}^2\dd s
 &\le C\mathcal J_2^4(1+t)^{-9/4}.
\end{aligned}
\end{equation}
\end{lemma}

\begin{proof}
\smallskip\noindent\emph{Bounds for $N$ and the horizontal viscous flux.}
Write $f=e^{-q}-1$. The vertical Sobolev inequality
\eqref{eq:mixedv}, together with
\eqref{eq:new-zero}, \eqref{eq:new-first}, and
\eqref{eq:new-compressible-second}, gives
\[
\begin{aligned}
 \|z\|_{L^\infty_3L^2_h}
 &\le C\|z\|_2^{1/2}\|\partial_3z\|_2^{1/2}
       \le C\mathcal J_2(1+s)^{-1},\\
 \|\nabla z\|_{L^\infty_3L^2_h}
 &\le C\|\nabla z\|_2^{1/2}\|\partial_3\nabla z\|_2^{1/2}
       \le C\mathcal J_2(1+s)^{-3/2},
\end{aligned}
\]
for $z=q,\tau,f,e^\tau-1$. The estimates for $f$ and $e^\tau-1$
follow by differentiating the exponential compositions and
using the uniform $H^2$ bound \eqref{eq:new-base-energy}.
The velocity estimates \eqref{eq:new-seed}, \eqref{eq:new-base-energy},
\eqref{eq:new-first}, and \eqref{eq:new-third-final},
together with \eqref{eq:mixedv}, give
\[
 \|u\cdot\nabla u\|_{L^2_3L^1_h}
 \le \|u_h\|_{L^\infty_3L^2_h}\|\nabla_hu\|_2
       +\|u_3\|_{L^\infty_3L^2_h}\|\partial_3u\|_2
 \le C\mathcal J_2^2(1+s)^{-5/4}.
\]
Here the horizontal and vertical products have decay exponents
$11/8$ and $5/4$, respectively.
For the vertical derivative, expand
\[
\partial_3(u\cdot\nabla u)
 ={}(\partial_3u_h)\cdot\nabla_hu
       +u_h\cdot\nabla_h\partial_3u
 +(\partial_3u_3)\partial_3u+u_3\partial_3^2u.
\]
By \eqref{eq:mixedv}, \eqref{eq:new-base-energy},
\eqref{eq:new-first}, \eqref{eq:new-S2-structure},
\eqref{eq:new-second-preliminary}, and \eqref{eq:new-third-final},
the mixed-norm decay exponents are $9/8$ for the first three products
and $1$ for the last.
For the third product, the identity
$\partial_3^2u_3=\partial_3\Div w-\Div_h\partial_3u_h$
and Lemma~\ref{lem:new-joint-second} give
$\|\partial_3^2u_3\|_2\le C\mathcal J_2(1+s)^{-3/4}$.
Thus \eqref{eq:mixedv} gives the mixed-norm exponent
$7/8$ for $\partial_3u_3$.
For $u_3\partial_3^2u$, \eqref{eq:new-base-energy} gives
only a uniform bound on $\partial_3^2u$:
\[
 \|u_3\partial_3^2u\|_{L^2_3L^1_h}
 \le \|u_3\|_{L^\infty_3L^2_h}\|\partial_3^2u\|_2
 \le C\mathcal J_2^2(1+s)^{-1}.
\]
The product $u_3\partial_3^2u$ therefore has the slowest
decay among the four differentiated convection terms.

For the viscosity remainder in \eqref{eq:new-N-definitions},
\eqref{eq:mixedv}, \eqref{eq:new-first},
\eqref{eq:new-S2-structure}, and \eqref{eq:new-second-preliminary}
give
\[
 \|(\nabla_hu,\Div u)\|_{L^\infty_3L^2_h}
       \le C\mathcal J_2(1+s)^{-7/8},
 \qquad \|\partial_3(\nabla_hu,\Div u)\|_2
       \le C\mathcal J_2(1+s)^{-3/4}.
\]
Hence
\[
\begin{aligned}
 \|\partial_3\Pi_v\|_{L^2_3L^1_h}
 &\le C\|\nabla\partial_3f\|_2\|(\nabla_hu,\Div u)\|_{L^\infty_3L^2_h}
       +C\|\nabla f\|_{L^\infty_3L^2_h}\|\partial_3(\nabla_hu,\Div u)\|_2\\
 &\le C\mathcal J_2^2\bigl((1+s)^{-21/8}+(1+s)^{-9/4}\bigr).
\end{aligned}
\]
The pressure formula in \eqref{eq:new-N-definitions},
combined with \eqref{eq:mixedv}, \eqref{eq:new-zero},
\eqref{eq:new-first}, and \eqref{eq:new-compressible-second}, gives
$\|\Pi_p\|_{L^2_3L^1_h}\le C\mathcal J_2^2(1+s)^{-9/4}$
and
$\|\partial_3\Pi_p\|_{L^2_3L^1_h}
\le C\mathcal J_2^2(1+s)^{-11/4}$.
The bounds on $\nabla f$ and $(\nabla_hu,\Div u)$ also give
$\|\Pi_v\|_{L^2_3L^1_h}\le C\mathcal J_2^2(1+s)^{-5/2}$.
Combining the convection, pressure, and viscosity remainder
estimates proves the mixed-norm bounds on $N$ and $\partial_3N$ in
\eqref{eq:new-profile-source-ledger}.

The mixed-norm estimates for the horizontal flux follow from
\eqref{eq:mixedv}, \eqref{eq:new-first},
\eqref{eq:new-compressible-second},
\eqref{eq:new-S2-structure}, \eqref{eq:new-second-preliminary},
and the product rule
$\partial_3(f\nabla_hu)
 =(\partial_3f)\nabla_hu+f\nabla_h\partial_3u$.
The undifferentiated product has exponent $2$; the two
differentiated products have exponents $5/2$ and $7/4$.
For the $L^2$ bound on $N$, use
$\|u_h\|_\infty\|\nabla_hu\|_2+
\|u_3\|_\infty\|\partial_3u\|_2$,
whose exponents are $15/8$ and $7/4$ by
\eqref{eq:new-uniform-basic}, \eqref{eq:new-third-final},
\eqref{eq:new-first}, and the vertical interpolation of
\eqref{eq:new-seed} with \eqref{eq:new-base-energy}.
By \eqref{eq:new-N-definitions}, \eqref{eq:new-first},
\eqref{eq:new-compressible-second}, and \eqref{eq:new-uniform-basic},
$\|\Pi_p\|_2\le C\mathcal J_2^2(1+s)^{-11/4}$ and
$\|\Pi_v\|_2\le C\mathcal J_2^2(1+s)^{-21/8}$.
Finally, the $L^2$ flux bounds follow from
\[
\begin{aligned}
 \|f\nabla_hu\|_2+\|\nabla_h(f\nabla_hu)\|_2
 &\le C\|f\|_\infty
       \bigl(\|\nabla_hu\|_2+\|\nabla_h^2u\|_2\bigr)
       +C\|\nabla_hf\|_6\|\nabla_hu\|_3\\
 &\le C\mathcal J_2^2(1+s)^{-5/2}.
\end{aligned}
\]
Here \eqref{eq:new-first}, \eqref{eq:new-compressible-second},
\eqref{eq:new-S2-structure}, \eqref{eq:new-second-preliminary},
and Sobolev interpolation give
$\|\nabla_hf\|_6\le C\mathcal J_2(1+s)^{-7/4}$
and $\|\nabla_hu\|_3\le C\mathcal J_2(1+s)^{-7/8}$.
This completes \eqref{eq:new-profile-source-ledger}.

\smallskip\noindent\emph{The $H^1$ bound for $F_2$.}
For convection, the product rule and Sobolev embedding give
\[
\begin{aligned}
 \|u\cdot\nabla u\|_{H^1}
 \le C\bigl(&\|u\|_\infty\|\nabla_hu\|_{H^1}
       +\|\nabla u\|_3\|\nabla_hu\|_6\\
       &+\|u_3\|_\infty\|\partial_3u\|_{H^1}
       +\|\nabla u_3\|_6\|\partial_3u\|_3\bigr)
 \le C\mathcal J_2^2(1+s)^{-9/8}.
\end{aligned}
\]
The products $\|\nabla u\|_3\|\nabla_hu\|_6$ and
$\|\nabla u_3\|_6\|\partial_3u\|_3$ use
\eqref{eq:new-vertical-L3}, \eqref{eq:new-G-structure},
\eqref{eq:new-second-preliminary}, and
\eqref{eq:H2-main-vertical-weak}; their exponent is $9/8$.
For variable viscosity, the definition
\eqref{eq:H2-dissipation-definition} controls
$\|\mathcal Lu\|_{H^1}$ by $C\cD_2(U)^{1/2}$, while
\eqref{eq:new-S2-structure} and \eqref{eq:new-second-preliminary}
give $\|\mathcal Lu\|_2\le C\mathcal J_2(1+s)^{-3/2}$.
Thus
\[
 \|\mathcal Lu\|_3
 \le C\|\mathcal Lu\|_2^{1/2}\|\nabla\mathcal Lu\|_2^{1/2}
 \le C\mathcal J_2^{1/2}(1+s)^{-3/4}\cD_2(U)^{1/4},
\]
and
\[
\begin{aligned}
 \|f\mathcal Lu\|_{H^1}
 &\le C\|f\|_\infty\|\mathcal Lu\|_{H^1}
       +C\|\nabla f\|_6\|\mathcal Lu\|_3\\
 &\le C\mathcal J_2(1+s)^{-3/2}\cD_2(U)^{1/2}
       +C\mathcal J_2^{3/2}(1+s)^{-5/2}\cD_2(U)^{1/4}.
\end{aligned}
\]
The pressure term satisfies
\[
 \|(e^\tau-1)\nabla(q+\tau)\|_{H^1}
       \le C\mathcal J_2^2(1+s)^{-11/4}
\]
by \eqref{eq:new-first}, \eqref{eq:new-compressible-second},
and \eqref{eq:new-uniform-basic}.
Combining the convection, viscosity, and pressure estimates
proves the pointwise $H^1$ bound in \eqref{eq:new-F2H1}.

Square the pointwise bound in \eqref{eq:new-F2H1} and apply
\eqref{eq:new-window-convolution} with \eqref{eq:new-base-energy}.
The convection term and the term containing $\cD_2(U)$ contribute
$C\mathcal J_2^4(1+t)^{-9/4}$ and
$C\mathcal J_2^4(1+t)^{-3}$, respectively.
For the remaining term, Cauchy--Schwarz and the total dissipation
bound \eqref{eq:new-base-energy} give
\[
 \int_0^t e^{-c_*(t-s)}(1+s)^{-5}\cD_2(U(s))^{1/2}\dd s
       \le C\mathcal J_2(1+t)^{-5}.
\]
Multiplication by $\mathcal J_2^3$ and summation prove the
exponential convolution bound in \eqref{eq:new-F2H1}.
\end{proof}

\subsection{Improved time-integrated dissipation estimates}

\begin{proposition}[Improved time-integrated estimates]\label{prop:r8-windows}
For every $t\ge0$,
\begin{align}\label{eq:r8-windows}
 \int_t^{t+1}\|\nabla_h^3u(s)\|_2^2\dd s
   &\le C\mathcal J_2^2(1+t)^{-4},\notag\\
 \int_t^{t+1}\|\nabla_h^2\partial_3u(s)\|_2^2\dd s
   &\le C\mathcal J_2^2(1+t)^{-5/2},\\
 \int_t^{t+1}(\|\nabla_hu\|_\infty^2+\|\nabla u_3\|_\infty^2)\dd s
   &\le C\mathcal J_2^2(1+t)^{-11/4}.\notag
\end{align}
Consequently,
\begin{equation}\label{eq:r8-lips-tail}
 \int_t^\infty(\|\nabla_hu\|_\infty+\|\nabla u_3\|_\infty)\dd s
       \le C\mathcal J_2(1+t)^{-3/8}.
\end{equation}
\end{proposition}

\begin{proof}
\Needspace{6\baselineskip}
\smallskip\noindent\emph{Step 1. Horizontally differentiated source estimates.}
Coordinate scaling
in the three-dimensional Sobolev inequality gives
\[
 \|f\|_6\le C\prod_{i=1}^3\|\partial_i f\|_2^{1/3}.
\]
The established directional second-order estimates
\eqref{eq:new-S2-structure}, \eqref{eq:new-second-preliminary},
and \eqref{eq:new-third-final} therefore imply
\begin{equation}\label{eq:r8-directional-L6}
 \|\nabla_hu_3\|_6\le C\mathcal J_2(1+t)^{-5/3},\qquad
 \|\nabla_hu\|_6+\|\partial_3u_3\|_6
       \le C\mathcal J_2(1+t)^{-5/4}.
\end{equation}
The $L^2$ decay exponents for
$\partial_1\nabla_hu_3$, $\partial_2\nabla_hu_3$, and
$\partial_3\nabla_hu_3$ are $7/4,7/4,3/2$, respectively.
To estimate the product
$\nabla_hu_3\,\partial_3u$, apply H\"older's inequality and use
\eqref{eq:r8-directional-L6} and \eqref{eq:new-vertical-L3}:
\[
\|\nabla_hu_3\,\partial_3u\|_2
 \le\|\nabla_hu_3\|_6\|\partial_3u\|_3
 \le C\mathcal J_2^2(1+t)^{-5/3}(1+t)^{-3/8}
 =C\mathcal J_2^2(1+t)^{-49/24}.
\]
The remaining convection products decay at least as $(1+t)^{-9/4}$,
by \eqref{eq:new-first}, \eqref{eq:new-second-preliminary}, and
\eqref{eq:new-third-final}. Expanding the scalar sources and the
viscosity coefficient gives
\begin{align}\label{eq:r8-horizontal-sources}
 \|\nabla_hF_1(t)\|_2&\le C\mathcal J_2^2(1+t)^{-21/8},\notag\\
 \int_t^{t+1}\|\nabla_hF_2(s)\|_2^2\dd s
       &\le C\mathcal J_2^4(1+t)^{-49/12},\\
 \int_t^{t+1}\|\nabla_hF_3(s)\|_2^2\dd s
       &\le C\mathcal J_2^4(1+t)^{-61/12}.\notag
\end{align}
We next estimate the horizontally differentiated variable-viscosity
term in $F_2$:
\[
 \bigl\|\nabla_h[(e^{-q}-1)\mathcal Lu]\bigr\|_2
 \le C\mathcal J_2(1+t)^{-3/2}\|\nabla_h\mathcal Lu\|_2
 +C\mathcal J_2^{3/2}(1+t)^{-5/2}\cD_2(U)^{1/4}.
\]
Since $\Div u=\Div w$, \eqref{eq:H2-final-windows} gives
$\int_t^{t+1}\|\nabla_h\mathcal Lu(s)\|_2^2\dd s
\le C\mathcal J_2^2(1+t)^{-3}$.
Thus the viscous term containing $\|\nabla_h\mathcal Lu\|_2$ has
squared unit-interval integral at most
$C\mathcal J_2^4(1+t)^{-6}$.
For the term containing $\cD_2(U)^{1/4}$, Cauchy--Schwarz and the total dissipation
bound \eqref{eq:new-base-energy} give
\[
\mathcal J_2^3\int_t^{t+1}(1+s)^{-5}\cD_2(U(s))^{1/2}\dd s
 \le C\mathcal J_2^3(1+t)^{-5}
       \left(\int_t^{t+1}\cD_2(U(s))\dd s\right)^{1/2}
 \le C\mathcal J_2^4(1+t)^{-5}.
\]

For the heating term in $F_3$, the product rule gives
\[
 \|\nabla_hQ(u)\|_2
 \le C(\|\nabla_hu\|_\infty+\|\nabla u_3\|_\infty)
          \|\nabla_h(\nabla_hu,\Div u)\|_2.
\]
By \eqref{eq:new-S2-structure} and \eqref{eq:new-second-preliminary},
$\|\nabla_h(\nabla_hu,\Div u)\|_2^2
\le C\mathcal J_2^2(1+s)^{-3}$.
Use \eqref{eq:H2-final-lipschitz-windows} for the unit-interval
integral of $\|\nabla_hu\|_\infty^2+\|\nabla u_3\|_\infty^2$.
Consequently,
\[
\begin{aligned}
 \int_t^{t+1}\|\nabla_hQ(u(s))\|_2^2\dd s
 &\le C\mathcal J_2^2(1+t)^{-3}
       \int_t^{t+1}(\|\nabla_hu\|_\infty^2+\|\nabla u_3\|_\infty^2)\dd s\\
 &\le C\mathcal J_2^4(1+t)^{-3-25/12}
  =C\mathcal J_2^4(1+t)^{-61/12}.
\end{aligned}
\]
The heating coefficient $e^{-q-\tau}$ is uniformly bounded by
\eqref{eq:new-base-energy}. Its differentiated contribution satisfies
$\|\nabla_h(e^{-q-\tau})Q(u)\|_2
\le C\mathcal J_2^3(1+s)^{-17/4}$ by
\eqref{eq:new-compressible-second} and \eqref{eq:r8-directional-L6},
and hence has faster decay after squaring and integration.

For variable thermal diffusion, write $f=e^{-q}-1$ and expand
\[
 \nabla_h(f\Delta\tau)
 =f\nabla_h\Delta\tau+(\nabla_hf)\Delta\tau.
\]
Use \eqref{eq:new-uniform-basic} for $\|f\|_\infty$,
\eqref{eq:new-compressible-second} for $\|\nabla_hf\|_6$, and
$\|\Delta\tau\|_3\le
C\|\nabla^2\tau\|_2^{1/2}\|\nabla^3\tau\|_2^{1/2}$
for the second product. Together with \eqref{eq:H2-final-windows}
and \eqref{eq:new-base-energy}, these bounds give
\[
 \int_t^{t+1}\|\nabla_h(f\Delta\tau)(s)\|_2^2\dd s
 \le C\mathcal J_2^4\bigl((1+t)^{-13/2}+(1+t)^{-21/4}\bigr).
\]
The remaining thermal production term
$\nabla_h(e^{-q}|\nabla\tau|^2)$ has squared unit-interval
integral at most $C\mathcal J_2^4(1+t)^{-7}$, by
\eqref{eq:new-compressible-second} and
\eqref{eq:H2-final-lipschitz-windows}.

For scalar transport, H\"older's inequality gives, for $z=q,\tau$,
\[
\|\nabla_h(u\cdot\nabla z)\|_2
 \le C\|u\|_\infty\|\nabla^2z\|_2
       +C\|\nabla_hu\|_6\|\nabla z\|_3
 \le C\mathcal J_2^2\bigl((1+s)^{-21/8}+(1+s)^{-11/4}\bigr).
\]
Here use \eqref{eq:new-uniform-basic},
\eqref{eq:new-compressible-second}, and \eqref{eq:r8-directional-L6};
the bound for $\|\nabla z\|_3$ follows by interpolation between
\eqref{eq:new-first} and \eqref{eq:new-compressible-second}.
The differentiated pressure term in $F_2$ is bounded by
$C\mathcal J_2^2(1+s)^{-11/4}$, as shown in the proof of
Lemma~\ref{lem:new-profile-sources}.
Together with the convection estimates above, these bounds prove
\eqref{eq:r8-horizontal-sources}.

\Needspace{6\baselineskip}
\smallskip\noindent\emph{Step 2. The horizontal third-order dissipation.}
We first estimate the shear part $v=\mathbb Pu$, using the horizontal
cutoffs $P_L^h,P_H^h$. At low horizontal frequencies, the Fourier
calculation in the proof of Lemma~\ref{lem:linear-semigroup} gives
\[
 \|\nabla_h^jP_L^hS_h(a)\mathbb PG\|_2
 \le C(1+a)^{-1/2-j/2}\|G\|_{L^2_3L^1_h},
 \qquad a\ge0,\quad j=3,4.
\]
Taking $j=3$ and $G=u_0$ gives the initial-data bound
\[
 \|\nabla_h^3P_L^hS_h(t)\mathbb Pu_0\|_2
 \le C\mathcal J_2(1+t)^{-2}.
\]

For $t\ge2$, split the Duhamel integral at $t/2$.
On $[0,t/2]$, use the decomposition \eqref{eq:new-projected-source}
and the mixed-norm bounds \eqref{eq:new-profile-source-ledger} for
$N$ and the horizontal viscous flux
$\nu_h(e^{-q}-1)\nabla_hu$.
The corresponding kernel decay exponents are $2$ and $5/2$,
respectively: the horizontal divergence of the viscous flux
contributes one additional derivative. Since the source decay
exponents $5/4$ and $2$ exceed one, we obtain
\[
 \left\|\int_0^{t/2}\nabla_h^3P_L^hS_h(t-s)\mathbb PF_2(s)\dd s\right\|_2
 \le C\mathcal J_2^2(1+t)^{-2}.
\]

On $[t/2,t]$, place one horizontal derivative on $F_2$.
The remaining two derivatives on the low-frequency heat operator
have $L^2$ operator norm at most $C(1+t-s)^{-1}$.
By Cauchy--Schwarz and \eqref{eq:r8-horizontal-sources}, each
subinterval $I\subset[t/2,t]$ of length at most one satisfies
$\int_I\|\nabla_hF_2(s)\|_2\dd s
\le C\mathcal J_2^2(1+t)^{-49/24}$.
Summing the kernel weights over these intervals gives
\[
 \left\|\int_{t/2}^t\nabla_h^3P_L^hS_h(t-s)\mathbb PF_2(s)\dd s\right\|_2
 \le C\mathcal J_2^2(1+t)^{-49/24}\log(e+t).
\]
Since $49/24>2$, adding the initial-data and the two Duhamel
contributions yields
\[
 \|\nabla_h^3P_L^hv(t)\|_2\le C\mathcal J_2(1+t)^{-2}.
\]
For $0\le t\le2$, this low-frequency bound follows from
\eqref{eq:new-base-energy} and the fixed horizontal cutoff.

At high horizontal frequencies, testing the equation for
$\nabla_h^2P_H^hv$ and moving one horizontal derivative from the
source to the test function gives
\[
 (\|\nabla_h^2P_H^hv\|_2^2)'
 +c_*\|\nabla_h^2P_H^hv\|_2^2
 +c_*\|\nabla_h^3P_H^hv\|_2^2
 \le C\|\nabla_hF_2\|_2^2.
\]
Apply \eqref{eq:new-window-convolution} with the unit-interval
source bound \eqref{eq:r8-horizontal-sources}, and then integrate
the energy inequality over $[t,t+1]$. Since the initial energy
is at most $C\mathcal J_2^2$ by \eqref{eq:new-base-energy},
\[
 \int_t^{t+1}\|\nabla_h^3P_H^hv(s)\|_2^2\dd s
 \le C\mathcal J_2^2e^{-c_*t}+C\mathcal J_2^4(1+t)^{-49/12}.
\]
Combining the low- and high-horizontal-frequency estimates gives
\[
 \int_t^{t+1}\|\nabla_h^3v(s)\|_2^2\dd s
 \le C\mathcal J_2^2(1+t)^{-4}.
\]

For the potential part $w=\mathbb Qu$, apply two horizontal
derivatives to the compensated energy for $W=(q,w,\tau)$,
using the multiplier $\nabla(1-\Delta)^{-1}$
in the density cross term. Let $E_{h,H}$ denote this energy on
a fixed high \emph{full} frequency region; it is equivalent to
$\|\nabla_h^2(q_H,w_H,\tau_H)\|_2^2$. The compensated energy
calculation in \eqref{eq:new-joint-energy}--\eqref{eq:new-joint-differential},
with density transport kept on the left, gives
\begin{align}\label{eq:r8-compressible-henergy}
 E_{h,H}'+c_*E_{h,H}
 +c_*\|\nabla\nabla_h^2(w_H,\tau_H)\|_2^2
 \le{}&C\|\Div u\|_\infty E_{h,H}
       +C\|\nabla_h(F_1,F_2,F_3)\|_2^2\notag\\
 &+C\mathcal J_2^4(1+t)^{-17/4}+C\mathcal R_h(t),
\end{align}
where
$\int_t^{t+1}\mathcal R_h\dd s\le
C\mathcal J_2^4(1+t)^{-5}$.
To check the commutators in this formula, write
\[
 [\partial_j\partial_k,u\cdot\nabla]q
 =\partial_ju\cdot\nabla\partial_kq
  +\partial_ku\cdot\nabla\partial_jq
  +\partial_j\partial_ku\cdot\nabla q,\qquad j,k=1,2.
\]
The first two products have squared time-integrated norms with decay exponent
$67/12$. The last uses
$\|\nabla q\|_6\le C\mathcal J_2(1+t)^{-7/4}$ and
\[
 \|\nabla_h^2u\|_3
 \le C[\mathcal J_2(1+t)^{-3/2}]^{1/2}\cD_2(U)^{1/4},
\]
giving a time-integrated estimate no larger than $C\mathcal J_2^4(1+t)^{-5}$.
The high-frequency density commutator from the fixed cutoff is
bounded by $C\mathcal J_2^2(1+t)^{-17/8}$ before squaring.
For a horizontal multi-index $\beta=\gamma+e_j$, $|\beta|=2$,
$j=1,2$, the density transport is paired with $\partial_h^\beta q_H$.
Integration by parts and the expansion of
$[\nabla_h^2,u\cdot\nabla]q$ use only two density derivatives.
The contribution of $F_1$ to the compensated cross term satisfies
\[
 \begin{aligned}
 &\left|\eta\ip{\partial_h^\beta w_H}
       {\nabla(1-\Delta)^{-1}\partial_h^\beta P_HF_1}\right|\\
 &\quad=\left|\eta\ip{\partial_h^\gamma P_HF_1}
       {\partial_j\Div(1-\Delta)^{-1}\partial_h^\beta w_H}\right|
 \le\varepsilon E_{h,H}
                   +C_\varepsilon\|\nabla_hF_1\|_2^2,
\end{aligned}
\]
since $\partial_j\Div(1-\Delta)^{-1}$ is bounded on $L^2$.
For $F_2,F_3$, integrate one horizontal derivative onto the
velocity and temperature test functions. The $F_2$ cross pairing
uses the bounded multiplier $\partial_j\nabla(1-\Delta)^{-1}$.
For every $\varepsilon>0$, Young's inequality bounds their sum by
\[
 \varepsilon\bigl(E_{h,H}
       +\|\nabla\nabla_h^2(w_H,\tau_H)\|_2^2\bigr)
       +C_\varepsilon\|\nabla_h(F_2,F_3)\|_2^2.
\]
Choosing $\varepsilon$ small gives the source bound in
\eqref{eq:r8-compressible-henergy}. The coefficient $\|\Div u\|_\infty$ is integrable
by $\Div u=\Div w$ and \eqref{eq:new-lipschitz-tails}, so its
integrating factor is bounded. The source bounds
\eqref{eq:r8-horizontal-sources} and the exponential convolution
estimate \eqref{eq:new-window-convolution} then give
\[
 E_{h,H}(t)+\int_t^{t+1}\|\nabla\nabla_h^2(w_H,\tau_H)\|_2^2\dd s
 \le C\mathcal J_2^2e^{-c_*t}+C\mathcal J_2^4(1+t)^{-49/12}.
\]
On the low spatial-frequency region, apply the Green-kernel
calculation in Lemma~\ref{lem:linear-semigroup} to the initial
term and split the source integral at $s=t/2$, using
\eqref{eq:new-second-low-ledger} and
\eqref{eq:r8-horizontal-sources}. The resulting bound is
$C\mathcal J_2(1+t)^{-9/4}
+C\mathcal J_2^2(1+t)^{-49/24}\log(e+t)=o(t^{-2})$.
Combining the two projections proves the first line of
\eqref{eq:r8-windows}.

\Needspace{6\baselineskip}
\smallskip\noindent\emph{Step 3. The mixed third-order dissipation.}
We prove the second line of \eqref{eq:r8-windows}, which
contains one vertical derivative. Expanding $F_2$ in
\eqref{eq:F2} and using \eqref{eq:H2-final-windows} for the
viscous derivatives gives
\begin{equation}\label{eq:r8-F2-H1-window}
 \int_b^{b+1}\|F_2(s)\|_{H^1}^2\dd s
       \le C\mathcal J_2^4(1+b)^{-3}.
\end{equation}
The slowest term is $u_3\partial_3^2u$, whose $L^2$ norm is bounded by
$C\mathcal J_2^2(1+t)^{-3/2}$ using
\eqref{eq:new-third-final} and \eqref{eq:new-base-energy}.
The differentiated horizontal transport terms have decay exponent
$13/8$. The $H^1$ norm of the variable-viscosity term
$(e^{-q}-1)\mathcal Lu$ is bounded by
\begin{equation}\label{eq:r8-visc-H1}
 C\mathcal J_2(1+t)^{-3/2}\cD_2(U)^{1/2}
 +C\mathcal J_2^{3/2}(1+t)^{-5/2}\cD_2(U)^{1/4}.
\end{equation}
At low horizontal frequencies Fourier Cauchy--Schwarz gives
\[
 \|\nabla_h^2\partial_3P_L^hS_h(t)\mathbb Pu_0\|_2^2
 \le\|\nabla_h^4P_L^hS_h(t)\mathbb Pu_0\|_2
      \|\partial_3^2P_L^hS_h(t)\mathbb Pu_0\|_2
 \le C\mathcal J_2^2(1+t)^{-5/2}.
\]
For the nonlinear integral over $s\in[0,t/2]$, use the
second and fourth lines of \eqref{eq:new-profile-source-ledger}:
$\|\partial_3N\|_{L^2_3L^1_h}\le C\mathcal J_2^2(1+s)^{-1}$,
while the vertically differentiated horizontal viscous flux
has integrable decay $(1+s)^{-7/4}$.
The integral over $[0,t/2]$ is therefore
$O(\mathcal J_2^2(1+t)^{-3/2}\log(e+t))$.
The integral over $s\in[t/2,t]$ satisfies the same estimate by
\eqref{eq:r8-F2-H1-window}. Both contributions are bounded by
$C\mathcal J_2(1+t)^{-5/4}$. The high-frequency energy of
$\nabla_h\partial_3P_H^hv$ has right-hand side
$C\|\partial_3F_2\|_2^2$, and hence its next horizontal
derivative has a time-integrated estimate
$C\mathcal J_2^2e^{-c_*t}+C\mathcal J_2^4(1+t)^{-3}$.
For the potential part, the third line of
\eqref{eq:H2-final-windows} bounds the unit-interval integral
of $\|\nabla_h^2\partial_3w\|_2^2$ by
$C\mathcal J_2^2(1+t)^{-7/2}$. Adding the shear and potential
parts proves the second line of \eqref{eq:r8-windows}.

\Needspace{6\baselineskip}
\smallskip\noindent\emph{Step 4. Gradient estimates and their time-integral tails.}
To prove the third line of \eqref{eq:r8-windows}, use the directional interpolation inequality
\eqref{eq:new-horizontal-lipschitz-interpolation}:
\begin{equation}\label{eq:r8-Agmon-directional}
 \|\nabla_hu\|_\infty
 \le C\|u\|_2^{1/12}\|\nabla_h^3u\|_2^{2/3}
                     \|\partial_3^2u\|_2^{1/4}.
\end{equation}
Squaring and integrating \eqref{eq:r8-Agmon-directional} over a unit
interval, then applying H\"older with \eqref{eq:new-seed},
\eqref{eq:new-base-energy}, and the first line of \eqref{eq:r8-windows},
gives the decay exponent $11/4$ for $\nabla_hu$.
For $\nabla u_3$, use
$\partial_3u_3=\Div w-\Div_hu_h$ and
\[
 \|\nabla u_3\|_\infty^2
 \le C\bigl(\|\nabla_hu\|_\infty^2+\|\nabla w\|_\infty^2\bigr).
\]
The first term has the unit-interval bound just proved; the
second is bounded by the second line of
\eqref{eq:H2-final-lipschitz-windows}. This proves the third
line of \eqref{eq:r8-windows}.

Apply \eqref{eq:r13-unit-interval-tail} to
$\|\nabla_hu\|_\infty+\|\nabla u_3\|_\infty$, using the third
line of \eqref{eq:r8-windows} with $\beta=11/4$ and
$A=C\mathcal J_2$. This gives \eqref{eq:r8-lips-tail}.
\end{proof}

\subsection{An integrable mixed norm for the potential velocity}

\begin{lemma}[Mixed-norm decay for the potential part of the velocity]
\label{lem:r8-acoustic}
Let $S_c(t)$ be the compressible semigroup defined by
\eqref{eq:fourierblock} and \eqref{eq:physical-to-block}.
For $f=(f_q,f_u,f_\tau)$, write
$S_c(t)f=(q_f(t),w_f(t),\tau_f(t))$.
Define the low-frequency velocity Green matrix $K_w(t,x)$ by
\[
 P_Lw_f(t)=K_w(t)*f.
\]
Here $f_q,f_\tau$ are scalar-valued and $f_u$ is $\R^3$-valued.
Thus $K_w$ is a $3\times5$ matrix-valued convolution kernel;
$P_L$ is the fixed cutoff with respect to $|\xi|$.
For $j=0,1$ and $t\ge0$,
\begin{equation}\label{eq:r8-acoustic-kernel}
 \|\nabla^jK_w(t,\cdot)\|_{L^\infty_3L^2_h}
       \le C(1+t)^{-5/4-j/2}.
\end{equation}
The potential velocity $w=\mathbb Qu$ of the nonlinear solution satisfies
\begin{equation}\label{eq:r8-w-mixed}
 \|w(t)\|_{L^\infty_3L^2_h}
 \le C\mathcal J_2(1+t)^{-5/4}
       +C\mathcal J_2^2(1+t)^{-9/8}.
\end{equation}
\end{lemma}
\begin{proof}

\medskip\noindent\emph{Step 1. Low-frequency spectral decomposition.}
Write $\varrho=|\xi|$, $r=|\xi_h|$, and
$c_s=\sqrt{c^2(1+\alpha)}$. The characteristic polynomial of
\eqref{eq:fourierblock} is
\[
 \lambda^3+(\nu+\chi)\varrho^2\lambda^2
 +(c_s^2\varrho^2+\nu\chi\varrho^4)\lambda+c^2\chi\varrho^4,
 \qquad \nu=\nu_c+\nu_hr^2/\varrho^2.
\]
As $\varrho\to0$, the acoustic eigenvalues have the expansions
\[
 \lambda_\pm(\xi)=\pm ic_s\varrho
 -\tfrac12\left(\nu+\frac{\alpha\chi}{1+\alpha}\right)\varrho^2
 +O(\varrho^3),
\]
uniformly in $\xi/|\xi|$. The entropy eigenvalue is
$-\chi\varrho^2/(1+\alpha)+O(\varrho^3)$, and the velocity component
of its spectral projection is $O(\varrho)$.
Consequently, before multiplication by the low-frequency cutoff,
the leading velocity multiplier applied to $f=(f_q,f_u,f_\tau)$ is
\[
 e^{-a_ht r^2-a_zt\varrho^2}
 \left[\mathbb Q(\xi)\cos(c_st\varrho)\widehat f_u
 -i\frac{c^2}{c_s}\frac\xi\varrho\sin(c_st\varrho)
                  (\widehat f_q+\widehat f_\tau)\right],
\]
where $a_h=\nu_h/2$ and
$a_z=(\nu_c+\alpha\chi/(1+\alpha))/2>0$.
On a sufficiently small neighborhood of $\xi=0$, the difference
between the exact velocity multiplier and the leading acoustic
multiplier displayed above is bounded by
$C(\varrho+t\varrho^3)e^{-c_*t\varrho^2}$.
Let $K_{\rm rem}$ be the inverse Fourier transform of this remainder.
Minkowski's inequality bounds its $L^\infty_3L^2_h$ norm by the
integral in $\xi_3$ of the horizontal $L^2$ norm. Plancherel's theorem
in $x_h$ and Gaussian scaling then give
\[
 \|\nabla^jK_{\rm rem}(t)\|_{L^\infty_3L^2_h}
       \le Ct^{-3/2-j/2},\qquad t\ge1,\quad j=0,1.
\]
The part of the cutoff supported away from $\xi=0$ decays exponentially
by \eqref{eq:blockmultiplier}. It remains to estimate the leading
acoustic contribution.

\medskip\noindent\emph{Step 2. Mixed-norm estimate for the acoustic Green matrix.}
Put $R=|x|$ and $b=c_st$. The inverse Fourier transform of
$e^{-a_zt\varrho^2}\sin(b\varrho)/\varrho$ is
\[
 H(t,b,R)=\frac{e^{-(R-b)^2/(4a_zt)}-e^{-(R+b)^2/(4a_zt)}}
                     {8\pi^{3/2}R\sqrt{a_zt}}.
\]
The scalar components $f_q,f_\tau$ produce the kernel
$-(c^2/c_s)\nabla H$, while the velocity component $f_u$ produces
$\mathbb Q\partial_bH$. The factor $e^{-a_ht r^2}$ will be treated
by horizontal heat convolution.
Set $H_b=\partial_bH$ and $M_b(R)=\int_0^R s^2H_b(s)\dd s$.
Here $\mathbb QH_b$ denotes the matrix kernel with Fourier symbol
$\mathbb Q(\xi)\widehat H_b(\xi)$.
The radial representation of the Newtonian potential gives
\[
 (\mathbb QH_b)_{ij}=H_b(R)\frac{x_ix_j}{R^2}
       +\left(\delta_{ij}-3\frac{x_ix_j}{R^2}\right)\frac{M_b(R)}{R^3},
 \qquad |M_b(R)|\le C\sqrt t.
\]
For $t\ge1$ and $R\ge b/2$, each kernel and its first derivatives satisfy
\[
 |\nabla^j(-\nabla H)|+|\nabla^j(\mathbb QH_b)|
 \le Ct^{-2-j/2}e^{-(R-b)^2/(Ct)}+C\sqrt t\,R^{-3-j},
 \qquad j=0,1.
\]
The first term is concentrated near the acoustic wave front $R=b$;
the second describes the algebraic decay of the Newtonian-potential
contribution. In $R<b/2$, the kernels and their first derivatives are
exponentially small. The difference of the two Gaussians in $H$
removes the apparent singularity at $R=0$.

For each fixed $x_3$, horizontal integration of a radial envelope
$\Phi(R)$ satisfies
\[
 \int_{\R^2}\Phi\bigl(\sqrt{|x_h|^2+x_3^2}\bigr)\dd x_h
       =2\pi\int_{|x_3|}^\infty R\Phi(R)\dd R.
\]
Applying this identity to the squared pointwise bounds above gives
\[
 \begin{aligned}
 &\|\nabla^j(-\nabla H)\|_{L^\infty_3L^2_h}^2
       +\|\nabla^j(\mathbb QH_b)\|_{L^\infty_3L^2_h}^2\\
 &\quad\le Ct^{-4-j}\int_{b/2}^\infty
                   R e^{-(R-b)^2/(Ct)}\dd R
       +Ct\int_{b/2}^\infty R^{-5-2j}\dd R+Ce^{-c_*t}
 \le C\bigl(t^{-5/2-j}+t^{-3-2j}\bigr).
\end{aligned}
\]
Taking square roots yields $Ct^{-5/4-j/2}$ for $j=0,1$.
Horizontal heat convolution is contractive in $L^\infty_3L^2_h$,
and the smooth frequency cutoff is bounded on this space because
its convolution kernel belongs to $L^1$.
Together with Step~1, these facts prove
\eqref{eq:r8-acoustic-kernel} for $t\ge1$.
For $0\le t\le1$, compact frequency support and
\eqref{eq:blockmultiplier} give a uniform bound on
$\|\nabla^jK_w(t)\|_{L^\infty_3L^2_h}$, $j=0,1$.
Since $(1+t)^{-5/4-j/2}$ is bounded below on this interval,
enlarging $C$ proves \eqref{eq:r8-acoustic-kernel} for all $t\ge0$.

\Needspace{8\baselineskip}
\medskip\noindent\emph{Step 3. Duhamel estimate for the nonlinear potential velocity.}
Write $F=\mathcal R-(0,\Div(u\otimes u),0)$,
with $F=(F_1,F_2,F_3)$ and the remainder from
\eqref{eq:new-second-low-ledger}.
Let $K_{wu}$ denote the $3\times3$ block of $K_w$ acting on the
velocity component. With $U_0=(q_0,u_0,\tau_0)$, Duhamel's formula and
integration by parts in the flux term give
\[
 w_L(t)={}K_w(t)*U_0
       +\int_0^tK_w(t-s)*\mathcal R(s)\dd s
 -\sum_{k=1}^3\int_0^t
       \partial_kK_{wu}(t-s)*(u_k u)(s)\dd s.
\]
The remainder estimate \eqref{eq:new-second-low-ledger}, with the
bootstrap bound established in Lemma~\ref{lem:new-joint-second}, and
the velocity bounds \eqref{eq:new-seed} and
\eqref{eq:new-uniform-basic} imply
\[
 \begin{aligned}
 \|\mathcal R(s)\|_1&\le C\mathcal J_2^2(1+s)^{-7/4},&
 \|\mathcal R(s)\|_2&\le C\mathcal J_2^2(1+s)^{-17/8},\\
 \|(u\otimes u)(s)\|_1&\le C\mathcal J_2^2(1+s)^{-1},&
 \|(u\otimes u)(s)\|_2&\le C\mathcal J_2^2(1+s)^{-11/8}.
 \end{aligned}
\]
The momentum flux bounds follow from
$\|u\otimes u\|_1\le\|u\|_2^2$ and
$\|u\otimes u\|_2\le\|u\|_\infty\|u\|_2$.

For $t\ge2$, split both Duhamel integrals at $s=t/2$.
On $[0,t/2]$, Young's inequality and
\eqref{eq:r8-acoustic-kernel} with $j=0$ for $\mathcal R$ and $j=1$
for $u\otimes u$ bound the combined contribution by
\[
 \begin{aligned}
 &C\mathcal J_2^2(1+t)^{-5/4}
                 \int_0^{t/2}(1+s)^{-7/4}\dd s
 +C\mathcal J_2^2(1+t)^{-7/4}
                 \int_0^{t/2}(1+s)^{-1}\dd s\\
 &\le C\mathcal J_2^2\bigl[(1+t)^{-5/4}
                     +(1+t)^{-7/4}\log(e+t)\bigr].
\end{aligned}
\]
For the Duhamel integrals over $s\in[t/2,t]$, use the following
bound for a source $G=(G_q,G_u,G_\tau)$, with vector-valued $G_u$.
Plancherel's theorem in $x_h$ and Cauchy--Schwarz in $\xi_3$,
applied to \eqref{eq:blockmultiplier}, give
\[
 \|\nabla^j(K_w(a)*G)\|_{L^\infty_3L^2_h}
 \le C(1+a)^{-1/4-j/2}\|G\|_2,
 \qquad a\ge0,\quad j=0,1.
\]
Setting $a=t-s\in[0,t/2]$ and applying this bound with $j=0$ to
$\mathcal R$ and $j=1$ to the momentum flux bounds the combined
contribution by
\[
 \begin{aligned}
 &C\mathcal J_2^2(1+t)^{-17/8}
                 \int_0^{t/2}(1+a)^{-1/4}\dd a
 +C\mathcal J_2^2(1+t)^{-11/8}
                 \int_0^{t/2}(1+a)^{-3/4}\dd a\\
 &\le C\mathcal J_2^2\bigl[(1+t)^{-11/8}+(1+t)^{-9/8}\bigr].
\end{aligned}
\]
The homogeneous term is bounded by
$C\mathcal J_2(1+t)^{-5/4}$ by \eqref{eq:r8-acoustic-kernel} and
$\|U_0\|_1\le\mathcal J_2$ from \eqref{eq:new-data-sizes}.
Combining the homogeneous term with the source contributions
on $[0,t/2]$ and $[t/2,t]$ gives
\[
 \|w_L(t)\|_{L^\infty_3L^2_h}
 \le C\mathcal J_2(1+t)^{-5/4}
       +C\mathcal J_2^2(1+t)^{-9/8}.
\]

For $w_H=P_Hw$, the one-dimensional Sobolev inequality in $x_3$ and
the high-frequency cutoff imply
\[
 \|w_H\|_{L^\infty_3L^2_h}
 \le C\|w_H\|_2^{1/2}\|\partial_3w_H\|_2^{1/2}
 \le C\|\nabla^2w_H\|_2
 \le C\mathcal J_2(1+t)^{-7/4},
\]
where the last inequality is \eqref{eq:new-compressible-second}.
Adding $w_L$ and $w_H$ proves \eqref{eq:r8-w-mixed} for $t\ge2$.
The remaining interval $0\le t\le2$ follows from the uniform $H^2$
bound \eqref{eq:new-base-energy}.
\end{proof}

\subsection{Decay of the full \texorpdfstring{$H^2$}{H2} norm}

The second-order bounds leave only
$\|\partial_3^2v(t)\|_2\to0$ to be proved for each fixed initial
datum. We obtain this limit from a transport--heat comparison and
the time-integral tails of the source and flux. The corresponding refined decay estimates
for first vertical and mixed derivatives then follow by Fourier
interpolation.

\begin{lemma}[Horizontal source--flux decomposition]
\label{lem:r13-horizontal-product-flux}
For $A\in H^1(\R^3)$ and $B\in L^2(\R^3)$, there are an
$L^2$ source $R$ and a horizontal $L^2$ flux $B_h$ such that,
in the sense of distributions,
\begin{equation}\label{eq:r8-product-flux}
 AB=R+\Div_hB_h,\qquad
 \|R\|_2+\|B_h\|_2\le C\|A\|_{H^1}\|B\|_2.
\end{equation}
For vector-valued products and contractions, the decomposition
is applied componentwise.
\end{lemma}
\begin{proof}
Indeed, $\|A\|_{L^3_hL^\infty_3}\le C\|A\|_{H^1}$ follows
by integrating $\partial_3|A|^3$ and using
$\|A\|_4^2\|\partial_3A\|_2\le C\|A\|_{H^1}^3$.
Hilbert-valued horizontal interpolation then gives
\[
|\langle AB,\phi\rangle|
 \le\|A\|_{L^3_hL^\infty_3}\|B\|_2
                                      \|\phi\|_{L^6_hL^2_3}
 \le C\|A\|_{H^1}\|B\|_2
                              (\|\phi\|_2+\|\nabla_h\phi\|_2).
\]
Set $R=P_L^h(AB)$ and
$B_h=\nabla_h\Delta_h^{-1}(I-P_L^h)(AB)$.
The adjoints of these two fixed horizontal multipliers map
$L^2$ to functions with bounded $L^2$ and horizontal-gradient
norms. Testing with those adjoints proves
\eqref{eq:r8-product-flux} by duality.
\end{proof}

\begin{proposition}[Qualitative decay at the initial regularity]
\label{prop:r8-H2-zero}
For each fixed initial datum with sufficiently small $\mathcal J_2$,
\begin{equation}\label{eq:r8-H2-zero}
 \lim_{t\to\infty}\|U(t)\|_{H^2}=0.
\end{equation}
Moreover,
\begin{equation}\label{eq:r8-weak-littleo}
 \|\partial_3u(t)\|_2=o(t^{-1/4}),\qquad
 \|\nabla_h\partial_3u(t)\|_2+\|\partial_3^2u_3(t)\|_2
       =o(t^{-3/4}).
\end{equation}
\end{proposition}

\begin{proof}
Set $g=\partial_3^2v$. By \eqref{eq:new-first},
\eqref{eq:new-second-preliminary}, and
\eqref{eq:new-compressible-second}, \eqref{eq:r8-H2-zero}
reduces to proving $\|g(t)\|_2\to0$.

\Needspace{6\baselineskip}
\smallskip\noindent\emph{Step 1. Stability and decay of the comparison equation.}
For each fixed $a\ge0$, consider
\begin{equation}\label{eq:r8-transport-comparison}
 F_t-\nu_h\Delta_hF+\mathbb P(u\cdot\nabla F)=0,
 \qquad F(a)=F_a,\quad\Div F_a=0.
\end{equation}

The vertical commutator requires the following mixed-norm
coefficient bound. The two-dimensional Agmon inequality
and Cauchy--Schwarz in $x_3$ give
\[
 \|f\|_{L^2_3L^\infty_h}^2
       \le C\|f\|_2\|\nabla_h^2f\|_2.
\]
Interpolation between \eqref{eq:new-seed} and
\eqref{eq:new-base-energy} gives
$\|\partial_3u(s)\|_2\le C\mathcal J_2(1+s)^{-1/4}$.
Using the second line of \eqref{eq:r8-windows}, we obtain
\[
 \int_b^{b+1}\|\partial_3u\|_{L^2_3L^\infty_h}^2\dd s
 \le C\mathcal J_2(1+b)^{-1/4}
       \left(\int_b^{b+1}\|\nabla_h^2\partial_3u\|_2^2\dd s\right)^{1/2}
 \le C\mathcal J_2^2(1+b)^{-3/2}.
\]
The same bound holds for $v$, by the $L^2$ boundedness of
$\mathbb P$ and its commutation with derivatives.
Summing over $b=a+k$, $k\ge0$, yields
\begin{equation}\label{eq:r8-Z-tail}
 \int_a^\infty
  (\|\partial_3u\|_{L^2_3L^\infty_h}^2
       +\|\partial_3v\|_{L^2_3L^\infty_h}^2)\dd s
       \le C\mathcal J_2^2(1+a)^{-1/2}.
\end{equation}

The $L^2$ energy estimate for \eqref{eq:r8-transport-comparison}
involves $\|\Div u\|_\infty$, which is integrable by $\Div u=\Div w$ and the second line of
\eqref{eq:new-lipschitz-tails}.
Set $H^1_\sigma=\{F\in H^1(\R^3;\R^3):\Div F=0\}$,
and define $L^2_\sigma$ analogously.
For $F_a\in H^1_\sigma$, horizontal commutators use
$\|\nabla_hu\|_\infty$. The vertical commutator splits into
$(\partial_3u_h)\nabla_hF+(\partial_3u_3)\partial_3F$;
the first product satisfies
\[
 \begin{aligned}
 |\langle(\partial_3u_h)\nabla_hF,\partial_3F\rangle|
 &\le C\|\partial_3u_h\|_{L^2_3L^\infty_h}
       \|\nabla_hF\|_2^{1/2}\|\nabla_h\partial_3F\|_2^{1/2}
       \|\partial_3F\|_2\\
 &\le \frac{\nu_h}{4}\|\nabla_hF\|_{H^1}^2
       +C\|\partial_3u_h\|_{L^2_3L^\infty_h}^2\|F\|_{H^1}^2.
 \end{aligned}
\]
The coefficients $\|\nabla_hu\|_\infty$ and
$\|\partial_3u_3\|_\infty$ are integrable by
\eqref{eq:r8-lips-tail}, and the squared mixed coefficient
by \eqref{eq:r8-Z-tail}. Gronwall gives
\[
\sup_{s\ge a}\|F(s)\|_2\le C\|F_a\|_2,
 \qquad \sup_{s\ge a}\|F(s)\|_{H^1}^2
       +\int_a^\infty\|\nabla_hF\|_{H^1}^2\dd s
       \le C\|F_a\|_{H^1}^2,
\]
with $C$ independent of $a$.
By \eqref{eq:new-uniform-basic} and
\eqref{eq:new-third-final}, the convection source is integrable:
\[
 \int_a^\infty\|u\cdot\nabla F\|_2\dd s
 \le\|u_h\|_{L^2(a,\infty;L^\infty)}
       \|\nabla_hF\|_{L^2(a,\infty;L^2)}
       +\|u_3\|_{L^1(a,\infty;L^\infty)}
       \|\partial_3F\|_{L^\infty(a,\infty;L^2)}<\infty.
\]
The Fourier multiplier $e^{-\nu_ht|\xi_h|^2}$ tends to zero
for almost every $\xi$ and is bounded by one. Plancherel's theorem
and dominated convergence therefore give $S_h(t)f\to0$ in $L^2$
for $f\in L^2$.
For a fixed $T>a$, the Duhamel integral of the convection source
over $[a,T]$ tends to zero by the same convergence and the
dominating function $\|u(s)\cdot\nabla F(s)\|_2$.
The integral over $[T,t]$ is bounded by
$\int_T^\infty\|u\cdot\nabla F\|_2\dd s$, which tends to zero
as $T\to\infty$ by the preceding time-integrability estimate.
Thus \eqref{eq:r8-transport-comparison} gives $F(t)\to0$ for
$F_a\in H^1_\sigma$.
The uniform $L^2$ evolution bound and density of
$H^1_\sigma$ in $L^2_\sigma$ extend the decay to every
$F_a\in L^2_\sigma$.

\Needspace{6\baselineskip}
\smallskip\noindent\emph{Step 2. The forcing of $g=\partial_3^2v$.}
To compare $g$ with a solution of \eqref{eq:r8-transport-comparison},
we write the forcing as an $L^1(a,\infty;L^2)$ source plus the
horizontal divergence of an $L^2(a,\infty;L^2)$ flux.
These spaces match the $L^2$ energy estimate in Step~3:
the source is paired with $g-F$, and the flux with
$\nabla_h(g-F)$. The uniform $L^2$ bound for $g$ and the
time-integrated $L^2$ bound for $\nabla_hg$ follow from
\eqref{eq:new-base-energy} and the dissipation definition
\eqref{eq:H2-dissipation-definition}.

We use Lemma~\ref{lem:r13-horizontal-product-flux} for the
pressure and variable-viscosity products below.

\Needspace{11\baselineskip}
\smallskip\noindent\emph{Equation for $g$.}
Write $f=e^{-q}-1$ and
$J=(\nu_h\Delta_h\Delta^{-1}+\nu_c)\Div w$, so that
$\mathcal Lw=\nabla J$.
Using \eqref{eq:new-projected-source}, differentiate the
projected momentum equation twice in $x_3$ and retain
$\mathbb P(u\cdot\nabla g)$ on the left:
\begin{align}\label{eq:r8-top-equation}
 g_t-\nu_h\Delta_hg+\mathbb P(u\cdot\nabla g)
 =\mathbb P\big[&-2\partial_3u\cdot\nabla\partial_3v
       -\partial_3^2u\cdot\nabla v-\partial_3^2(u\cdot\nabla w)\notag\\
 &+c^2(\bigl((q+\tau)e^\tau\bigr)\nabla\partial_3^2\tau
       +2\partial_3\bigl((q+\tau)e^\tau\bigr)\nabla\partial_3\tau+\partial_3^2\bigl((q+\tau)e^\tau\bigr)\nabla\tau)\notag\\
 &+\nu_h(f\Delta_hg+2\partial_3f\Delta_h\partial_3v
                     +\partial_3^2f\Delta_hv)\\
 &-(\nabla f)\partial_3^2J+2\partial_3f\nabla\partial_3J
                     +\partial_3^2f\nabla J\big].\notag
\end{align}
In deriving the last line of \eqref{eq:r8-top-equation}, we used
\[
 \mathbb P(f\nabla\partial_3^2J)
 =-\mathbb P((\nabla f)\partial_3^2J).
\]
Thus the viscosity terms involving $J$ require at most two
derivatives of $J$.

\smallskip\noindent\emph{Convection terms.}
We estimate every convection term in \eqref{eq:r8-top-equation}
directly in $L^1(a,\infty;L^2)$.
For the two commutator products requiring mixed-norm estimates, we have
\[
\begin{aligned}
 &\|(\partial_3u_h)\cdot\nabla_h\partial_3v\|_2
       +\|(\partial_3^2u_3)\partial_3v\|_2\\
 &\qquad\le
 C\mathcal J_2^{1/2}(1+s)^{-3/8}\cD_2(U(s))^{1/4}
   \bigl(\|\partial_3u\|_{L^2_3L^\infty_h}
       +\|\partial_3v\|_{L^2_3L^\infty_h}\bigr).
\end{aligned}
\]
This follows from vertical Sobolev \eqref{eq:mixedv},
the structure bound \eqref{eq:new-S2-structure},
\eqref{eq:new-second-preliminary}, and
\eqref{eq:new-base-energy}.
For $(\partial_3^2u_3)\partial_3v$, the vertical derivative
needed in \eqref{eq:mixedv} is controlled through
$\partial_3^3u_3=\partial_3^2\Div u-\Div_h\partial_3^2u_h$.
H\"older in time with exponents $2,4,4$, using
\eqref{eq:r8-Z-tail} and \eqref{eq:new-base-energy}, gives
\[
\begin{aligned}
 &\int_a^\infty
 \bigl(\|(\partial_3u_h)\cdot\nabla_h\partial_3v\|_2
       +\|(\partial_3^2u_3)\partial_3v\|_2\bigr)\dd s\\
 &\quad\le C\mathcal J_2^{1/2}
 \left(\int_a^\infty
    (\|\partial_3u\|_{L^2_3L^\infty_h}^2
       +\|\partial_3v\|_{L^2_3L^\infty_h}^2)\dd s\right)^{1/2}\\
 &\qquad\quad{}\times
 \left(\int_a^\infty\cD_2(U(s))\dd s\right)^{1/4}
 \left(\int_a^\infty(1+s)^{-3/2}\dd s\right)^{1/4}\\
 &\quad\le C\mathcal J_2^2(1+a)^{-1/4-1/8}
       =C\mathcal J_2^2(1+a)^{-3/8}.
\end{aligned}
\]
For the other two commutator products, the uniform $H^2$
bound in \eqref{eq:new-base-energy}, the tail estimate
\eqref{eq:r8-lips-tail}, and the second line of
\eqref{eq:new-lipschitz-tails} yield
\[
 \|(\partial_3u_3)g\|_{L^1(a,\infty;L^2)}
 +\|(\partial_3^2u_h)\cdot\nabla_hv\|_{L^1(a,\infty;L^2)}
 \le C\mathcal J_2^2(1+a)^{-3/8},
\]
where $\nabla_hv=\nabla_hu-\nabla_hw$.

Expanding $\partial_3^2(u\cdot\nabla w)$ gives three
potential-velocity products, with the bounds
\[
\begin{aligned}
 \|(\partial_3^2u)\cdot\nabla w\|_{L^1(a,\infty;L^2)}
       &\le C\mathcal J_2^2(1+a)^{-3/4},\\
 \|2\partial_3u\cdot\nabla\partial_3w\|_{L^1(a,\infty;L^2)}
       &\le C\mathcal J_2^2(1+a)^{-9/8},
 \qquad \|u\cdot\nabla\partial_3^2w\|_{L^1(a,\infty;L^2)}
       \le C\mathcal J_2^2(1+a)^{-13/8}.
\end{aligned}
\]
The bound for $(\partial_3^2u)\cdot\nabla w$ follows from
\eqref{eq:new-base-energy} and the second line of
\eqref{eq:new-lipschitz-tails}.
For $2\partial_3u\cdot\nabla\partial_3w$ and
$u\cdot\nabla\partial_3^2w$, sum the unit-interval bounds
obtained from the third line of \eqref{eq:H2-final-windows},
using \eqref{eq:new-vertical-L3} and
\eqref{eq:new-uniform-basic}, respectively.

\smallskip\noindent\emph{Pressure and variable-viscosity terms.}
First separate the horizontal divergence in the leading
variable-viscosity term:
\[
 f\Delta_hg=\Div_h(f\nabla_hg)-\nabla_hf\cdot\nabla_hg.
\]
The flux $f\nabla_hg$ and the leading pressure product have
the direct bounds
\[
\|f\nabla_hg\|_{L^2(a,\infty;L^2)}
       \le C\mathcal J_2^2(1+a)^{-3/2},
 \qquad \|((q+\tau)e^\tau)\nabla\partial_3^2\tau\|_{L^1(a,\infty;L^2)}
       \le C\mathcal J_2^2(1+a)^{-1}.
\]
Both estimates use the coefficient decay in
\eqref{eq:new-uniform-basic} and the dissipation bounds
\eqref{eq:new-base-energy},
\eqref{eq:H2-dissipation-definition}.

Apply \eqref{eq:r8-product-flux} to the remaining eight products.
For the five products with one derivative on a coefficient,
choose the $H^1$ factor $A$ and the $L^2$ factor $B$ as
\[
(A,B)\in\{
 (\nabla_hf,\nabla_hg),\
 (\partial_3f,\Delta_h\partial_3v),\
 (\nabla f,\partial_3^2J),
 \qquad (\partial_3f,\nabla\partial_3J),\
 (\partial_3((q+\tau)e^\tau),\nabla\partial_3\tau)\}.
\]
For the three products with two derivatives on a coefficient,
put that coefficient in $L^2$ and the other factor in $H^1$:
\[
 (A,B)\in\{
 (\Delta_hv,\partial_3^2f),\
 (\nabla J,\partial_3^2f),\
 (\nabla\tau,\partial_3^2((q+\tau)e^\tau))\}.
\]
\Needspace{7\baselineskip}
By the composition estimates
\eqref{eq:Mosercomp}--\eqref{eq:Mosercompgrad},
\eqref{eq:new-first}, \eqref{eq:new-compressible-second}, and
\eqref{eq:new-uniform-basic}, the coefficients differentiated once
have $H^1$ norm bounded by $C\mathcal J_2(1+s)^{-5/4}$,
and those differentiated twice have $L^2$ norm bounded by
$C\mathcal J_2(1+s)^{-7/4}$.
The other factor in each pair is controlled by
$\cD_2(U(s))^{1/2}$. Hence, for each of these choices,
\[
 \|A\|_{H^1}\|B\|_2
       \le C\mathcal J_2(1+s)^{-5/4}\cD_2(U(s))^{1/2}.
\]
Consequently, each source $R$ and flux $B_h$ supplied by
\eqref{eq:r8-product-flux} satisfies
\[
\begin{aligned}
 \|R\|_{L^1(a,\infty;L^2)}
 &\le C\mathcal J_2
 \left(\int_a^\infty(1+s)^{-5/2}\dd s\right)^{1/2}
 \left(\int_a^\infty\cD_2(U(s))\dd s\right)^{1/2}
 \le C\mathcal J_2^2(1+a)^{-3/4},\\
 \|B_h\|_{L^2(a,\infty;L^2)}
 &\le C\mathcal J_2(1+a)^{-5/4}
 \left(\int_a^\infty\cD_2(U(s))\dd s\right)^{1/2}
 \le C\mathcal J_2^2(1+a)^{-5/4},
\end{aligned}
\]
where the dissipation integrals are bounded by
\eqref{eq:new-base-energy}.

Collect the directly estimated sources and the remainders $R$
from \eqref{eq:r8-product-flux} into $S$, and all horizontal
fluxes into $B_h$, with the coefficients in
\eqref{eq:r8-top-equation}. This gives
\begin{equation}\label{eq:r8-top-forcing}
 \begin{aligned}
 &g_t-\nu_h\Delta_hg+\mathbb P(u\cdot\nabla g)
       =\mathbb PS+\mathbb P\Div_hB_h,\\
 &\|S\|_{L^1(a,\infty;L^2)}+\|B_h\|_{L^2(a,\infty;L^2)}
       \le C\mathcal J_2^2(1+a)^{-3/8}.
 \end{aligned}
\end{equation}
The exponent $3/8$ in \eqref{eq:r8-top-forcing} is determined
by the convection bounds; the pressure and viscosity
contributions have faster tails.

\Needspace{6\baselineskip}
\smallskip\noindent\emph{Step 3. Decay of $\partial_3^2v$ and the full $H^2$ norm.}
Fix $a\ge0$ and let $F$ solve
\eqref{eq:r8-transport-comparison} with $F(a)=g(a)\in L^2_\sigma$.
The difference $g-F$ has zero initial value at $a$.
Testing its equation in $L^2$ places the flux $B_h$ on
$\nabla_h(g-F)$. Young's inequality and the integrability of
$\|\Div u\|_\infty$ give
\[
 \sup_{t\ge a}\|g(t)-F(t)\|_2
 \le C\bigl(\|S\|_{L^1(a,\infty;L^2)}
                   +\|B_h\|_{L^2(a,\infty;L^2)}\bigr)
 \le C\mathcal J_2^2(1+a)^{-3/8},
\]
where the last inequality is \eqref{eq:r8-top-forcing}.
Step~1 gives $F(t)\to0$ for this fixed $a$. Hence
\[
 \limsup_{t\to\infty}\|g(t)\|_2
       \le C\mathcal J_2^2(1+a)^{-3/8}.
\]
Letting $a\to\infty$ proves $\|\partial_3^2v(t)\|_2\to0$.
Fourier interpolation gives
\[
 \|U(t)\|_{H^2}
 \le C\bigl(\|W(t)\|_{H^2}+\|v(t)\|_2
              +\|\nabla_h^2v(t)\|_2+\|g(t)\|_2\bigr).
\]
Together with \eqref{eq:new-seed}, \eqref{eq:new-first},
\eqref{eq:new-compressible-second}, and
\eqref{eq:new-second-preliminary}, this proves
\eqref{eq:r8-H2-zero}.
Spatial regularization, local integrability of
$\|\nabla u\|_\infty$, and the energy bounds for
\eqref{eq:r8-transport-comparison} and \eqref{eq:r8-top-forcing}
justify the comparison for the strong solution.
The $H^1$ regularity in Step~1 is required only for the dense
comparison data; $g(a)\in L^2_\sigma$ suffices.

\Needspace{6\baselineskip}
\smallskip\noindent\emph{Step 4. The fixed-data derivative estimates.}
Step~3 gives $\|U(t)\|_{H^2}\to0$. Fourier interpolation and
the zeroth- and second-order horizontal decay estimates give
\[
\begin{aligned}
 \|\partial_3u(t)\|_2^2
 &\le\|u(t)\|_2\|\partial_3^2u(t)\|_2=o(t^{-1/2}),\\
 \|\nabla_h\partial_3u(t)\|_2^2
 &\le\|\nabla_h^2u(t)\|_2\|\partial_3^2u(t)\|_2
       =o(t^{-3/2}).
\end{aligned}
\]
Moreover, $\partial_3^2u_3=\partial_3\Div w-
\Div_h\partial_3u_h$. The second-order decay of $w$ therefore yields
\[
 \|\partial_3^2u_3(t)\|_2
 \le C\|\nabla^2w(t)\|_2
      +C\|\nabla_h\partial_3u_h(t)\|_2=o(t^{-3/4}).
\]
This proves \eqref{eq:r8-weak-littleo} and Theorem~\ref{thm:decay-H2}(ii).
\end{proof}

\subsection{Energy and dissipation tails}
\begin{proposition}[Energy and dissipation tails]
\label{prop:r8-energy-tail}
The energy $\mathfrak E_2$ of Proposition~\ref{prop:H2-closed}
is locally absolutely continuous and satisfies almost everywhere
\begin{equation}\label{eq:r8-energy-derivative}
 c_*\cD_2(U)\le-\mathfrak E_2'\le C\cD_2(U).
\end{equation}
Thus $\mathfrak E_2(a)\simeq\int_a^\infty\cD_2(U(s))\dd s$.
\end{proposition}
\begin{proof}
The lower bound in \eqref{eq:r8-energy-derivative} is
\eqref{eq:H2-closed}. For the upper bound, use the energy
identities in Lemmas~\ref{lem:zeroth-order-energy} and
\ref{prop:H2-nonlinear} before absorption: the absolute
derivatives of $\cE_{\rm th}$ and $\cE_{+,2}$ are bounded by
$C(1+\|U\|_{H^2})\cD_2(U)$. Summing
\eqref{eq:crosscalculation} for $|\beta|\le1$ gives the compensation identity
\begin{align*}
 \cX_2'=\sum_{|\beta|\le1}\big[&-c^2\|\nabla\partial^\beta q\|_2^2
 -c^2\langle\nabla\partial^\beta\tau,\nabla\partial^\beta q\rangle
 +\langle\partial^\beta \mathcal Lu,\nabla\partial^\beta q\rangle\\
 &+\|\Div\partial^\beta u\|_2^2
 +\langle\partial^\beta F_2,\nabla\partial^\beta q\rangle
 -\langle\Div\partial^\beta u,\partial^\beta F_1\rangle\big].
\end{align*}
The dissipation definition \eqref{eq:H2-dissipation-definition}
and the source bound \eqref{eq:H2-NLsource} imply
$|\cX_2'|\le C(1+\|U\|_{H^2})\cD_2(U)$.
Combining the three energy derivatives proves
\eqref{eq:r8-energy-derivative}. The regularization argument used for
\eqref{eq:H2-global-exact} gives local
absolute continuity of $\mathfrak E_2$.

By \eqref{eq:r8-H2-zero}, $\mathfrak E_2(t)\to0$.
Integrating \eqref{eq:r8-energy-derivative} from $a$ to infinity gives
\[
 c_*\int_a^\infty\cD_2(U(s))\dd s
 \le\mathfrak E_2(a)
 \le C\int_a^\infty\cD_2(U(s))\dd s.\qedhere
\]
\end{proof}

\begin{corollary}[Energy balance for the second vertical shear derivative]

Assume that $\mathcal J_2$ is sufficiently small.
With $g=\partial_3^2v$, for every $t\ge0$,
\begin{equation}\label{eq:r8-top-balance}
 \left|2\nu_h\int_t^\infty\|\nabla_hg(s)\|_2^2\dd s
                                  -\|g(t)\|_2^2\right|
       \le C\mathcal J_2(1+t)^{-3/8}\mathfrak E_2(t).
\end{equation}
Also
$\mathfrak E_2(t)\le C[\|g(t)\|_2^2+
\mathcal J_2^2(1+t)^{-1}]$.
\end{corollary}
\begin{proof}
\emph{Step 1. Source tails with the remaining energy.}
Propositions~\ref{prop:H2-closed} and \ref{prop:r8-energy-tail} give
$\|U(s)\|_{H^2}^2\le C\mathfrak E_2(s)$,
$\mathfrak E_2(s)\le\mathfrak E_2(a)$ for $s\ge a$, and
$\int_a^\infty\cD_2(U)\dd s\le C\mathfrak E_2(a)$.
Retaining this energy in the vertical interpolations used in
\eqref{eq:H2-main-vertical-weak} yields
\[
 \begin{aligned}
 \|\partial_3u(s)\|_2
       &\le C\mathcal J_2^{1/2}\mathfrak E_2(s)^{1/4}(1+s)^{-1/4},\\
 \|\nabla_h\partial_3u(s)\|_2+\|\partial_3^2u_3(s)\|_2
       &\le C\mathcal J_2^{1/2}\mathfrak E_2(s)^{1/4}(1+s)^{-3/4}.
 \end{aligned}
\]
The proof of \eqref{eq:r8-Z-tail}, with the first interpolation
and \eqref{eq:r8-windows}, then gives
\[
 \|\partial_3u\|_{L^2(a,\infty;L^2_3L^\infty_h)}
 +\|\partial_3v\|_{L^2(a,\infty;L^2_3L^\infty_h)}
       \le C\mathcal J_2^{3/4}\mathfrak E_2(a)^{1/8}(1+a)^{-1/4}.
\]
For $f=\nabla_h\partial_3v$, $\partial_3^2u_3$, or
$\nabla\partial_3w$, vertical Sobolev and
\eqref{eq:H2-dissipation-definition} give
\[
 \|f(s)\|_{L^\infty_3L^2_h}
       \le C\mathcal J_2^{1/4}\mathfrak E_2(s)^{1/8}
                  (1+s)^{-3/8}\cD_2(U(s))^{1/4}.
\]
Repeat the source estimates for \eqref{eq:r8-top-forcing} in Step~2
of the proof of Proposition~\ref{prop:r8-H2-zero}.
H\"older's inequality in time with exponents $2,4,4$ uses the preceding
two bounds and the remaining dissipation.
The pressure and viscosity terms use the same dissipation tail;
the product of $u$ with a third potential derivative uses
\eqref{eq:new-H2-main-Lp} with $p=\infty$ and time Cauchy--Schwarz.
Thus the sources $S,B_h$ in \eqref{eq:r8-top-forcing} satisfy
\[
 \|S\|_{L^1(a,\infty;L^2)}+\|B_h\|_{L^2(a,\infty;L^2)}
       \le C\mathcal J_2\mathfrak E_2(a)^{1/2}(1+a)^{-3/8}.
\]

\smallskip\noindent
\emph{Step 2. Integration of the shear energy identity.}
Testing the regularized equation \eqref{eq:r8-top-forcing} with the
solenoidal field $g$ and passing to the limit gives
\[
 \frac12\frac{\dd}{\dd s}\|g\|_2^2+\nu_h\|\nabla_hg\|_2^2
       =\frac12\int_{\R^3}\Div u\,|g|^2\dd x
          +\langle S,g\rangle-\langle B_h,\nabla_hg\rangle.
\]
Proposition~\ref{prop:r8-energy-tail} implies
\[
 \sup_{s\ge t}\|g(s)\|_2^2
             +\int_t^\infty\|\nabla_hg(s)\|_2^2\dd s
       \le C\mathfrak E_2(t).
\]
The refined source tails from Step~1 control the two source pairings.
For the transport term, \eqref{eq:r8-lips-tail} gives
$\int_t^\infty\|\Div u\|_\infty\dd s
\le C\mathcal J_2(1+t)^{-3/8}$.
Hence the absolute time integral of the right side of the energy
identity is bounded by
$C\mathcal J_2(1+t)^{-3/8}\mathfrak E_2(t)$.
Since $g(s)\to0$ by \eqref{eq:r8-H2-zero}, integration to infinity
proves \eqref{eq:r8-top-balance}.

Finally, the directional norm equivalence and $u=v+w$ give
\[
 \|u\|_{H^2}^2\le C\big(\|u\|_2^2+\|\nabla_h^2u\|_2^2
                         +\|g\|_2^2+\|w\|_{H^2}^2\big).
\]
The energy equivalence in Proposition~\ref{prop:H2-closed} and
\eqref{eq:new-H2-main-rates} with $k=0,1,2$ therefore yield
$\mathfrak E_2(t)\le C[\|g(t)\|_2^2+\mathcal J_2^2(1+t)^{-1}]$.
\end{proof}

\subsection{The vertical derivative of the nonlinear remainder}
Let $r_v=v-S_h(t)\mathbb Pu_0$ denote the nonlinear correction to
the shear velocity $v=\mathbb Pu$.
To estimate its vertical derivative using horizontal diffusion,
we rewrite the contribution of $v_3\partial_3^2u$ with an auxiliary
horizontal vector field satisfying $\Div_hb_h=v_3$.
The contribution of $w_3\partial_3^2u$ is controlled by
\eqref{eq:r8-w-mixed}.

\begin{proposition}[Vertical derivative of the nonlinear remainder]
\label{prop:r8-vertical-correction}
If $\mathcal J_2$ is sufficiently small, then for $j=0,1$ and $t\ge0$,
\begin{equation}\label{eq:r8-sqrtlog}
 \|\nabla_h^j\partial_3r_v(t)\|_2
 \le C\mathcal J_2^2(1+t)^{-(1+j)/2}\sqrt{\log(e+t)}.
\end{equation}
For each fixed solution the left side is
$o(t^{-(1+j)/2}\sqrt{\log t})$.
\end{proposition}
\begin{proof}
\Needspace{6\baselineskip}
\smallskip\noindent\emph{Step 1. An auxiliary horizontal vector field.}
Define $b_h=\nabla_h\Delta_h^{-1}v_3$, where
$\Delta_h=\partial_1^2+\partial_2^2$ and the inverse is defined
by horizontal frequency truncation. We claim that
\begin{equation}\label{eq:r8-lifting}
 \Div_hb_h=v_3,\qquad
 \|b_h(t)\|_{L^\infty_3L^2_h}\le C\mathcal J_2(1+t)^{-1/2}.
\end{equation}
To estimate $b_h$, fix $r=|\xi_h|>0$. The vertical kernels of
$\nabla_h\Delta_h^{-1}\mathbb P_{3\bullet}$ are bounded
combinations of $\operatorname{sgn}(x_3)e^{-r|x_3|}$ and
$e^{-r|x_3|}$. Their vertical derivatives are finite measures
with uniformly bounded total variation, including the Dirac
mass at zero. After horizontal heat evolution, an $L^1$
source therefore contributes at most $Ca^{-1/2}\|f\|_1$ to
the $L^\infty_3L^2_h$ norm. A source of the form $\partial_3f$
contributes at most $Ca^{-1/2}\|f\|_{L^\infty_3L^1_h}$,
by placing the vertical derivative on the kernel.

Apply these bounds to the projected flux decomposition
\eqref{eq:new-third-flux-decomposition}. The first bound in
\eqref{eq:new-third-source-ledger} controls $\mathcal R_u$.
The other two source bounds are
\[
\begin{aligned}
 \|\Div_h(u\otimes u_h)\|_1
 &\le C\|u\|_2\|\nabla_hu\|_2
       \le C\mathcal J_2^2(1+s)^{-3/2},\\
 \|u_3u\|_{L^\infty_3L^1_h}
 &\le \|u_3\|_{L^\infty_3L^2_h}\|u\|_{L^\infty_3L^2_h}
       \le C\mathcal J_2^2(1+s)^{-11/8}.
\end{aligned}
\]
The first line follows from \eqref{eq:new-seed} and
\eqref{eq:new-first}; the second follows from \eqref{eq:mixedv},
\eqref{eq:new-seed}, \eqref{eq:new-base-energy}, and
\eqref{eq:new-third-final}.
The initial-data contribution is
$Ct^{-1/2}\|u_0\|_1\le C\mathcal J_2t^{-1/2}$ by
\eqref{eq:new-data-sizes}. Hence the Duhamel formula yields,
for $t\ge2$,
\[
\begin{aligned}
 \|b_h(t)\|_{L^\infty_3L^2_h}
 &\le C\mathcal J_2t^{-1/2}
   +C\mathcal J_2^2\int_0^t(t-s)^{-1/2}\\
 &\qquad\times\bigl[(1+s)^{-7/4}+(1+s)^{-3/2}
                         +(1+s)^{-11/8}\bigr]\dd s
 \le C\mathcal J_2(1+t)^{-1/2}.
\end{aligned}
\]
Indeed, each source exponent is greater than $1$, and splitting
the convolution at $s=t/2$ gives the last inequality.
For $0\le t\le2$, the low horizontal frequencies of $b_h(0)$
are controlled by $\|u_0\|_1$, and the high horizontal frequencies
by $\|u_0\|_{H^2}$ and vertical Sobolev. The nonlinear integral
uses the integrable kernel $a^{-1/2}$.
These bounds are uniform in the frequency truncation and
justify its removal. The identity $\Div_hb_h=v_3$ follows
from the Fourier multiplier, completing \eqref{eq:r8-lifting}.

\Needspace{6\baselineskip}
\smallskip\noindent\emph{Step 2. Decomposition of the differentiated source.}
Using $\Div_hb_h=v_3$ and $u_3=v_3+w_3$ in
\eqref{eq:new-projected-source}, we obtain
\[
\partial_3\mathbb PF_2
 =\mathbb P\bigl(E+b_h\cdot\nabla_h\partial_3^2u\bigr)
 +\mathbb P\Div_h\bigl(
    \partial_3[\nu_h(e^{-q}-1)\nabla_hu]
                      -\partial_3^2u\otimes b_h\bigr),
\]
where, with $\Pi_p,\Pi_v$ defined in
\eqref{eq:new-N-definitions},
\[
E={}-\partial_3u_h\cdot\nabla_hu
       -u_h\cdot\nabla_h\partial_3u
       -\partial_3u_3\partial_3u
 -w_3\partial_3^2u-c^2\partial_3\Pi_p+\partial_3\Pi_v.
\]
For $E$, the three differentiated convection products and
$\partial_3\Pi_p,\partial_3\Pi_v$ have the bounds established
in the proof of Lemma~\ref{lem:new-profile-sources} for
\eqref{eq:new-profile-source-ledger}.
The remaining product $w_3\partial_3^2u$ is controlled by
\eqref{eq:r8-w-mixed} and the uniform $H^2$ bound
\eqref{eq:new-base-energy}. H\"older's inequality in $x_h$, together with
\eqref{eq:r8-lifting}, therefore gives
\[
\begin{aligned}
 \|E(s)\|_{L^2_3L^1_h}
 &\le C\mathcal J_2^2(1+s)^{-9/8},
 \qquad \|\partial_3^2u\otimes b_h\|_{L^2_3L^1_h}
 \le C\mathcal J_2^2(1+s)^{-1/2},\\
 \|b_h\cdot\nabla_h\partial_3^2u\|_{L^2_3L^1_h}
 &\le C\mathcal J_2(1+s)^{-1/2}
                         \|\nabla_h\partial_3^2u(s)\|_2.
\end{aligned}
\]
The fourth line of \eqref{eq:new-profile-source-ledger} also gives
\[
 \|\partial_3[\nu_h(e^{-q}-1)\nabla_hu]\|_{L^2_3L^1_h}
                  \le C\mathcal J_2^2(1+s)^{-7/4}.
\]

\Needspace{6\baselineskip}
\smallskip\noindent\emph{Step 3. Estimates for the Duhamel integral.}
Since $r_v(0)=0$, we have
\[
 \partial_3r_v(t)=\int_0^t S_h(t-s)\mathbb P\partial_3F_2(s)\dd s.
\]
Use the fixed horizontal cutoffs $P_L^h,P_H^h$ from
Lemma~\ref{lem:linear-semigroup}, and first assume $t\ge2$.
At low horizontal frequencies and $s\in[0,t/2]$, apply
\eqref{eq:Sh-low-all-orders} with $m=0,1,2$ to the decomposition in Step~2.
The norms of $E$ and the differentiated viscous flux are
integrable in time. For the other horizontal divergence term,
the extra derivative yields
\[
 C\mathcal J_2^2(1+t)^{-(2+j)/2}
       \int_0^{t/2}(1+s)^{-1/2}\dd s
       \le C\mathcal J_2^2(1+t)^{-(1+j)/2}.
\]
Consequently, for $j=0,1$,
\[
\begin{aligned}
 &\left\|\nabla_h^jP_L^h\int_0^{t/2}
      S_h(t-s)\mathbb P\partial_3F_2(s)\dd s\right\|_2\\
 &\quad\le C\mathcal J_2^2(1+t)^{-(1+j)/2}
       +C\mathcal J_2(1+t)^{-(1+j)/2}
          \int_0^{t/2}(1+s)^{-1/2}
                       \|\nabla_h\partial_3^2u(s)\|_2\dd s.
\end{aligned}
\]

For the low horizontal frequencies and $s\in[t/2,t]$, use the
complete source $\mathbb P\partial_3F_2$, with
\[
 \|\nabla_h^jP_L^hS_h(a)\mathbb P\partial_3F_2(s)\|_2
       \le C(1+a)^{-j/2}\|F_2(s)\|_{H^1}.
\]
Partition $[t/2,t]$ into intervals of length at most one.
Cauchy--Schwarz on each interval and
\eqref{eq:r8-F2-H1-window} give
\[
\begin{aligned}
 &\left\|\nabla_h^jP_L^h\int_{t/2}^t
         S_h(t-s)\mathbb P\partial_3F_2(s)\dd s\right\|_2\\
 &\quad\le C\mathcal J_2^2(1+t)^{-3/2}
         \sum_{k=0}^{\lceil t/2\rceil}(1+k)^{-j/2}
       \le C\mathcal J_2^2(1+t)^{-(1+j)/2}.
\end{aligned}
\]
Here $k$ indexes the intervals on which $t-s$ lies between
$k$ and $k+1$; throughout this time range, $1+s$ is comparable
to $1+t$.

For the high horizontal frequencies, use the exponential
heat bound when $j=0$, and \eqref{eq:new-maximal-horizontal}
with $k=1$ and $G=\partial_3F_2$ when $j=1$.
The unit-interval estimate \eqref{eq:r8-F2-H1-window}, summed
with the exponential weight by \eqref{eq:new-window-convolution},
then gives
\[
\|\nabla_h^jP_H^h\partial_3r_v(t)\|_2^2
 \le C\int_0^t e^{-c_*(t-s)}\|F_2(s)\|_{H^1}^2\dd s
 \le C\mathcal J_2^4(1+t)^{-3},\qquad j=0,1.
\]
Adding the low and high horizontal frequencies yields
\[
\|\nabla_h^j\partial_3r_v(t)\|_2
 \le C\mathcal J_2^2(1+t)^{-(1+j)/2}
 +C\mathcal J_2(1+t)^{-(1+j)/2}
       \int_0^{t/2}(1+s)^{-1/2}
                       \|\nabla_h\partial_3^2u(s)\|_2\dd s.
\]
For $0\le t\le2$, \eqref{eq:H2-NLsource} and
\eqref{eq:new-base-energy} give
$\|F_2\|_{L^2(0,2;H^1)}\le C\mathcal J_2^2$.
The $L^2$ boundedness of $S_h(t)\mathbb P$ handles $j=0$,
and \eqref{eq:r8-horizontal-maximal} handles $j=1$.
This proves \eqref{eq:r8-sqrtlog} on the bounded time interval.

\Needspace{6\baselineskip}
\smallskip\noindent\emph{Step 4. The logarithmic bound and the fixed-data limit.}
The remaining weighted integral is bounded by
\[
\begin{aligned}
 &\int_0^{t/2}(1+s)^{-1/2}\|\nabla_h\partial_3^2u(s)\|_2\dd s\\
 &\quad\le
 \left[\int_0^{t/2}(1+s)^{-1}\dd s\right]^{1/2}
 \left[\int_0^{t/2}\|\nabla_h\partial_3^2u(s)\|_2^2\dd s\right]^{1/2}
 \le C\mathcal J_2\sqrt{\log(e+t)},
\end{aligned}
\]
where the second factor is controlled by
\eqref{eq:new-base-energy}. Substitution into the bound for
$\|\nabla_h^j\partial_3r_v(t)\|_2$ obtained in Step~3 proves
\eqref{eq:r8-sqrtlog} for $t\ge2$.

For a fixed solution, split the weighted integral into
$[0,a]$ and $[a,t/2]$, with $a$ fixed and $t>2a$.
The integral over $[0,a]$, divided by $\sqrt{\log t}$, tends
to zero. Cauchy--Schwarz on $[a,t/2]$ gives
\[
 \limsup_{t\to\infty}
 \frac{\displaystyle\int_0^{t/2}(1+s)^{-1/2}
                    \|\nabla_h\partial_3^2u(s)\|_2\dd s}
      {\sqrt{\log t}}
 \le C\left[\int_a^\infty
                    \|\nabla_h\partial_3^2u(s)\|_2^2\dd s\right]^{1/2}.
\]
The right side tends to zero as $a\to\infty$ by
\eqref{eq:new-base-energy}. The term
$C\mathcal J_2^2(1+t)^{-(1+j)/2}$ in Step~3, divided by
$t^{-(1+j)/2}\sqrt{\log t}$, also tends to zero.
This proves the asymptotic refinement for each fixed initial datum following
\eqref{eq:r8-sqrtlog}.
\end{proof}

\section{Third-order decay and the nonlinear thermal contribution}
\label{sec:third-order-decay}

\subsection{Horizontal and mixed derivative estimates}
\label{sec:r9-pointwise}

We work with the global strong solution constructed in Section~\ref{sec:r10-foundations}
and assume that $\mathcal J_2$ is sufficiently small.
The time-integrated estimates \eqref{eq:H2-final-windows} and
\eqref{eq:r8-windows} first yield pointwise bounds for the time
derivatives. Elliptic estimates then give preliminary third-order
spatial bounds. These bounds are used in the separate proofs of the horizontal
and mixed third-order decay in \eqref{eq:r9-horizontal-third}
and \eqref{eq:r9-mixed-third}, respectively.
The energy arguments start at an averaged positive time and require
only the $H^2$ initial regularity. We state decay estimates for $t\ge1$;
the estimates in Lemmas~\ref{lem:r9-time-derivatives} and
\ref{lem:r9-third-preliminary} also hold on any fixed interval
$[t_0,\infty)$, $t_0>0$, with constants depending on $t_0$.

In Sections~\ref{sec:r9-pointwise} and \ref{sec:r9-thermal-upper},
we write $V=u_t$, $\Theta=\tau_t$, and $f=e^{-q}-1$.
The potential viscous scalar
\[
 J=(\nu_h\Delta_h\Delta^{-1}+\nu_c)\Div u
\]
satisfies $\mathcal Lw=\nabla J$.

\subsubsection{Trace estimates and a uniform Gronwall lemma}

\begin{lemma}[Trace and interpolation estimates]
Let $a\in H^1(\R^3)$ and $b\in H^2(\R^3)$. Then
\begin{align}
 \|a\|_{L^\infty_3L^4_h}&\le C\|\nabla a\|_2,&
 \|b\|_{L^2_3L^4_h}&\le C\|b\|_2^{1/2}\|\nabla_hb\|_2^{1/2},
                                                    \label{eq:r9-trace}\\
 \|ab\|_2&\le C\|\nabla a\|_2\|b\|_2^{1/2}
                                      \|\nabla_hb\|_2^{1/2},&
 \|b\|_\infty&\le C\|\nabla b\|_2^{1/2}\|\nabla^2b\|_2^{1/2}.
                                                    \label{eq:r9-agmon}
\end{align}
The mixed-norm estimate for $b$ and the product estimate remain valid
under the weaker assumption $b,\nabla_hb\in L^2(\R^3)$.
\end{lemma}
\begin{proof}
The two estimates in \eqref{eq:r9-trace} were proved after
Lemma~\ref{lem:plane-trace}. H\"older in $x_h$ and $x_3$ gives
$\|ab\|_2\le\|a\|_{L^\infty_3L^4_h}\|b\|_{L^2_3L^4_h}$.
The last inequality follows by splitting the Fourier integral
at a radius and optimizing that radius. Approximation by smooth,
compactly supported functions in the respective norms extends
these estimates to the stated spaces.
\end{proof}

The following uniform Gronwall lemma converts the time-integrated
energy bounds into pointwise estimates.

\begin{lemma}[Uniform Gronwall estimate]
\label{lem:r9-uniform-gronwall}
Let $E\ge0$ be locally absolutely continuous on $(0,\infty)$,
and let $D,a,R\ge0$ be locally integrable. Suppose that
$E'+c_*D\le aE+R$ almost everywhere, with $c_*>0$, and that
for every $b\ge0$,
\[
 \int_b^{b+1}E(s)\dd s\le C_0(1+b)^{-p},\qquad
 \int_b^{b+1}R(s)\dd s\le C_1(1+b)^{-p},\qquad
 \int_b^{b+1}a(s)\dd s\le C_2,
\]
where $p\ge0$. Then
\begin{equation}\label{eq:r9-average-principle}
 E(t)+\int_t^{t+1}D(s)\dd s
       \le C(C_0+C_1)(1+t)^{-p},\qquad t\ge1,
\end{equation}
where $C$ depends only on $c_*,p,C_2$.
If the first two time integrals are $o((1+b)^{-p})$ as
$b\to\infty$, then the left side of
\eqref{eq:r9-average-principle} is $o((1+t)^{-p})$.
\end{lemma}
\begin{proof}
Fix $t\ge1$. By averaging, choose
$s_t\in[t-\tfrac12,t-\tfrac14]$ such that
\[
 E(s_t)\le4\int_{t-1/2}^{t-1/4}E(s)\dd s
        \le CC_0(1+t)^{-p}.
\]
Since
$[s_t,t+1]\subset[t-\tfrac12,t+\tfrac12]\cup[t,t+1]$,
the hypotheses with $b=t-\tfrac12$ and $b=t$ imply
$\int_{s_t}^{t+1}a(s)\dd s\le2C_2$ and
$\int_{s_t}^{t+1}R(s)\dd s\le CC_1(1+t)^{-p}$.
Applying Gronwall's inequality to $E'\le aE+R$ yields
\[
 \sup_{s_t\le s\le t+1}E(s)
 \le e^{2C_2}\left(E(s_t)+\int_{s_t}^{t+1}R(s)\dd s\right)
 \le C(C_0+C_1)(1+t)^{-p}.
\]

Integrating $E'+c_*D\le aE+R$ over $[s_t,t+1]$ and using
$E(t+1)\ge0$ and $D\ge0$ yields
\[
 c_*\int_t^{t+1}D(s)\dd s
 \le E(s_t)+2C_2\sup_{s_t\le s\le t+1}E(s)+\int_{s_t}^{t+1}R(s)\dd s
 \le C(C_0+C_1)(1+t)^{-p}.
\]
Together with $E(t)\le\sup_{s_t\le s\le t+1}E(s)$, this proves
\eqref{eq:r9-average-principle}.

For the refined asymptotic assertion, the averaging argument and
the source estimate give
$E(s_t)+\int_{s_t}^{t+1}R(s)\dd s=o((1+t)^{-p})$.
The Gronwall estimate and the integrated energy inequality then yield
$\sup_{s_t\le s\le t+1}E(s)+\int_t^{t+1}D(s)\dd s=o((1+t)^{-p})$.
Since $s_t\ge\tfrac12$, the argument requires no value of $E(0)$.
\end{proof}

\subsubsection{Estimates for the first time derivatives}

\begin{lemma}[First time derivatives]
\label{lem:r9-time-derivatives}
For $t\ge1$, the solution satisfies
\begin{equation}\label{eq:r9-V-point}
 \begin{aligned}
 &\|\nabla_hu_t(t)\|_2+\|\Div u_t(t)\|_2\\
 &\quad+\|\Delta_hu_t\|_{L^2(t,t+1;L^2)}
 +\|\nabla\Div u_t\|_{L^2(t,t+1;L^2)}
 \le C\mathcal J_2(1+t)^{-7/4}.
 \end{aligned}
\end{equation}
\begin{equation}\label{eq:r9-Theta-point}
 \|\nabla\tau_t(t)\|_2+
       \|\Delta\tau_t\|_{L^2(t,t+1;L^2)}
       \le C\mathcal J_2(1+t)^{-7/4}.
\end{equation}
\begin{equation}\label{eq:r9-vt-point}
 \|\nabla_hv_t(t)\|_2+
       \|\Delta_hv_t\|_{L^2(t,t+1;L^2)}
       \le C\mathcal J_2(1+t)^{-2}.
\end{equation}
\end{lemma}
\begin{proof}
\emph{Step 1. Time difference quotients and preliminary estimates.}
By the regularity \eqref{eq:solspace3}, we use
forward time difference quotients
\[
 V_\sigma(t)=\frac{u(t+\sigma)-u(t)}{\sigma},\qquad
 \Theta_\sigma(t)=\frac{\tau(t+\sigma)-\tau(t)}{\sigma},
 \qquad 0<\sigma<1/8.
\]
Define $q_\sigma$ and $f_\sigma$ analogously.
The energy estimates in Steps~2--4 are first
proved for these quotients; the equations are displayed with
$V=u_t$ and $\Theta=\tau_t$. In particular,
\[
 \frac{f(t+\sigma)\mathcal Lu(t+\sigma)-f(t)\mathcal Lu(t)}{\sigma}
       =f(t+\sigma)\mathcal LV_\sigma+f_\sigma(t)\mathcal Lu(t).
\]

\smallskip\noindent\emph{Preliminary pointwise estimates.}
\begin{align}
 \|\nabla q_t\|_2+\|q_t\|_6+\|\nabla f_t\|_2
       &\le C\mathcal J_2(1+t)^{-7/4},\label{eq:r9-qt}\\
 \|V\|_2+\|\Theta\|_2&\le C\mathcal J_2(1+t)^{-5/4},
                                                    \label{eq:r9-time-L2}\\
 \|f\|_\infty&\le C\mathcal J_2(1+t)^{-3/2},\qquad
 \|\nabla^2f\|_2\le C\mathcal J_2(1+t)^{-7/4}.\notag
\end{align}
The mean value theorem and \eqref{eq:new-uniform-basic} give
$\|f\|_\infty\le C\|q\|_\infty$.
For the second derivative, use
$\nabla^2f=-e^{-q}\nabla^2q+e^{-q}\nabla q\otimes\nabla q$ and
$\|\nabla q\|_4^2\le C\|q\|_\infty\|\nabla^2q\|_2$.
Thus \eqref{eq:new-compressible-second} yields
$\|\nabla^2f\|_2\le C(1+\|q\|_\infty)\|\nabla^2q\|_2$.

For \eqref{eq:r9-qt}, write $q_t=-\Div w+F_1$.
By \eqref{eq:new-base-energy}, \eqref{eq:new-first}, and
\eqref{eq:new-vertical-L3},
$\|\nabla u\|_3\le C\mathcal J_2(1+t)^{-3/8}$.
Together with \eqref{eq:new-compressible-second} and
\eqref{eq:new-uniform-basic}, this gives
\[
 \|\nabla F_1\|_2
 \le C\bigl(\|\nabla u\|_3\|\nabla q\|_6
             +\|u\|_\infty\|\nabla^2q\|_2\bigr)
 \le C\mathcal J_2^2(1+t)^{-17/8}.
\]
The bound for $\nabla q_t$ follows from
\eqref{eq:new-compressible-second}; the Sobolev inequality then
bounds $q_t$ in $L^6$. Finally,
$\nabla f_t=-e^{-q}\nabla q_t+e^{-q}q_t\nabla q$, and
\eqref{eq:new-base-energy} controls $\|\nabla q\|_3$.

For \eqref{eq:r9-time-L2}, the momentum and temperature equations
in \eqref{eq:perturbed} give
\[
 V=\mathcal Lu-c^2\nabla(q+\tau)+F_2,\qquad
 \Theta=\chi\Delta\tau-\alpha\Div w+F_3.
\]
The linear terms are controlled by \eqref{eq:new-first},
\eqref{eq:new-compressible-second}, and
$\|\mathcal Lu\|_2\le C\mathcal J_2(1+t)^{-3/2}$ from
\eqref{eq:new-S2-structure} and \eqref{eq:new-second-preliminary}.
For the sources in \eqref{eq:F2}--\eqref{eq:F3},
\eqref{eq:new-uniform-basic}, \eqref{eq:new-third-final}, and
\eqref{eq:new-production-L2} imply
\[
 \|F_2\|_2\le C\mathcal J_2^2(1+t)^{-7/4},\qquad
 \|F_3\|_2\le C\mathcal J_2^2(1+t)^{-17/8}.
\]
For the vertical convection term, \eqref{eq:new-seed},
\eqref{eq:new-base-energy}, and \eqref{eq:new-third-final} give
\[
 \|u_3\partial_3u\|_2
 \le\|u_3\|_\infty\|u\|_2^{1/2}\|\partial_3^2u\|_2^{1/2}
 \le C\mathcal J_2^2(1+t)^{-7/4}.
\]
The pressure gradient $\nabla(q+\tau)$ and the divergence term
$\Div w$ have decay exponent $5/4$ by \eqref{eq:new-first};
together with the preceding bounds, this proves \eqref{eq:r9-time-L2}.

\smallskip\noindent\emph{Time-integrated source estimates.}
\begin{align}
 \|\Div F_2\|_{L^2(b,b+1;L^2)}
 +\|(F_2)_3\|_{L^2(b,b+1;H^1)}
       &\le C\mathcal J_2^2(1+b)^{-49/24},\label{eq:r9-source-special}\\
 \|F_3\|_{L^2(b,b+1;H^1)}
       &\le C\mathcal J_2^2(1+b)^{-17/8}.\label{eq:r9-source-thermal-old}
\end{align}
For the momentum source in \eqref{eq:r9-source-special}, use
\[
 \Div(u\cdot\nabla u)
   =u\cdot\nabla\Div u+\sum_{i,j}\partial_i u_j\,\partial_j u_i .
\]
The mixed product is bounded by
$\|\nabla_hu_3\|_6\|\partial_3u\|_3$, with exponent
$49/24$ by \eqref{eq:r8-directional-L6} and
\eqref{eq:new-vertical-L3}.
The horizontal convection terms have exponent at least $19/8$
by \eqref{eq:new-first}, \eqref{eq:r8-directional-L6}, and
\eqref{eq:new-uniform-basic}.
The variable-viscosity contribution is
$f(\nu_h\Delta_h+\nu_c\Delta)\Div u+\nabla f\cdot\mathcal Lu$;
its time-integrated norm is controlled by the bounds for
$\|f\|_\infty$ and $\|\nabla^2f\|_2$,
\eqref{eq:H2-final-windows}, and the product inequality
\eqref{eq:r9-agmon}.

For $\partial_3(F_2)_3$, the identity
\begin{equation}\label{eq:r9-third-identity}
 \partial_3\Delta_hu_3
       =\Delta_h\Div w-\Delta_h\Div_hu_h
\end{equation}
reduces the viscous term to the third derivatives in
\eqref{eq:H2-final-windows} and \eqref{eq:r8-windows}.
The convection product $\partial_3u_h\,\nabla_hu_3$ has the same
bound as in the divergence estimate.
For $u_h\nabla_h\partial_3u_3$, \eqref{eq:new-uniform-basic},
\eqref{eq:new-S2-structure}, and \eqref{eq:new-second-preliminary}
give exponent $19/8$.
For $(\partial_3u_3)^2$, interpolation between
\eqref{eq:new-first} and \eqref{eq:r8-directional-L6} gives $19/8$.
For $u_3\partial_3^2u_3$, \eqref{eq:new-third-final} and
\eqref{eq:H2-main-vertical-weak} give $9/4$.
The horizontal derivatives are already controlled by
\eqref{eq:r8-horizontal-sources}. For the zeroth-order component,
$\|(F_2)_3\|_2\le C\mathcal J_2^2(1+t)^{-17/8}$ by
\eqref{eq:new-first} and \eqref{eq:new-uniform-basic}; the slowest
product is $u_h\cdot\nabla_hu_3$. The pressure terms follow from
\eqref{eq:new-first}, \eqref{eq:new-compressible-second}, and
\eqref{eq:new-uniform-basic}. This proves
\eqref{eq:r9-source-special}.

For the temperature source in \eqref{eq:r9-source-thermal-old},
the differentiated vertical transport term is bounded by
$\|\partial_3u\|_3\|\nabla\tau\|_6$, with exponent
$17/8$ by \eqref{eq:new-vertical-L3} and
\eqref{eq:new-compressible-second}.
For the gradient of the viscous heating term, use the third line of
\eqref{eq:r8-windows} for $\nabla_hu$ in $L^2(b,b+1;L^\infty)$,
and \eqref{eq:new-H2-main-rates} and \eqref{eq:H2-main-vertical-weak}
for $\nabla\nabla_hu$ in $L^\infty(b,b+1;L^2)$.
Their respective exponents $11/8$ and $3/4$ give $17/8$.
The bounds for $\|f\|_\infty$ and $\|\nabla^2f\|_2$, together
with the third line of \eqref{eq:H2-final-windows}, control
the variable-coefficient thermal diffusion term.
This proves \eqref{eq:r9-source-thermal-old}.

\smallskip\noindent\emph{Time-integrated derivative estimates.}
The momentum and temperature equations in \eqref{eq:perturbed},
\eqref{eq:new-compressible-second}, \eqref{eq:H2-final-windows},
the first line of \eqref{eq:r8-windows}, and
\eqref{eq:r8-horizontal-sources} give
\begin{equation}\label{eq:r9-time-old}
\begin{aligned}
 &\|\nabla_hV\|_{L^2(b,b+1;L^2)}
 +\|\Div V\|_{L^2(b,b+1;L^2)}\\
 &\qquad+\|\nabla V_3\|_{L^2(b,b+1;L^2)}
 +\|\nabla\Theta\|_{L^2(b,b+1;L^2)}
       \le C\mathcal J_2(1+b)^{-7/4}.
\end{aligned}
\end{equation}
Here $\Div F_2$ and $(F_2)_3$ are controlled by
\eqref{eq:r9-source-special}, and $F_3$ by
\eqref{eq:r9-source-thermal-old}; the vertical viscous derivative
in $\nabla V_3$ is reduced using \eqref{eq:r9-third-identity}.
The momentum equation, the second line of \eqref{eq:r8-windows},
the third line of \eqref{eq:H2-final-windows},
\eqref{eq:new-compressible-second}, and \eqref{eq:r8-F2-H1-window}
yield
\begin{equation}\label{eq:r9-time-z-old}
 \|\partial_3V\|_{L^2(b,b+1;L^2)}
 +\|\nabla\mathcal Lu\|_{L^2(b,b+1;L^2)}
       \le C\mathcal J_2(1+b)^{-5/4}.
\end{equation}
Finally, \eqref{eq:projected-P-intro}, the first line of
\eqref{eq:r8-windows}, and \eqref{eq:r8-horizontal-sources} give
\begin{equation}\label{eq:r9-shear-time-old}
 \|\nabla_hv_t\|_{L^2(b,b+1;L^2)}
       \le C\mathcal J_2(1+b)^{-2}.
\end{equation}
The anisotropic Sobolev inequality
$\|V\|_6\le C\|\nabla_hV\|_2^{2/3}\|\partial_3V\|_2^{1/3}$,
H\"older's inequality in time, and
\eqref{eq:r9-time-old}--\eqref{eq:r9-time-z-old}
give exponent $19/12$ for $V$ in $L^2(b,b+1;L^6)$.
The standard Sobolev inequality and \eqref{eq:r9-time-old}
give exponent $7/4$ for $V_3$ and $\Theta$ in $L^2(b,b+1;L^6)$.

For $z=u,\tau,q,f$ and a Banach space $X$, Jensen's inequality gives
\[
 \|z_\sigma\|_{L^2(b,b+1;X)}^2
       \le\int_b^{b+1+\sigma}\|z_t(s)\|_X^2\dd s.
\]
Since $[b,b+1+\sigma]\subset[b,b+1]\cup[b+1,b+2]$,
\eqref{eq:r9-time-old}--\eqref{eq:r9-shear-time-old} hold for the
quotients with constants independent of $\sigma$.
The averaging formula also transfers the pointwise bounds
\eqref{eq:r9-qt} and \eqref{eq:r9-time-L2} to $q_\sigma,f_\sigma$,
$V_\sigma,\Theta_\sigma$.

\emph{Step 2. Energy estimate for $u_t$.}
Differentiating the momentum equation in time gives
\begin{equation}\label{eq:r9-V-equation}
 V_t-\mathcal LV+c^2\nabla(q_t+\Theta)=f\mathcal LV+G,
\end{equation}
where
\begin{equation}\label{eq:r9-G}
 G={}-V\cdot\nabla u-u\cdot\nabla V
       -c^2e^\tau\Theta\nabla(q+\tau)
 -c^2(e^\tau-1)\nabla(q_t+\Theta)+f_t\mathcal Lu .
\end{equation}
In $L^2(b,b+1;L^2)$, the terms in \eqref{eq:r9-G}
have decay exponents at least
$65/24$ for $V_h\nabla_hu$, $17/8$ for $V_3\partial_3u$,
$21/8$ for $u_h\nabla_hV$, $11/4$ for $u_3\partial_3V$,
$13/4$ for each pressure product, and $25/8$ for
$f_t\mathcal Lu$. Thus
$\|G\|_{L^2(b,b+1;L^2)}
\le C\mathcal J_2^2(1+b)^{-17/8}$.
Each product is estimated by H\"older's inequality, with one factor
controlled pointwise in time and the other in $L^2(b,b+1)$.
For the coefficient difference, \eqref{eq:r9-agmon},
\eqref{eq:r9-qt}, \eqref{eq:H2-final-windows}, and
\eqref{eq:r9-time-z-old} give
\[
 \begin{aligned}
 \|f_\sigma\mathcal Lu\|_{L^2(b,b+1;L^2)}
 &\le C\|\nabla f_\sigma\|_{L^\infty(b,b+1;L^2)}
       \|\mathcal Lu\|_{L^2(b,b+1;L^2)}^{1/2}
       \|\nabla_h\mathcal Lu\|_{L^2(b,b+1;L^2)}^{1/2}\\
 &\le C\mathcal J_2^2(1+b)^{-25/8}.
 \end{aligned}
\]
The remaining quotient products corresponding to \eqref{eq:r9-G}
use \eqref{eq:r9-qt}--\eqref{eq:r9-time-L2} and
\eqref{eq:r9-time-old}--\eqref{eq:r9-time-z-old}.

Test \eqref{eq:r9-V-equation} with
$-\Delta_hV-\nabla\Div V$.
The mixed viscous term is nonnegative. The contribution of $f\mathcal LV$
is absorbed into the dissipation by the smallness of $\|f\|_\infty$,
without differentiating $f$. Thus
\[
 (\|\nabla_hV\|_2^2+\|\Div V\|_2^2)'
 +c_*(\|\Delta_hV\|_2^2+\|\nabla\Div V\|_2^2)
 \le C(\|\nabla q_t\|_2^2+\|\nabla\Theta\|_2^2+\|G\|_2^2).
\]
Apply Lemma~\ref{lem:r9-uniform-gronwall} with
\eqref{eq:r9-time-old}, \eqref{eq:r9-qt}, and the
$L^2(b,b+1;L^2)$ source estimate for $G$.
This proves \eqref{eq:r9-V-point}, first uniformly for $V_\sigma$.
Applying Lemma~\ref{lem:r9-uniform-gronwall} on shorter intervals
gives the same quotient bounds for $t\ge t_0>0$, with constants
depending on $t_0$.

\emph{Step 3. Energy estimate for $\tau_t$.}
Differentiating the equation for $\tau$ in time gives
\begin{equation}\label{eq:r9-Theta-equation}
 \Theta_t-\chi\Delta\Theta+\alpha\Div V
 =\chi f\Delta\Theta-u\cdot\nabla\Theta
       +2\chi e^{-q}\nabla\tau\cdot\nabla\Theta+S,
\end{equation}
with
\begin{equation}\label{eq:r9-S}
 \begin{aligned}
 S={}&-V\cdot\nabla\tau+\chi f_t\Delta\tau
       -\chi e^{-q}q_t|\nabla\tau|^2\\
 &+\frac{e^{-q-\tau}}{\bar\rho\bar\theta c_v}
 \left[-(q_t+\Theta)Q(u)+2\mu\nabla_hu:\nabla_hV
                  +2(\lambda+\mu)\Div u\,\Div V\right].
 \end{aligned}
\end{equation}
For the transport term $-V\cdot\nabla\tau$ in \eqref{eq:r9-S},
\eqref{eq:r9-time-L2} and the second line of
\eqref{eq:H2-final-lipschitz-windows} give the time-integrated
decay exponent $3$.
For the horizontal heating difference, the quotient bound
\eqref{eq:r9-V-point} established in Step~2 and the third line of
\eqref{eq:r8-windows} yield, for $b\ge1/4$,
\[
 \begin{aligned}
 \left\|\frac{|\nabla_hu(t+\sigma)|^2-|\nabla_hu(t)|^2}{\sigma}
                 \right\|_{L^2(b,b+1;L^2)}
 &\le2\|\nabla_hu\|_{L^2(b,b+1+\sigma;L^\infty)}
       \|\nabla_hV_\sigma\|_{L^\infty(b,b+1;L^2)}\\
 &\le C\mathcal J_2^2(1+b)^{-25/8}.
 \end{aligned}
\]
The divergence heating difference uses
\eqref{eq:H2-final-lipschitz-windows} and the quotient estimate
\eqref{eq:r9-V-point}.
Differences of the exponential coefficients use the mean value
formula, \eqref{eq:r9-qt}--\eqref{eq:r9-time-L2}, and
\eqref{eq:r9-time-old}.
The remaining terms in \eqref{eq:r9-S} have larger decay exponents;
here $\|Q(u)\|_3\le C\mathcal J_2^2(1+t)^{-5/2}$.
Consequently $\|S\|_{L^2(b,b+1;L^2)}
\le C\mathcal J_2^2(1+b)^{-3}$ after a positive time.
Testing \eqref{eq:r9-Theta-equation} with $-\Delta\Theta$ gives
\[
 (\|\nabla\Theta\|_2^2)'+c_*\|\Delta\Theta\|_2^2
 \le C(\|u\|_\infty^2+\|\nabla\tau\|_6^4)\|\nabla\Theta\|_2^2
       +C\|\Div V\|_2^2+C\|S\|_2^2 .
\]
Apply Lemma~\ref{lem:r9-uniform-gronwall} with
\eqref{eq:r9-time-old}, \eqref{eq:r9-V-point}, and the
$L^2(b,b+1;L^2)$ source estimate for $S$.
This proves \eqref{eq:r9-Theta-point}.

\emph{Step 4. Energy estimate for $v_t$.}
Since $\mathcal Lw=\nabla J$,
$\mathbb P(f\mathcal Lw_t)=-\mathbb P(J_t\nabla f)$.
Applying the Leray projection to \eqref{eq:r9-V-equation} gives
\[
 (v_t)_t-\nu_h\Delta_hv_t
   =\nu_h\mathbb P(f\Delta_hv_t)+\mathbb P(G-J_t\nabla f).
\]
By \eqref{eq:r9-trace} and \eqref{eq:r9-V-point},
$\|J_t\nabla f\|_{L^2(b,b+1;L^2)}
\le C\mathcal J_2^2(1+b)^{-7/2}$.
Testing the projected equation with $-\Delta_hv_t$ and applying
Lemma~\ref{lem:r9-uniform-gronwall} with \eqref{eq:r9-shear-time-old}
yield \eqref{eq:r9-vt-point}.

For $t\ge1$, choose the averaging interval $[t-1/2,t-1/4]$
independently of $0<\sigma<1/8$. Jensen's inequality bounds the
quotient averages by the corresponding derivative bounds on the enlarged
interval, which is covered by two unit intervals. The mean-value starting
time may depend on $\sigma$, but its energy bound and all subsequent
constants do not. Weak lower semicontinuity as $\sigma\downarrow0$ gives
\eqref{eq:r9-V-point}--\eqref{eq:r9-vt-point} almost everywhere.
The differentiated equations give local $L^2$ continuity of $V,\Theta$;
together with the uniform spatial bounds, this yields weak continuity
in the displayed energy spaces and extends the estimates to every $t>0$.
For $t\ge t_0>0$, use shorter averaging intervals, with constants depending
on $t_0$; no value of $V$ or $\Theta$ at time zero is needed.
\end{proof}

\subsubsection{Preliminary third-order spatial estimates}

\begin{lemma}[Preliminary third-order spatial decay]
\label{lem:r9-third-preliminary}
For $t\ge1$, the solution satisfies
\begin{equation}\label{eq:r9-third-preliminary}
 \|\nabla_h^3v(t)\|_2\le C\mathcal J_2(1+t)^{-2},\qquad
 \|\nabla^3(w,\tau)(t)\|_2+
       \|\nabla(w,\tau)(t)\|_\infty
       \le C\mathcal J_2(1+t)^{-7/4}.
\end{equation}
The velocity Lipschitz estimate \eqref{eq:r9-lips} also holds, and
\begin{equation}\label{eq:r9-point-sources-pre}
 \|\nabla_hF_2\|_2+\|\Div F_2\|_2
       \le C\mathcal J_2^2(1+t)^{-49/24},\qquad
 \|\nabla_hF_3\|_2\le C\mathcal J_2^2(1+t)^{-21/8}.
\end{equation}
\end{lemma}
\begin{proof}
\emph{Step 1. Coupled elliptic estimates for $v$ and $w$.}
We use \eqref{eq:r9-V-point} and \eqref{eq:r9-vt-point} to
recover spatial derivatives from the momentum equation.
The projected equation is
\begin{equation}\label{eq:r9-shear-elliptic}
 v_t-\nu_h\Delta_hv
 =\mathbb P[-u\cdot\nabla u+c^2(q+\tau)e^\tau\nabla\tau
                     +\nu_h f\Delta_hv-J\nabla f].
\end{equation}
Here the pressure identity follows by subtracting the gradient
of $(q+\tau)(e^\tau-1)$.
Take one horizontal derivative. Use \eqref{eq:r9-trace}
on $\nabla_hf\,\Delta_hv$ and $\nabla_hJ\,\nabla f$,
and Agmon on $J\nabla_h\nabla f$.
Since $\|\nabla J\|_2\le C\mathcal J_2(1+t)^{-7/4}$,
one obtains
\begin{align}
 \|\nabla_h^3v\|_2\le{}&
 C\mathcal J_2(1+t)^{-2}
 +C\mathcal J_2^2(1+t)^{-49/24}\notag\\
 &+C\mathcal J_2^{3/2}(1+t)^{-5/2}\|\nabla_h^3v\|_2^{1/2}
 +C\mathcal J_2^{3/2}(1+t)^{-21/8}\|\nabla^3w\|_2^{1/2}.
                                                   \label{eq:r9-joint-shear}
\end{align}
Taking divergence of the original momentum equation gives
\begin{equation}\label{eq:r9-div-elliptic}
 (\nu_h\Delta_h+\nu_c\Delta)\Div u
       =\Div V+c^2\Delta(q+\tau)-\Div F_2.
\end{equation}
The principal variable-viscosity term is
$f(\nu_h\Delta_h+\nu_c\Delta)\Div u$ and is absorbed.
The remaining coefficient term satisfies
\[
 \|\nabla f\cdot\mathcal Lu\|_2
 \le C\|\nabla^2f\|_2\|\mathcal Lu\|_2^{1/2}
                                \|\nabla_h\mathcal Lu\|_2^{1/2},
 \quad
 \|\nabla_h\mathcal Lu\|_2
       \le C(\|\nabla_h^3v\|_2+\|\nabla^3w\|_2).
\]
Full ellipticity follows from $\nu_c>0$. Thus
\begin{equation}\label{eq:r9-joint-potential}
 \|\nabla^3w\|_2
 \le C\mathcal J_2(1+t)^{-7/4}
   +C\mathcal J_2^2(1+t)^{-49/24}
   +C\mathcal J_2^{3/2}(1+t)^{-5/2}
                       (\|\nabla_h^3v\|_2+\|\nabla^3w\|_2)^{1/2}.
\end{equation}
Multiply \eqref{eq:r9-joint-shear} by $(1+t)^2$ and
\eqref{eq:r9-joint-potential} by $(1+t)^{7/4}$.
Young's inequality absorbs one half of each weighted unknown and gives
\[
 \|\nabla_h^3v(t)\|_2\le C\mathcal J_2(1+t)^{-2},\qquad
 \|\nabla^3w(t)\|_2\le C\mathcal J_2(1+t)^{-7/4}.
\]

\emph{Step 2. Velocity Lipschitz decay.}
The directional Agmon inequality gives
\[
 \|\nabla_hv\|_\infty
 \le C\|v\|_2^{1/12}\|\nabla_h^3v\|_2^{2/3}
                              \|\partial_3^2v\|_2^{1/4}
 \le C\mathcal J_2(1+t)^{-11/8},
\]
using \eqref{eq:new-seed} and the uniform $H^2$ bound
\eqref{eq:new-base-energy}. The estimate for $\nabla^3w$
in Step~1, together with \eqref{eq:new-compressible-second}
and the isotropic Agmon inequality, yields
$\|\nabla w\|_\infty\le C\mathcal J_2(1+t)^{-7/4}$.
Since $v$ is divergence-free,
$\partial_3u_3=\Div w-\Div_hu_h$. The bounds for
$\nabla_hu$ and $\nabla w$, together with this identity,
prove \eqref{eq:r9-lips}.

\emph{Step 3. Temperature and nonlinear source estimates.}
Differentiate $\chi\Delta\tau=\Theta+\alpha\Div u-F_3$ once.
The coefficient term $f\nabla\Delta\tau$ is absorbed by smallness.
The horizontal heating gradient is bounded by
$\|\nabla_hu\|_\infty\|\nabla\nabla_hu\|_2$,
using the velocity Lipschitz estimate proved in Step~2.
The divergence heating gradient is bounded by
$\|\Div u\|_\infty\|\nabla\Div u\|_2$ using the potential
velocity estimates in Steps~1--2.
The other terms containing the new temperature derivative are bounded by
$C\mathcal J_2^{3/2}(1+t)^{-21/8}\|\nabla^3\tau\|_2^{1/2}$
and are absorbed by Young's inequality.
Together with \eqref{eq:r9-Theta-point}, this proves the
third-order temperature bound in \eqref{eq:r9-third-preliminary}.
The isotropic Agmon inequality in \eqref{eq:r9-agmon} and
\eqref{eq:new-compressible-second} then give the asserted
$L^\infty$ bound for $\nabla\tau$.

Apply the product estimates \eqref{eq:r9-trace}--\eqref{eq:r9-agmon}
to \eqref{eq:F2} and \eqref{eq:F3}, using
\eqref{eq:r9-lips} and \eqref{eq:r9-third-preliminary}.
The convection term $\nabla_hu_3\,\partial_3u$ has decay
exponent $49/24$, as in the proof of
\eqref{eq:r9-source-special}; the horizontal scalar transport
term $u\,\nabla_h\nabla\tau$ has exponent $21/8$.
All other terms decay at least as fast. This proves
\eqref{eq:r9-point-sources-pre}.
The elliptic estimates are first justified at times with finite
dissipation; $H^2$ continuity and weak compactness extend them
to every $t\ge1$. Applying the same elliptic estimates on
$[t_0,\infty)$ gives the positive-time version described at the
start of this subsection.
\end{proof}

\subsubsection{Pointwise horizontal third-order decay}

Lemma~\ref{lem:r9-third-preliminary} already gives the
$(1+t)^{-2}$ bound for $\nabla_h^3v$.
We now improve the horizontal third derivatives of $(w,\tau)$
and estimate the high horizontal frequencies of $v$.
The two refined bounds below are also used in the proof of
Proposition~\ref{prop:gaussian}.

\begin{proposition}[Horizontal third-order decay]
\label{prop:r9-horizontal-third}
For $t\ge1$, the horizontal third-order estimate
\eqref{eq:r9-horizontal-third} holds. More precisely,
\begin{equation}\label{eq:r9-horizontal-compressible-small}
 \|\nabla_h^3(w,\tau)(t)\|_2
       \le C\mathcal J_2(1+t)^{-49/24}\log(e+t)=o(t^{-2}).
\end{equation}
\begin{equation}\label{eq:r9-high-shear-small}
 \|\nabla_h^3\mathbf1_{\{|D_h|\ge1\}}v(t)\|_2
       \le C\mathcal J_2 e^{-c_*t}
                    +C\mathcal J_2^2(1+t)^{-49/24}.
\end{equation}
\end{proposition}
\begin{proof}
\emph{Step 1. Time averages at high spatial frequencies.}
Localization to high spatial frequencies does not give a spectral gap for the shear
semigroup. To estimate the compressible variables, use the operator
$\mathcal A_h=-\Delta_h(-\Delta)^{-1/2}$, whose symbol is
$|\xi_h|^2/|\xi|$. It satisfies
$\|\mathcal A_hg\|_2\le\|\nabla_hg\|_2$ and
$\|\nabla\mathcal A_hg\|_2=\|\nabla_h^2g\|_2$.
Let $P_H$ be a fixed smooth high spatial-frequency cutoff.
The compensated energy estimate at high spatial frequencies
\eqref{eq:r8-compressible-henergy} in the proof of
Proposition~\ref{prop:r8-windows}, together with the horizontal
source bounds \eqref{eq:r8-horizontal-sources}, gives
\begin{align}
 \|\mathcal A_hV_H\|_{L^2(b,b+1;L^2)}
 &\le C\mathcal J_2 e^{-c_*b}+C\mathcal J_2(1+b)^{-5/2}
                         +C\mathcal J_2^2(1+b)^{-49/24},\notag\\
 \|\mathcal A_h(q_t)_H(b)\|_2+
 \|\mathcal A_h\Theta_H\|_{L^2(b,b+1;L^2)}
 &\le C\mathcal J_2 e^{-c_*b}
                         +C\mathcal J_2^2(1+b)^{-49/24}.
                                                     \label{eq:r9-A-average}
\end{align}
Indeed, $\mathcal A_h\Div u$ has the norm of $\nabla_h^2w$,
and $\mathcal A_h\Delta\tau$ that of $\nabla\nabla_h^2\tau$.
For the shear term the symbol is $|\xi_h|^4/|\xi|$.
On $|\xi_h|\ge1$, use the high-horizontal-frequency energy bound
for $\nabla_h^3P_H^hv$ in the proof of
Proposition~\ref{prop:r8-windows}; on $|\xi_h|\le1$ and the support
of $P_H$, the symbol is bounded by $C|\xi_h|^4$.
The low-horizontal heat estimate in the proof of
Proposition~\ref{prop:r8-windows} gives decay exponent $5/2$
for the initial term. On $s\in[0,t/2]$, use the integrable
mixed source bounds \eqref{eq:new-profile-source-ledger}.
On $s\in[t/2,t]$, place one horizontal derivative on $F_2$,
use \eqref{eq:r8-horizontal-sources}, and integrate the kernel
$(1+t-s)^{-3/2}$. This proves \eqref{eq:r9-A-average}.

\emph{Step 2. Frequency-localized time energy.}
The frequency-localized energy estimates are justified by the time
difference quotients in the proof of
Lemma~\ref{lem:r9-time-derivatives}.
Take $P_H$ real with $0\le P_H\le1$, and write $B=\mathcal A_hP_H$,
$E=\|BV\|_2^2$ and
$D=\nu_h\|\nabla_hBV\|_2^2+\nu_c\|B\Div V\|_2^2$.
Testing \eqref{eq:r9-V-equation} with $B^2V$ gives the exact identity
\[
 \tfrac12E'+D
 =c^2\ip{B(q_t+\Theta)}{B\Div V}+\ip{f\mathcal LV+G}{B^2V}.
\]
Split the last test using $V=v_t+w_t$. Since
$\|B^2v_t\|_2\le\|\Delta_hv_t\|_2$ and
$\|B^2w_t\|_2\le\|B\Div V\|_2$, Young's inequality yields
\[
 \begin{aligned}
 E'+c_*D&\le C\mathcal R_H,\\
 \mathcal R_H&=\|B(q_t+\Theta)\|_2^2+\|G\|_2^2
       +\|f\|_\infty^2\|\mathcal LV\|_2^2
       +(\|G\|_2+\|f\|_\infty\|\mathcal LV\|_2)\|\Delta_hv_t\|_2.
 \end{aligned}
\]
No commutator of $B$ with $f$ occurs. On each unit interval, the last
two products have bounds $C\mathcal J_2^3(1+b)^{-33/8}$ and
$C\mathcal J_2^3(1+b)^{-21/4}$, respectively, by
\eqref{eq:r9-V-point} and \eqref{eq:r9-vt-point}.
The squared $G$ and $f\mathcal LV$ terms have exponents $17/4$ and
$13/2$. Together with \eqref{eq:r9-A-average}, this gives
$\int_b^{b+1}\mathcal R_H\dd s\le C\mathcal J_2^2(1+b)^{-49/12}$.
Lemma~\ref{lem:r9-uniform-gronwall} and the average of $E$ therefore imply
\begin{equation}
 \|\mathcal A_hV_H(t)\|_2+
 \|\mathcal A_h\Div V_H\|_{L^2(t,t+1;L^2)}
       \le C\mathcal J_2(1+t)^{-49/24}.
\end{equation}
For the heat equation, put
$H_\Theta=\chi f\Delta\Theta-u\cdot\nabla\Theta
 +2\chi e^{-q}\nabla\tau\cdot\nabla\Theta+S$.
Testing \eqref{eq:r9-Theta-equation} with $B^2\Theta$ and using
$\|B^2\Theta\|_2\le\|\nabla B\Theta\|_2$ gives
\[
 (\|B\Theta\|_2^2)'+c_*\|\nabla B\Theta\|_2^2
 \le C\|B\Theta\|_2^2+C\|B\Div V\|_2^2+C\|H_\Theta\|_2^2.
\]
In $L^2(b,b+1;L^2)$, convection has exponent $21/8$ by
\eqref{eq:new-uniform-basic} and \eqref{eq:r9-Theta-point}.
For the gradient product, use
$\|\nabla\tau\|_6\|\nabla\Theta\|_2^{1/2}\|\nabla^2\Theta\|_2^{1/2}$
and time H\"older to obtain exponent $7/2$.
The terms $f\Delta\Theta$ and $S$ have exponents $13/4$ and $3$.
Thus $\|H_\Theta\|_{L^2(b,b+1;L^2)}\le C\mathcal J_2^2(1+b)^{-21/8}$.
The preceding divergence bound, \eqref{eq:r9-A-average}, and the same
Gronwall lemma yield
$\|\mathcal A_h\Theta_H(t)\|_2\le C\mathcal J_2(1+t)^{-49/24}$.

\emph{Step 3. Elliptic estimates for the compressible derivatives.}
Apply the high-frequency multiplier to
\eqref{eq:r9-div-elliptic}. Its symbol gives
\[
 \|\nabla_h^3w_H\|_2
 \le C\bigl(\|\mathcal A_hw_{t,H}\|_2
       +\|\nabla_h^2(q_H,\tau_H)\|_2+\|\Div F_2\|_2\bigr).
\]
The respective multiplier bounds are
$C|\xi_h|^3/|\xi|^3$ for the divergence source and
$C|\xi_h|^3/|\xi|\le C|\xi_h|^2$ for the scalar second derivatives.
The heat equation similarly gives
\[
 \|\nabla_h^3\tau_H\|_2
 \le C(\|\mathcal A_h\Theta_H\|_2+
                  \|\nabla_h^2w_H\|_2+\|\nabla_hF_3\|_2).
\]
At low spatial frequencies, the Duhamel estimates in Step~2 of the proof of
Proposition~\ref{prop:r8-windows}, with
\eqref{eq:r9-point-sources-pre}, give
$C\mathcal J_2(1+t)^{-9/4}
+C\mathcal J_2^2(1+t)^{-49/24}\log(e+t)$.
Adding the high- and low-frequency bounds proves
\eqref{eq:r9-horizontal-compressible-small}.
Together with the shear bound in
\eqref{eq:r9-third-preliminary}, this proves
\eqref{eq:r9-horizontal-third}.

\emph{Step 4. High horizontal frequencies of the shear velocity.}
To prove \eqref{eq:r9-high-shear-small}, split $f\Delta_hv_t$
in the projected equation obtained from \eqref{eq:r9-V-equation}
into high and low horizontal frequencies. The high part is
absorbed in the energy of $\nabla_hv_t$; the low part is bounded
by $C\mathcal J_2^2(1+t)^{-3}$ since
$\|v_t\|_2\le C\mathcal J_2(1+t)^{-3/2}$.
The spectral gap on $|\xi_h|\ge1$ gives
$\|\nabla_h\mathbf1_{\{|D_h|\ge1\}}v_t\|_2
\le C\mathcal J_2 e^{-c_*t}+C\mathcal J_2^2(1+t)^{-17/8}$.
Insert this bound and \eqref{eq:r9-point-sources-pre} in the
projected spatial equation \eqref{eq:r9-shear-elliptic} to prove
\eqref{eq:r9-high-shear-small}.
\end{proof}

\subsubsection{Pointwise mixed third-order decay}

\begin{proposition}[Mixed third-order decay]
\label{prop:r9-mixed-third}
For $t\ge1$, the mixed third-order estimate
\eqref{eq:r9-mixed-third} holds, and
\begin{equation}\label{eq:r9-Y-point}
 \|\partial_3u_t(t)\|_2+
 \|\nabla_h\partial_3u_t\|_{L^2(t,t+1;L^2)}
 +\|\Div\partial_3u_t\|_{L^2(t,t+1;L^2)}
       \le C\mathcal J_2(1+t)^{-5/4}.
\end{equation}
\begin{equation}\label{eq:r9-F2-H1}
 \|F_2(t)\|_{H^1}\le C\mathcal J_2^2(1+t)^{-3/2}.
\end{equation}
For each fixed initial datum,
$\|\nabla_h^2\partial_3u(t)\|_2=o(t^{-5/4})$ and
$\|F_2(t)\|_{H^1}=o(t^{-3/2})$.
\end{proposition}
\begin{proof}
\emph{Step 1. The vertically differentiated time equation.}
Set $Z=\partial_3u$ and $Y=\partial_3V$, with $V=u_t$.
Differentiate \eqref{eq:r9-V-equation} vertically to obtain
\begin{equation}
 \begin{aligned}
 Y_t-(1+f)\mathcal LY+u\cdot\nabla Y
       +c^2e^\tau\nabla(\partial_3q_t+\partial_3\Theta)
 ={}& \partial_3f\mathcal LV-Z\cdot\nabla V-Y\cdot\nabla u\\
 &-V\cdot\nabla Z+\partial_3(f_t\mathcal Lu)+R_p ,
 \end{aligned}
\end{equation}
where
\[
 R_p=-c^2e^\tau\partial_3\tau\nabla(q_t+\Theta)
             -c^2\partial_3(e^\tau\Theta\nabla(q+\tau)),\qquad
 \|R_p\|_{L^2(b,b+1;L^2)}
       \le C\mathcal J_2^2(1+b)^{-13/4}.
\]
For $\Theta\nabla\partial_3(q+\tau)$, the isotropic Agmon
inequality and \eqref{eq:r9-Theta-point} give the time-integrated
$L^\infty$ bound for $\Theta$; use
\eqref{eq:new-compressible-second} for the other factor.
For $\partial_3\Theta\nabla(q+\tau)$, use H\"older's inequality
in $L^6\times L^3$ with \eqref{eq:r9-Theta-point} and
\eqref{eq:new-compressible-second}.
Testing the principal pressure with $Y$ and integrating by
parts requires only $\partial_3q_t,\partial_3\Theta$ in $L^2$.
Its bound is
$\eta\|\Div Y\|_2^2+C_\eta\|Y\|_2^2+
C_\eta\mathcal J_2^2(1+t)^{-7/2}$.

\emph{Step 2. Energy estimates for transport and viscosity.}
For this calculation, write
$E=\|Y\|_2^2$,
$\mathscr D_Y=(\|\nabla_hY\|_2^2+\|\Div Y\|_2^2)^{1/2}$,
and $K=\|\Delta_hV\|_2+\|\nabla\Div V\|_2$.
The following bounds specify the available derivatives:
\begin{align}
 \|\nabla_hV\|_6+\|\nabla_hV\|_{L^\infty_3L^4_h}
                     +\|Y_3\|_6&\le C(K+\mathscr D_Y),\notag\\
 \|V\|_6&\le C\mathcal J_2^{2/3}(1+t)^{-7/6}E^{1/6},\notag\\
 \|V_3\|_\infty&\le
       C\mathcal J_2^{1/2}(1+t)^{-7/8}(K+\mathscr D_Y)^{1/2},\notag\\
 \|Y\|_{L^2_3L^4_h}&\le CE^{1/4}\mathscr D_Y^{1/2}.
\end{align}
Here $Y_3=\Div V-\Div_hV_h$ and
$\partial_3^2V_3=\partial_3\Div V-\Div_hY_h$.
Thus no full $\partial_3Y$ dissipation is used.

The principal transport and the factors
$Z_3Y$, $Y_h\nabla_hu$ are bounded by
$C\mathcal J_2(1+t)^{-11/8}E$.
The remaining products $Z_h\nabla_hV$ and $Y_3Z$ satisfy
\[
 C\mathcal J_2(1+t)^{-3/8}(K+\mathscr D_Y)E^{1/2}
 \le\eta \mathscr D_Y^2+C_\eta E+
                 C_\eta\mathcal J_2^2(1+t)^{-3/4}K^2.
\]
Integrate $V_h\nabla_hZ$ horizontally by parts.
The term with $\Div_hV_h$ is controlled by the preceding
Young estimate for $(K+\mathscr D_Y)E^{1/2}$, and the
other term is at most
$C\mathcal J_2^{5/3}(1+t)^{-37/24}E^{1/6}\mathscr D_Y$.
Young's inequality leaves
$\eta \mathscr D_Y^2+C_\eta E+C_\eta\mathcal J_2^5(1+t)^{-37/8}$.
Keep the vertical transport in its original form:
\[
 \left|\int V_3\partial_3Z\cdot Y\right|
 \le C\mathcal J_2^{3/2}(1+t)^{-7/8}(K+\mathscr D_Y)^{1/2}E^{1/2}.
\]
Its bound is
$\eta \mathscr D_Y^2+C_\eta E+
C_\eta\mathcal J_2^3(1+t)^{-7/4}K+
C_\eta\mathcal J_2^6(1+t)^{-7/2}$.
Only the bounded second vertical derivative of $u$ was used.

Coefficient gradients from the principal viscosity are bounded
by $\|\nabla^2f\|_2E^{1/4}\mathscr D_Y^{3/2}$.
For $\partial_3f\mathcal LV$, move one viscous derivative onto $Y$
or $\partial_3f$; \eqref{eq:r9-trace} gives
\[
 C\mathcal J_2^{3/2}(1+t)^{-21/8}K^{1/2}\mathscr D_Y,\qquad
 C\mathcal J_2(1+t)^{-7/4}(K+\mathscr D_Y)E^{1/4}\mathscr D_Y^{1/2}.
\]
Their remainders are bounded by
$C_\eta E+C_\eta\mathcal J_2^3(1+t)^{-21/4}K+
C_\eta\mathcal J_2^2(1+t)^{-7/2}K^2$.
Finally, the exact split
$\partial_3(f_t\mathcal Lu)=\partial_3f_t\mathcal Lu+f_t\mathcal LZ$
and the trace inequalities yield
\[
 \left|\langle\partial_3(f_t\mathcal Lu),Y\rangle\right|
 \le C\mathcal J_2(1+t)^{-7/4}
   (\|\nabla\mathcal Lu\|_2+\|\mathcal LZ\|_2)E^{1/4}\mathscr D_Y^{1/2}.
\]
After absorption, its remainder is bounded by $C_\eta E$ plus
$C_\eta\mathcal J_2^2(1+t)^{-7/2}
(\|\nabla\mathcal Lu\|_2^2+\|\mathcal LZ\|_2^2)$.
The unit-interval integral of
$\|\nabla\mathcal Lu\|_2^2+\|\mathcal LZ\|_2^2$
is controlled by \eqref{eq:r9-time-z-old} and
\eqref{eq:H2-final-windows}.

\emph{Step 3. Pointwise estimate for $\partial_3u_t$.}
Adding the pressure, transport, and viscosity estimates with
small fixed $\eta$ gives
\begin{equation}\label{eq:r9-Y-energy}
 E'+c_*\mathscr D_Y^2\le CE+R,\qquad
       \int_b^{b+1}R(s)\dd s\le C\mathcal J_2^2(1+b)^{-7/2}.
\end{equation}
The average of $Y=\partial_3u_t$ in
\eqref{eq:r9-time-z-old} has squared decay exponent $5/2$.
Lemma~\ref{lem:r9-uniform-gronwall}, applied to
\eqref{eq:r9-Y-energy} with \eqref{eq:r9-time-z-old},
therefore proves \eqref{eq:r9-Y-point}.
The uniform Gronwall argument uses the coefficient $C$ in
\eqref{eq:r9-Y-energy} only on unit intervals.
To justify the energy estimate, apply it to $\partial_3V_\sigma$.
Its time derivative is in
$H^{-1}$, and all coefficient quotients are controlled by
\eqref{eq:r9-qt}. Jensen controls the $K$ terms.
An averaged positive initial time and weak lower semicontinuity
give \eqref{eq:r9-Y-point} without $\partial_3u_t(0)$.

\emph{Step 4. Estimate for the mixed spatial derivative.}
The vertical derivative of the momentum equation gives
\begin{equation}\label{eq:r9-Z-recovery}
 \begin{aligned}
 (1+f)\nu_h\Delta_hZ={}&Y+u\cdot\nabla Z+Z\cdot\nabla u
       +c^2e^\tau\nabla\partial_3(q+\tau)\\
 &+c^2e^\tau\partial_3\tau\nabla(q+\tau)
       -(1+f)\nu_c\nabla\Div Z-\partial_3f\mathcal Lu .
 \end{aligned}
\end{equation}
The potential third derivative has decay exponent $7/4$
by \eqref{eq:r9-third-preliminary}.
The convective products have exponent at least $13/8$ except
$u_3\partial_3^2u$, whose exponent is $3/2$.
The last coefficient product has exponent $27/8$ by
\eqref{eq:r9-trace}. Thus
\[
 \|\nabla_h^2Z(t)\|_2
 \le C\|Y(t)\|_2+C\mathcal J_2(1+t)^{-7/4}
                         +C\mathcal J_2^2(1+t)^{-3/2},
\]
which, together with \eqref{eq:r9-Y-point}, proves
\eqref{eq:r9-mixed-third}. Expanding \eqref{eq:F2} and
using \eqref{eq:r9-third-preliminary} then gives
\eqref{eq:r9-F2-H1}; the slowest product is $u_3\partial_3^2u$.

\emph{Step 5. Refined decay estimates for fixed initial data.}
The full-energy convergence \eqref{eq:r8-H2-zero} implies
$\|\partial_3^2u(t)\|_2\to0$, so \eqref{eq:r9-F2-H1}
is $o(t^{-3/2})$ for fixed data.
For the mixed third derivative, the low-horizontal linear part satisfies
\[
 \|\nabla_h^2\partial_3P_L^hS_h(t)\mathbb Pu_0\|_2^2
 \le \|\nabla_h^4P_L^hS_h(t)\mathbb Pu_0\|_2
      \|\partial_3^2P_L^hS_h(t)\mathbb Pu_0\|_2=o(t^{-5/2}).
\]
The first factor is $O(t^{-5/2})$ and the second tends to zero.
In the proof of the second line of \eqref{eq:r8-windows},
the low-horizontal nonlinear contribution is
$O(t^{-3/2}\log(e+t))$, the high-horizontal contribution is
$O(e^{-c_*t}+t^{-3/2})$, and the potential contribution is
$O(t^{-7/4})$. Together with the linear estimate above, these
bounds give
$\|\nabla_h^2Z\|_{L^2(t,t+1;L^2)}=o(t^{-5/4})$.
Rearranging the spatial identity \eqref{eq:r9-Z-recovery},
use \eqref{eq:r9-third-preliminary} and the faster convective
and coefficient terms estimated in Step~4 to obtain the same
time-integrated asymptotic estimate for $Y$.
Apply the asymptotic refinement in
Lemma~\ref{lem:r9-uniform-gronwall} to \eqref{eq:r9-Y-energy};
its remainder has the strictly faster squared exponent $7/2$.
This proves $\|Y(t)\|_2=o(t^{-5/4})$.
Finally, \eqref{eq:r9-Z-recovery} transfers this limit to
$\nabla_h^2Z$.
\end{proof}

\subsection{Third-order estimates for the potential velocity and temperature}
\label{sec:r9-thermal-upper}

The preliminary third-order estimates in
Lemma~\ref{lem:r9-third-preliminary} and the mixed derivative estimate in
Proposition~\ref{prop:r9-mixed-third} provide the regularity needed below.
We improve the decay of $w$ and $\tau$ by combining a damped
divergence--temperature system with the low-frequency Green matrix.

\subsubsection{Scalar source estimates}

\begin{lemma}[Pointwise scalar source estimates]
\label{lem:r10-scalar-sources}
Under the assumptions of Theorem~\ref{thm:decay-H2}, the solution satisfies
\begin{equation}\label{eq:r10-scalar-sources}
       \|(F_1,F_3)(t)\|_{H^1}
                  \le C\mathcal J_2^2(1+t)^{-17/8},\qquad t\ge1 .
\end{equation}
\end{lemma}
\begin{proof}

For transport, \eqref{eq:new-uniform-basic} and
\eqref{eq:new-first} bound
$\|u\|_\infty\|\nabla(q,\tau)\|_2$ with decay exponent $17/8$.
Their differentiated terms use
$\|\partial_3u\|_3\|\nabla(q,\tau)\|_6$,
whose exponent is $17/8$ by
\eqref{eq:new-vertical-L3} and \eqref{eq:new-compressible-second}.
The product $\|u\|_\infty\|\nabla^2(q,\tau)\|_2$ decays
faster by \eqref{eq:new-uniform-basic} and
\eqref{eq:new-compressible-second}.

For heating, \eqref{eq:r9-lips}, \eqref{eq:new-G-structure},
and \eqref{eq:new-second-preliminary} control
$\|\nabla_hu\|_\infty\|\nabla\nabla_hu\|_2$,
with power $17/8$.
The divergence heating is faster because the potential
velocity has the full third-order bound
\eqref{eq:r9-third-preliminary}.
For diffusion, use \eqref{eq:new-uniform-basic},
\eqref{eq:new-compressible-second}, and
\eqref{eq:r9-third-preliminary}. The differentiated conductivity terms are
$(e^{-q}-1)\nabla\Delta\tau$ and
$\nabla e^{-q}\Delta\tau$; the first has power $13/4$,
and the second uses $L^6$ and $L^3$ with power at least $7/2$.
The gradient of $e^{-q}|\nabla\tau|^2$ uses
\eqref{eq:r9-third-preliminary} for $\|\nabla\tau\|_\infty$
and \eqref{eq:new-compressible-second} for
$\|\nabla^2\tau\|_2$, and decays faster than $(1+t)^{-17/8}$.
Coefficient derivatives on the heating terms use
$\|\nabla(q,\tau)\|_6$ and
$\|Q(u)\|_3\le C\mathcal J_2^2(1+t)^{-5/2}$.
Combining the transport, heating, and diffusion bounds proves
\eqref{eq:r10-scalar-sources}.
\end{proof}

\subsubsection{Potential velocity and thermal comparison}

\begin{proposition}[Third-order decay and thermal comparison]
\label{prop:r9-potential-thermal-third}
For the global solution with sufficiently small $\mathcal J_2$, the
potential velocity and temperature satisfy
\eqref{eq:r9-w-third}--\eqref{eq:r9-thermal-profile}, and the gradient bound
\eqref{eq:r9-compressible-lips} holds, for every $t\ge1$.
For each fixed initial datum,
$\|\nabla^3w(t)\|_2=o(t^{-2})$ and
$\|\nabla^3\tau(t)\|_2=o(t^{-17/8})$.
The thermal comparison \eqref{eq:r9-thermal-profile} also holds for
$\theta-\bar\theta$ with coefficient $\mu/(\bar\rho c_v\chi)$.
\end{proposition}

\begin{proof}
Throughout the proof, $t\ge1$ and unit intervals have left endpoint
$b\ge1$, unless an earlier integration interval is specified.
Write $Z=\partial_3u$ and $\mathfrak d=\Div u$; then $\Div V=\partial_t \mathfrak d$.
Step~1 gives preliminary high-frequency time estimates; Step~2
uses them to close the estimate for $V_3$.
In Step~3, the low-frequency gradient estimates and density relaxation
precede the closure for $\nabla_hu_3$.
The improved sources in Step~4 then give the third-order potential
velocity estimate, followed by the thermal comparison in Step~5.

\Needspace{6\baselineskip}
\smallskip\noindent\emph{Step 1. The damped divergence--temperature system.}

Eliminate $q_t=-\mathfrak d+F_1$ from \eqref{eq:r9-V-equation} and use
\eqref{eq:r9-Theta-equation}. This gives
\begin{equation}\label{eq:r9-div-time-system}
 \left\{\begin{aligned}
 \partial_t^2\mathfrak d-(\nu_h\Delta_h+\nu_c\Delta)\partial_t \mathfrak d-c^2\Delta \mathfrak d+c^2\Delta\Theta
       &=\Div(F_{2,t}-c^2\nabla F_1),\\
 \Theta_t-\chi\Delta\Theta+\alpha\partial_t \mathfrak d&=F_{3,t}.
 \end{aligned}\right.
\end{equation}
Let $\widetilde P_L$ be a fixed smooth cutoff supported in a sufficiently small ball in frequency space and equal to one near zero,
and put $\widetilde P_H=I-\widetilde P_L$.
Throughout this proof, subscripts $L,H$ refer to these cutoffs.
They differ from those in Section~\ref{sec:H2-decay}; on the compact
intermediate annulus, \eqref{eq:blockmultiplier} gives exponential
semigroup decay, so the change affects only the fixed constants.
For the high-frequency part of \eqref{eq:r9-div-time-system}, the energy
\[
 \|\partial_t \mathfrak d_H\|_2^2+c^2\|\nabla \mathfrak d_H\|_2^2
                    +(c^2/\alpha)\|\nabla\Theta_H\|_2^2
\]
cancels the divergence--temperature coupling.
Add a sufficiently small multiple of
\[
 \langle \partial_t \mathfrak d_H,\mathfrak d_H\rangle+
       \tfrac12(\nu_h\|\nabla_h\mathfrak d_H\|_2^2+\nu_c\|\nabla \mathfrak d_H\|_2^2).
\]
The derivative of this correction is
$\|\partial_t \mathfrak d_H\|_2^2-c^2\|\nabla \mathfrak d_H\|_2^2+
c^2\langle\nabla\Theta_H,\nabla \mathfrak d_H\rangle
-\langle(F_{2,t}-c^2\nabla F_1)_H,\nabla \mathfrak d_H\rangle$.
The Poincar\'e inequality on the support of $\widetilde P_H$ and Young's
inequality yield an energy $E_H$ equivalent to the sum above, with
\begin{equation}\label{eq:r9-damped-time}
 E_H'+c_*E_H+c_*(\|\nabla \partial_t \mathfrak d_H\|_2^2+\|\Delta\Theta_H\|_2^2)
 \le C(\|F_{2,t}-c^2\nabla F_1\|_2^2+\|F_{3,t}\|_2^2).
\end{equation}
The damping in \eqref{eq:r9-damped-time} is specific to the
compressible variables. We use the time difference quotients from
Lemma~\ref{lem:r9-time-derivatives} and start at a positive time
selected by the averaging argument in
Lemma~\ref{lem:r9-uniform-gronwall}; no additional derivatives of the
initial data are needed.

To improve the right-hand side of \eqref{eq:r9-damped-time}, separate
the vertical convection term from $F_{2,t}$:
\begin{equation}\label{eq:r9-F2t-split}
 F_{2,t}=-V_3Z+\mathcal E,\qquad
 \begin{aligned}
 \mathcal E={}&-V_h\cdot\nabla_hu-u\cdot\nabla V
       -c^2e^\tau\Theta\nabla(q+\tau)\\
 &-c^2(e^\tau-1)\nabla(q_t+\Theta)+f_t\mathcal Lu+f\mathcal LV
 \end{aligned}
\end{equation}
Using \eqref{eq:r9-qt}, \eqref{eq:r9-time-L2},
\eqref{eq:r9-V-point}, \eqref{eq:r9-Theta-point}, and
\eqref{eq:r9-Y-point} in this decomposition gives
\begin{equation}\label{eq:r9-time-remainders}
 \|\mathcal E\|_{L^2(b,b+1;L^2)}
 +\|(F_{1,t},F_{3,t})\|_{L^2(b,b+1;L^2)}
       \le C\mathcal J_2^2(1+b)^{-21/8}.
\end{equation}
The slowest term in $\mathcal E$ is $u_h\cdot\nabla_hV$:
\eqref{eq:new-uniform-basic} and \eqref{eq:r9-V-point} give the rate
$21/8$ in $L^2(b,b+1;L^2)$.
The product $V_h\cdot\nabla_hu$ decays with exponent $65/24$.
For the scalar sources, use
$F_{1,t}=-V\cdot\nabla q-u\cdot\nabla q_t$ and
\eqref{eq:r9-Theta-equation}--\eqref{eq:r9-S}.
Since $\partial_3V_3=\Div V-\Div_hV_h$, the Sobolev inequality
and \eqref{eq:r9-V-point} give
$\|V_3(t)\|_6\le C\mathcal J_2(1+t)^{-7/4}$;
\eqref{eq:new-vertical-L3} gives
$\|Z(t)\|_3\le C\mathcal J_2(1+t)^{-3/8}$.
Thus $V_3Z$ has decay exponent $17/8$, and
\eqref{eq:r9-F2t-split}--\eqref{eq:r9-time-remainders} control
$F_{2,t}$ in the same time-integrated norm.
Exponential convolution in \eqref{eq:r9-damped-time} now yields
\begin{equation}\label{eq:r9-pre-high-time}
 \begin{aligned}
 &\|\partial_t \mathfrak d_H(t)\|_2+\|\nabla \mathfrak d_H(t)\|_2+\|\nabla\Theta_H(t)\|_2\\
 &\quad+\|\nabla \partial_t \mathfrak d_H\|_{L^2(t,t+1;L^2)}
            +\|\Delta\Theta_H\|_{L^2(t,t+1;L^2)}
 \le C\mathcal J_2 e^{-c_*t}+C\mathcal J_2^2(1+t)^{-17/8}.
 \end{aligned}
\end{equation}
Since $(w_t)_H=\nabla\Delta^{-1}\partial_t \mathfrak d_H$, the Agmon inequality
\eqref{eq:r9-agmon} applied to \eqref{eq:r9-pre-high-time} gives
\eqref{eq:r9-wt-high} below. The thermal equation in
\eqref{eq:perturbed}, together with \eqref{eq:r10-scalar-sources}, gives
\eqref{eq:r9-tau-high-pre}:
\begin{align}
 \|(w_t)_H\|_{L^2(t,t+1;L^\infty)}
       &\le C\mathcal J_2 e^{-c_*t}+C\mathcal J_2^2(1+t)^{-17/8},
                                                        \label{eq:r9-wt-high}\\
 \|\tau_H\|_{H^2}+\|\nabla^3\tau_H\|_2
       &\le C\mathcal J_2 e^{-c_*t}+C\mathcal J_2^2(1+t)^{-17/8}.
                                                        \label{eq:r9-tau-high-pre}
\end{align}

\Needspace{6\baselineskip}
\smallskip\noindent\emph{Step 2. The time derivative of the third velocity component.}

The second-order estimates \eqref{eq:new-first},
\eqref{eq:new-compressible-second}, \eqref{eq:new-uniform-basic}, and
\eqref{eq:new-vertical-L3}, applied to
\eqref{eq:F1}--\eqref{eq:F3}, give
\[
 \|F_2\|_1\le C\mathcal J_2^2(1+t)^{-1},\quad
 \|(F_1,F_3)\|_1\le C\mathcal J_2^2(1+t)^{-7/4},\quad
 \|F_2\|_2\le C\mathcal J_2^2(1+t)^{-7/4}.
\]
In the Duhamel integral over $s\in[0,t/2]$, we use the decomposition
\begin{equation}\label{eq:r9-early-flux}
 \begin{aligned}
 F_2&=N_0+\partial_3B_0,\qquad B_0=-u_3u,\\
 N_0&=-u_h\cdot\nabla_hu+(\partial_3u_3)u
       -c^2(e^\tau-1)\nabla(q+\tau)+f\mathcal Lu,\\
 \|N_0(t)\|_1&\le C\mathcal J_2^2(1+t)^{-3/2},\qquad
 \|B_0(t)\|_1\le C\mathcal J_2^2(1+t)^{-5/4}.
 \end{aligned}
\end{equation}
Both $\|N_0(s)\|_1$ and $\|B_0(s)\|_1$ are integrable in time.
This decomposition avoids the logarithm that would result from
integrating the $(1+s)^{-1}$ bound for $F_2$.

Let $\mathscr S(t)$ denote the homogeneous semigroup for the full
linearized system \eqref{eq:perturbed}.
Differentiate its Duhamel formula and integrate by parts in time
on $s\in[t/2,t]$ to obtain
\begin{equation}\label{eq:r9-time-duhamel}
 \begin{aligned}
 U_t(t)={}&\partial_t\mathscr S(t)U_0
    +\int_0^{t/2}\partial_t\mathscr S(t-s)F(s)\dd s\\
 &+\mathscr S(t/2)F(t/2)
    +\int_{t/2}^t\mathscr S(t-s)F_t(s)\dd s .
 \end{aligned}
\end{equation}
The third-component estimates needed for the shear part of
\eqref{eq:r9-time-duhamel} are
\begin{equation}\label{eq:r9-third-row}
 \|\nabla_hS_h(a)(\mathbb Pg)_3\|_\infty\le Ca^{-2}\|g\|_1,
 \qquad \|S_h(a)(\mathbb Pg)_3\|_\infty\le Ca^{-3/4}\|g\|_2.
\end{equation}
Indeed, the vertical kernels of its Fourier row are multiples of
$|\xi_h|e^{-|\xi_h||x_3|}$, and the vertical square integral of
the row symbol is $C|\xi_h|$.
For the potential part at low spatial frequencies, the
$L^2$--$L^\infty$ semigroup bound is $C(1+a)^{-3/4}$ by
\eqref{eq:blockmultiplier}. In particular, the shear kernel
$a^{-3/4}$ and the low-frequency potential kernel are integrable
near $a=0$.

We estimate the third velocity component in
\eqref{eq:r9-time-duhamel}. The homogeneous term and the boundary
term $\mathscr S(t/2)F(t/2)$ contribute at most
$C\mathcal J_2(1+t)^{-2}$.
The shear integral over $[0,t/2]$ is bounded by
$C\mathcal J_2^2t^{-5/2}\log(e+t)$.
For the potential integral over $[0,t/2]$, use
\eqref{eq:r9-early-flux}: the time-differentiated $L^1$ kernel
decays with exponent $2$, and the kernel acting on $\partial_3B_0$
with exponent $5/2$.

It remains to control the integral over $[t/2,t]$, which contains
$F_{2,t}=-V_3Z+\mathcal E$.
Define the finite bootstrap maximum
$M_N=\max_{1\le n\le N}(1+n)^2\|V_3\|_{L^2(n,n+1;L^\infty)}$.
The anisotropic Sobolev estimate
$\|g\|_\infty\le C(\|g\|_{H^1}+\|\nabla_hg\|_{H^1})$
and \eqref{eq:r9-V-point}, \eqref{eq:r9-Y-point} ensure $M_N<\infty$.
Split the source interval $[t/2,t]$ into unit intervals. Summing the
kernel $(1+a)^{-3/4}$ gives $C(1+n)^{1/4}$ on the time interval
$[n,n+1]$. By \eqref{eq:r9-time-remainders} and
$\|Z(s)\|_2\le C\mathcal J_2(1+s)^{-1/4}$, the source is bounded by
\[
 C\mathcal J_2^2(1+n)^{-21/8}
       +C\mathcal J_2 M_N(1+n)^{-9/4}
\]
in $L^2_sL^2_x$ on each relevant unit interval.
For the convolution near $s=t$, apply Young's inequality in time.
Together with the high-frequency estimate \eqref{eq:r9-wt-high},
these bounds give
\[
 (1+n)^2\|V_3\|_{L^2(n,n+1;L^\infty)}
       \le C\mathcal J_2+C\mathcal J_2M_N,
       \qquad 1\le n\le N.
\]
All source intervals are covered by the intervals in $M_N$ or
by the fixed initial interval. Taking the maximum and absorbing
$C\mathcal J_2M_N$ gives a bound independent of $N$.
Let $N\to\infty$ and cover any unit interval by two consecutive
integer intervals. We obtain
\begin{equation}\label{eq:r9-V3-window}
 \|V_3\|_{L^2(t,t+1;L^\infty)}
       \le C\mathcal J_2(1+t)^{-2}.
\end{equation}
Insert \eqref{eq:r9-V3-window} and the $L^2$ bound for $Z$ into
\eqref{eq:r9-F2t-split}. The product $V_3Z$ has exponent
$9/4$, while $\mathcal E$ decays faster by
\eqref{eq:r9-time-remainders}. Hence
\begin{equation}\label{eq:r9-F2t-improved}
 \|F_{2,t}\|_{L^2(t,t+1;L^2)}
       \le C\mathcal J_2^2(1+t)^{-9/4}.
\end{equation}
For a fixed datum, \eqref{eq:r8-weak-littleo} gives
$\|Z(t)\|_2=o(t^{-1/4})$, so the left-hand side of
\eqref{eq:r9-F2t-improved} is $o(t^{-9/4})$.
For use in the Green-matrix primitive, retain the divergence form
\begin{equation}\label{eq:r9-late-flux}
 \begin{aligned}
 F_{2,t}&=N_1+\partial_3B_1,\quad
 B_1=-V_3u,\quad N_1=\mathcal E+(\partial_3V_3)u,\\
 \|B_1\|_{L^2(b,b+1;L^2)}
       &\le C\mathcal J_2^2(1+b)^{-5/2},\qquad
 \|N_1\|_{L^2(b,b+1;L^2)}
       \le C\mathcal J_2^2(1+b)^{-21/8}.
 \end{aligned}
\end{equation}

\Needspace{6\baselineskip}
\smallskip\noindent\emph{Step 3. Low-frequency gradients and density relaxation.}

\emph{(a) Low-frequency gradients.}
We first estimate the instantaneous low-frequency momentum forcing.
The Hilbert-valued horizontal Agmon inequality gives
\[
 \|g\|_{L^\infty_hL^2_3}
       \le C\|g\|_2^{1/2}\|\nabla_h^2g\|_2^{1/2}.
\]
By \eqref{eq:new-seed}, \eqref{eq:new-zero}, \eqref{eq:new-first},
\eqref{eq:new-second-preliminary}, \eqref{eq:new-third-final},
\eqref{eq:r9-third-preliminary}, and \eqref{eq:r9-mixed-third},
the corresponding $L^\infty_hL^2_3$ decay exponents for
$u,\nabla_hu,u_3,Z$ are $1,3/2,5/4,3/4$.
The kernel of $\nabla^j\widetilde P_L\mathbb Q$ satisfies
$|K_j(x)|\le C_j(1+|x|)^{-3-j}$ and hence
$\int_{\R^2}\sup_{x_3}|K_j(x)|\dd x_h<\infty$.
Near the origin this follows from the compact Fourier
support; away from it, differentiate the Coulomb potential
convolved with the cutoff kernel.
It follows that
\[
 \|\nabla^j\widetilde P_L\mathbb Qg\|_\infty
       \le C_j\|g\|_{L^\infty_hL^1_3},\qquad j=0,1,2.
\]
H\"older's inequality in $x_3$ bounds
$u_h\cdot\nabla_hu$ and $u_3Z$ in $L^\infty_hL^1_3$
with decay exponents $5/2$ and $2$, respectively.
The pressure and viscosity terms use the faster $L^2$ bounds from
\eqref{eq:new-compressible-second} and \eqref{eq:new-uniform-basic}. Thus
\begin{equation}\label{eq:r9-instantaneous-Q}
 \|\nabla^j\widetilde P_L\mathbb QF_2(t)\|_\infty
       \le C_j\mathcal J_2^2(1+t)^{-2},\qquad j=0,1,2.
\end{equation}
The cutoff $\widetilde P_L$ is essential in this kernel estimate.

For the columns of the compressible Green matrix corresponding
to momentum forcing, define
$\Psi_{\ell2}(a,\xi)=\int_a^\infty S_{\ell2}(s,\xi)\dd s$,
where $\ell=q,w,\tau$ and $|\xi|>0$.
At low spatial frequencies,
\begin{equation}\label{eq:r9-Green-primitive}
 |\Psi_{\ell2}(a,\xi)|\le C|\xi|^{-1}e^{-c_*|\xi|^2a},
 \quad
 \Psi_{w2}(0)=\Psi_{\tau2}(0)=0,\quad
 \nabla\Psi_{q2}(0)=c^{-2}\mathbb Q .
\end{equation}
To verify the primitive bound, set $r=|\xi|$, $\sigma=\xi/r$,
and $\widehat w_\parallel=i\sigma\cdot\widehat w$.
The generator on $(\widehat q,\widehat w_\parallel,\widehat\tau)$ is
\[
 \mathsf L(r,\sigma)=
 \begin{pmatrix}
 0&-r&0\\
 c^2r&-\nu(\sigma)r^2&c^2r\\
 0&-\alpha r&-\chi r^2
 \end{pmatrix},\qquad
 \nu(\sigma)=\nu_c+\nu_h|\sigma_h|^2.
\]
The matrix $\mathsf L_0=\lim_{r\downarrow0}r^{-1}\mathsf L(r,\sigma)$
has the distinct eigenvalues $0,\,\pm ic\sqrt{1+\alpha}$.
For its zero eigenvalue, right and left eigenvectors are
$(1,0,-1)^{\mathsf T}$ and $(-\alpha,0,1)$, respectively.
The characteristic polynomial and the spectral projectors therefore give,
for sufficiently small $r>0$,
\[
 \lambda_0(r,\sigma)=-\frac{\chi}{1+\alpha}r^2+O(r^3),\qquad
 \lambda_\pm(r,\sigma)=\pm ic\sqrt{1+\alpha}\,r+O(r^2).
\]
These expansions are uniform in $\sigma\in\mathbb S^2$, since
$\nu_c\le\nu(\sigma)\le\nu_c+\nu_h$ and the three eigenvalues
of $\mathsf L_0$ are separated. If $\mathsf E_j(r,\sigma)$ denotes
the corresponding spectral projector and
$e_\parallel=(0,1,0)^{\mathsf T}$, then
\[
 \|\mathsf E_0(r,\sigma)e_\parallel\|\le Cr,
 \qquad \|\mathsf E_\pm(r,\sigma)e_\parallel\|\le C.
\]
Indeed, $\mathsf E_0(0,\sigma)e_\parallel=0$ by the left
eigenvector above, and the projectors are smooth in $r$.
Equation~\eqref{eq:blockmultiplier} gives
$\Re\lambda_j(r,\sigma)\le-c_0r^2$ after decreasing the fixed
low-frequency radius. Consequently,
\[
 \left\|\int_a^\infty e^{s\mathsf L(r,\sigma)}
                         e_\parallel\dd s\right\|
 =\left\|-\sum_{j=0,+,-}
       \frac{e^{a\lambda_j(r,\sigma)}}{\lambda_j(r,\sigma)}
          \mathsf E_j(r,\sigma)e_\parallel\right\|
 \le Cr^{-1}e^{-c_0r^2a}.
\]
The passage between the longitudinal scalar and the physical velocity
has uniformly bounded multipliers, so this proves the bound for
$\Psi_{\ell2}$ in \eqref{eq:r9-Green-primitive}.
At $a=0$, integrate the linear equations for initial data supported in the momentum component. The resulting stationary equations give $\Div w=0$,
$\tau=0$, and $c^2\nabla q=\mathbb QF_2$.
Since $w$ is potential, this yields
\[
 \Psi_{q2}(0)=c^{-2}\Delta^{-1}\Div,\qquad
 \Psi_{w2}(0)=\Psi_{\tau2}(0)=0,
\]
and hence the endpoint identities in \eqref{eq:r9-Green-primitive}.
Since $\partial_a\Psi_{\ell2}(a)=-S_{\ell2}(a)$,
integration by parts over $[t/2,t]$ gives
\begin{equation}\label{eq:r9-primitive-ibp}
 \begin{aligned}
 \int_{t/2}^tS_{\ell2}(t-s)F_2(s)\dd s
 ={}&\Psi_{\ell2}(0)F_2(t)-\Psi_{\ell2}(t/2)F_2(t/2)\\
 &-\int_{t/2}^t\Psi_{\ell2}(t-s)F_{2,t}(s)\dd s .
 \end{aligned}
\end{equation}
After one, two, or three spatial derivatives, the low-frequency
$L^2$--$L^\infty$ norms of $\Psi_{\ell2}(a)$ decay with
exponents $3/4,5/4,7/4$, respectively.
For the endpoint $s=t$ in \eqref{eq:r9-primitive-ibp}, use
\eqref{eq:r9-Green-primitive} and \eqref{eq:r9-instantaneous-Q}.
For the integral containing $F_{2,t}$, use
\eqref{eq:r9-F2t-improved}. The one-derivative kernel has
unit-interval sum $O(t^{1/4})$, converting the source exponent
$9/4$ into $2$; the kernels with two and three derivatives are
integrable at large times. The integral over $[0,t/2]$ uses
\eqref{eq:r9-early-flux}, and the scalar source terms use the
$L^2$ bound \eqref{eq:r10-scalar-sources} with exponent $17/8$.
We obtain
\begin{equation}\label{eq:r9-low-gradients}
 \|\nabla q_L\|_{W^{2,\infty}}
       +\|\nabla w_L\|_\infty+\|\nabla\tau_L\|_\infty
       \le C\mathcal J_2(1+t)^{-2}.
\end{equation}

\smallskip\noindent\emph{(b) Density relaxation and high-frequency gradients.}
The positive order-zero operator
$\mathcal N=\nu_c+\nu_h\Delta_h\Delta^{-1}$ satisfies
$\nu_c\le\mathcal N(\xi)\le\nu_c+\nu_h$.
The divergence equation and continuity give the exact relaxation law
\begin{equation}\label{eq:r9-density-relaxation}
 q_t+u\cdot\nabla q+c^2\mathcal N^{-1}q
 =-c^2\mathcal N^{-1}\tau
       +\mathcal N^{-1}(-\Delta)^{-1}(\partial_t \mathfrak d-\Div F_2).
\end{equation}
To estimate $q_H$ in $H^2$, we use the transport bound
\begin{equation}\label{eq:r9-H2-transport}
 \left|\sum_{|\beta|\le2}
       \langle\partial^\beta(u\cdot\nabla\varphi),
                                  \partial^\beta\varphi\rangle\right|
 \le C(\|u\|_{H^2}+\|\nabla_hu\|_{H^2})\|\varphi\|_{H^2}^2.
\end{equation}
In the terms with two derivatives on $u$, the mixed-norm
H\"older inequality gives
\[
 \|\nabla^2u\|_{L^2_3L^4_h}
 \|\nabla\varphi\|_{L^\infty_3L^4_h}\|\nabla^2\varphi\|_2
 \le C\|u\|_{H^2}^{1/2}\|\nabla_hu\|_{H^2}^{1/2}
                                      \|\varphi\|_{H^2}^2.
\]
The terms with one velocity derivative follow from
\eqref{eq:H2-Lipschitz}.
Applying $\widetilde P_H$ to \eqref{eq:r9-density-relaxation}
gives
\[
 (q_H)_t+u\cdot\nabla q_H+c^2\mathcal N^{-1}q_H=R_q,
\]
where, in this argument,
\[
 R_q={}-c^2\mathcal N^{-1}\tau_H
 +\mathcal N^{-1}(-\Delta)^{-1}
                  (\partial_t \mathfrak d-\Div F_2)_H+\mathcal C_q,
 \qquad \mathcal C_q={}\widetilde P_L(u\cdot\nabla q)-u\cdot\nabla q_L.
\]
The fixed low-frequency cutoff and the product rule imply
\[
 \|\mathcal C_q(t)\|_{H^2}
 \le C\|F_1(t)\|_2
      +C\|u(t)\|_{H^2}\|\nabla q_L(t)\|_{W^{2,\infty}}
 \le C\mathcal J_2^2(1+t)^{-2},
\]
by \eqref{eq:r10-scalar-sources}, \eqref{eq:r9-low-gradients},
and \eqref{eq:new-base-energy}.
On the support of $\widetilde P_H$, the operator
$(-\Delta)^{-1}$ maps $L^2$ to $H^2$.
Thus \eqref{eq:r9-pre-high-time}, \eqref{eq:r9-tau-high-pre},
and the divergence bound in \eqref{eq:r9-point-sources-pre} give
\[
 \|R_q(t)\|_{H^2}
 \le C\mathcal J_2e^{-c_*t}+C\mathcal J_2^2(1+t)^{-2},
 \qquad t\ge1.
\]
Set $Y_q(t)=\|q_H(t)\|_{H^2}^2$ and
$\mathscr A(t)=\|u(t)\|_{H^2}+\|\nabla_hu(t)\|_{H^2}$.
The self-adjoint operator $c^2\mathcal N^{-1}$ has lower bound
$\kappa_q=c^2/(\nu_c+\nu_h)>0$.
Using \eqref{eq:r9-H2-transport} and Young's inequality in the
$H^2$ energy identity yields
\[
 Y_q'+\kappa_qY_q
 \le C\mathscr A(t)Y_q+C\|R_q(t)\|_{H^2}^2.
\]
Friedrichs regularization and \eqref{eq:r9-H2-transport}
justify this inequality for the strong solution.
For $t\ge s\ge1$, \eqref{eq:new-base-energy} gives
\[
 \int_s^t\mathscr A(r)\dd r
 \le C\mathcal J_2\bigl[(t-s)+(t-s)^{1/2}\bigr]
 \le C\mathcal J_2(1+t-s).
\]
Consequently, if $\mathcal J_2$ is sufficiently small,
\[
 \exp\left[-\kappa_q(t-s)+C\int_s^t\mathscr A(r)\dd r\right]
 \le C e^{-\kappa_q(t-s)/2}.
\]
Since $Y_q(1)\le C\mathcal J_2^2$, integration of the energy
inequality gives
\[
 Y_q(t)\le C\mathcal J_2^2e^{-c_*t}
                        +C\mathcal J_2^4(1+t)^{-4}.
\]
Taking square roots, we obtain
\begin{equation}\label{eq:r9-q-high}
 \|q_H(t)\|_{H^2}
       \le C\mathcal J_2 e^{-c_*t}+C\mathcal J_2^2(1+t)^{-2}.
\end{equation}
The elliptic equation \eqref{eq:r9-div-elliptic} gives the same
bound for $\|\nabla^3w_H\|_2$. Also,
$\|\nabla^2w_H\|_2$ has decay exponent $17/8$ by
\eqref{eq:r9-pre-high-time}. Interpolation in
\eqref{eq:r9-agmon} therefore gives exponent
$33/16>2$ for $\|\nabla w_H\|_\infty$.
Combining this with \eqref{eq:r9-low-gradients} and
\eqref{eq:r9-tau-high-pre} proves the $w,\tau$ part of
\eqref{eq:r9-compressible-lips}.

\smallskip\noindent\emph{(c) Closing the horizontal gradient of $u_3$.}
To complete \eqref{eq:r9-compressible-lips}, split the Duhamel
integral at $t/2$. The part over $[0,t/2]$ uses
\eqref{eq:r9-early-flux} and \eqref{eq:r9-third-row}.
In the flux term use
$\partial_3(\mathbb PB_0)_3=-\Div_h(\mathbb PB_0)_h$.
The two horizontal heat derivatives use
\[
 \|B_0\|_{L^2_3L^1_h}\le C\mathcal J_2^2(1+t)^{-3/2},\qquad
 \|\partial_3B_0\|_{L^2_3L^1_h}
       \le C\mathcal J_2^2(1+t)^{-5/4}.
\]
The geometric mean of these two mixed source bounds decays
with exponent $11/8>1$. The third-row semigroup estimate
\eqref{eq:new-mixed-kernel-all}, with one additional horizontal derivative,
then gives decay exponent $2$.
On $[t/2,t]$, expand $F_2$ as in \eqref{eq:F2}:
\[
 \|\nabla_hF_2(t)\|_2
 \le C\mathcal J_2^2(1+t)^{-9/4}
       +C\mathcal J_2(1+t)^{-1/4}\|\nabla_hu_3(t)\|_\infty.
\]
Set $M_h(T)=\sup_{1\le s\le T}(1+s)^2\|\nabla_hu_3(s)\|_\infty$,
finite by \eqref{eq:r9-lips}. The kernel in \eqref{eq:r9-third-row}
has exponent $3/4$ and time integral $O(t^{1/4})$, so
$M_h(T)\le C\mathcal J_2+C\mathcal J_2M_h(T)$.
Choosing $C\mathcal J_2\le1/2$ closes \eqref{eq:r9-compressible-lips}
uniformly in $T$.
The source integrals over $[0,1]$ are controlled by
\eqref{eq:new-base-energy}; no third derivative at time zero is used.

\Needspace{6\baselineskip}
\smallskip\noindent\emph{Step 4. Third-order decay of the potential velocity.}

Expand \eqref{eq:F1}--\eqref{eq:F3} using
\eqref{eq:r9-low-gradients}, \eqref{eq:r9-q-high}, and
\eqref{eq:r9-compressible-lips} to obtain
\begin{equation}\label{eq:r9-improved-sources}
 \begin{aligned}
 &\|F_1(t)\|_2\le C\mathcal J_2^2(1+t)^{-5/2},\qquad
 \|F_3(t)\|_2\le C\mathcal J_2^2(1+t)^{-19/8},\\
 &\|\nabla F_1(t)\|_2+\|\nabla_hF_2(t)\|_2+\|\Div F_2(t)\|_2
 \le C\mathcal J_2^2(1+t)^{-9/4}.
\end{aligned}
\end{equation}
For $\nabla F_1$, split $q=q_L+q_H$. The slowest differentiated
product is $\|Z\|_2\|\nabla q_L\|_\infty$, with exponent
$9/4$. For $F_3$, the horizontal heating term uses
\eqref{eq:r9-lips} and \eqref{eq:new-first}, with exponent
$19/8$.
In $\nabla_hF_2$ and $\Div F_2$, the slowest products are
$\|\nabla_hu_3\|_\infty\|Z\|_2$ and
$\|u_3\|_\infty\|\nabla_hZ\|_2$; both have exponent
$9/4$ by
\eqref{eq:r9-compressible-lips}, \eqref{eq:new-third-final},
and \eqref{eq:new-second-preliminary}.
Insert \eqref{eq:r9-F2t-improved} and
\eqref{eq:r9-improved-sources} in \eqref{eq:r9-damped-time}.
We obtain
\begin{equation}\label{eq:r9-high-time-final}
 \begin{aligned}
 &\|\partial_t \mathfrak d_H(t)\|_2+\|\nabla \mathfrak d_H(t)\|_2+\|\nabla\Theta_H(t)\|_2\\
 &\quad+\|\nabla \partial_t \mathfrak d_H\|_{L^2(t,t+1;L^2)}
          +\|\Delta\Theta_H\|_{L^2(t,t+1;L^2)}
 \le C\mathcal J_2 e^{-c_*t}+C\mathcal J_2^2(1+t)^{-9/4}.
 \end{aligned}
\end{equation}
The thermal equation in \eqref{eq:perturbed} and
\eqref{eq:r9-improved-sources} give the same bound for
$\|\tau_H\|_{H^2}$.

For $\nabla^3(w_L,\tau_L)$, use \eqref{eq:r9-primitive-ibp}.
The endpoint $s=t$ vanishes because
$\Psi_{w2}(0)=\Psi_{\tau2}(0)=0$ in
\eqref{eq:r9-Green-primitive}.
Insert the decomposition $F_{2,t}=N_1+\partial_3B_1$ from
\eqref{eq:r9-late-flux}. The term $N_1$ contributes at most
$C\mathcal J_2^2t^{-21/8}\log(e+t)$.
The kernel acting on the vertical flux $B_1$ decays as
$C(1+a)^{-3/2}$, so its contribution is
$C\mathcal J_2^2t^{-5/2}$. The boundary term at $s=t/2$ has
decay exponent $11/4$.
Use \eqref{eq:r9-early-flux} on $[0,t/2]$ and
\eqref{eq:r9-improved-sources} for the scalar source terms. These bounds give
\begin{equation}\label{eq:r9-low-third}
 \|\nabla^3(w_L,\tau_L)(t)\|_2
       \le C\mathcal J_2(1+t)^{-9/4}.
\end{equation}
Combine \eqref{eq:r9-low-third}, \eqref{eq:r9-q-high}, and
\eqref{eq:r9-div-elliptic} to prove \eqref{eq:r9-w-third}.
For each fixed datum, \eqref{eq:r8-H2-zero} gives
$\|u(t)\|_{H^2}\to0$. In the density relaxation argument, this
improves the cutoff commutator to $o(t^{-2})$.
The source bounds \eqref{eq:r9-pre-high-time},
\eqref{eq:r9-tau-high-pre}, and \eqref{eq:r9-low-third} decay
strictly faster than $t^{-2}$. The same damped estimate and
elliptic estimate therefore yield
$\|\nabla^3w(t)\|_2=o(t^{-2})$.

\Needspace{6\baselineskip}
\smallskip\noindent\emph{Step 5. Thermal comparison and third-order temperature decay.}
The horizontal viscous heating determines the thermal contribution.
Keep this term separate in the temperature equation:

\begin{equation}\label{eq:r9-thermal-elliptic}
 \begin{aligned}
 -\chi\Delta\tau-\frac{\mu}{\bar\rho\bar\theta c_v}|\nabla_hu|^2
 ={}&-\Theta-\alpha \mathfrak d-u\cdot\nabla\tau\\
 &+\frac{\mu(e^{-q-\tau}-1)}{\bar\rho\bar\theta c_v}|\nabla_hu|^2
 +\frac{(\lambda+\mu)e^{-q-\tau}}{\bar\rho\bar\theta c_v}\mathfrak d^2\\
 &+\chi f\Delta\tau+\chi e^{-q}|\nabla\tau|^2 .
 \end{aligned}
\end{equation}
Apply one gradient and $\widetilde P_H$ to \eqref{eq:r9-thermal-elliptic}.
The terms $\nabla\Theta_H$ and $\nabla \mathfrak d_H$ are controlled by
\eqref{eq:r9-high-time-final}. For the differentiated transport term,
\[
 \|\nabla u\|_2\|\nabla\tau\|_\infty
       +\|u\|_\infty\|\nabla^2\tau\|_2
       \le C\mathcal J_2^2(1+t)^{-9/4}.
\]
The remaining coefficient-weighted heating, divergence heating,
variable conductivity, and temperature-gradient-square terms
decay faster than $(1+t)^{-9/4}$.
The conductivity term is controlled by the preliminary bound
\eqref{eq:r9-third-preliminary} for $\nabla^3\tau$.
Elliptic inversion of \eqref{eq:r9-thermal-elliptic} gives
\eqref{eq:r9-thermal-profile} for the high-frequency part.
For the low-frequency part, use \eqref{eq:r9-low-third} and
\[
 \|\nabla^3\widetilde P_L(-\Delta)^{-1}|\nabla_hu|^2\|_2
 \le C\||\nabla_hu|^2\|_2
 \le C\mathcal J_2^2(1+t)^{-19/8}.
\]
Combining the two frequency regions proves
\eqref{eq:r9-thermal-profile} with the full differentiated Newton potential.

For the differentiated heating term,
\[
 \|\partial_3|\nabla_hu|^2\|_2
       \le2\|\nabla_hu\|_\infty\|\nabla_hZ\|_2
       \le C\mathcal J_2^2(1+t)^{-17/8},\qquad
 \|\nabla_h|\nabla_hu|^2\|_2
       \le C\mathcal J_2^2(1+t)^{-23/8}.
\]
The Fourier identity
$\|\nabla^3(-\Delta)^{-1}g\|_2=\|\nabla g\|_2$
for the ordered derivative norms, together with
\eqref{eq:r9-thermal-profile}, proves \eqref{eq:r9-tau-third}.
For fixed data, \eqref{eq:r8-weak-littleo} gives
$\|\nabla_hZ\|_2=o(t^{-3/4})$. The differentiated heating term
is therefore $o(t^{-17/8})$, and the remainder in
\eqref{eq:r9-thermal-profile} decays faster.

For the physical temperature, differentiate $\theta=\bar\theta e^\tau$ three times.
The coefficient correction
$\bar\theta(e^\tau-1)\nabla^3\tau$ and the products containing
$\nabla\tau\nabla^2\tau$ or $(\nabla\tau)^3$ decay faster than
$t^{-9/4}$. Here use \eqref{eq:new-uniform-basic} for
$\|\tau\|_\infty=O(t^{-3/2})$,
\eqref{eq:r9-compressible-lips} for
$\|\nabla\tau\|_\infty=O(t^{-2})$, and
\eqref{eq:new-compressible-second} for $\nabla^2\tau$.
Multiplying the leading term by $\bar\theta$ gives the
physical-temperature coefficient $\mu/(\bar\rho c_v\chi)$.
Together with Propositions~\ref{prop:r9-horizontal-third}
and~\ref{prop:r9-mixed-third}, this proves Theorem~\ref{thm:decay-H2}(iii).
\end{proof}

\section{Enhanced vertical decay}
\label{sec:r13-enhanced}

Throughout this section $\mathcal J_{2,1}$ is sufficiently small.
Starting from Theorem~\ref{thm:decay-H2}, we obtain vertical decay
without logarithmic loss and the enhanced third-order bounds.
Write $v=\mathbb Pu$, $w=\mathbb Qu$, and
$r_v=v-S_h(t)\mathbb Pu_0$.

\subsection{Preliminary estimates under first-derivative localization}
\begin{proposition}[Preliminary vertical estimates]\label{prop:new-intermediate}
If $\mathcal J_{2,1}$ is sufficiently small, then
\begin{align}\label{eq:new-intermediate}
 \|u\|_\infty&\le C\mathcal J_{2,1}(1+t)^{-1},\notag\\
 \|\partial_3u\|_2&\le C\mathcal J_{2,1}(1+t)^{-1/2}\log(e+t),\qquad
 \|\nabla_h\partial_3u\|_2\le C\mathcal J_{2,1}(1+t)^{-1}\log(e+t).
\end{align}
Thus $\|u\|_p\le C_p\mathcal J_{2,1}(1+t)^{-1+1/p}$ for
$2\le p\le\infty$, and
$\|U\|_{H^1}\le C\mathcal J_{2,1}(1+t)^{-1/2}\log(e+t)$.
Moreover,
\[
 \|w(t)\|_{L^\infty_3L^2_h}
       \le C\mathcal J_{2,1}(1+t)^{-5/4}.
\]
\end{proposition}
\begin{proof}
The estimates for $0\le t\le2$ follow from
\eqref{eq:new-base-energy}, so assume $t\ge2$.

\emph{Step 1. $L^\infty$ velocity decay and the mixed-norm bound for $w$.}
For a source $g$, the projected horizontal heat semigroup satisfies
\[
 \|\nabla_h^jS_h(a)\mathbb Pg\|_\infty
 \le Ca^{-1-j/2}
       \bigl(\|g\|_{L^2_3L^1_h}\|\partial_3g\|_{L^2_3L^1_h}\bigr)^{1/2},
 \qquad j=0,1.
\]
For $0<a<1$, Fourier Cauchy--Schwarz also gives
\[
 \|S_h(a)\mathbb Pg\|_\infty\le Ca^{-1/4}\|g\|_{H^1}.
\]
Indeed, its squared Fourier kernel factor is bounded by
$C\int_0^\infty r e^{-2\nu_har^2}(1+r^2)^{-1/2}\dd r\le Ca^{-1/2}$.

Split the Duhamel integral for $v$ into $[0,t/2]$,
$[t/2,t-1]$, and $[t-1,t]$.
On the first two intervals, use \eqref{eq:new-projected-source}:
apply the $a^{-1}$ kernel to $N$ and the $a^{-3/2}$ kernel to
$\nu_h(e^{-q}-1)\nabla_hu$.
By \eqref{eq:new-profile-source-ledger}, their source factors
decay as $(1+s)^{-9/8}$ and $(1+s)^{-15/8}$.
The contributions of $[0,t/2]$ are bounded by
$C\mathcal J_2^2(t^{-1}+t^{-3/2})$, and those of
$[t/2,t-1]$ by
$C\mathcal J_2^2[t^{-9/8}\log(e+t)+t^{-15/8}]$.
On $[t-1,t]$, use the $a^{-1/4}$ kernel with $F_2$.
Cauchy--Schwarz and the second line of \eqref{eq:new-F2H1} give
\[
 \int_{t-1}^t(t-s)^{-1/4}\|F_2(s)\|_{H^1}\dd s
 \le C\left(\int_{t-1}^t\|F_2(s)\|_{H^1}^2\dd s\right)^{1/2}
 \le C\mathcal J_2^2(1+t)^{-9/8}.
\]
The homogeneous term is bounded by
$C\mathcal J_{2,1}t^{-1}$ using the two mixed initial norms
in \eqref{eq:new-data-sizes}. Adding the bound for $w$ in
\eqref{eq:new-uniform-basic} proves the first line of
\eqref{eq:new-intermediate}.

The basic $L^2$ decay then gives
$\|u\otimes u\|_2\le C\mathcal J_{2,1}^2(1+t)^{-3/2}$.
Repeat Step~3 of the proof of Lemma~\ref{lem:r8-acoustic}.
On $[0,t/2]$, use \eqref{eq:r8-acoustic-kernel} with $j=0$
for the remainder and $j=1$ for the momentum flux.
On $[t/2,t]$, use the $L^2\to L^\infty_3L^2_h$ kernels
$(1+t-s)^{-1/4}$ and $(1+t-s)^{-3/4}$, respectively,
with the improved product bound above. This yields
\[
 \|w_L(t)\|_{L^\infty_3L^2_h}
 \le C\mathcal J_{2,1}(1+t)^{-5/4}
 +C\mathcal J_{2,1}^2
 \bigl[(1+t)^{-7/4}\log(e+t)+(1+t)^{-11/8}
                         +(1+t)^{-5/4}\bigr].
\]
For $w_H$, vertical Sobolev embedding and
\eqref{eq:new-compressible-second} give
$\|w_H\|_{L^\infty_3L^2_h}
\le C\|\nabla^2w_H\|_2
\le C\mathcal J_2(1+t)^{-7/4}$.
This proves the asserted mixed-norm bound.

\emph{Step 2. Preliminary vertical derivative decay.}
For $m=0,1$, \eqref{eq:Shestimate} and the initial norms in
\eqref{eq:new-data-sizes} give
\[
 \|\nabla_h^m\partial_3S_h(t)\mathbb Pu_0\|_2
       \le C\mathcal J_{2,1}(1+t)^{-(1+m)/2}.
\]
Using $u=S_h(t)\mathbb Pu_0+r_v+w$, combine this bound with
\eqref{eq:r8-sqrtlog} and the faster potential bounds
\eqref{eq:new-first}, \eqref{eq:new-compressible-second}
to obtain the vertical estimates in \eqref{eq:new-intermediate},
using the weaker logarithmic form of the square-root logarithm here.
Interpolation with the basic $L^2$ decay gives the stated
$L^p$ estimates; the basic compressible estimates and the
vertical derivative bounds give the asserted $H^1$ bound.
\end{proof}

\subsection{Vertical transport and first derivative decay}
\label{sec:r8-endpoint}

The preliminary bounds \eqref{eq:new-intermediate} enter the
verification of the coefficient and source hypotheses of
Lemma~\ref{lem:r8-forced-endpoint}.
To obtain the heat decay rates for the first vertical derivative
without a logarithmic loss, we retain the full drift
$u_3\partial_3$ in the evolution. Lemma~\ref{lem:r8-graph} controls
horizontal Fourier sections under vertical composition, and
Lemma~\ref{lem:r8-forced-endpoint} gives an abstract estimate for the
projected transport--diffusion equation. Lemma~\ref{lem:r8-real-drift}
verifies its coefficient hypotheses for $a=u_3$;
Lemma~\ref{lem:r13-first-vertical-decay} verifies the source hypotheses
and applies the estimate to obtain first vertical derivative decay.

\subsubsection{Vertical changes of variables}
\begin{lemma}[Horizontal Fourier estimates under vertical changes of variables]
\label{lem:r8-graph}
Let $\{\phi_{x_h}\}_{x_h\in\R^2}$ be a measurable family of increasing
bi-Lipschitz maps $\phi_{x_h}:\R\to\R$, with inverses $\psi_{x_h}$.
Assume that, for some $0<c_0\le1$ and a nonnegative function $d$,
\[
 c_0\le\partial_z\psi_{x_h}\le c_0^{-1},\qquad
 |\psi_{x_h}(z)-z|\le d(x_h),\qquad
 \|d\|_{L^{2,\infty}_h}<\infty.
\]
For $f\in L^2\cap L^2_3L^1_h$ with $\partial_3f\in L^2$, set
$f_\phi(x_h,z)=f(x_h,\phi_{x_h}(z))$. Then $\mathcal F_hf_\phi$
has a continuous $L^2_z$-valued representative and satisfies
\begin{equation}\label{eq:r8-graph}
 \sup_{\xi_h\in\R^2}
 \|\mathcal F_hf_\phi(\xi_h,\cdot)\|_{2,z}
       \le C(\|f\|_{L^2_3L^1_h}+\|d\|_{L^{2,\infty}_h}\|\partial_3f\|_2).
\end{equation}
If, in addition, $\|\partial_z\psi-1\|_\infty\le\varepsilon_J$, then
\begin{equation}\label{eq:r8-graph-difference}
 \begin{aligned}
 &\sup_{\xi_h\in\R^2}
 \|\mathcal F_hf_\phi(\xi_h,\cdot)-\mathcal F_hf(\xi_h,\cdot)\|_{2,z}\\
 &\qquad\le C\left[\varepsilon_J\|f\|_{L^2_3L^1_h}
 +\|d\|_{L^{2,\infty}_h}\|\partial_3f\|_2
 \sqrt{\log\left(e+\frac{\|f\|_{L^2_3L^1_h}}
                              {\|d\|_{L^{2,\infty}_h}\|\partial_3f\|_2}\right)}\right].
 \end{aligned}
\end{equation}
Here $C$ depends only on $c_0$.
The second summand is defined by continuous extension to be zero
when $\|d\|_{L^{2,\infty}_h}\|\partial_3f\|_2=0$.
\end{lemma}
\begin{proof}
\emph{Step 1. Duality and low vertical frequencies.}
First assume that $f$ has compact horizontal support.
Fix $\xi_h\in\R^2$ and test the Fourier difference against
$g\in L^2_z$, $\|g\|_2=1$. Changing the vertical variable pairs
$f$ with $\psi'g(\psi)-g$, where
$\psi=\psi_{x_h}$ and $\psi'=\partial_z\psi_{x_h}$.
Fix $\chi_0\in C_c^\infty(\R)$, equal to one on $[-1,1]$ and
supported in $[-2,2]$. For $0<\varepsilon\le1/2$, write
$g_{\rm low}=\chi_0((d(x_h)/\varepsilon)D_3)g$ and
$g_{\rm high}=g-g_{\rm low}$.
The displacement and Jacobian bounds give
\[
 |\psi'g_{\rm low}(\psi)-g_{\rm low}|
       \le C(\varepsilon_J+\varepsilon)M_3g,
\]
where $M_3$ is the one-dimensional Hardy--Littlewood maximal
operator. Its $L^2$ boundedness therefore bounds the
low-frequency pairing by
$C(\varepsilon_J+\varepsilon)\|f\|_{L^2_3L^1_h}$, uniformly in $\xi_h$.

\smallskip\noindent
\emph{Step 2. High vertical frequencies.}
Set $H_{x_h}=\partial_3^{-1}g_{\rm high}$, taking $H_{x_h}=0$
where $d(x_h)=0$. The identity
\[
 \psi'g_{\rm high}(\psi)-g_{\rm high}
       =\partial_3[H_{x_h}(\psi)-H_{x_h}]
\]
allows integration by parts, placing one vertical derivative on $f$.
For fixed $x_h$ with $d=d(x_h)>0$, split
$H=H_{\le d^{-1}}+H_{>d^{-1}}$ at vertical frequency $d^{-1}$.
The maps
$\psi_\vartheta(z)=(1-\vartheta)z+\vartheta\psi(z)$,
$0\le\vartheta\le1$, satisfy
$\partial_z\psi_\vartheta\ge\min\{1,c_0\}$.
The fundamental theorem of calculus therefore gives
\[
 \|H_{\le d^{-1}}\circ\psi-H_{\le d^{-1}}\|_2
 \le d\int_0^1
       \|\partial_3H_{\le d^{-1}}\circ\psi_\vartheta\|_2
           \dd\vartheta
 \le C d\|\partial_3H_{\le d^{-1}}\|_2.
\]
For the high-frequency part, the Jacobian bound gives
$\|H_{>d^{-1}}\circ\psi-H_{>d^{-1}}\|_2
\le C\|H_{>d^{-1}}\|_2$.
Combining the two estimates and using Plancherel yields
\[
 \|H\circ\psi-H\|_2^2
       \le C\int_\R\min(d^2\zeta^2,1)
                                     |\widehat H(\zeta)|^2\dd\zeta.
\]
Since $\chi_0=1$ on $[-1,1]$,
\[
 |\widehat H_{x_h}(\zeta)|
 \le C|\zeta|^{-1}|\widehat g(\zeta)|
                    \mathbf1_{\{|\zeta|d(x_h)\ge\varepsilon\}}.
\]
For $\varrho=|\zeta|>0$, the distribution estimate
$|\{d>s\}|\le\|d\|_{L^{2,\infty}_h}^2s^{-2}$ gives
\[
 \begin{aligned}
 &\int_{\R^2}\min(d(x_h)^2,\varrho^{-2})
                  \mathbf1_{\{\varrho d(x_h)\ge\varepsilon\}}\dd x_h\\
 &\quad\le (\varepsilon/\varrho)^2
                  |\{d\ge\varepsilon/\varrho\}|
       +2\int_{\varepsilon/\varrho}^{1/\varrho}s|\{d>s\}|\dd s
 \le C\|d\|_{L^{2,\infty}_h}^2\log(e/\varepsilon).
\end{aligned}
\]
Tonelli's theorem and the preceding two estimates now imply
\[
 \int_{\R^2}\|H_{x_h}\circ\psi_{x_h}-H_{x_h}\|_{2,z}^2\dd x_h
 \le C\|d\|_{L^{2,\infty}_h}^2
                          \log(e/\varepsilon)\|g\|_2^2.
\]
After the vertical integration by parts, Cauchy--Schwarz in
$(x_h,z)$ bounds the high-frequency pairing by
$C\|d\|_{L^{2,\infty}_h}\|\partial_3f\|_2
\sqrt{\log(e/\varepsilon)}$, uniformly in $\xi_h$.

\smallskip\noindent
\emph{Step 3. Choice of the cutoff and continuity.}
Combining Steps~1 and~2 with
\[
 \varepsilon=\min\left\{\frac12,\,
 \frac{\|d\|_{L^{2,\infty}_h}\|\partial_3f\|_2}
      {\|f\|_{L^2_3L^1_h}}\right\}
\]
gives \eqref{eq:r8-graph-difference} when the numerator and
denominator are positive.
If $\|d\|_{L^{2,\infty}_h}\|\partial_3f\|_2=0$, let $\varepsilon\downarrow0$;
if $\|f\|_{L^2_3L^1_h}=0$, then $f=0$.
The Jacobian bounds give $\|\partial_z\psi-1\|_\infty\le1+c_0^{-1}$.
Using this bound, $B\sqrt{\log(e+A/B)}\le C(A+B)$ for $A,B>0$
(with its continuous extension at $B=0$), and the triangle inequality
in \eqref{eq:r8-graph-difference} yields \eqref{eq:r8-graph}.

For $f_R=\mathbf1_{\{|x_h|\le R\}}f$, the Jacobian bounds imply
$(f_R)_\phi\in L^1_hL^2_z$, so its horizontal Fourier transform
is continuous in $\xi_h$.
Applying \eqref{eq:r8-graph} to $f-f_R$ gives uniform convergence
in $\xi_h$ as $R\to\infty$. Since $(f_R)_\phi\to f_\phi$ also in
$L^2$, the limit is a continuous representative of $\mathcal F_hf_\phi$.
\end{proof}

\subsubsection{A forced transport--heat estimate}
For a real coefficient $a(t,x)$ and $\delta\ge0$, consider the
displacement condition
\begin{equation}\label{eq:r8-drift-displacement}
 \sup_{b>0}\sup_\gamma
 \left\|\int_0^b\sup_z|a(s,x_h+\gamma(s),z)|\dd s
                                      \right\|_{L^{2,\infty}_h}\le\delta,
\end{equation}
where $\gamma$ ranges over measurable horizontal translation
paths. The path supremum is outside the spatial norm.

\Needspace{8\baselineskip}
\begin{lemma}[Decay estimates for the projected transport--diffusion equation]
\label{lem:r8-forced-endpoint}
Consider
\[
 G_t-\nu_h\Delta_hG+\mathbb P(a\partial_3G)
       =\mathbb P(S+\Div_hB),\qquad G(0)=\mathbb Pf_0,
\]
where $a$ is real-valued, $f_0,S$ are $\R^3$-valued, and
$B$ is a $\R^{3\times2}$-valued horizontal flux.
Suppose $a\in L^1_{\rm loc}([0,\infty);W^{1,\infty})$ satisfies
\eqref{eq:r8-drift-displacement} and, for all $t,b\ge0$,
\begin{subequations}\label{eq:r8-drift-assumptions}
\begin{align}
 \begin{split}
 \int_0^\infty\|\nabla a\|_\infty\dd t&\le\delta,\\
 \int_0^\infty\|\nabla_ha\|_{L^2_hL^\infty_3}\dd t&\le\delta
 \end{split}
 &\label{eq:r13-drift-characteristics}\\
 \|a(t)\|_\infty\le\delta(1+t)^{-3/2}
 &\label{eq:r13-drift-gradient}\\
 \begin{split}
 \|\nabla a(t)\|_{L^\infty_3L^2_h}&\le\delta(1+t)^{-7/8},\\
 \|\nabla a\|_{L^2(b,b+1;L^\infty)}
       &\le\delta(1+b)^{-11/8},\\
 \int_0^\infty(1+t)^{-2/3}\|\partial_3\nabla a\|_3\dd t
       &\le\delta
 \end{split}
 &\label{eq:r13-drift-projection}
\end{align}
\end{subequations}
The characteristic and flux changes of variables use
\eqref{eq:r13-drift-characteristics}.
The horizontal-gradient estimate uses \eqref{eq:r13-drift-gradient},
and the projection correction uses \eqref{eq:r13-drift-projection}.
Assume $f_0\in H^1\cap L^2_3L^1_h$ and
\begin{equation}\label{eq:r8-source-space}
\begin{aligned}
 &\Big[\|f_0\|_{H^1}+\|f_0\|_{L^2_3L^1_h}
       +\int_0^\infty
          \big(\|S\|_{L^2_3L^1_h}+\|\partial_3S\|_2\big)\dd t\Big]\\[2pt]
 &\quad+\Big[\int_0^\infty
          \big(\|B\|_{L^2_3L^1_h}+\|\partial_3B\|_2\big)\dd t
       +\|\partial_3B\|_{L^2_tL^2}
       +\sup_{t\ge0}(1+t)\|B(t)\|_2\Big]\\[2pt]
 &\quad+\sup_{b\ge0}\Big[(1+b)^{1/2}\int_b^\infty\|S\|_2\dd t\\
 &\hspace{5em}+(1+b)\big(\|S\|_{L^2(b,\infty;L^2)}
                     +\|\Div_hB\|_{L^2(b,\infty;L^2)}\big)\Big]
       \le N.
\end{aligned}
\end{equation}
In \eqref{eq:r8-source-space}, the mixed-norm and time-integrability
bounds control the low horizontal frequencies and vertical energy.
The weighted pointwise bound for $B$ gives the $L^2$ decay;
the tails of $S$ and $\Div_hB$ control the late-time heat integral.
For a vector field $F$, define
\[
 \|F\|_{\mathcal X}
 :=\sup_{t\ge0}\big[(1+t)^{1/2}\|F\|_2
       +(1+t)\|\nabla_hF\|_2+\|\partial_3F\|_2\big]
       +\|\nabla_h\partial_3F\|_{L^2_tL^2}.
\]
For sufficiently small $\delta$, the solution satisfies
\begin{equation}\label{eq:r8-forced-endpoint}
 \|G\|_{\mathcal X}\le CN.
\end{equation}
\end{lemma}
\begin{proof}
\emph{Step 1. Time splitting and vertical characteristics.}
First consider the unprojected equation
\[
 F_t-\nu_h\Delta_hF+a\partial_3F=S+\Div_hB,\qquad F(0)=f_0,
\]
componentwise. Steps~1--3 estimate this equation; Step~4 adds
the correction needed to obtain $G=\mathbb PF$.
Fix $s<t$ and an initial field $f\in H^1\cap L^2_3L^1_h$ at time $s$.
For a partition $s=t_0<\cdots<t_m=t$, let
$\mathcal T_a(t_{k+1},t_k)$ denote the vertical transport propagator
associated with $\partial_r+a(r,x)\partial_3$, and set
\[
 F_{k+1}=S_h(t_{k+1}-t_k)\mathcal T_a(t_{k+1},t_k)F_k,
 \qquad F_0=f.
\]
Expand each heat operator as a Gaussian convolution. For every
fixed set of horizontal convolution variables, the product is a
horizontal translation followed by a vertical composition
$f(x_h,\phi_{x_h}(z))$. After absorbing the terminal translation,
the coefficient is evaluated along a piecewise constant horizontal
translation path $\gamma$. Both the forward and inverse vertical
characteristics have displacement bounded by
\[
 d_\gamma(x_h)=\int_s^t\sup_z|a(r,x_h+\gamma(r),z)|\dd r.
\]
Extending $\gamma$ to $[0,t]$, \eqref{eq:r8-drift-displacement}
gives $\|d_\gamma\|_{L^{2,\infty}_h}\le\delta$.
Writing $\psi=\phi^{-1}$, differentiation of the characteristics
and \eqref{eq:r13-drift-characteristics} give
\[
 e^{-\delta}\le\partial_z\psi\le e^\delta,\qquad
 \left\|\sup_z|\nabla_h\phi_{x_h}(z)|\right\|_{2,h}
 \le C\int_s^t\|\nabla_ha\|_{L^2_hL^\infty_3}\dd r\le C\delta.
\]
These bounds are independent of the partition and the convolution
variables.

\smallskip\noindent
\emph{Step 2. Low horizontal frequencies and passage to the limit.}
Apply \eqref{eq:r8-graph} to the entire vertical composition and use
the displacement bound $\|d_\gamma\|_{L^{2,\infty}_h}\le\delta$ from Step~1.
The terminal translation gives a Fourier factor of modulus one.
Since the product of the Gaussian kernels has total mass one,
Minkowski's inequality yields
\[
 \sup_{\xi_h}\|\mathcal F_hF_m(\xi_h,\cdot)\|_{2,3}
 \le C(\|f\|_{L^2_3L^1_h}+\delta\|\partial_3f\|_2).
\]
We next identify the limit of these products.
Write $h_k=t_{k+1}-t_k$ and $\mathfrak h=\max_k h_k$.
For $r\in[t_k,t_{k+1}]$ and $0\le\theta\le h_k$, set
\[
 X_k(r)=\mathcal T_a(r,t_k)F_k,\qquad
 Y_k(\theta)=S_h(\theta)X_k(t_{k+1}).
\]
The transport $H^1$ estimate and heat contractivity give
\[
 \sup_k\left(\sup_{r\in[t_k,t_{k+1}]}\|X_k(r)\|_{H^1}
       +\sup_{0\le\theta\le h_k}\|Y_k(\theta)\|_{H^1}\right)
 \le M,
 \quad
 M=C e^{C\int_s^t\|\nabla a(r)\|_\infty\dd r}\|f\|_{H^1}.
\]
Put
\[
 A_k=\int_{t_k}^{t_{k+1}}\|a(r)\|_\infty\dd r,
 \qquad
 B_{s,t}=\int_s^t\|a(r)\|_{W^{1,\infty}}\dd r.
\]
The transport equation and the Fourier estimate
$\|(S_h(\theta)-I)g\|_2\le C\theta^{1/2}\|\nabla_hg\|_2$
give
\[
 \|X_k(r)-F_k\|_2\le M A_k,
 \qquad
 \|Y_k(\theta)-F_k\|_2\le CM(A_k+\theta^{1/2}).
\]
Let $\varphi\in C_c^\infty([s,t]\times\R^3)$ and
$\varphi_k=\varphi(t_k,\cdot)$.
Integration of the two substep equations gives the exact identity
\[
 \langle F_{k+1}-F_k,\varphi_k\rangle
 ={}\int_{t_k}^{t_{k+1}}
    \langle X_k(r),a(r)\partial_3\varphi_k
                       +\partial_3a(r)\varphi_k\rangle\dd r
 +\nu_h\int_0^{h_k}
                      \langle Y_k(\theta),\Delta_h\varphi_k\rangle
                         \dd\theta.
\]
Define $\overline F_{\mathfrak h}(r)=F_k$ for
$t_k\le r<t_{k+1}$, with value $F_m$ at $r=t$.
Replacing $X_k(r)$ by $F_k$ in the transport integral contributes,
after summation, at most
$C_\varphi M(\max_k A_k)B_{s,t}$.
The corresponding replacement in the heat integral is bounded by
\[
 C_\varphi M\sum_k\int_0^{h_k}(A_k+\theta^{1/2})\dd\theta
 \le C_\varphi M\bigl[\mathfrak h B_{s,t}
                             +\mathfrak h^{1/2}(t-s)\bigr].
\]
Replacing $\varphi_k$ by $\varphi(r)$ in these integrals adds at most
$C_\varphi M\mathfrak h(B_{s,t}+t-s)$.
Finally, summing
\[
 \langle F_{k+1},\varphi(t_{k+1})\rangle
       -\langle F_k,\varphi(t_k)\rangle
 =\langle F_{k+1}-F_k,\varphi_k\rangle
       +\int_{t_k}^{t_{k+1}}
                    \langle F_{k+1},\partial_r\varphi(r)\rangle\dd r
\]
and replacing $F_{k+1}$ by $F_k$ in the last integral incurs
another $C_\varphi M[\mathfrak h B_{s,t}
+\mathfrak h^{1/2}(t-s)]$.
We obtain
\[
 \begin{aligned}
 &\langle F_m,\varphi(t)\rangle-\langle f,\varphi(s)\rangle
 =\int_s^t\bigl\langle\overline F_{\mathfrak h}(r),
       \partial_r\varphi+\nu_h\Delta_h\varphi
                   +a\partial_3\varphi+(\partial_3a)\varphi
                                      \bigr\rangle\dd r
       +\mathfrak r_{\mathfrak h}(\varphi),\\
 &|\mathfrak r_{\mathfrak h}(\varphi)|
 \le C_\varphi M\bigl[
       (\max_k A_k+\mathfrak h)B_{s,t}
                   +\mathfrak h^{1/2}(t-s)\bigr],\qquad
                   0<\mathfrak h\le1.
\end{aligned}
\]
Absolute continuity of the integral of $\|a(r)\|_\infty$ gives
$\max_kA_k\to0$ as $\mathfrak h\to0$, so the full residual tends
to zero.
For every $\vartheta\in C_c^\infty(\R^3)$, the substep equations also give,
at partition nodes,
\[
 |\langle F_j-F_i,\vartheta\rangle|
 \le C_\vartheta M\left[(t_j-t_i)
                   +\int_{t_i}^{t_j}\|a(r)\|_\infty\dd r\right].
\]
The piecewise affine interpolants are therefore equicontinuous in
the weak $L^2$ topology, and their difference from
$\overline F_{\mathfrak h}$ tends to zero in the same topology,
uniformly on $[s,t]$.
The uniform $H^1$ bound and a diagonal weak compactness argument
give a limit in $C([s,t];L^2_{\rm weak})\cap L^\infty(s,t;H^1)$.
Since
$a\partial_3\varphi+(\partial_3a)\varphi\in L^1(s,t;L^2)$,
the displayed weak identity passes to the limit and gives
the scalar transport--diffusion equation with initial value $f$.
Its equation gives a time derivative in $L^1(s,t;H^{-1})$;
interpolation between $H^{-1}$ and the uniform $H^1$ bound
therefore gives $L^2$ continuity, including at the initial time.
Its $L^2$ energy identity, obtained by spatial regularization and
time averaging, yields uniqueness by
$\partial_3a\in L^1(s,t;L^\infty)$.
Denote this evolution by $\mathcal U_a(t,s)$.
At each terminal time, the Fourier bound passes to this limit by
weak-star compactness in $L^\infty_{\xi_h} L^2_3$.

For a horizontal flux, the chain rule gives
\begin{equation}\label{eq:r8-flux-chain}
 (\Div_hB)(x_h,\phi(z))
       =\Div_h[B(x_h,\phi(z))]
             -\nabla_h\phi(z)\cdot\partial_3B(x_h,\phi(z)).
\end{equation}
The first term has a Fourier factor $i\xi_h$ and is bounded using
\eqref{eq:r8-graph}. Minkowski's inequality, horizontal
Cauchy--Schwarz, and the Jacobian bounds give
\[
 \|\nabla_h\phi\cdot\partial_3B(x_h,\phi)\|_{L^2_3L^1_h}
 \le\|\sup_z|\nabla_h\phi|\|_{2,h}
       \|\partial_3B(x_h,\phi)\|_2
 \le C\delta\|\partial_3B\|_2.
\]
Applying these estimates at each source time in the Duhamel formula,
\eqref{eq:r8-source-space} yields
\[
 \sup_{t\ge0}\operatorname*{ess\,sup}_{|\xi_h|\le1}
       \|\mathcal F_hF(t,\xi_h,\cdot)\|_{2,3}\le CN.
\]
To extend the estimates to \eqref{eq:r8-source-space}, convolve
$f_0,S,B$ in space, truncate the sources smoothly in time, and
mollify them in time after extension by zero to $t<0$, using a
nonnegative one-sided mollifier supported in $[0,n^{-1}]$.
The approximations satisfy \eqref{eq:r8-source-space} with $CN$
uniformly in $n$. They converge in the initial $H^1$ and mixed norms
and in every unweighted time-integrable norm in
\eqref{eq:r8-source-space}. In particular,
$(S_n,\partial_3S_n)$ converge in $L^1_tL^2$, and
$(B_n,\partial_3B_n)$ converge in $L^2_tL^2$.
The $L^2$ energy estimates in Step~3, applied to the difference of
two approximating solutions and its vertical derivative, give
convergence of $F_n,\partial_3F_n$ in $C([0,T];L^2)$ for each
finite $T$. Weak-star compactness preserves the horizontal Fourier
bound, while weak lower semicontinuity and Fatou's lemma preserve
the energy and weighted bounds of Step~3.

\smallskip\noindent
\emph{Step 3. Energy estimates and decay.}
Testing the scalar equation and its vertical derivative in $L^2$,
then applying Gronwall's inequality with
\eqref{eq:r8-drift-assumptions} and \eqref{eq:r8-source-space}, gives
\[
 \sup_{t\ge0}(\|F(t)\|_2^2+\|\partial_3F(t)\|_2^2)
 +\int_0^\infty
       (\|\nabla_hF\|_2^2+\|\nabla_h\partial_3F\|_2^2)\dd t
       \le CN^2.
\]
Here $B\in L^2_tL^2$ follows from its pointwise bound in
\eqref{eq:r8-source-space}.
First take $S=0$. The low-frequency estimate in Step~2 gives
\[
 \|F\|_2^2\le CN^2r^2+r^{-2}\|\nabla_hF\|_2^2,
 \qquad 0<r\le1.
\]
For $N>0$, choose $r^2$ proportional to $\|F\|_2^2/N^2$;
the uniform energy bound permits $r\le1$. The energy inequality
and $\|B(t)\|_2\le N(1+t)^{-1}$ imply
\[
 (\|F\|_2^2)'+c_*N^{-2}\|F\|_2^4
       \le\|\partial_3a\|_\infty\|F\|_2^2+CN^2(1+t)^{-2}.
\]
An integrating factor and comparison with $CN^2/(1+t)$ yield
$\|F(t)\|_2\le CN(1+t)^{-1/2}$.
For the homogeneous evolution, the same argument on $[s,t]$
gives this decay with elapsed time $t-s$ and data size
$\|f\|_2+\|f\|_{L^2_3L^1_h}+\delta\|\partial_3f\|_2$.
Consequently, the source term is estimated over $[0,t/2]$
by this homogeneous bound and over $[t/2,t]$ by the $L^2$
stability of $\mathcal U_a$ and the $L^1_tL^2$ tail in
\eqref{eq:r8-source-space}. This proves the same $L^2$ decay
for general $S$.

For $t\ge1$, the vertical energy bound and
$\|a(s)\|_\infty\le\delta(1+s)^{-3/2}$ give
\[
 \|S+\Div_hB-a\partial_3F\|_{L^2(t/2,t;L^2)}
       \le CN(1+t)^{-1}.
\]
Writing the horizontal heat Duhamel formula from $t/2$ and
applying \eqref{eq:r8-horizontal-maximal}, we obtain
\[
 \|\nabla_hF(t)\|_2
 \le Ct^{-1/2}\|F(t/2)\|_2
       +C\|S+\Div_hB-a\partial_3F\|_{L^2(t/2,t;L^2)}
 \le CN(1+t)^{-1}.
\]
For $0\le t\le1$, the formula starting at zero uses
$\|f_0\|_{H^1}\le N$ and the same smoothing estimate.
Together with the vertical energy estimate, this proves
\eqref{eq:r8-forced-endpoint} for the unprojected equation.

\smallskip\noindent
\emph{Step 4. The projection correction.}
Define the correction
$H(F)=-(\nabla a)\partial_3\Delta^{-1}\Div F$.
The $L^2$ multiplier bound and \eqref{eq:r13-drift-projection}
imply
\[
 \|H(F)(t)\|_{L^2_3L^1_h}
       \le C\delta\|F\|_{\mathcal X}(1+t)^{-11/8},
 \qquad \|H(F)\|_{L^2(b,b+1;L^2)}
       \le C\delta\|F\|_{\mathcal X}(1+b)^{-15/8}.
\]
The first bound is time integrable. Summing the second bound over
unit intervals gives both time tails required of the source term
in \eqref{eq:r8-source-space}.
For the vertical derivative, the anisotropic Sobolev inequality gives
\[
 \|\partial_3\Delta^{-1}\Div F\|_6
       \le C\|\nabla_hF\|_2^{2/3}\|\partial_3F\|_2^{1/3}
       \le C\|F\|_{\mathcal X}(1+t)^{-2/3}.
\]
The weighted integral in \eqref{eq:r13-drift-projection} controls
the product with $\partial_3\nabla a$; the $L^1_tL^\infty$
bound in \eqref{eq:r13-drift-characteristics} controls the
other differentiated product. Thus
$\int_0^\infty\|\partial_3H(F)\|_2\dd t
\le C\delta\|F\|_{\mathcal X}$.
Therefore $H(F)$ satisfies the source conditions in
\eqref{eq:r8-source-space} with size $C\delta\|F\|_{\mathcal X}$.
Steps~1--3 and a contraction for small $\delta$ solve
\[
 F_t-\nu_h\Delta_hF+a\partial_3F=S+\Div_hB+H(F),
 \qquad F(0)=f_0,
\]
with $\|F\|_{\mathcal X}\le CN$.
Finally,
\[
 \mathbb P(a\partial_3\mathbb QF)=\mathbb PH(F)
\]
shows that $G=\mathbb PF$ solves the stated projected equation.
The $L^2$ boundedness of $\mathbb P$ gives
\eqref{eq:r8-forced-endpoint}, and the projected energy estimate
proves uniqueness.
\end{proof}

\subsubsection{Estimates for the vertical transport coefficient}
\begin{lemma}[Bounds for the vertical transport coefficient]\label{lem:r8-real-drift}
For sufficiently small $\mathcal J_{2,1}$, the coefficient $a=u_3$
satisfies \eqref{eq:r8-drift-displacement} and
\eqref{eq:r8-drift-assumptions} with $\delta=C\mathcal J_2$.
\end{lemma}
\begin{proof}
\emph{Step 1. Time integrability of the unprojected sources.}
The unprojected sources $(F_1,F_2,F_3)$ are defined in
\eqref{eq:F1}--\eqref{eq:F3}. Theorem~\ref{thm:decay-H2} and
\eqref{eq:new-intermediate} give
\[
 \|u_3\partial_3u\|_1
       \le\|u_3\|_2\|\partial_3u\|_2
       \le C\mathcal J_2\mathcal J_{2,1}
                       (1+t)^{-5/4}\log(e+t).
\]
The other $L^1$ products in $(F_1,F_2,F_3)$ are time integrable by
Theorem~\ref{thm:decay-H2}.
For the $H^1$ integral on $[1,\infty)$, use
\eqref{eq:new-F2H1} for $F_2$ and
\eqref{eq:r10-scalar-sources} for $F_1,F_3$, with
Cauchy--Schwarz and H\"older in time for the terms containing
$\cD_2(U)$. On $[0,1]$, the uniform $H^2$ bound in
\eqref{eq:new-base-energy} controls $F_1$, while
\eqref{eq:new-F2H1} controls $F_2$ and
\eqref{eq:r9-source-thermal-old} with $b=0$ controls $F_3$. Hence
\[
 \int_0^\infty\|(F_1,F_2,F_3)(t)\|_1\dd t
       \le C\mathcal J_2\mathcal J_{2,1},\qquad
 \int_0^\infty\|(F_1,F_2,F_3)(t)\|_{H^1}\dd t\le C\mathcal J_2^2.
\]

\emph{Step 2. Displacement bound for the shear velocity.}
Write $u_3=v_3+w_3$, where $v=\mathbb Pu$ and $w=\mathbb Qu$.
Let $K_v(t,x_h,z)$ be the third row of the convolution kernel
of $S_h(t)\mathbb P$. With the horizontal heat kernel $K_h$
and the normalized Poisson kernel
$P_r(x_h)=cr(r^2+|x_h|^2)^{-3/2}$, $r>0$,
vertical Fourier inversion gives
\[
 \begin{aligned}
 (K_v)_j(t,x_h,z)
       &=-\tfrac12\operatorname{sgn}(z)
                   \partial_j(K_h(t)*_hP_{|z|})(x_h),\quad j=1,2,\\
 (K_v)_3(t,x_h,z)
       &=-\tfrac12\left.\partial_r(K_h(t)*_hP_r)(x_h)
                                                    \right|_{r=|z|}.
 \end{aligned}
\]
The evaluation at $r=0$ is understood as a right limit.
Consequently,
$k_v(t,x_h):=\sup_z|K_v(t,x_h,z)|$
satisfies $k_v(t,x_h)\le C(\nu_ht+|x_h|^2)^{-3/2}$.
For every measurable $E\subset\R^2$ with $0<|E|<\infty$ and
every measurable horizontal translation path $\gamma$,
\[
 \int_0^\infty\int_E k_v(t,x_h+\gamma(t))\dd x_h\dd t
 \le C\int_0^\infty
          \min\{t^{-1/2},|E|t^{-3/2}\}\dd t
 \le C|E|^{1/2},
\]
where the last integral is split at $t=|E|$.
For a measurable function $d:\R^2\to\R$, the equivalent characterization
\[
 \|d\|_{L^{2,\infty}_h}
       \simeq\sup_{0<|E|<\infty}|E|^{-1/2}\int_E|d(x_h)|\dd x_h
\]
therefore gives a bound uniform in $\gamma$.
Duhamel's formula for $v$, applied at each source time, and
Tonelli's theorem prove \eqref{eq:r8-drift-displacement} with
$a=v_3$ and
\[
 \delta\le C\left(\|u_0\|_1
                    +\int_0^\infty\|F_2(t)\|_1\dd t\right)
       \le C(\mathcal J_2+\mathcal J_2\mathcal J_{2,1})
       \le C\mathcal J_2.
\]
Here the data bound is \eqref{eq:new-data-sizes} and the source
integral was established in Step~1.

\emph{Step 3. Displacement bound for the potential velocity.}
Let $K_{\rm ac}$ denote the leading acoustic part of the
low-frequency velocity Green matrix in the proof of
Lemma~\ref{lem:r8-acoustic}. The radial pointwise bounds in
Step~2 of the proof of Lemma~\ref{lem:r8-acoustic}, including the Newtonian-potential
contribution, give, for $k_{\rm ac}(t,x_h)=\sup_z|K_{\rm ac}(t,x_h,z)|$,
\[
 \|k_{\rm ac}(t)\|_{L^1_h}\le C,\qquad
 \|k_{\rm ac}(t)\|_{L^\infty_h}\le Ct^{-2},\qquad t\ge1.
\]
Horizontal heat convolution and the fixed smooth frequency
cutoff preserve these bounds. Thus
\[
 \int_1^\infty\int_E
     k_{\rm ac}(t,x_h+\gamma(t))\dd x_h\dd t
       \le C\int_1^\infty\min\{1,|E|t^{-2}\}\dd t
       \le C|E|^{1/2}.
\]
For the remainder kernel $K_{\rm rem}$, the multiplier bound in
Step~1 of the proof of Lemma~\ref{lem:r8-acoustic} gives
$\|K_{\rm rem}(t)\|_2\le Ct^{-5/4}$ and
$\|\partial_3K_{\rm rem}(t)\|_2\le Ct^{-7/4}$ for $t\ge1$.
The vertical Sobolev inequality \eqref{eq:mixedv} yields
\[
 \|K_{\rm rem}(t)\|_{L^2_hL^\infty_3}
       \le C\|K_{\rm rem}(t)\|_2^{1/2}
                   \|\partial_3K_{\rm rem}(t)\|_2^{1/2}
       \le Ct^{-3/2}.
\]
On $0\le t\le1$, compact frequency support and
\eqref{eq:blockmultiplier} give a uniform
$L^2_hL^\infty_3$ bound for the entire low-frequency kernel.
For the high-frequency semigroup, \eqref{eq:blockmultiplier} gives
$\|P_HS_c(t)f\|_{H^1}\le Ce^{-c_*t}\|f\|_{H^1}$.
The embedding $H^1\hookrightarrow L^2_hL^\infty_3$ follows
from \eqref{eq:mixedv}.
Horizontal translations preserve $L^2_hL^\infty_3$, so
Minkowski's inequality bounds the displacement of the short-time
low-frequency term and the spectral remainder in $L^2_h$ by
$C\|f\|_1$, and that of the high-frequency evolution by
$C\|f\|_{H^1}$.
Since $L^2_h\hookrightarrow L^{2,\infty}_h$, these bounds and
the acoustic estimate, applied in Duhamel's formula, give
\eqref{eq:r8-drift-displacement} with
$a=w_3$ and
\[
 \begin{aligned}
 \delta&\le C\left[\|U_0\|_1+\|U_0\|_{H^1}
       +\int_0^\infty\big(\|(F_1,F_2,F_3)(t)\|_1
                 +\|(F_1,F_2,F_3)(t)\|_{H^1}\big)\dd t\right]\\
       &\le C(\mathcal J_2+\mathcal J_2\mathcal J_{2,1}
                                      +\mathcal J_2^2)
        \le C\mathcal J_2.
 \end{aligned}
\]
Since $\mathcal J_2\le\mathcal J_{2,1}$ and
$\mathcal J_{2,1}$ is sufficiently small, the last inequalities
in Steps~2--3 hold with a fixed constant.
The bounds for $v_3$ and $w_3$ therefore prove
\eqref{eq:r8-drift-displacement} for $a=u_3$ with
$\delta=C\mathcal J_2$.

\emph{Step 4. Verification of the coefficient bounds.}
Equation~\eqref{eq:new-H2-main-Lp} with $p=\infty$ gives
$\|u_3(t)\|_\infty\le C\mathcal J_2(1+t)^{-3/2}$.
The third line of \eqref{eq:r8-windows} gives
$\|\nabla u_3\|_{L^2(b,b+1;L^\infty)}
\le C\mathcal J_2(1+b)^{-11/8}$, and
\eqref{eq:r8-lips-tail} at $t=0$ gives
$\int_0^\infty\|\nabla u_3\|_\infty\dd t\le C\mathcal J_2$.
In particular, $u_3\in L^1_{\rm loc}([0,\infty);W^{1,\infty})$.

Apply \eqref{eq:mixedv} and Minkowski's inequality with the
$L^2$ decay bounds \eqref{eq:new-H2-main-rates},
\eqref{eq:H2-main-vertical-weak}, \eqref{eq:new-S2-structure},
and \eqref{eq:new-second-preliminary}. This gives
\[
 \begin{aligned}
 \|\nabla u_3(t)\|_{L^\infty_3L^2_h}
 &\le C\|\nabla u_3(t)\|_2^{1/2}
                \|\partial_3\nabla u_3(t)\|_2^{1/2}
       \le C\mathcal J_2(1+t)^{-7/8},\\
 \|\nabla_hu_3(t)\|_{L^2_hL^\infty_3}
 &\le C\|\nabla_hu_3(t)\|_2^{1/2}
                \|\partial_3\nabla_hu_3(t)\|_2^{1/2}
       \le C\mathcal J_2(1+t)^{-11/8}.
 \end{aligned}
\]
The second bound is time integrable. Finally,
\eqref{eq:verticalcomp} and \eqref{eq:H2-dissipation-definition}
control $\|\nabla\partial_3\nabla u_3\|_2$ by
$C\cD_2(U)^{1/2}$. Sobolev interpolation and H\"older in time,
using $\|\partial_3\nabla u_3\|_2\le C\mathcal J_2(1+t)^{-3/4}$,
give
\[
 \int_0^\infty(1+t)^{-2/3}
                     \|\partial_3\nabla u_3\|_3\dd t
 \le C\mathcal J_2^{1/2}
       \left(\int_0^\infty\cD_2(U)\dd t\right)^{1/4}
       \left(\int_0^\infty(1+t)^{-25/18}\dd t\right)^{3/4}
       \le C\mathcal J_2,
\]
where \eqref{eq:new-base-energy} bounds the dissipation integral.
Thus all six conditions in \eqref{eq:r8-drift-assumptions}
hold with $\delta=C\mathcal J_2$.
\end{proof}

\subsubsection{Source estimates and first vertical derivative decay}

\begin{lemma}[Decay of the first vertical derivative]
\label{lem:r13-first-vertical-decay}
If $\mathcal J_{2,1}$ is sufficiently small, the solution of
Theorem~\ref{thm:global} satisfies, for $j=0,1$ and $t\ge0$,
\[
 \|\nabla_h^j\partial_3u(t)\|_2
       \le C\mathcal J_{2,1}(1+t)^{-(1+j)/2}.
\]
\end{lemma}
\begin{proof}
Differentiating $v_t-\nu_h\Delta_hv=\mathbb PF_2$ and using
\eqref{eq:new-projected-source}, we obtain
\[
 (\partial_t-\nu_h\Delta_h)\partial_3v
       +\mathbb P(u_3\partial_3^2v)
       =\mathbb P(S_{\mathrm{reg}}+\mathcal V),
\]
where the convection and pressure source $S_{\mathrm{reg}}$ and the viscous
source $\mathcal V$ are defined below.
We verify the source condition \eqref{eq:r8-source-space} of
Lemma~\ref{lem:r8-forced-endpoint}; its coefficient conditions
were proved in Lemma~\ref{lem:r8-real-drift}.

\emph{Step 1. Convection and pressure.}
The convection and pressure terms give
\begin{equation}\label{eq:r8-regular-source}
 S_{\mathrm{reg}}=-\partial_3(u_h\cdot\nabla_hv)-\partial_3u_3\,\partial_3v
 -\partial_3(u\cdot\nabla w)
       +c^2\partial_3[(q+\tau)e^\tau\nabla\tau].
\end{equation}
By \eqref{eq:new-H2-main-rates}, \eqref{eq:H2-main-vertical-weak},
\eqref{eq:r9-lips}, and \eqref{eq:r9-compressible-lips},
\[
 \|S_{\mathrm{reg}}(t)\|_{L^2_3L^1_h}
       \le C\mathcal J_2^2(1+t)^{-9/8},
 \qquad \|S_{\mathrm{reg}}(t)\|_2
       \le C\mathcal J_2^2(1+t)^{-13/8}.
\]
Integration gives the mixed-norm and $L^2$ tail bounds in
\eqref{eq:r8-source-space}. The vertically differentiated convection
terms are bounded by sums of
\[
 C\mathcal J_2(\|\nabla_hv\|_\infty+\|\partial_3u_3\|_\infty
                                      +\|\nabla w\|_\infty),\qquad
 C(\|u\|_\infty+\|\partial_3u\|_3+\|\partial_3v\|_3)
                                                    \cD_2(U)^{1/2}.
\]
The first expression is integrable by \eqref{eq:r9-lips} and
\eqref{eq:r9-compressible-lips} for $t\ge1$.
On $[0,1]$, use \eqref{eq:r8-lips-tail},
$\|\nabla w\|_\infty\le C\cD_2(U)^{1/2}$, and
\eqref{eq:new-base-energy}. The directional Sobolev inequality
\eqref{eq:coordtriple} gives
\[
 \|\partial_3u\|_3
       \le C\|\partial_3u\|_2^{1/2}
             \|\nabla_h\partial_3u\|_2^{1/3}
             \|\partial_3^2u\|_2^{1/6},
 \qquad \|\partial_3v\|_3
       \le C\|\partial_3v\|_2^{1/2}
             \|\nabla_h\partial_3v\|_2^{1/3}
             \|\partial_3^2v\|_2^{1/6}.
\]
Thus \eqref{eq:new-intermediate} and the $H^2$ bound in
\eqref{eq:new-base-energy} imply
$\|\partial_3u\|_{L^2_tL^3}+\|\partial_3v\|_{L^2_tL^3}
\le C\mathcal J_{2,1}$.
Cauchy--Schwarz in time and \eqref{eq:new-base-energy} control
the second expression. The pressure terms use
$\nabla\partial_3^2\tau$ in \eqref{eq:H2-dissipation-definition}.
Consequently,
\[
 \int_0^\infty\|\partial_3S_{\mathrm{reg}}\|_2\dd t
       \le C\mathcal J_2\mathcal J_{2,1}.
\]

\emph{Step 2. The viscous source and horizontal flux.}
The viscous source is
\begin{equation}\label{eq:r8-viscosity-representative}
 \mathcal V={}\nu_h\partial_3[(e^{-q}-1)\Delta_hu]
 +\nu_c\big[\partial_3(e^{-q}-1)\nabla\Div u
              -\nabla(e^{-q}-1)\partial_3\Div u\big].
\end{equation}
Since
\[
 \partial_3[(e^{-q}-1)\mathcal Lu]-\mathcal V
       =\nu_c\nabla[(e^{-q}-1)\partial_3\Div u],
\]
we have $\mathbb P\partial_3[(e^{-q}-1)\mathcal Lu]=\mathbb P\mathcal V$.
We first estimate $\mathcal V$ and $\partial_3\mathcal V$ in a horizontal negative
Sobolev norm, which permits a low-frequency source term
and a high-frequency horizontal flux without requiring an
additional derivative of the data. We use the norm
\[
 \|h\|_{H^{-1}_hL^2_3}^2
       =(2\pi)^{-3}\int(1+|\xi_h|^2)^{-1}|\widehat h(\xi)|^2\dd\xi.
\]
The product decomposition \eqref{eq:r8-product-flux} implies
$\|AB\|_{H^{-1}_hL^2_3}\le C\|A\|_{H^1}\|B\|_2$.
Expanding \eqref{eq:r8-viscosity-representative} and using
\eqref{eq:new-compressible-second} gives
\begin{align}\label{eq:r8-viscosity-bounds}
 \|\mathcal V\|_{L^2_3L^1_h}&\le C\mathcal J_2^2(1+t)^{-3}
                  +C\mathcal J_2(1+t)^{-1}\cD_2(U)^{1/2},\notag\\
 \|\mathcal V\|_2&\le C\mathcal J_2(1+t)^{-3/2}\cD_2(U)^{1/2}
               +C\mathcal J_2^{3/2}(1+t)^{-5/2}\cD_2(U)^{1/4},\notag\\
 \|\mathcal V\|_{H^{-1}_hL^2_3}&\le C\mathcal J_2^2(1+t)^{-2},\\
 \|\partial_3\mathcal V\|_{H^{-1}_hL^2_3}
       &\le C\mathcal J_2(1+t)^{-5/4}\cD_2(U)^{1/2}.\notag
\end{align}
For the last inequality, the horizontal terms are
\[
 (e^{-q}-1)\Delta_h\partial_3^2u
 +2\partial_3(e^{-q}-1)\Delta_h\partial_3u
 +\partial_3^2(e^{-q}-1)\Delta_hu.
\]
Write
\[
 (e^{-q}-1)\Delta_h\partial_3^2u
  ={}\Div_h[(e^{-q}-1)\nabla_h\partial_3^2u]
 -\nabla_h(e^{-q}-1)\cdot\nabla_h\partial_3^2u.
\]
The bulk-viscosity terms are
\[
 \partial_3^2(e^{-q}-1)\nabla\Div u
   +\partial_3(e^{-q}-1)\nabla\partial_3\Div u
 -\nabla\partial_3(e^{-q}-1)\partial_3\Div u
        -\nabla(e^{-q}-1)\partial_3^2\Div u.
\]
Apply the product bound from \eqref{eq:r8-product-flux}, with the once
differentiated coefficient in $H^1$, or with the velocity factor in
$H^1$ when the coefficient has two derivatives.
The remaining factors are controlled by
\eqref{eq:H2-dissipation-definition} and
\eqref{eq:new-compressible-second}.

For the fixed horizontal cutoff $P_L^h$, set
\[
 S_\nu=P_L^h\mathcal V,\qquad
 B_\nu=\nabla_h\Delta_h^{-1}(I-P_L^h)\mathcal V.
\]
Then $\mathcal V=S_\nu+\Div_hB_\nu$ is the source--flux decomposition
required in \eqref{eq:r8-source-space}.
Both operators are bounded from $H^{-1}_hL^2_3$ to $L^2$
and on $L^2_3L^1_h$.
Indeed their horizontal kernels belong to $L^1$; the high-frequency
inverse-derivative kernel has summable dyadic bounds $C2^{-j}$.
The first and last lines of \eqref{eq:r8-viscosity-bounds},
Cauchy--Schwarz, and \eqref{eq:new-base-energy} give
\[
 \int_0^\infty
 \bigl(\|(S_\nu,B_\nu)\|_{L^2_3L^1_h}
       +\|\partial_3(S_\nu,B_\nu)\|_2\bigr)\dd t
       +\|\partial_3B_\nu\|_{L^2_tL^2}
       \le C\mathcal J_2^2.
\]
The second and third lines of \eqref{eq:r8-viscosity-bounds} yield
\[
 \begin{aligned}
 \|S_\nu(t)\|_2+\|B_\nu(t)\|_2
       &\le C\mathcal J_2^2(1+t)^{-2},\\
 \|S_\nu\|_{L^2(b,\infty;L^2)}
       +\|\Div_hB_\nu\|_{L^2(b,\infty;L^2)}
       &\le C\mathcal J_2^2(1+b)^{-3/2}.
 \end{aligned}
\]
Integrating the pointwise bound for $S_\nu$ also gives
$\int_b^\infty\|S_\nu(t)\|_2\dd t
\le C\mathcal J_2^2(1+b)^{-1}$.
Thus $(S_\nu,B_\nu)$ satisfies \eqref{eq:r8-source-space}
with $f_0=0$ and $N=C\mathcal J_2^2$.

\emph{Step 3. Application of the vertical transport estimate.}
Apply Lemma~\ref{lem:r8-forced-endpoint} with
\[
 f_0=\partial_3u_0,\qquad a=u_3,\qquad
 S=S_{\mathrm{reg}}+S_\nu,\qquad B=B_\nu.
\]
Lemma~\ref{lem:r8-real-drift} gives $\delta\le C\mathcal J_2$,
and Steps~1--2 verify \eqref{eq:r8-source-space} with
\[
 N\le C\big(\mathcal J_{2,1}
             +\mathcal J_2\mathcal J_{2,1}+\mathcal J_2^2\big)
       \le C\mathcal J_{2,1}.
\]
The smallness of $\mathcal J_{2,1}$ ensures the required smallness
of $\delta$. By uniqueness, the solution in the lemma is
$\partial_3v$. Adding the potential-velocity bounds in
\eqref{eq:new-H2-main-rates} to \eqref{eq:r8-forced-endpoint} yields
\begin{equation}\label{eq:r8-logfree}
 \|\partial_3u(t)\|_2+(1+t)^{1/2}\|\nabla_h\partial_3u(t)\|_2
       \le C\mathcal J_{2,1}(1+t)^{-1/2}.
\end{equation}
\end{proof}

\begin{corollary}[Nonlinear correction and dissipation tails]
If $\mathcal J_{2,1}$ is sufficiently small, then for $j=0,1$,
\begin{equation}\label{eq:r8-amplitude}
 \|\nabla_h^j\partial_3r_v(t)\|_2
 \le C\left[\mathcal J_2\mathcal J_{2,1}+
       \mathcal J_2^2\sqrt{\log(e+\mathcal J_{2,1}/\mathcal J_2^2)}\right]
                    (1+t)^{-(1+j)/2}.
\end{equation}
For $U_0=0$, the right-hand side is defined to be zero
by continuous extension.
Also, for $t\ge1$,
\begin{equation}\label{eq:r8-first-diss-tails}
 \int_t^\infty\|\nabla_h\partial_3u\|_2^2\dd s
       \le C\mathcal J_{2,1}^2(1+t)^{-1},\qquad
 \int_t^\infty\|\nabla_h^2\partial_3u\|_2^2\dd s
       \le C\mathcal J_{2,1}^2(1+t)^{-2}.
\end{equation}
\end{corollary}
\begin{proof}
\emph{Step 1. The nonlinear correction.}
Let $F$ solve the homogeneous unprojected equation, componentwise,
\[
 F_t-\nu_h\Delta_hF+u_3\partial_3F=0,
 \qquad F(0)=\partial_3u_0.
\]
The vertical energy estimate in Step~3 of the proof of
Lemma~\ref{lem:r8-forced-endpoint}, together with
$\|\partial_3^2u_0\|_2\le\mathcal J_2$, gives
\[
 \sup_{t\ge0}\|\partial_3F(t)\|_2\le C\mathcal J_2.
\]
By Lemma~\ref{lem:r8-real-drift}, the characteristic maps used
in Steps~1--2 of the proof of Lemma~\ref{lem:r8-forced-endpoint}
have displacement satisfying
$\|d_\gamma\|_{L^{2,\infty}_h}\le C\mathcal J_2$.
Their inverse Jacobians satisfy
\[
 \|\partial_z\psi-1\|_\infty
       \le e^{C\mathcal J_2}-1\le C\mathcal J_2.
\]
Apply \eqref{eq:r8-graph-difference}
to $\partial_3u_0$, using
$\|\partial_3u_0\|_{L^2_3L^1_h}\le\mathcal J_{2,1}$ and
$\|\partial_3^2u_0\|_2\le\mathcal J_2$.
Subtracting the corresponding pure-heat splitting approximation
and passing to the limit as in Step~2 of the proof of Lemma~\ref{lem:r8-forced-endpoint} yields
\[
 \sup_{t\ge0}\operatorname*{ess\,sup}_{\xi_h\in\R^2}
   \|\mathcal F_h[F(t)-S_h(t)\partial_3u_0](\xi_h,\cdot)\|_2
 \le C\left[\mathcal J_2\mathcal J_{2,1}
       +\mathcal J_2^2
           \sqrt{\log(e+\mathcal J_{2,1}/\mathcal J_2^2)}\right].
\]
The difference $F-S_h(t)\partial_3u_0$ has zero initial datum
and solves a horizontal heat equation with source
$-u_3\partial_3F$, for which
\[
 \|u_3\partial_3F(t)\|_2
       \le C\mathcal J_2^2(1+t)^{-3/2}.
\]
The horizontal Nash inequality with the Fourier bound above,
followed by the energy comparison in Step~3 of the proof of
Lemma~\ref{lem:r8-forced-endpoint}, gives the $(1+t)^{-1/2}$
decay of $F-S_h(t)\partial_3u_0$ with the bracketed coefficient
in \eqref{eq:r8-amplitude}. Indeed, a sufficiently large constant
times its square divided by $1+t$ is an energy supersolution,
since that coefficient is at least $\mathcal J_2^2$.
Duhamel's formula for the horizontal heat equation on $[t/2,t]$,
together with \eqref{eq:r8-horizontal-maximal}, gives the
$(1+t)^{-1}$ gradient decay with the same coefficient.
For $0\le t\le1$, both estimates follow from the zero-initial-data
heat formula; the zero-data case follows by continuous extension.

Let $\widetilde F$ be the unprojected lift constructed in
Step~4 of the proof of Lemma~\ref{lem:r8-forced-endpoint},
with initial datum $\partial_3u_0$ and the source pair
$(S_{\mathrm{reg}}+S_\nu,B_\nu)$ from the proof of
Lemma~\ref{lem:r13-first-vertical-decay}.
Thus $\mathbb P\widetilde F=\partial_3v$ and
$\|\widetilde F\|_{\mathcal X}\le C\mathcal J_{2,1}$,
with $\mathcal X$ as in Lemma~\ref{lem:r8-forced-endpoint}.
The difference $\widetilde F-F$ has zero initial datum and satisfies
\[
 (\partial_t-\nu_h\Delta_h)(\widetilde F-F)
                  +u_3\partial_3(\widetilde F-F)
 =S_{\mathrm{reg}}+S_\nu+\Div_hB_\nu
          -(\nabla u_3)\partial_3\Delta^{-1}\Div\widetilde F.
\]
Steps~1--2 of the proof of
Lemma~\ref{lem:r13-first-vertical-decay} bound the physical
source pair in \eqref{eq:r8-source-space}, with zero initial
datum, by $C\mathcal J_2\mathcal J_{2,1}$.
The projection correction has the same source bound by
Step~4 of the proof of Lemma~\ref{lem:r8-forced-endpoint},
since $\delta=C\mathcal J_2$ and
$\|\widetilde F\|_{\mathcal X}\le C\mathcal J_{2,1}$.
The unprojected estimate in Steps~1--3 of the proof of
Lemma~\ref{lem:r8-forced-endpoint} gives
\[
 \|\widetilde F-F\|_{\mathcal X}
       \le C\mathcal J_2\mathcal J_{2,1}.
\]
Finally,
\[
 \partial_3r_v
       =\mathbb P\big[F-S_h(t)\partial_3u_0
                                      +(\widetilde F-F)\big].
\]
The $L^2$ boundedness of $\mathbb P$ proves
\eqref{eq:r8-amplitude}.

\emph{Step 2. The dissipation tails.}
Squaring the horizontal-gradient bound in
\eqref{eq:r8-logfree} and integrating proves the first
inequality in \eqref{eq:r8-first-diss-tails}.
For the second, test
\[
 (\partial_t-\nu_h\Delta_h)\partial_3v
       =\mathbb P\partial_3F_2
\]
with $-\Delta_h\partial_3v$.
Integration and Young's inequality give
\[
 \int_t^\infty\|\nabla_h^2\partial_3v(s)\|_2^2\dd s
 \le C\|\nabla_h\partial_3v(t)\|_2^2
           +C\int_t^\infty\|F_2(s)\|_{H^1}^2\dd s
 \le C\mathcal J_{2,1}^2(1+t)^{-2}.
\]
The last inequality uses \eqref{eq:r8-logfree} and the sum
over unit intervals of \eqref{eq:r8-F2-H1-window}.
Adding the potential-velocity tail obtained by squaring
and integrating \eqref{eq:r9-w-third} proves the second
inequality in \eqref{eq:r8-first-diss-tails}.
\end{proof}

\subsection{Enhanced mixed third-order and thermal decay}
\label{sec:r9-vertical-third}

Multiply the momentum equation in \eqref{eq:perturbed} by $e^q$
and apply the Leray projection to eliminate the pressure. This gives
\begin{equation}\label{eq:r9-weighted-shear}
 v_t-\nu_h\Delta_hv=-\mathbb P\mathcal R_{\rm sh},\qquad
 \mathcal R_{\rm sh}=(e^q-1)u_t+e^q(u\cdot\nabla u).
\end{equation}

\begin{lemma}[Enhanced third-order decay]
\label{lem:r13-enhanced-third}
If $\mathcal J_{2,1}$ is sufficiently small, the solution of
Theorem~\ref{thm:global} satisfies, for $t\ge1$,
\[
 \|\nabla_h^2\partial_3u(t)\|_2
       \le C\mathcal J_{2,1}(1+t)^{-3/2},\qquad
 \|\nabla^3\tau(t)\|_2
       \le C\mathcal J_{2,1}(1+t)^{-9/4}.
\]
Moreover,
\begin{equation}\label{eq:r9-strong-time-third}
 \|\partial_3u_t(t)\|_2+
 \|\nabla_h\partial_3u_t\|_{L^2(t,t+1;L^2)}
 +\|\Div\partial_3u_t\|_{L^2(t,t+1;L^2)}
       \le C\mathcal J_{2,1}(1+t)^{-3/2}.
\end{equation}
\end{lemma}
\begin{proof}
\emph{Step 1. Mixed third-order velocity decay.}
Equation~\eqref{eq:r9-weighted-shear} gives $\mathbb PF_2=-\mathbb P\mathcal R_{\rm sh}$.
Lemma~\ref{lem:r9-source-tail} in Appendix~\ref{app:r10-sources}
gives \eqref{eq:r9-source-tail} using only
Theorem~\ref{thm:decay-H2}.
For $t\ge2$, apply the Duhamel formula on $[t/2,t]$.
The horizontal heat semigroup estimate
\eqref{eq:r8-horizontal-maximal} with $m=0$,
\eqref{eq:r8-logfree}, and \eqref{eq:r9-source-tail} yield
\[
 \|\nabla_h^2\partial_3v(t)\|_2
 \le Ct^{-1/2}\|\nabla_h\partial_3v(t/2)\|_2
       +C\|\nabla_h\partial_3\mathcal R_{\rm sh}\|_{L^2(t/2,t;L^2)}
 \le C\mathcal J_{2,1}(1+t)^{-3/2}.
\]
Adding the potential-velocity bound \eqref{eq:r9-third-preliminary},
which decays as $(1+t)^{-7/4}$, proves the mixed third-order
estimate for $t\ge2$. Estimate~\eqref{eq:r9-mixed-third}
covers $1\le t\le2$.

\emph{Step 2. Third-order temperature decay.}
By \eqref{eq:r9-lips}, \eqref{eq:r8-logfree}, and
\eqref{eq:new-second-preliminary},
\[
 \begin{aligned}
 \|\partial_3|\nabla_hu|^2\|_2
 &\le2\|\nabla_hu\|_\infty\|\nabla_h\partial_3u\|_2
       \le C\mathcal J_2\mathcal J_{2,1}(1+t)^{-19/8},\\
 \|\nabla_h|\nabla_hu|^2\|_2
 &\le2\|\nabla_hu\|_\infty\|\nabla_h^2u\|_2
       \le C\mathcal J_2^2(1+t)^{-23/8}.
 \end{aligned}
\]
Using the Fourier identity
$\|\nabla^3(-\Delta)^{-1}h\|_2=\|\nabla h\|_2$
in the thermal comparison \eqref{eq:r9-thermal-profile} gives
\[
 \|\nabla^3\tau(t)\|_2
 \le C\mathcal J_2(1+t)^{-9/4}
       +C\mathcal J_2\mathcal J_{2,1}(1+t)^{-19/8}
 \le C\mathcal J_{2,1}(1+t)^{-9/4}.
\]

\emph{Step 3. Energy estimate for $\partial_3u_t$.}
Rearranging the vertically differentiated momentum equation
\eqref{eq:r9-Z-recovery}, and using Step~1,
\eqref{eq:new-compressible-second}, \eqref{eq:r9-third-preliminary},
and the product estimates in Step~4 of the proof of
Proposition~\ref{prop:r9-mixed-third}, we obtain
\[
 \int_b^{b+1}\|\partial_3u_t(s)\|_2^2\dd s
       \le C\mathcal J_{2,1}^2(1+b)^{-3},\qquad b\ge2.
\]
In the energy inequality \eqref{eq:r9-Y-energy} for
$\|\partial_3u_t\|_2^2$, the remainder has integral at most
$C\mathcal J_2^2(1+b)^{-7/2}$ over $[b,b+1]$.
Applying Lemma~\ref{lem:r9-uniform-gronwall} with $p=3$ on
$[2,\infty)$ proves \eqref{eq:r9-strong-time-third} for $t\ge3$.
Estimate~\eqref{eq:r9-Y-point} covers $1\le t\le3$.
\end{proof}

\subsection{Enhanced decay and estimates for the physical variables}
\label{sec:r10-strong-sources}

\begin{proof}[Completion of Theorem~\ref{thm:H2-intermediate-main}]
Proposition~\ref{prop:new-intermediate} gives the $L^p$ bounds in
\eqref{eq:H2-intermediate-main}.
Lemma~\ref{lem:r13-first-vertical-decay} and
\eqref{eq:new-H2-main-rates} give
\[
 \|\partial_3u\|_2+\|U\|_{H^1}
       \le C\mathcal J_{2,1}(1+t)^{-1/2}.
\]
Since $\Div u=\Div w$, we have
\[
 \|\partial_3^2u_3\|_2
 \le\|\partial_3\Div w\|_2+\|\Div_h\partial_3u_h\|_2
 \le C\mathcal J_{2,1}(1+t)^{-1},
\]
by the $k=2$ case of \eqref{eq:new-H2-main-rates} and
\eqref{eq:r8-logfree}. These estimates prove
\eqref{eq:H2-intermediate-main}.
Lemma~\ref{lem:r13-enhanced-third} gives
\eqref{eq:r9-strong-third}, completing the theorem.
\end{proof}
\begin{proof}[Proof of Corollary~\ref{cor:physicaldecay}]
Recall that $\rho-\bar\rho=\bar\rho(e^q-1)$ and
$\theta-\bar\theta=\bar\theta(e^\tau-1)$.
For $z=q,\tau$, the small $H^2$ bound and
\begin{equation}
 e^z-1=z\int_0^1e^{sz}\dd s,\quad
 \partial_i(e^z-1)=e^z\partial_i z,\quad
 \partial_i\partial_j(e^z-1)
       =e^z\partial_i\partial_jz+e^z\partial_i z\partial_jz
\end{equation}
transfer the $L^p$ and derivative bounds through order two in
Theorems~\ref{thm:decay-H2} and \ref{thm:H2-intermediate-main}
to the physical variables, since
$\|\nabla z\|_4^2\le C\|z\|_\infty\|\nabla^2z\|_2$.
The inverse logarithmic composition gives the corresponding
equivalence of the small $H^2$ norms.
\end{proof}

\section{Gaussian profiles}\label{sec:r10-profiles}

For sufficiently small $\mathcal J_2$, we construct the effective
datum and obtain horizontal Gaussian profiles. Small
$\mathcal J_{2,1}$ also yields the vertical derivative profiles.

\subsection{Effective datum and comparison with the horizontal heat flow}

We use the unprojected momentum source $F_2$ defined in \eqref{eq:F2}.

\begin{lemma}[Momentum-source decay]\label{lem:r10-momentum-source}
For sufficiently small $\mathcal J_2$ and all $t\ge0$,
\begin{equation}\label{eq:r10-momentum-source}
 \|F_2(t)\|_{L^2_3L^1_h}\le C\mathcal J_2^2(1+t)^{-5/4},
 \qquad
 \|F_2(t)\|_2\le C\mathcal J_2^2(1+t)^{-7/4}.
\end{equation}
\end{lemma}
\begin{proof}
By \eqref{eq:H2-main-vertical-weak}, \eqref{eq:new-H2-main-Lp}, and
\eqref{eq:H2-main-mixed-and-full}, the vertical convection term satisfies
\[
 \begin{aligned}
 \|u_3\partial_3u\|_2
 &\le\|u_3\|_\infty\|\partial_3u\|_2
       \le C\mathcal J_2^2(1+t)^{-7/4},\\
 \|u_3\partial_3u\|_{L^2_3L^1_h}
 &\le\|u_3\|_{L^\infty_3L^2_h}\|\partial_3u\|_2
       \le C\mathcal J_2^2(1+t)^{-5/4}.
 \end{aligned}
\]
For the horizontal convection term, Sobolev's inequality in $x_3$ gives
\[
 \|u\|_{L^\infty_3L^2_h}
 \le C\|u\|_2^{1/2}\|\partial_3u\|_2^{1/2}
 \le C\mathcal J_2(1+t)^{-3/8}.
\]
Together with \eqref{eq:new-H2-main-rates} and \eqref{eq:new-H2-main-Lp},
this yields decay $(1+t)^{-11/8}$ in $L^2_3L^1_h$ and
$(1+t)^{-15/8}$ in $L^2$.
The pressure and variable-viscosity terms decay faster in both norms
by \eqref{eq:new-H2-main-rates}, \eqref{eq:new-H2-main-Lp}, and
\eqref{eq:new-compressible-second}. Summing the terms of \eqref{eq:F2}
proves \eqref{eq:r10-momentum-source}.
\end{proof}

\begin{proposition}[Effective datum and heat-flow comparison]\label{prop:r8-H1-scattering}
Assume $\mathcal J_2$ is sufficiently small. The effective velocity datum defined in
\eqref{eq:r8-physical-effective} belongs to
$H^1\cap L^2_3L^1_h$ and satisfies
\begin{equation}\label{eq:new-effective-datum}
 \|u_\infty-u_0\|_{H^1}+\|u_\infty-u_0\|_{L^2_3L^1_h}
       \le C\mathcal J_2^2.
\end{equation}
For $t\ge1$,
\begin{align}\label{eq:r8-H1-scattering}
 &\|v-S_h(t)\mathbb Pu_\infty\|_{H^1}
       +t^{1/2}\|\nabla_h(v-S_h(t)\mathbb Pu_\infty)\|_{H^1}\notag\\
 &\qquad\le Ct^{-1/2}
       [\mathcal J_2\mathfrak E_2(\sqrt t)^{1/2}
                         +\mathcal J_2^2t^{-1/8}]
       =o(t^{-1/2}).
\end{align}
In particular, the left side is at most $C\mathcal J_2^2t^{-1/2}$.
The same comparison for the full tuple $U$ and
$(0,S_h(t)\mathbb Pu_\infty,0)$ has the additional bound
$C\mathcal J_2 t^{-3/4}$.
For $t\ge1$,
\begin{align}
 \label{eq:new-effective-errors}
 \|v-S_h(t)\mathbb Pu_\infty\|_2&\le C\mathcal J_2^2t^{-3/4},\notag
 \qquad \|\nabla_h(v-S_h(t)\mathbb Pu_\infty)\|_2\le C\mathcal J_2^2t^{-5/4},\notag\\*
 \|\nabla_h^2(v-S_h(t)\mathbb Pu_\infty)\|_2
       &\le C\mathcal J_2^2t^{-7/4}\sqrt{\log(e+t)}.
\end{align}
\end{proposition}
The $H^1$ comparison \eqref{eq:r8-H1-scattering} will be used in
Proposition~\ref{prop:r8-vertical-profile}; the directional $L^2$
bounds \eqref{eq:new-effective-errors} give the horizontal profiles
in Proposition~\ref{prop:gaussian}.
\begin{proof}
\emph{Step 1. Source integrability and the effective datum.}
In the estimate of $u_3\partial_3^2u$, retain the remaining
$H^2$ energy $\mathfrak E_2(b)$ of Proposition~\ref{prop:r8-energy-tail}.
The proof of \eqref{eq:r8-F2-H1-window}, together with
\eqref{eq:r8-visc-H1}, then gives
\begin{equation}\label{eq:r8-refined-source}
 \|F_2\|_{L^2(b,b+1;H^1)}
 \le C\mathcal J_2\mathfrak E_2(b)^{1/2}(1+b)^{-3/2}
       +C\mathcal J_2^2(1+b)^{-13/8}.
\end{equation}
The dissipative viscosity term uses
$\int_b^\infty\cD_2(U)\dd s\le C\mathfrak E_2(b)$ from
Proposition~\ref{prop:r8-energy-tail}.
Summation over unit intervals gives
\begin{equation}\label{eq:r13-effective-source-tails}
 \begin{aligned}
 \|F_2\|_{L^1(b,\infty;H^1)}
 &\le C\mathcal J_2\mathfrak E_2(b)^{1/2}(1+b)^{-1/2}
       +C\mathcal J_2^2(1+b)^{-5/8},\\
 \|F_2\|_{L^2(b,\infty;H^1)}
 &\le C\mathcal J_2\mathfrak E_2(b)^{1/2}(1+b)^{-1}
       +C\mathcal J_2^2(1+b)^{-9/8}.
 \end{aligned}
\end{equation}
Thus the integral in \eqref{eq:r8-physical-effective} converges in $H^1$.
Its convergence in $L^2_3L^1_h$ follows from
Lemma~\ref{lem:r10-momentum-source}. The source integrals on $[0,1]$
are bounded by \eqref{eq:new-base-energy} and \eqref{eq:H2-NLsource}.
These estimates prove \eqref{eq:new-effective-datum}.

\smallskip\noindent
\emph{Step 2. Comparison in $H^1$.}
The projected equation \eqref{eq:projected-P-intro} and Duhamel's formula give
\begin{equation}\label{eq:r13-effective-duhamel}
 \begin{aligned}
 v(t)-S_h(t)\mathbb Pu_\infty
 ={}&\int_0^{t/2}[S_h(t-s)-S_h(t)]\mathbb PF_2(s)\dd s\\
 &+\int_{t/2}^tS_h(t-s)\mathbb PF_2(s)\dd s
       -S_h(t)\int_{t/2}^\infty\mathbb PF_2(s)\dd s.
 \end{aligned}
\end{equation}
First assume $t\ge4$. For $j=0,1$, the horizontal heat bound and
the $H^1$ contractivity of $\mathbb P$ yield for the first Duhamel term
\[
 \left\|\nabla_h^j\int_0^{t/2}
       [S_h(t-s)-S_h(t)]\mathbb PF_2(s)\dd s\right\|_{H^1}
       \le Ct^{-1-j/2}\int_0^{t/2}s\|F_2(s)\|_{H^1}\dd s.
\]
For $b\ge1$, Cauchy--Schwarz and \eqref{eq:r8-refined-source} give
\[
 \int_b^{b+1}s\|F_2(s)\|_{H^1}\dd s
 \le(1+b)\|F_2\|_{L^2(b,b+1;H^1)}
 \le C\mathcal J_2\mathfrak E_2(b)^{1/2}(1+b)^{-1/2}
       +C\mathcal J_2^2(1+b)^{-5/8}.
\]
The contribution of $[0,1]$ to this time moment is bounded by
$C\mathcal J_2^2$ by \eqref{eq:new-base-energy} and \eqref{eq:H2-NLsource}.
Partition at $\sqrt t$ and sum the unit-interval estimate, using
$\mathfrak E_2(b)^{1/2}\le C\mathcal J_2$ on $[1,\sqrt t]$ and
$\mathfrak E_2(b)\le\mathfrak E_2(\sqrt t)$ on $[\sqrt t,t/2]$
from Proposition~\ref{prop:r8-energy-tail}.
The sums of $(1+b)^{-1/2}$ over these two ranges are bounded by
$Ct^{1/4}$ and $Ct^{1/2}$, respectively; summing $(1+b)^{-5/8}$
over $[1,t/2]$ gives $Ct^{3/8}$. Hence
\[
 \int_0^{t/2}s\|F_2(s)\|_{H^1}\dd s
 \le C\mathcal J_2^2t^{1/4}
       +C\mathcal J_2\mathfrak E_2(\sqrt t)^{1/2}t^{1/2}
       +C\mathcal J_2^2t^{3/8}.
\]
For the Duhamel integral over $[t/2,t]$, contractivity of $S_h$
on $H^1$ gives, for $j=0$,
\[
 \left\|\int_{t/2}^tS_h(t-s)\mathbb PF_2(s)\dd s\right\|_{H^1}
       \le\|F_2\|_{L^1(t/2,t;H^1)}.
\]
For $j=1$, applying \eqref{eq:r8-horizontal-maximal} with $m=1$ gives
\[
 \left\|\nabla_h\int_{t/2}^tS_h(t-s)\mathbb PF_2(s)\dd s\right\|_{H^1}
       \le C\|F_2\|_{L^2(t/2,t;H^1)}.
\]
For the source tail over $[t/2,\infty)$, the homogeneous heat bound gives
\[
 \left\|\nabla_h^jS_h(t)\int_{t/2}^\infty\mathbb PF_2(s)\dd s\right\|_{H^1}
       \le Ct^{-j/2}\|F_2\|_{L^1(t/2,\infty;H^1)}.
\]
Substituting \eqref{eq:r13-effective-source-tails} at $b=t/2$
yields, for $j=0,1$,
{\small
\[
 \left\|\nabla_h^j\int_{t/2}^tS_h(t-s)\mathbb PF_2(s)\dd s\right\|_{H^1}
 +\left\|\nabla_h^jS_h(t)\int_{t/2}^\infty\mathbb PF_2(s)\dd s\right\|_{H^1}
 \le Ct^{-(1+j)/2}
       [\mathcal J_2\mathfrak E_2(t/2)^{1/2}+\mathcal J_2^2t^{-1/8}].
\]
}
Proposition~\ref{prop:r8-energy-tail} gives
$\mathfrak E_2(t/2)\le\mathfrak E_2(\sqrt t)$ for $t\ge4$ and
$\mathfrak E_2(t)\to0$.
The estimates for the three Duhamel terms therefore yield
\eqref{eq:r8-H1-scattering}, including the vanishing of the
corresponding time-weighted error.
For $1\le t\le4$, \eqref{eq:new-effective-datum},
\eqref{eq:r13-effective-source-tails} at $b=1$, and the source estimates on
$[0,1]$ from \eqref{eq:new-base-energy} and \eqref{eq:H2-NLsource}
bound the left side of \eqref{eq:r8-H1-scattering} by $C\mathcal J_2^2$,
using \eqref{eq:r8-horizontal-maximal} for the horizontal gradient.
Since $t^{-5/8}\ge4^{-5/8}$ on $[1,4]$, increasing $C$ proves
\eqref{eq:r8-H1-scattering} on this interval.
Finally, \eqref{eq:new-H2-main-rates} with $k=0,1,2$ gives
\[
 \|W(t)\|_{H^1}+t^{1/2}\|\nabla_hW(t)\|_{H^1}
       \le C\mathcal J_2t^{-3/4},\qquad W=(q,w,\tau).
\]
Since $u=v+w$, the $W$ estimate supplies the additional
$C\mathcal J_2t^{-3/4}$ term in the comparison of $U$ with
$(0,S_h(t)\mathbb Pu_\infty,0)$.

\smallskip\noindent
\emph{Step 3. Horizontal derivative estimates.}
For $t\ge1$, apply $\nabla_h^k$, $k=0,1,2$, to
\eqref{eq:r13-effective-duhamel}. The horizontal heat bounds of
Lemma~\ref{lem:linear-semigroup} and \eqref{eq:r10-momentum-source} give
\[
 \left\|\nabla_h^k\int_0^{t/2} [S_h(t-s)-S_h(t)]\mathbb PF_2(s)\dd s\right\|_2\le C\mathcal J_2^2t^{-(3+k)/2} \int_0^{t/2}s(1+s)^{-5/4}\dd s \le C\mathcal J_2^2t^{-3/4-k/2}.
\]
On $[t/2,t]$, use the $L^2_3L^1_h\to L^2$ heat bound for $k=0$
and the $L^2\to L^2$ gradient bound for $k=1$.
With the respective source estimates in \eqref{eq:r10-momentum-source},
we obtain
\[
 \begin{aligned}
 \left\|\int_{t/2}^tS_h(t-s)\mathbb PF_2(s)\dd s\right\|_2&\le C\mathcal J_2^2t^{-5/4} \int_{t/2}^t(t-s)^{-1/2}\dd s \le C\mathcal J_2^2t^{-3/4},\\
 \left\|\nabla_h\int_{t/2}^tS_h(t-s)\mathbb PF_2(s)\dd s\right\|_2&\le C\mathcal J_2^2t^{-7/4} \int_{t/2}^t(t-s)^{-1/2}\dd s \le C\mathcal J_2^2t^{-5/4}.
 \end{aligned}
\]

For $k=2$, use the horizontal frequency cutoffs $P_L^h,P_H^h$
of Lemma~\ref{lem:linear-semigroup}.
At low horizontal frequencies, weighted Cauchy--Schwarz in time at each
frequency, followed by Plancherel's theorem, gives
\[
 \begin{aligned}
 \left\|\nabla_h^2P_L^h \int_{t/2}^tS_h(t-s)\mathbb PF_2(s)\dd s\right\|_2^2&\le C\int_{t/2}^t(1+t-s)^{-1}\|F_2(s)\|_2^2\dd s\\
 &\le C\mathcal J_2^4t^{-7/2} \int_{t/2}^t(1+t-s)^{-1}\dd s \le C\mathcal J_2^4t^{-7/2}\log(e+t).
 \end{aligned}
\]
Here the time kernel follows from
$|\xi_h|^2e^{-\nu_h(t-s)|\xi_h|^2}\le C(1+t-s)^{-1}$ on the support
of $P_L^h$; the source bound is \eqref{eq:r10-momentum-source}.
At high horizontal frequencies, \eqref{eq:new-maximal-horizontal}
with $k=2$, \eqref{eq:r8-horizontal-sources}, and
\eqref{eq:new-window-convolution} yield
\[
 \left\|\nabla_h^2P_H^h \int_{t/2}^tS_h(t-s)\mathbb PF_2(s)\dd s\right\|_2^2\le C\int_{t/2}^te^{-c_*(t-s)}\|\nabla_hF_2(s)\|_2^2\dd s \le C\mathcal J_2^4t^{-49/12}.
\]
The source tail is bounded by
\[
 \left\|\nabla_h^kS_h(t) \int_{t/2}^\infty\mathbb PF_2(s)\dd s\right\|_2\le C\mathcal J_2^2t^{-(1+k)/2} \int_{t/2}^\infty(1+s)^{-5/4}\dd s \le C\mathcal J_2^2t^{-3/4-k/2}.
\]

Taking square roots in the $k=2$ estimates and adding the three
Duhamel terms proves \eqref{eq:new-effective-errors}.
\end{proof}

\subsection{Horizontal Gaussian profiles}
The following heat-flow limit identifies the leading term determined by
the horizontal mean of the effective datum.

\begin{lemma}[Linear Gaussian profiles]
\label{lem:r13-linear-gaussian}
Let $g\in L^2\cap L^2_3L^1_h$ be a vector field, and set
$m_g(x_3)=\int_{\R^2}g_h(x_h,x_3)\dd x_h$.
Then $m_g\in L^2_3$, with
$\|m_g\|_{L^2_3}\le\|g\|_{L^2_3L^1_h}$, and
\begin{equation}\label{eq:r13-linear-gaussian}
 \|\nabla_h^kS_h(t)\mathbb Pg-(\nabla_h^kK_h(t)m_g,0)\|_2
       =o(t^{-(1+k)/2}),\qquad k=0,1,2,3.
\end{equation}
The same profile holds in $L^2(t,t+1;L^2)$.
For the tensor of all ordered horizontal derivatives,
\[
 \|\nabla_h^kK_h(t)\|_{L^2_h}^2
       =\frac{k!}{8\pi\nu_h(2\nu_h)^k}t^{-1-k}.
\]
\end{lemma}
\begin{proof}
Dominated convergence gives
$\mathcal F_hg(\xi_h,\cdot)\to\int_{\R^2}g\dd x_h$
in $L^2_3$ as $\xi_h\to0$: the Fourier difference is bounded
pointwise by $2\int_{\R^2}|g|\dd x_h$, which belongs to $L^2_3$.
The Leray symbol converges to $\operatorname{diag}(1,1,0)$ for almost
every $\xi_3\ne0$ and has operator norm at most one.
Vertical Plancherel and dominated convergence therefore yield
\[
 \mathbb P(\xi_h,D_3)\mathcal F_hg(\xi_h,\cdot)
       \longrightarrow(m_g,0)\quad\text{in }L^2_3.
\]
Apply horizontal Plancherel and set $\upsilon_h=\sqrt t\,\xi_h$.
The integrand of the squared profile error multiplied by $t^{1+k}$
tends to zero and is bounded by
$C|\upsilon_h|^{2k}e^{-2\nu_h|\upsilon_h|^2}\|g\|_{L^2_3L^1_h}^2$.
Dominated convergence proves \eqref{eq:r13-linear-gaussian}.
The supremum over $s\ge t$ of the scaled profile error tends to zero,
which also gives the unit-interval profile.
Finally, polar integration gives
\[
 \|\nabla_h^kK_h(t)\|_{L^2_h}^2
 =\frac{1}{(2\pi)^2}\int_{\R^2}|\xi_h|^{2k}
                   e^{-2\nu_ht|\xi_h|^2}\dd\xi_h
 =\frac{k!}{8\pi\nu_h(2\nu_h)^k}t^{-1-k}.\qedhere
\]
\end{proof}

\begin{proposition}[Gaussian profiles and exact constants]\label{prop:gaussian}
Assume small $\mathcal J_2$, with the mean $m_\infty$ and
kernel $K_h$ defined in Section~\ref{sec:r12-leading-terms}.
For $k=0,1,2,3$,
\begin{equation}\label{eq:new-gaussian}
 \|\nabla_h^ku(t)-(\nabla_h^kK_h(t)m_\infty,0)\|_2
                         =o(t^{-(1+k)/2}),
\end{equation}
and
\begin{equation}\label{eq:new-gaussian-constants}
 \lim_{t\to\infty}t^{(1+k)/2}\|\nabla_h^ku(t)\|_2
 =\left(\frac{k!}{8\pi\nu_h(2\nu_h)^k}\right)^{1/2}\|m_\infty\|_2.
\end{equation}
The profile also holds in $L^2(t,t+1;L^2)$, and
\[
 \lim_{t\to\infty}t^{1+k}\int_t^{t+1}\|\nabla_h^ku(s)\|_2^2\dd s
       =\frac{k!}{8\pi\nu_h(2\nu_h)^k}\|m_\infty\|_2^2.
\]
Here $\nabla_h^k$ contains all ordered horizontal derivatives.
If $m_\infty\ne0$, these are matching upper and lower bounds
for the fixed nonlinear solution. Also
\begin{equation}\label{eq:new-mass-correction}
 \|m_\infty-\int_{\R^2}(u_0)_h\dd x_h\|_2\le C\mathcal J_2^2.
\end{equation}
\end{proposition}
\begin{proof}
\emph{Step 1. Comparison through second order and mean correction.}
Proposition~\ref{prop:r8-H1-scattering} gives
$u_\infty\in H^1\cap L^2_3L^1_h$.
Since $u=v+w$, \eqref{eq:new-effective-errors} and
\eqref{eq:new-H2-main-rates} imply
\[
 \|\nabla_h^k[u(t)-S_h(t)\mathbb Pu_\infty]\|_2
       =o(t^{-(1+k)/2}),\qquad k=0,1,2.
\]
The horizontal mean bound in Lemma~\ref{lem:r13-linear-gaussian}
and \eqref{eq:new-effective-datum} give \eqref{eq:new-mass-correction}.

\smallskip\noindent
\emph{Step 2. The third-order comparison.}
We prove
\begin{equation}\label{eq:r9-third-effective}
 \|\nabla_h^3(v(t)-S_h(t)\mathbb Pu_\infty)\|_2=o(t^{-2}).
\end{equation}
For $t\ge2$, apply $\nabla_h^3$ to
\eqref{eq:r13-effective-duhamel} and split at $|\xi_h|=1$.
At low horizontal frequencies, \eqref{eq:Sh-low-all-orders} with five
horizontal derivatives and \eqref{eq:r10-momentum-source}
bound the semigroup difference over $[0,t/2]$ by
\[
 C\mathcal J_2^2t^{-3}\int_0^{t/2}s(1+s)^{-5/4}\dd s
       \le C\mathcal J_2^2t^{-9/4}.
\]
For the convolution over $[t/2,t]$, place one horizontal
derivative on $F_2$. The low-frequency bound
$|\xi_h|^2e^{-\nu_ha|\xi_h|^2}\le C(1+a)^{-1}$ and
\eqref{eq:r9-point-sources-pre} give
\[
 \begin{aligned}
 \left\|\nabla_h^3\mathbf1_{\{|D_h|\le1\}}
       \int_{t/2}^tS_h(t-s)\mathbb PF_2(s)\dd s\right\|_2
 &\le C\int_{t/2}^t(1+t-s)^{-1}\|\nabla_hF_2(s)\|_2\dd s\\
 &\le C\mathcal J_2^2t^{-49/24}\log(e+t)=o(t^{-2}).
 \end{aligned}
\]
For the source tail over $[t/2,\infty)$,
\eqref{eq:Sh-low-all-orders} with three horizontal derivatives and
\eqref{eq:r10-momentum-source} give
\[
 C\mathcal J_2^2t^{-2}\int_{t/2}^\infty(1+s)^{-5/4}\dd s
       \le C\mathcal J_2^2t^{-9/4}.
\]
At high horizontal frequencies, \eqref{eq:r9-high-shear-small}
bounds $\|\nabla_h^3\mathbf1_{\{|D_h|\ge1\}}v(t)\|_2$ by
$C\mathcal J_2e^{-c_*t}+C\mathcal J_2^2(1+t)^{-49/24}$.
The horizontal spectral gap and
$\|u_\infty\|_{H^1}\le C\mathcal J_2$ from \eqref{eq:new-effective-datum}
and \eqref{eq:new-data-sizes} bound
$\|\nabla_h^3\mathbf1_{\{|D_h|\ge1\}}S_h(t)\mathbb Pu_\infty\|_2$
by $C\mathcal J_2e^{-c_*t}$.
Combining these estimates proves \eqref{eq:r9-third-effective}.
Since $u=v+w$, \eqref{eq:r9-horizontal-compressible-small}
gives the same $o(t^{-2})$ comparison for $u$.

\smallskip\noindent
\emph{Step 3. Gaussian profiles and exact limits.}
Apply Lemma~\ref{lem:r13-linear-gaussian} with $g=u_\infty$.
The comparisons in Steps~1--2 and the triangle inequality prove
\eqref{eq:new-gaussian} for $k=0,1,2,3$.
The Gaussian norm identity in that lemma and the reverse triangle
inequality give \eqref{eq:new-gaussian-constants}.
The scaled remainder in \eqref{eq:new-gaussian} tends uniformly
to zero on $s\ge t$, which gives the unit-interval profile.
Since $t^{1+k}\int_t^{t+1}s^{-1-k}\dd s\to1$, the reverse triangle
inequality in $L^2(t,t+1;L^2)$ gives the unit-interval norm limit.
\end{proof}

\begin{corollary}[Sharpness of the horizontal decay rates]

There are arbitrarily small $H^2$ data in the class $\mathcal J_2$
whose nonlinear solutions satisfy
$\|\nabla_h^ku(t)\|_2\simeq t^{-(1+k)/2}$ for $k=0,1,2,3$.
\end{corollary}
\begin{proof}
Fix real functions $\phi\in C_c^\infty(\R^2)$ and
$\psi\in C_c^\infty(\R)$ with $\int_{\R^2}\phi\dd x_h\ne0$ and
$\psi\ne0$, and set
\[
 q_0=\tau_0=0,\qquad
 u_0^{(\delta)}(x)=\delta(\phi(x_h)\psi(x_3),0,0),\qquad\delta>0.
\]
By \eqref{eq:new-data-sizes}, these data satisfy $\mathcal J_2\le C\delta$.
The mean correction estimate \eqref{eq:new-mass-correction}
and the reverse triangle inequality give, for sufficiently small $\delta$,
\[
 \|m_\infty\|_2
 \ge\delta\left|\int_{\R^2}\phi\dd x_h\right|\|\psi\|_2-C\delta^2
 \ge\frac{\delta}{2}\left|\int_{\R^2}\phi\dd x_h\right|\|\psi\|_2>0.
\]
Here $C$ depends only on the fixed functions $\phi,\psi$.
For each fixed sufficiently small $\delta$, the positive limits in
\eqref{eq:new-gaussian-constants} of Proposition~\ref{prop:gaussian} imply
$\|\nabla_h^ku(t)\|_2\simeq t^{-(1+k)/2}$ as $t\to\infty$,
for $k=0,1,2,3$, with comparison constants depending on the chosen datum.
\end{proof}

Without a horizontal moment condition, \eqref{eq:new-gaussian}
does not provide a uniform gain in the algebraic decay exponent
of the remainder.
The quantitative bounds \eqref{eq:new-effective-errors} concern
the comparison of $v$ with the effective heat flow.

\subsection{Vertical Gaussian profiles}

Assume in this subsection that $\mathcal J_{2,1}$ is sufficiently small.
Since $u_\infty\in H^1\cap L^2_3L^1_h$, its vertical derivative belongs to $L^2$, but the required horizontal mixed-norm bound is obtained only after vertical frequency truncation. We first control the vertical high-frequency tails of the shear solution at the horizontal heat decay scale. The resulting uniform cutoff estimate gives $m_\infty\in H^1(\R)$ and allows us to remove the vertical cutoff from the differentiated Gaussian profile.

\subsubsection{Control of vertical frequency tails and the derivative profile}
\begin{proposition}[Vertical Gaussian profiles and dissipation tails]\label{prop:r8-vertical-profile}
For each fixed datum with small $\mathcal J_{2,1}$ the effective
horizontal mean $m_\infty=\int_{\R^2}(u_\infty)_h\dd x_h$
belongs to $H^1(\R)$, and for $j=0,1$,
\begin{equation}\label{eq:r8-vertical-profile}
 \|\nabla_h^j\partial_3u(t)
       -(\nabla_h^jK_h(t)\partial_3m_\infty,0)\|_2
       =o(t^{-(1+j)/2}).
\end{equation}
The corresponding norm constants are
$[j!/(8\pi\nu_h(2\nu_h)^j)]^{1/2}\|\partial_3m_\infty\|_2$.
The profile also holds in $L^2(t,t+1;L^2)$, and
\[
 \lim_{t\to\infty}t^{1+j}\int_t^{t+1}
        \|\nabla_h^j\partial_3u(s)\|_2^2\dd s
       =\frac{j!}{8\pi\nu_h(2\nu_h)^j}\|\partial_3m_\infty\|_2^2.
\]
The dissipation tails satisfy
\begin{align}\label{eq:r8-profile-tails}
 \lim_{t\to\infty}t\int_t^\infty\|\nabla_h\partial_3u(s)\|_2^2\dd s
       &=\frac{\|\partial_3m_\infty\|_2^2}{16\pi\nu_h^2},\notag\\
 \lim_{t\to\infty}t^2\int_t^\infty\|\nabla_h^2\partial_3u(s)\|_2^2\dd s
       &=\frac{\|\partial_3m_\infty\|_2^2}{32\pi\nu_h^3}.
\end{align}
\end{proposition}
\begin{proof}
\emph{Step 1. Vertical continuity of the characteristic maps.}
Use the backward vertical characteristic maps
$\phi_{x_h}$ from Step~1 of the proof of
Lemma~\ref{lem:r8-forced-endpoint}, along the horizontal translation paths
in its time-splitting construction.
The supremum below is over the maps on all intervals $[s,t]$,
$0\le s\le t$, and all horizontal translation paths in that construction.
Lemma~\ref{lem:r8-real-drift} and \eqref{eq:r8-drift-assumptions}
bound their Jacobians above and below uniformly.
We first prove the vertical translation estimate
\begin{equation}\label{eq:r8-graph-modulus}
 \sup_\phi\left\|\sup_z
 |\nabla_h\phi(x_h,z+h)-\nabla_h\phi(x_h,z)|\right\|_{2,h}
       \le C\mathcal J_{2,1}|h|^{3/14},\qquad 0<|h|\le1.
\end{equation}
By \eqref{eq:verticalcomp} and \eqref{eq:H2-dissipation-definition},
$\|\partial_3u_3\|_{H^2}\le C\cD_2(U)^{1/2}$.
The identity $\partial_3u_3=\Div u-\Div_hu_h$ and
\eqref{eq:new-H2-main-rates} with $k=2$ also give
$\int_0^\infty\|\nabla_h\partial_3u_3\|_2\dd s\le C\mathcal J_2$.
For $0<r\le1$, set
$\omega_3(a,r)=\sup_{|\eta|\le r}\|a(x_h,z+\eta)-a(x_h,z)\|_\infty$.
Fourier Cauchy--Schwarz gives
$\omega_3(\partial_3u_3,r)\le Cr^{1/2}\|\partial_3u_3\|_{H^2}$.
Splitting time at $r^{-4/7}$ and using \eqref{eq:new-base-energy} and
\eqref{eq:r8-lips-tail}, we obtain
\[
 \int_0^\infty\omega_3(\partial_3u_3,r)\dd s
 \le C\mathcal J_2\bigl(r^{1/2}r^{-2/7}+r^{3/14}\bigr)
 \le C\mathcal J_2r^{3/14}.
\]
In elapsed backward time,
$\frac{\dd}{\dd s}\partial_i\phi
=-\partial_i u_3-(\partial_3u_3)\partial_i\phi$, $i=1,2$,
with coefficients evaluated along the characteristics.
Their separation is at most $C|h|$ by the Jacobian bounds.
For any measurable displacement $\eta(x_h,z)$ with $|\eta|\le r$,
the fundamental theorem of calculus gives the variable-shift bound
\[
 \sup_z|\partial_i u_3(x_h,z+\eta)-\partial_i u_3(x_h,z)|
 \le r^{1/2}\|\partial_3\partial_i u_3(x_h,\cdot)\|_{L^2_3}.
\]
Let $M_\phi=\sup_\phi\|\nabla_h\phi\|_{L^2_hL^\infty_3}$, bounded by
Step~1 of Lemma~\ref{lem:r8-forced-endpoint}.
Subtracting the characteristic equations and applying Gronwall bounds
the left-hand side of \eqref{eq:r8-graph-modulus} by
\[
 C\left(|h|^{1/2}\int_0^\infty\|\partial_3\nabla_hu_3\|_2\dd s
       +M_\phi\int_0^\infty\omega_3(\partial_3u_3,C|h|)\dd s\right).
\]
Horizontal translations preserve these norms. The preceding bounds prove
\eqref{eq:r8-graph-modulus}; for $C|h|>1$, use the bound for $M_\phi$.

\smallskip\noindent
\emph{Step 2. Uniform vertical frequency tails of the auxiliary solution.}
Let $F$ be the auxiliary solution constructed in Step~4 of the proof of
Lemma~\ref{lem:r8-forced-endpoint}, which satisfies
\[
 \begin{gathered}
 F_t-\nu_h\Delta_hF+u_3\partial_3F
       =S_{\mathrm{reg}}+S_\nu+\Div_hB_\nu+H(F),\qquad F(0)=\partial_3u_0,\\
 H(F)=-(\nabla u_3)\partial_3\Delta^{-1}\Div F.
 \end{gathered}
\]
Here $S_{\mathrm{reg}}$ is given by \eqref{eq:r8-regular-source},
and $(S_\nu,B_\nu)$ is the viscosity decomposition in Step~2 of the proof
of Lemma~\ref{lem:r13-first-vertical-decay}.
That lemma verifies the source bounds and identifies
$\partial_3v=\mathbb PF$ by uniqueness.
Set $H_A^v=\mathbf1_{\{|D_3|\ge A\}}$.
For fixed $f\in L^2\cap L^2_3L^1_h$ with $\partial_3f\in L^2$,
truncate $f(x_h,\phi_{x_h}(z))$ to $|x_h|\le R$.
By \eqref{eq:r8-graph}, the discarded horizontal Fourier sections
tend uniformly to zero as $R\to\infty$, uniformly in $\phi$.
On $|x_h|\le R$, the Jacobian bounds give a vertical derivative
of each section bounded by $CR\|\partial_3f\|_2$.
Its $H_A^v$ tail is therefore at most $CRA^{-1}\|\partial_3f\|_2$.
Choose $R$ first and then $A$.

For the correction in \eqref{eq:r8-flux-chain}, fix a source time
and write
$g=\partial_3(B_\nu)_i$,
$g_M=(I-H_M^v)g$, and
$Q_\phi(x_h,z)=\partial_i\phi_{x_h}(z)g(x_h,\phi_{x_h}(z))$.
Horizontal Cauchy--Schwarz, the Jacobian bounds, and
\eqref{eq:r8-graph-modulus} give, for $0<|h|\le1$,
\[
 \sup_\phi\|Q_\phi(\cdot,\cdot+h)-Q_\phi\|_{L^2_3L^1_h}
 \le C\mathcal J_{2,1}
 \bigl(\|g-g_M\|_2+|h|\|\partial_3g_M\|_2
                         +|h|^{3/14}\|g\|_2\bigr).
\]
The three terms come from truncating $g$, translating $g_M\circ\phi$,
and translating $\partial_i\phi$, respectively;
$|\phi(z+h)-\phi(z)|\le C|h|$ controls the second term.
Since
$A\int_{-A^{-1}}^{A^{-1}}|e^{ih\zeta}-1|^2\dd h\ge c_*>0$
for $|\zeta|\ge A$, Plancherel yields
\[
 \sup_{\phi,\xi_h}\|H_A^v\mathcal F_hQ_\phi(\xi_h,\cdot)\|_{2,z}^2
 \le CA\int_{-A^{-1}}^{A^{-1}}
       \sup_\phi\|Q_\phi(\cdot,\cdot+h)-Q_\phi\|_{L^2_3L^1_h}^2\dd h
       \longrightarrow0.
\]
Take $A\to\infty$ with $M$ fixed, then $M\to\infty$.

The composition and flux bounds are uniform in the time interval,
horizontal translation path, and partition mesh.
The Gaussian averages have mass one, so weak-star lower
semicontinuity in $L^\infty_{\xi_h} L^2_3$ preserves these bounds
in the splitting limit. For the flux divergence in
\eqref{eq:r8-flux-chain}, the Fourier factor is bounded by
$|\xi_h|\le1$.
Apply the bounds at each source time to
$S_{\mathrm{reg}}+S_\nu+H(F)$ and $B_\nu$.
The $L^1_t$ bounds in \eqref{eq:r8-source-space}, verified in
Lemma~\ref{lem:r13-first-vertical-decay} and Step~4 of the proof of
Lemma~\ref{lem:r8-forced-endpoint}, dominate the source integrals
uniformly in the terminal time. Dominated convergence gives
\begin{equation}\label{eq:r8-lift-tight}
 \lim_{A\to\infty}\sup_{t\ge0}\sup_{|\xi_h|\le1}
                  \|H_A^v\mathcal F_hF(t)(\xi_h,\cdot)\|_{2,z}=0.
\end{equation}
The suprema use the continuous horizontal Fourier representatives
obtained by the truncation argument in Lemma~\ref{lem:r8-graph};
for $Q_\phi$, the discarded part is bounded by
$C\mathcal J_{2,1}\|\mathbf1_{\{|x_h|>R\}}g\|_2$ uniformly in $\phi$.

\smallskip\noindent
\emph{Step 3. The effective datum and the vertical Gaussian profile.}
We first transfer the vertical frequency-tail control from the
auxiliary solution to the effective datum.
Proposition~\ref{prop:r8-H1-scattering} gives
$\partial_3\mathbb Pu_\infty\in L^2$.
At $t=r^{-2}$, the heat multiplier is bounded below on $|\xi_h|\le r$.
Combining \eqref{eq:r8-H1-scattering} with \eqref{eq:r8-lift-tight} gives
\begin{equation}\label{eq:r8-effective-tight}
 \lim_{A\to\infty}\sup_{0<r\le1}
 r^{-1}\|\mathbf1_{\{|D_h|\le r\}}H_A^v\partial_3\mathbb Pu_\infty\|_2=0.
\end{equation}
Indeed, fix $r_0>0$.
For $0<r\le r_0$, the scaled comparison error from
\eqref{eq:r8-H1-scattering} tends uniformly to zero as $r_0\to0$,
while \eqref{eq:r8-lift-tight} controls the auxiliary solution.
For $r\ge r_0$, use
$r_0^{-1}\|H_A^v\partial_3\mathbb Pu_\infty\|_2$.
Let $A\to\infty$ first and then $r_0\to0$.

We next prove $m_\infty\in H^1$ by keeping the vertical cutoff fixed.
Choose $\chi_v\in C_c^\infty(\R)$ with $0\le\chi_v\le1$,
$\chi_v=1$ on $[-1,1]$, and
$\operatorname{supp}\chi_v\subset[-2,2]$, and write
$\chi_A=\chi_v(D_3/A)$.
For fixed $A>0$, the vertical convolution kernel of
$\partial_3\chi_A$ has $L^1$ norm at most $CA$.
Thus $u_\infty\in H^1\cap L^2_3L^1_h$ implies
$\partial_3\chi_Au_\infty\in L^2\cap L^2_3L^1_h$.
Its horizontal mean is $\chi_A\partial_3m_\infty$, where the
derivative of $m_\infty$ is initially understood in distributions.
Lemma~\ref{lem:r13-linear-gaussian}, applied at this fixed cutoff,
and commutation of $\chi_A$ with $S_h$ and $\mathbb P$ give
\[
 \begin{aligned}
 \frac{\|\chi_A\partial_3m_\infty\|_2}{\sqrt{8\pi\nu_h}}
 &=\lim_{t\to\infty}t^{1/2}
             \|S_h(t)\mathbb P\partial_3\chi_Au_\infty\|_2\\
 &\le\limsup_{t\to\infty}t^{1/2}\|\chi_A\partial_3v(t)\|_2
 +\limsup_{t\to\infty}t^{1/2}
       \|\partial_3[v(t)-S_h(t)\mathbb Pu_\infty]\|_2
 \le C\mathcal J_{2,1}.
\end{aligned}
\]
Here \eqref{eq:H2-intermediate-main} bounds the first term,
\eqref{eq:r8-H1-scattering} makes the second limit zero, and
$\|\chi_A\|_{L^2\to L^2}\le1$ makes the bound independent of $A$.
Since $\chi_v(\zeta/A)\to1$ for every $\zeta$, Fatou's lemma
in vertical Fourier space gives
$\|\partial_3m_\infty\|_2\le C\mathcal J_{2,1}$.
Together with $m_\infty\in L^2$, this proves $m_\infty\in H^1$.

We now remove the vertical cutoff from the heat-flow profile.
For $t\ge1$ and $j=0,1$, horizontal dyadic decomposition gives
\[
 t^{1+j}\|\nabla_h^jS_h(t)H_A^v\partial_3\mathbb Pu_\infty\|_2^2
 \le C_j\sup_{0<r\le1}r^{-2}
       \|\mathbf1_{\{|D_h|\le r\}}H_A^v\partial_3\mathbb Pu_\infty\|_2^2
       +Ct^{1+j}e^{-c_*t}\|\partial_3\mathbb Pu_\infty\|_2^2.
\]
Indeed, on $2^{-k-1}<|\xi_h|\le2^{-k}$, the squared $L^2$ mass
is bounded by $4^{-k}$ times the supremum on the right.
The sum $\sum_{k\ge0}(t4^{-k})^{1+j}e^{-c_*t4^{-k}}$ is uniformly
bounded by geometric summation when $t4^{-k}\le1$ and exponential
decay when $t4^{-k}>1$; the last term comes from $|\xi_h|>1$.
By \eqref{eq:r8-effective-tight}, the supremum tends to zero as $A\to\infty$.
Consequently,
\[
 \lim_{A\to\infty}\limsup_{t\to\infty}
 t^{(1+j)/2}\|\nabla_h^jS_h(t)H_A^v\partial_3\mathbb Pu_\infty\|_2=0,
 \qquad j=0,1.
\]
Since $1-\chi_v(D_3/A)$ vanishes on $|\xi_3|\le A$, this removes the
vertical cutoff from the heat flow.
The Gaussian amplitude converges because
$\chi_v(D_3/A)\partial_3m_\infty\to\partial_3m_\infty$ in $L^2_3$.
Equation~\eqref{eq:r8-H1-scattering} and the potential-velocity bounds
in \eqref{eq:new-H2-main-rates} with $k=1,2$ prove
\eqref{eq:r8-vertical-profile}.
The reverse triangle inequality and the Gaussian norm identity in
Lemma~\ref{lem:r13-linear-gaussian} give the norm constants. The same argument on $[t,t+1]$, using the tail supremum
of the scaled profile error, gives the unit-interval limits.

\smallskip\noindent
\emph{Step 4. Dissipation tails.}
Integrate the $j=1$ profile in \eqref{eq:r8-vertical-profile} over
$[t,\infty)$ and apply the reverse triangle inequality in
$L^2(t,\infty;L^2)$.
Since $t\int_t^\infty s^{-2}\dd s=1$, this proves the first limit in
\eqref{eq:r8-profile-tails}.
For the second limit, the projected momentum equation
\eqref{eq:projected-P-intro} gives
\[
 (\partial_s-\nu_h\Delta_h)
       \partial_3[v(s)-S_h(s)\mathbb Pu_\infty]
       =\mathbb P\partial_3F_2(s),
\]
where $F_2$ is the unprojected momentum source in \eqref{eq:F2}.
Testing with $-\Delta_h\partial_3[v(s)-S_h(s)\mathbb Pu_\infty]$,
integrating from $t$ to $T$, and discarding the nonnegative terminal
energy give, after $T\to\infty$,
\[
 \begin{aligned}
 &\int_t^\infty\|\nabla_h^2\partial_3
       [v(s)-S_h(s)\mathbb Pu_\infty]\|_2^2\dd s\\
 &\quad\le C\|\nabla_h\partial_3
       [v(t)-S_h(t)\mathbb Pu_\infty]\|_2^2
       +C\|F_2\|_{L^2(t,\infty;H^1)}^2=o(t^{-2}),
 \end{aligned}
\]
by \eqref{eq:r8-H1-scattering}, \eqref{eq:r13-effective-source-tails},
and $\mathfrak E_2(t)\to0$ from Proposition~\ref{prop:r8-energy-tail}.
The potential-velocity contribution is $O(t^{-5/2})$ by
\eqref{eq:r9-third-preliminary}.
The horizontal heat energy identity is
\[
 2\nu_h\int_t^\infty
    \|\nabla_h^2\partial_3S_h(s)\mathbb Pu_\infty\|_2^2\dd s
       =\|\nabla_h\partial_3S_h(t)\mathbb Pu_\infty\|_2^2.
\]
Its Gaussian profile from Step~3 gives the second constant in
\eqref{eq:r8-profile-tails}; the comparison and potential errors
vanish after multiplication by $t^2$.
\end{proof}

\subsubsection{The mixed third derivative and completion of the proof}
\label{sec:r13-mixed-profile}
\begin{proof}[Completion of the proof of Theorem~\ref{thm:r9-gaussians}]
Proposition~\ref{prop:r8-H1-scattering} constructs the effective datum.
Proposition~\ref{prop:gaussian} proves \eqref{eq:r9-gaussian-all} for $k=0,1,2,3$.
Proposition~\ref{prop:r8-vertical-profile} gives
\eqref{eq:r9-gaussian-vertical-all} for $j=0,1$ and the dissipation
tail constants. It remains to prove the vertical profile for $j=2$.

By \eqref{eq:r8-vertical-profile},
$\|\nabla_h[\partial_3v(t)-(K_h(t)\partial_3m_\infty,0)]\|_2=o(t^{-1})$.
Since the Gaussian term solves the horizontal heat equation,
\eqref{eq:r9-weighted-shear} gives
\[
 (\partial_t-\nu_h\Delta_h)
       [\partial_3v-(K_h(t)\partial_3m_\infty,0)]
       =-\mathbb P\partial_3\mathcal R_{\rm sh}.
\]
Duhamel's formula on $[t/2,t]$, \eqref{eq:r8-horizontal-maximal}, and
Lemma~\ref{lem:r9-source-tail} yield
\[
 \begin{aligned}
 &\|\nabla_h^2[\partial_3v(t)-(K_h(t)\partial_3m_\infty,0)]\|_2
 \le Ct^{-1/2}
       \|\nabla_h[\partial_3v(t/2)-(K_h(t/2)\partial_3m_\infty,0)]\|_2\\
 &\qquad+C\|\nabla_h\partial_3\mathcal R_{\rm sh}\|_{L^2(t/2,t;L^2)}
       =o(t^{-3/2}).
\end{aligned}
\]
Adding the potential-velocity estimate \eqref{eq:r9-third-preliminary}
proves \eqref{eq:r9-gaussian-vertical-all} for $j=2$. Gaussian integration gives
\[
 \lim_{t\to\infty}t^3\|\nabla_h^2\partial_3u(t)\|_2^2
       =\frac{\|\partial_3m_\infty\|_2^2}{16\pi\nu_h^3}.
\]
The profile remainder is also $o(t^{-3/2})$ in $L^2(t,t+1;L^2)$;
thus the same constant is the limit of
$t^3\int_t^{t+1}\|\nabla_h^2\partial_3u(s)\|_2^2\dd s$.
\end{proof}

\subsection{Conservation-law representation of the effective mean}

The following proposition expresses the Gaussian amplitude
$m_\infty$ through the conservation laws.

\begin{proposition}[Conservation-law formula for the effective mean]
Under small $\mathcal J_2$, the normalized flux
\begin{equation}\label{eq:new-physical-flux}
 \mathscr J_\infty=\int_0^\infty\int_{\R^2}e^qu_3u_h\dd x_h\dd t
       \in H^1(\R;\R^2),\qquad
 \|\mathscr J_\infty\|_{H^1}\le C\mathcal J_2^2
\end{equation}
satisfies
\begin{equation}\label{eq:new-conservation-mean}
 m_\infty=\int_{\R^2}e^{q_0}(u_0)_h\dd x_h
                            -\partial_3\mathscr J_\infty.
\end{equation}
\end{proposition}
\begin{proof}
\emph{Step 1. The cutoff horizontal mean.}
The definitions of $F_2$ and $N$ in \eqref{eq:F2} and
\eqref{eq:new-projected-source}--\eqref{eq:new-N-definitions} give
\[
 (F_2)_h-N_h
 =\Div_h\bigl(\nu_h(e^{-q}-1)\nabla_hu_h\bigr)
 +\nabla_h\bigl[\nu_c(e^{-q}-1)\Div u
                         -c^2(e^\tau-1)(q+\tau)\bigr].
\]
The products on the right belong to $L^2_3L^1_h$, so their horizontal
cutoff means vanish. Hence $(F_2)_h$ and $N_h$ have the same horizontal mean.
The horizontal mean of $u(t)$ is understood as a smooth-cutoff limit;
no membership of $u(t)$ in $L^2_3L^1_h$ is required.
The shear Duhamel formula and \eqref{eq:new-profile-source-ledger} give
\begin{equation}\label{eq:new-velocity-mean}
 m_u(t)=\int_{\R^2}(u_0)_h\dd x_h+\int_0^t\int_{\R^2}N_h(s)\dd x_h\dd s.
\end{equation}
To justify this limit, choose $\eta\in C_c^\infty(\R^2)$ equal to one
near zero and set $\eta_L(x_h)=\eta(x_h/L)$.
The bounds \eqref{eq:new-base-energy} and the source formulas
\eqref{eq:F1}--\eqref{eq:F3} give local time integrability in
$L^1\cap L^2$. Since $L^1\cap L^2\subset\dot H^{-1}(\R^3)$, the
compressible Duhamel formula and Lemma~\ref{lem:linear-semigroup}
give $w\in\dot H^{-1}$ on each finite time interval.
Thus $w=\nabla\Psi$ with $\Psi\in L^2$, and integration by parts gives
\[
 \left\|\int_{\R^2}\eta_Lw_h\dd x_h\right\|_{L^2_3}
 \le C\|\nabla_h\eta_L\|_{L^2_h}
           \|\Psi\|_{L^2(\{L\lesssim|x_h|\lesssim CL\})}\longrightarrow0.
\]
For the shear, dominated convergence in its Duhamel formula gives the
horizontal Fourier trace; the horizontal diffusion flux has trace zero.
The Fourier kernel of $\eta_L$ approximates this trace and proves
\eqref{eq:new-velocity-mean} in $L^2_3$.

\smallskip\noindent
\emph{Step 2. Time integrability of the vertical momentum flux.}
Set $B(t)=\int_{\R^2}e^qu_3u_h\dd x_h$.
Sobolev's inequality in $x_3$, \eqref{eq:new-H2-main-rates} with
$k=0,1,2$, and \eqref{eq:H2-main-vertical-weak} give
\[
 \|u\|_{L^\infty_3L^2_h}\le C\mathcal J_2(1+t)^{-3/8},\qquad
 \|\partial_3q\|_{L^\infty_3L^2_h}\le C\mathcal J_2(1+t)^{-3/2}.
\]
Horizontal H\"older's inequality, \eqref{eq:H2-main-mixed-and-full},
\eqref{eq:H2-main-vertical-weak}, and \eqref{eq:new-H2-main-Lp}
with $p=\infty$ then yield
\[
 \|B(t)\|_2\le C\|u_3\|_{L^\infty_3L^2_h}\|u\|_2
                  \le C\mathcal J_2^2(1+t)^{-3/2}.
\]
Differentiating $B$ gives the horizontal integrals of
$e^qu_3\partial_3u_h$, $e^q(\partial_3u_3)u_h$, and
$e^qu_3u_h\partial_3q$. Their $L^2_3$ norms are bounded, respectively, by
\[
 \begin{aligned}
 C\|u_3\|_{L^\infty_3L^2_h}\|\partial_3u\|_2
       &\le C\mathcal J_2^2(1+t)^{-5/4},
 \qquad C\|\partial_3u_3\|_2\|u\|_{L^\infty_3L^2_h}
       \le C\mathcal J_2^2(1+t)^{-11/8},\\
 C\|u_3\|_\infty\|\partial_3q\|_{L^\infty_3L^2_h}\|u\|_2
       &\le C\mathcal J_2^3(1+t)^{-7/2}.
\end{aligned}
\]
All four time weights are integrable. Together with the $H^2$
continuity in \eqref{eq:solspace1}, these estimates show that
$\mathscr J_\infty=\int_0^\infty B(t)\dd t$ converges in $H^1_3$
and prove \eqref{eq:new-physical-flux}.

\smallskip\noindent
\emph{Step 3. Horizontal momentum balance and the effective mean.}
Define the normalized horizontal momentum mean by
$M(t)=m_u(t)+\int_{\R^2}(e^q-1)u_h\dd x_h$,
where $m_u$ is the cutoff mean in \eqref{eq:new-velocity-mean}.
The conservative momentum equation is
\[
 \partial_t(\rho u)+\Div(\rho u\otimes u)+\nabla(R\rho\theta)
       =\mu\Delta_hu+(\lambda+\mu)\nabla\Div u.
\]
Average its horizontal components with the cutoff $\eta_L$ from Step~1.
The pressure perturbation $\rho\theta-\bar\rho\bar\theta$ and
$\Div u$ have vanishing $L^2$ tails, while
$\|\nabla_h\eta_L\|_{L^2_h}$ is bounded.
Thus the horizontal pressure and bulk-viscosity terms vanish as
$L\to\infty$; the constant pressure contributes zero.
Two integrations by parts in the shear term use
$\|\Delta_h\eta_L\|_{L^2_h}=O(L^{-1})$.
The horizontal convective flux is in $L^1_x$ and uses
$\|\nabla_h\eta_L\|_\infty=O(L^{-1})$.
These limits hold after integration on compact time intervals.
Dividing by $\bar\rho$ gives $\partial_tM+\partial_3B=0$
in distributions. Since $M\in C([0,\infty);L^2_3)$ by
\eqref{eq:new-velocity-mean} and \eqref{eq:solspace1}, integration gives
\[
 M(t)=\int_{\R^2}e^{q_0}(u_0)_h\dd x_h
                    -\partial_3\int_0^tB(s)\dd s
 \quad\text{in }L^2_3.
\]
By the source-mean identity in Step~1,
\eqref{eq:r8-physical-effective},
\eqref{eq:new-velocity-mean}, and the bound for $N$ in
\eqref{eq:new-profile-source-ledger} imply $m_u(t)\to m_\infty$ in $L^2_3$.
Moreover, Sobolev's inequality in $x_3$ and
\eqref{eq:new-H2-main-rates} with $k=0,1$ give
\[
 \|M(t)-m_u(t)\|_2
       \le C\|q\|_{L^\infty_3L^2_h}\|u\|_2
       \le C\mathcal J_2^2(1+t)^{-3/2}.
\]
Hence $M(t)\to m_\infty$ in $L^2_3$.
Passing to the limit in the integrated balance, using the $H^1_3$
convergence of the flux from Step~2, proves \eqref{eq:new-conservation-mean}.
\end{proof}

\appendix

\section{Stability and comparison estimates}\label{app:construction-tools}

Throughout the appendices, the coefficients satisfy \eqref{eq:coeff}.
The difference estimate below gives uniqueness in the small $H^2$
energy class. Weak stability and relative energy give the comparison
principles for the nonlinear lower bounds in
Appendix~\ref{sec:r10-counterexamples}.

\subsection{Time regularity and uniqueness at \texorpdfstring{$H^2$}{H2} regularity}

\begin{lemma}[Difference energy]\label{lem:r10-difference}
There are constants $\varepsilon_*,c_*,C>0$, depending only on the
fixed coefficients and equilibrium, with the following property.
Let $T>0$ and let $U^{(1)},U^{(2)}$ be distributional solutions of
\eqref{eq:perturbed} on $[0,T]$ such that, for $i=1,2$,
\[
 \|U^{(i)}\|_{L^\infty(0,T;H^2)}\le\varepsilon_*,
 \qquad \int_0^T\cD_2(U^{(i)}(t))\dd t<\infty.
\]
Then $U^{(i)}\in C([0,T];H^1)\cap C_w([0,T];H^2)$.
For $\delta U=U^{(1)}-U^{(2)}$, the energy
\begin{equation}
 E_\Delta=\frac12\int_{\R^3}
 \left(c^2|\delta q|^2+e^{-\tau^{(1)}}|\delta u|^2
                         +\frac{c^2}{\alpha}|\delta\tau|^2\right)\dd x
\end{equation}
belongs to $AC([0,T])$, satisfies
$C^{-1}\|\delta U\|_2^2\le E_\Delta\le C\|\delta U\|_2^2$, and obeys
\begin{equation}\label{eq:r10-difference-inequality}
 \begin{aligned}
 E_\Delta' +c_*\bigl(\|\nabla_h\delta u\|_2^2
       +\|\Div\delta u\|_2^2+\|\nabla\delta\tau\|_2^2\bigr)
 \le C\bigl(1+\cD_2(U^{(1)})+\cD_2(U^{(2)})\bigr)E_\Delta
 \end{aligned}
\end{equation}
for almost every $t\in(0,T)$. In particular,
\[
 \sup_{0\le t\le T}\|\delta U(t)\|_2
 \le C\exp\!\left(C T+C\sum_{i=1}^2
                   \int_0^T\cD_2(U^{(i)}(t))\dd t\right)
                        \|\delta U(0)\|_2.
\]
Thus the initial trace determines at most one solution in this class.
\end{lemma}
\begin{proof}
The equations first give $\partial_tU^{(i)}\in L^2(0,T;H^1)$.
For example, $\mathcal Lu,\Delta\tau\in L^2H^1$, and multiplication
by a bounded $H^2$ perturbation preserves $H^1$.
For the temperature equation, the product estimate
\eqref{eq:H2-H1-multiplier} gives
\begin{align*}
 \|u\cdot\nabla\tau\|_{H^1}
 &\le C\|U\|_{H^2}\|\nabla\tau\|_{H^2},\\
 \|(e^{-q}-1)\Delta\tau\|_{H^1}
 &\le C\|U\|_{H^2}\|\nabla\tau\|_{H^2},\\
 \|e^{-q-\tau}Q(u)\|_{H^1}
 +\|e^{-q}|\nabla\tau|^2\|_{H^1}
 &\le C\|U\|_{H^2}
       \bigl(\|(\nabla_hu,\Div u)\|_{H^2}+\|\nabla\tau\|_{H^2}\bigr).
\end{align*}
For the quadratic terms, use
\[
 \|gh\|_{H^1}
 \le C\bigl(\|g\|_6\|h\|_3
       +\|\nabla g\|_3\|h\|_6+\|g\|_6\|\nabla h\|_3\bigr)
 \le C\bigl(\|g\|_{H^1}\|h\|_{H^2}
                   +\|g\|_{H^2}\|h\|_{H^1}\bigr).
\]
Apply it to components of $\nabla_hu$, $\Div u$ or $\nabla\tau$.
Their $H^1$ norms are bounded by $C\|U\|_{H^2}$, and their
$H^2$ norms are controlled by $C\cD_2(U)^{1/2}$.
The exponential factors are bounded multipliers on $H^1$ by
\eqref{eq:H2-H1-multiplier}. Hence
\begin{equation}\label{eq:H2-thermal-source}
 \|F_3\|_{H^1}
 \le C\|U\|_{H^2}\cD_2(U)^{1/2}.
\end{equation}
This gives $F_3\in L^2(0,T;H^1)$ and, through the third equation
of \eqref{eq:perturbed}, $\tau_t\in L^2(0,T;H^1)$.
The first two equations are treated by the same product estimate
and \eqref{eq:H2-NLsource}.
Thus each solution has a strongly continuous $H^1$ representative
and a weakly continuous $H^2$ representative.

The product estimate needed below is
\begin{equation}\label{eq:r10-difference-triple}
 \begin{aligned}
 \left|\int A f g\dd x\right|
 &\le C\|A\|_{H^1}\|f\|_2^{1/2}
                         \|\nabla_hf\|_2^{1/2}\|g\|_2\\
 &\le \eta\|\nabla_hf\|_2^2+
       C_\eta\|A\|_{H^1}^{4/3}(\|f\|_2^2+\|g\|_2^2).
 \end{aligned}
\end{equation}
Use the plane trace and horizontal interpolation proved after
Lemma~\ref{lem:plane-trace}, followed by Young's inequality.

The exact difference equations are
\begin{align}
 &(\partial_t+u^{(1)}\cdot\nabla)\delta q+\Div\delta u
                    =-\delta u\cdot\nabla q^{(2)},\notag\\*
 &(\partial_t+u^{(1)}\cdot\nabla)\delta u
       -e^{-q^{(1)}}\mathcal L\delta u
       +c^2e^{\tau^{(1)}}\nabla(\delta q+\delta\tau)\notag\\*
 &\quad=-\delta u\cdot\nabla u^{(2)}
       +(e^{-q^{(1)}}-e^{-q^{(2)}})\mathcal Lu^{(2)}
       -c^2(e^{\tau^{(1)}}-e^{\tau^{(2)}})
                                      \nabla(q^{(2)}+\tau^{(2)}),\notag\\
 &(\partial_t+u^{(1)}\cdot\nabla)\delta\tau
       -\chi e^{-q^{(1)}}\Delta\delta\tau+\alpha\Div\delta u\notag\\*
 &\quad=-\delta u\cdot\nabla\tau^{(2)}
       +\chi(e^{-q^{(1)}}-e^{-q^{(2)}})\Delta\tau^{(2)}\notag\\*
 &\qquad+\chi e^{-q^{(1)}}(\nabla\tau^{(1)}+\nabla\tau^{(2)})
                                                    \cdot\nabla\delta\tau\notag
 +\chi(e^{-q^{(1)}}-e^{-q^{(2)}})|\nabla\tau^{(2)}|^2\notag\\*
 &\qquad+\frac{e^{-q^{(1)}-\tau^{(1)}}}{\bar\rho\bar\theta c_v}
                                      [Q(u^{(1)})-Q(u^{(2)})]
 +\frac{e^{-q^{(1)}-\tau^{(1)}}-e^{-q^{(2)}-\tau^{(2)}}}
                  {\bar\rho\bar\theta c_v}Q(u^{(2)}).
\label{eq:r10-difference-system}
\end{align}
Test the three equations with
$c^2\delta q$, $e^{-\tau^{(1)}}\delta u$, and
$(c^2/\alpha)\delta\tau$, respectively.
The principal pressure has constant coefficient after multiplication
by the velocity weight; integration by parts cancels exactly the
two scalar-divergence pairings.

\paragraph{Principal diffusion and the transported weight.}
The velocity diffusion gives
\[
 \int e^{-q^{(1)}-\tau^{(1)}}
          (\nu_h|\nabla_h\delta u|^2+\nu_c|\Div\delta u|^2)\dd x
\]
and coefficient-gradient terms. By \eqref{eq:r10-difference-triple}, their
absolute sum is bounded by
\[
 C\|\nabla e^{-q^{(1)}-\tau^{(1)}}\|_{H^1}
       \|\delta u\|_2^{1/2}\|\nabla_h\delta u\|_2^{1/2}
       (\|\nabla_h\delta u\|_2+\|\Div\delta u\|_2).
\]
The coefficient norm is uniformly bounded by the $H^2$ size.
Young absorbs an arbitrarily small fraction of the displayed
diffusion and leaves $C_\eta\|\delta u\|_2^2$.
The gradient term from temperature diffusion has the identical
bound with $\delta\tau$ and its full gradient.

The contribution from differentiating the velocity weight uses
the exact identity
\begin{equation}\label{eq:r10-weight-transport}
 \begin{aligned}
 &(\partial_t+u^{(1)}\cdot\nabla)e^{-\tau^{(1)}}
                 +e^{-\tau^{(1)}}\Div u^{(1)}\\
 &\quad=e^{-\tau^{(1)}}\left[
 (1+\alpha)\Div u^{(1)}
 -\chi e^{-q^{(1)}}(\Delta\tau^{(1)}+|\nabla\tau^{(1)}|^2)
 -\frac{e^{-q^{(1)}-\tau^{(1)}}}{\bar\rho\bar\theta c_v}
                                                   Q(u^{(1)})\right].
 \end{aligned}
\end{equation}
The $L^\infty$ norms of $\Div u^{(i)}$, $\nabla\tau^{(i)}$,
and $(\nabla_hu^{(i)},\Div u^{(i)})$ are bounded by
$C\cD_2(U^{(i)})^{1/2}$ at almost every time.
All terms in \eqref{eq:r10-weight-transport} other than
$\Delta\tau^{(1)}$ therefore contribute at most
$C(1+\cD_2(U^{(1)}))E_\Delta$.
For that remaining term,
\[
 \|e^{-q^{(1)}-\tau^{(1)}}\Delta\tau^{(1)}\|_{H^1}
                         \le C\cD_2(U^{(1)})^{1/2}.
\]
Use \eqref{eq:r10-difference-triple} with both last factors
equal to $\delta u$. The scalar transport terms only produce
$\Div u^{(1)}|\delta q|^2$ and
$\Div u^{(1)}|\delta\tau|^2$ and obey the same integrable bound.

\paragraph{Transport and coefficient differences.}
The velocity transport remainder uses
\eqref{eq:H2-Lipschitz}:
$\|\nabla u^{(2)}\|_\infty
\le C(1+\cD_2(U^{(2)})^{1/2})$.
The temperature transport remainder uses
$\|\nabla\tau^{(2)}\|_\infty$.
For $\delta u\cdot\nabla q^{(2)}$, use
\eqref{eq:r10-difference-triple} with
$A=\nabla q^{(2)}$, $f=\delta u$, $g=\delta q$.

By the integral mean-value formula, the momentum viscosity
coefficient difference is a product $\delta q\,A\,\delta u$
in the energy pairing, with
\[
 A=-e^{-\tau^{(1)}}\left(\int_0^1
             e^{-q^{(2)}-s\delta q}\dd s\right)\mathcal Lu^{(2)},
 \qquad \|A\|_{H^1}\le C\cD_2(U^{(2)})^{1/2}.
\]
Apply \eqref{eq:r10-difference-triple}.
For the pressure coefficient difference, apply
\eqref{eq:r10-difference-triple} with $\delta\tau$ replacing
$\delta q$ and the coefficient containing
$\nabla(q^{(2)}+\tau^{(2)})$.
For the temperature diffusion coefficient difference use
$A=(\int_0^1-e^{-q^{(2)}-s\delta q}\dd s)\Delta\tau^{(2)}$,
$f=\delta\tau$, and $g=\delta q$.
Its $H^1$ coefficient is again bounded by
$C\cD_2(U^{(2)})^{1/2}$.
All exponential multipliers in this paragraph are bounded
on $H^1$ by their uniformly bounded $H^2$ perturbation norms.
Since $x^{2/3}\le1+x$ for $x\ge0$, the coefficients left by
Young are integrable with the claimed bound.

\paragraph{Temperature-gradient and heating differences.}
The term with $\nabla\delta\tau$ is bounded by
\[
 C(\|\nabla\tau^{(1)}\|_\infty+\|\nabla\tau^{(2)}\|_\infty)
                  \|\nabla\delta\tau\|_2\|\delta\tau\|_2.
\]
Absorb its gradient and use the squared coefficient.
The coefficient difference multiplying
$|\nabla\tau^{(2)}|^2$ is bounded by
$C\cD_2(U^{(2)})(\|\delta q\|_2^2+\|\delta\tau\|_2^2)$.
The viscous heating terms satisfy the identity
\[
 Q(u^{(1)})-Q(u^{(2)})
 ={}\mu(\nabla_hu^{(1)}+\nabla_hu^{(2)}):\nabla_h\delta u
 +(\lambda+\mu)(\Div u^{(1)}+\Div u^{(2)})\Div\delta u .
\]
Its pairing with $\delta\tau$ is a small fraction of the
velocity diffusion plus
$C(\cD_2(U^{(1)})+\cD_2(U^{(2)}))\|\delta\tau\|_2^2$.
The last coefficient difference in
\eqref{eq:r10-difference-system} uses
$\|Q(u^{(2)})\|_\infty\le C\cD_2(U^{(2)})$ and the
pointwise bound by $C(|\delta q|+|\delta\tau|)$.

Choose the finitely many absorption constants $\eta$ sufficiently
small and add the pairings. Uniform positivity gives
\eqref{eq:r10-difference-inequality}.
Spatial and time mollification justify the tests: all the pairings are integrable in the stated spaces, and the transport commutator is controlled by
$\|\nabla u^{(1)}\|_\infty\in L^1(0,T)$.
The $H^1$ time continuity fixes the initial energy.
Gronwall proves the final assertions.
\end{proof}

\subsection{Weak stability and lower semicontinuity}

\begin{lemma}[Finite-time weak stability]\label{lem:r10-weak-stability}
Let $U_{0,n},U_0\in H^2$ satisfy
$U_{0,n}\rightharpoonup U_0$ in $H^2$ and
\[
 \sup_n\|U_{0,n}\|_{H^2}\le\delta_* ,
\]
where $\delta_*>0$ is small enough for Theorem~\ref{thm:global}
and Lemma~\ref{lem:r10-difference}. Let $U_n,U$ be the corresponding
global solutions. For every $T>0$,
\begin{align}
 U_n&\longrightarrow U
                  &&\hbox{in }C([0,T];H^1_{\rm loc}),\notag\\
 U_n(t)&\rightharpoonup U(t)
                  &&\hbox{in }H^2\hbox{ for every }t\in[0,T],\notag\\
 (\nabla_hu_n,\Div u_n,\nabla\tau_n)
       &\rightharpoonup(\nabla_hu,\Div u,\nabla\tau)
                  &&\hbox{in }L^2(0,T;H^2),\notag\\
 \nabla q_n&\rightharpoonup\nabla q
                  &&\hbox{in }L^2(0,T;H^1).
                                                    \label{eq:r10-weak-stability}
\end{align}
These convergences hold for the full sequence.
For any constant-coefficient derivative, with or without a Hodge
projection, its $L^2$ norm at a fixed time and its $L^1_tL^2_x$ or
$L^2_tL^2_x$ norm on a fixed finite interval are lower semicontinuous
under this weak convergence of the initial data. The same assertion
holds with a fixed Fourier cutoff and with bounded nonnegative time
weights. Untruncated derivative norms are understood as extended
nonnegative quantities and may equal $+\infty$.
\end{lemma}
\begin{proof}
The global energy estimate bounds $U_n$ in $L^\infty H^2$
and all dissipative derivatives in the stated time $L^2$ spaces.
The calculation at the beginning of
Lemma~\ref{lem:r10-difference} also bounds
$\partial_tU_n$ in $L^2(0,T;H^1)$.
On each bounded spatial domain, Rellich and time
equicontinuity give subsequential convergence in
$C([0,T];H^1_{\rm loc})$.
Interpolation gives strong local convergence at any fixed
Sobolev order below two, hence uniform local convergence of the
thermodynamic coefficients at an order above $3/2$.
The gradients converge strongly locally in $L^2$, so the
heating squares converge in $L^1_{\rm loc}$.
Transport and pressure pass by the same local strong convergence;
variable diffusion passes by weak $L^2$ convergence and strong
coefficient convergence. The unprojected limiting equations
and their initial trace are therefore the original ones.

The limit has the energy regularity by lower semicontinuity.
The $H^2$ continuity in Lemma~\ref{lem:local}
and uniqueness in Lemma~\ref{lem:r10-difference} identify it with $U$.
Every subsequence has a further subsequence with this same
limit; consequently the compactness conclusions hold for the
full sequence. Uniform global $H^2$ boundedness and local
distributional convergence identify the weak $H^2$ limit at
each fixed time. The same argument identifies all weak
dissipation limits on the full space-time cylinder.

For a fixed radius $R$, a derivative multiplier restricted to
$|\xi|\le R$ is a bounded map from $H^2$ to $L^2$.
Its norm is therefore weakly lower semicontinuous at every
fixed time. Space-time $L^2$ norms have the same property.
For a bounded nonnegative weight $w$ on a finite interval,
the integral of the norm has the dual formula
\[
 \int w(t)\|g(t)\|_2\dd t
 =\sup_{\|\phi(t)\|_2\le w(t)}
                  \int\langle g(t),\phi(t)\rangle\dd t.
\]
It too is lower semicontinuous under weak space-time $L^2$
convergence. Weighted integrals of the squared norm follow
by multiplication by $\sqrt w$.
Finally let $R$ increase to infinity. Plancherel and monotone
convergence identify the extended derivative norms and integrals.
A supremum of lower semicontinuous functions is lower
semicontinuous, proving all the stated lower semicontinuity assertions.
\end{proof}

Let $X$ denote the real space of initial data for which
\begin{equation}\label{eq:r14-comparison-class}
 \|U_0\|_X=\mathcal J_{2,1}(U_0)
                    +\|\partial_3^2u_0\|_{L^2_3L^1_h}
\end{equation}
is finite. Here $\mathcal J_2(U_0)$ and $\mathcal J_{2,1}(U_0)$
denote the sizes in \eqref{eq:new-data-sizes} evaluated at $U_0$.

\begin{corollary}[Translation separation and the initial-data topology]
\label{cor:r10-translation}
The spaces of initial data with finite $\mathcal J_2$,
$\mathcal J_{2,1}$, and $\|\cdot\|_X$, equipped with their respective
norms, are real Banach spaces. Every open ball in a closed linear
subspace is a Baire space. On a sufficiently small ball, each strict
superlevel set of a derivative norm in
Lemma~\ref{lem:r10-weak-stability} is open in the corresponding norm
topology.

Let $F,P\in H^2$ satisfy $\|F\|_{H^2}+\|P\|_{H^2}\le\delta_*$,
with $\delta_*$ as in Lemma~\ref{lem:r10-weak-stability}, and set
$P_y(x_h,z)=P(x_h,z-y)$. If $U^{F+P_y}$ and $U^P$ denote the
solutions with these initial data, then
\[
 U^{F+P_y}(t,x_h,z+y)\longrightarrow U^P(t,x_h,z)
\]
in all the local strong and global weak senses of
\eqref{eq:r10-weak-stability} on every fixed finite time interval
as $|y|\to\infty$. Consequently, every finite collection of strict
lower bounds for translation-invariant derivative norms of $U^P$
also holds for $U^{F+P_y}$ when $|y|$ is sufficiently large.
\end{corollary}
\begin{proof}
A Cauchy sequence in a displayed sum norm converges in each
component Banach space. The limits agree as distributions;
the mixed-norm derivative limits equal the distributional
derivatives of the $H^2$ limit. This proves completeness.
The symmetry, zero-component, divergence, and fixed Fourier-support
conditions used below are closed under that convergence.
Open subsets of Banach spaces are Baire spaces.
Strong convergence in any of the three data norms implies
weak $H^2$ convergence, so the preceding lower semicontinuity
makes every strict superlevel set open.

Translations of a fixed $H^2$ function tend weakly to zero
as their displacement tends to infinity; prove this first
against compactly supported smooth test functions and then
use density for each derivative. In the translated frame
the initial datum is $P+F(\cdot,\cdot+y)\rightharpoonup P$.
Translation invariance of the original equations and
Lemma~\ref{lem:r10-weak-stability} give the assertion.
Lower semicontinuity preserves the finitely many strict lower bounds.
\end{proof}

\paragraph{Countable lower bounds.}
Let $Y$ be a closed linear subspace of one of the three data spaces
above, invariant under vertical translations, and let
$B_\delta=\{F\in Y:\|F\|_Y<\delta\}$, with $\delta$
sufficiently small. Consider a countable family of sets of initial data. Each set is the union, over finite times, of the sets on which finitely many norms in Lemma~\ref{lem:r10-weak-stability} satisfy prescribed strict lower bounds. These sets are open.
Suppose that each set contains data in $Y$ with arbitrarily small norm. For any $F\in B_\delta$, choose such a datum $P$ small enough
that $F+P_y$ remains in a prescribed neighborhood of $F$ for every
$y$. Corollary~\ref{cor:r10-translation} places $F+P_y$ in the prescribed set
for sufficiently large $|y|$. Each set is therefore dense, and
Baire's theorem gives a datum in their intersection.
Lower bounds for time-integral tails can be imposed on finite subintervals
of those tails. Thus the argument applies to pointwise norms,
finite-time integral norms, and tail lower bounds.

For a fixed solution, the truncated derivative norm is
continuous in time because $U\in C_tH^2$.
Its increasing supremum is therefore a lower semicontinuous,
extended-real-valued function of time. In particular, an
infinite norm at one time forces every finite threshold to be
exceeded in a neighborhood of that time.

\subsection{Comparison with a smooth reference}

\begin{lemma}[Relative energy for the eddy stress]\label{lem:r10-relative}
Let $T>0$, let $(\rho,u,\theta)$ be a solution in the class of
Theorem~\ref{thm:global}, and let $(r,V,\Theta)$ be a smooth reference
field with the same far-field state. Assume that
\[
 \widetilde U=(\log(r/\bar\rho),V,\log(\Theta/\bar\theta))
                       \in C([0,T];H^2)
\]
and that, for fixed $0<m_*<M_*<\infty$,
\[
 m_*\le\rho,r,\theta,\Theta\le M_*
                  \quad\hbox{on }[0,T]\times\R^3.
\]
Define the reference residuals in the nonconservative physical
equations by
\begin{equation}
 \begin{aligned}
 f_\rho&=r_t+\Div(rV),
 \qquad f_m=r(V_t+V\cdot\nabla V)+\nabla(Rr\Theta)
                                      -\Div\mathbb S^{\rm ed}(V),\\
 f_\theta&=rc_v(\Theta_t+V\cdot\nabla\Theta)
       +Rr\Theta\Div V-Q(V)-\kappa\Delta\Theta .
\end{aligned}
\end{equation}
Assume $(f_\rho,f_m,f_\theta)\in L^1(0,T;L^2)$ and
$A_{\rm ref}\in L^1(0,T)$, where
\begin{equation}\label{eq:r10-relative-coefficient}
 A_{\rm ref}={}\|\nabla V\|_\infty
 +\|\Div\mathbb S^{\rm ed}(V)\|_\infty
 +\|\nabla\Theta\|_\infty+\|\Delta\Theta\|_\infty
 +\|\nabla_hV\|_\infty^2+\|\Div V\|_\infty^2 .
\end{equation}
Set $z=u-V$, $b=\rho/r$, $d=\theta/\Theta$,
$H_\rho=\rho\log b-\rho+r$, and
$K_\theta=\theta-\Theta-\Theta\log d$.
Then
\begin{equation}
 E_{\rm rel}=\int_{\R^3}
       \left(\tfrac12\rho|z|^2+c_v\rho K_\theta+R\Theta H_\rho\right)\dd x
\end{equation}
satisfies
$C^{-1}\|U-\widetilde U\|_2^2\le E_{\rm rel}
\le C\|U-\widetilde U\|_2^2$, and
\begin{equation}\label{eq:r8-relative-stability}
 \sup_{0\le t\le T}\|U(t)-\widetilde U(t)\|_2
 \le C e^{C\int_0^T A_{\rm ref}(s)\dd s}
       \left(\|U_0-\widetilde U_0\|_2+
                    \int_0^T\|(f_\rho,f_m,f_\theta)\|_2\dd s\right).
\end{equation}
Here $C$ depends only on $m_*,M_*$ and the fixed physical
coefficients.
\end{lemma}
\begin{proof}
The identity \eqref{eq:r8-relative-identity}, derived in
Appendix~\ref{app:r10-relative}, gives the estimate after the following
quadratic bounds. The eddy stress satisfies
\[
 Q(V+z)=Q(V)+2\mathbb S^{\rm ed}(V):\nabla z+Q(z),\qquad
 |\mathbb S^{\rm ed}(V):\nabla z|\le Q(V)^{1/2}Q(z)^{1/2}.
\]
Thus the only viscous cross term is absorbed in half of
$d^{-1}Q(z)$, leaving $CQ(V)|d-1|^2$.
The scalar coefficients
$d^{-1}-1+b\log d$ and $(b-1)\log d$ are bounded by
$C(|b-1|^2+|d-1|^2)$ on the fixed compact range.
The relative potentials are quadratic and positive there.
Every unforced remainder is consequently bounded by
$CA_{\rm ref}E_{\rm rel}$.
The three residual coefficients in
\eqref{eq:r8-relative-identity} are bounded pointwise by
$C(|z|+|\rho-r|+|\theta-\Theta|)$.
Their integrals are bounded by
$C\|(f_\rho,f_m,f_\theta)\|_2 E_{\rm rel}^{1/2}$.
Apply Gronwall to $(E_{\rm rel}+\eta)^{1/2}$ and let
$\eta\downarrow0$.

For solutions in the stated class, justify the identities by
spatial regularization and a cutoff on the relative energy.
The energy bounds and \eqref{eq:r10-relative-coefficient} ensure
integrability of the pairings; $L^2$ flux estimates allow removal
of the cutoff. Approximation and the strong continuity in
Lemma~\ref{lem:local} give the identity at every endpoint.
See \cite{FeireislNovotnyRelative} for the classical relative-energy method.
\end{proof}

If a reference is written in logarithmic variables and has
residuals $R_q,R_u,R_\tau$ in the advective form of
\eqref{eq:perturbed}, then the residuals just defined are exactly
\begin{equation}
       f_\rho=rR_q,\qquad f_m=rR_u,\qquad
                                  f_\theta=rc_v\Theta R_\tau.
\end{equation}
The corresponding conservative momentum residual is $f_m+Vf_\rho$.

\section{Nonlinear obstructions and sharpness}\label{sec:r10-counterexamples}

The thermal and potential-velocity constructions give the two
fixed-data optimality assertions in Theorem~\ref{thm:decay-H2}(iv).
The lower bounds are first obtained for smooth families and then
transferred to fixed data by the countable-intersection argument
following Corollary~\ref{cor:r10-translation}. The choices of scale
and amplitude are summarized below.
\begin{center}\small
\begin{tabularx}{\linewidth}{@{}>{\RaggedRight\arraybackslash}p{.2\linewidth}>{\RaggedRight\arraybackslash}X>{\RaggedRight\arraybackslash}X@{}}
\toprule
Construction & Order of choices & Error relative to the leading term\\ \midrule
Thermal heating &
Fix small $\delta$, then let $\ell\to\infty$ with
$\eta_\ell=\delta/\log(e+\ell)$. &
The heating comparison and the normalized
$\eta_\ell^{-1}\ell^{-1/4}$ remainder tend to zero.\\
Potential-velocity layer &
For each fixed $\ell,\eta_\ell$, shrink $a_\ell$ subject to
\eqref{eq:r9-potential-thickness}; then let $\ell\to\infty$. &
The density-response error ratio is bounded by
$C_\ell \ell^4(a_\ell/\eta_\ell+\sqrt{a_\ell})\le \ell^{-1}$.\\
\bottomrule
\end{tabularx}
\end{center}
In the potential-velocity construction, $C_\ell$ is independent of
$a_\ell$ but may depend on the already fixed horizontal scale $\ell$.

\subsection{Initial localization and Gaussian remainders}

We first record the uniform lower bounds for the linearized system
and the absence of a uniform modulus for the nonlinear Gaussian
remainders.

\begin{proposition}[Linear scaling lower bounds]
For every $M>0$ there are smooth initial data
$U_0^{(\ell)}=(0,u_0^{(\ell)},0)$, $\ell\ge2$, such that
$\Div u_0^{(\ell)}=0$, $\mathcal J_2(U_0^{(\ell)})\le M$, and the
corresponding homogeneous linearized solutions satisfy
\[
 \begin{aligned}
 \|u^{\rm lin}(\ell^2)\|_p
       &\ge c_p M\ell^{-7/4+3/(2p)},&&2\le p\le\infty,\\
 \|\partial_3u^{\rm lin}(\ell^2)\|_2
       &\ge c_*M\ell^{-1/2},\\
 \|\nabla_h\partial_3u^{\rm lin}(\ell^2)\|_2
       &\ge c_*M\ell^{-3/2},\\
 \|\partial_3^2u^{\rm lin}(\ell^2)\|_2&\ge c_*M.
 \end{aligned}
\]
The constants are positive and independent of $M$ and $\ell$.
Consequently, on a $\mathcal J_2$-bounded class the corresponding
uniform decay powers $7/8-3/(4p)$, $1/4$, and $3/4$ cannot be
increased, and the full linear $H^2$ norm has no uniform positive
algebraic decay power.
\end{proposition}
\begin{proof}
Choose nonzero Schwartz $\phi:\R^2\to\R^2$ with
$\Div_h\phi=0$ and nonzero Schwartz $\psi:\R\to\R$.
For $\ell_0\ge2$ take
\[
 q_0=\tau_0=0,\qquad
 u_0^{(\ell_0)}=\varepsilon \ell_0^{-7/4}(\phi(x_h/\ell_0),0)
                                             \psi(\ell_0^{1/2}x_3).
\]
Direct scaling gives $\|u_0^{(\ell_0)}\|_{H^2}\le C\varepsilon$,
$\|u_0^{(\ell_0)}\|_1=C\varepsilon \ell_0^{-1/4}$ and
$\|u_0^{(\ell_0)}\|_{L^2_3L^1_h}=C\varepsilon$.
The linear solution is horizontal heat flow. At $t=\ell_0^2$,
\[
 \|u_{\rm lin}(\ell_0^2)\|_p=c_p\varepsilon \ell_0^{-7/4+3/(2p)},
 \quad \|\partial_3u_{\rm lin}(\ell_0^2)\|_2=c_0\varepsilon \ell_0^{-1/2},
\]
\[
 \|\nabla_h\partial_3u_{\rm lin}(\ell_0^2)\|_2=c_1\varepsilon \ell_0^{-3/2},
 \qquad \|\partial_3^2u_{\rm lin}(\ell_0^2)\|_2=c_2\varepsilon.
\]
All constants are positive and independent of $\ell_0$.
Choose $\varepsilon=M/C_0$, with $C_0$ large enough to make
$\mathcal J_2(U_0^{(\ell_0)})\le M$. This proves the lower bounds
and their consequences. The initial vertical mixed norms satisfy
$\|\partial_3^ju_0^{(\ell_0)}\|_{L^2_3L^1_h}\simeq\varepsilon \ell_0^{j/2}$.
Thus these mixed derivative norms are not uniformly bounded as
$\ell_0\to\infty$.
\end{proof}

\begin{proposition}[Linear $H^2$ decay under $X$ localization]

Let $U_0\in X$, where $X$ is defined by
\eqref{eq:r14-comparison-class}, and let $U^{\rm lin}$ be the
homogeneous linearized solution with initial datum $U_0$.
For every $t\ge0$,
\begin{equation}
 \|U^{\rm lin}(t)\|_{H^2}+\|\partial_3^2u^{\rm lin}(t)\|_2
       \le C\|U_0\|_X(1+t)^{-1/2},\qquad
 \|\nabla_h\partial_3u^{\rm lin}(t)\|_2
       \le C\|U_0\|_X(1+t)^{-1}.
\end{equation}
Here $C$ depends only on the fixed coefficients.
\end{proposition}
\begin{proof}
Apply \eqref{eq:Shestimate} with $m+n\le2$, using the three
mixed norms in $X$. The compressible contribution follows from
\eqref{eq:Scestimate} for $k=0,1,2$. Sum all derivatives through
order two. Every initial derivative used here belongs to the
stated $H^2$ or mixed norm.
\end{proof}

\begin{proposition}[Arbitrarily slow Gaussian remainders]

For every sufficiently small $\delta>0$ and every countable family
of functions $\omega_l:[1,\infty)\to(0,\infty)$ with
$\omega_l(t)\to0$, there is a datum
$U_0\in X\cap\bigcap_{m\ge0}H^m$ satisfying
\[
 \|U_0\|_X\le\delta,\qquad q_0=\tau_0=(u_0)_3=0,
 \qquad\Div_h(u_0)_h=0,
\]
whose global solution has $m_\infty=0$ and the following property.
Each of the four functions
\[
 t^{1/2}\|\partial_3u(t)\|_2,\qquad
 t\|\nabla_h\partial_3u(t)\|_2,
 \qquad t\int_t^\infty\|\nabla_h\partial_3u(s)\|_2^2\dd s,\qquad
 t^2\int_t^\infty\|\nabla_h^2\partial_3u(s)\|_2^2\dd s.
\]
tends to zero as $t\to\infty$, whereas its ratio to $\omega_l(t)$
has infinite limsup for every $l\ge1$. In particular, the normalized
remainders in \eqref{eq:r8-vertical-profile} and
\eqref{eq:r8-profile-tails} have no uniform vanishing modulus on the
small $X$ ball with $m_\infty=0$.
\end{proposition}
\begin{proof}
Choose a nonzero real radial Schwartz function $\varphi$ with
Fourier support in $1<|\xi_h|<2$, and set
$V=\nabla_h^\perp\varphi$. Let $\psi$ be a nonzero real Schwartz
function with Fourier support in $1<|\zeta|<2$. The physical
velocity data
\[
 p_{\delta,\ell}=\delta \ell^{-2}(V(x_h/\ell),0)\psi(z)
\]
have $X$ size comparable to $\delta$, uniformly in $\ell\ge1$.
Their Fourier support lies in
$\ell^{-1}<|\xi_h|<2\ell^{-1}$, $1<|\zeta|<2$.
Under $x_h\mapsto-x_h$, the horizontal velocity is odd, whereas
$q,\tau,u_3$ are even. The equations and uniqueness preserve this
half-turn symmetry, so $(F_2)_h$ is odd. Hence the horizontal mean
of $(u_\infty)_h$ vanishes: $m_\infty=0$.
At $t\in[\ell^2,2\ell^2]$ the horizontal heat-flow norm of
$\nabla_h^j\partial_3p_{\delta,\ell}$ is bounded below by
$c_*\delta \ell^{-1-j}$. Equation~\eqref{eq:r8-amplitude} bounds
the nonlinear error by
$C\delta^2\sqrt{\log(e/\delta)}\ell^{-1-j}$.
A fixed small $\delta$ therefore gives the same lower bound
for the nonlinear solution, uniformly in $\ell$.
The horizontal band also gives the two squared tail lower
bounds $c_*\delta^2\ell^{-2}$ and $c_*\delta^2\ell^{-4}$.
Taking $\ell\to\infty$ proves the uniform-modulus assertion.

For the fixed-datum construction use the lower semicontinuity
and translation separation of Lemma~\ref{lem:r10-weak-stability}
and Corollary~\ref{cor:r10-translation}. The symmetry and fixed
Fourier-support conditions define a closed subspace of $X$.
Recursively choose a summable amplitude $\delta_n$ inside a
nested $X$ ball preserving all previously established strict lower bounds;
next choose $\ell_n$ so large that
$\max_{l\le n}\omega_l(\ell_n^2)<c_*\delta_n^2/n$; finally choose
a vertical translation preserving the four new lower bounds.
Lower semicontinuity allows a smaller closed ball about this
new datum contained in all the required open sets. Choose its
radius at most $2^{-n}$ and keep it inside the prescribed
$\delta$ ball. Completeness gives one
limit datum satisfying all the prescribed lower bounds.
Its initial Fourier support remains in
$|\xi_h|\le2$, $1\le|\zeta|\le2$, so it belongs to every
integer Sobolev space. Proposition~\ref{prop:r8-vertical-profile} gives the four zero limits,
while the preserved lower bounds give the claimed ratios.
\end{proof}

\subsection{Sharpness of the thermal exponent}
\label{sec:r9-thermal-counter}

Horizontal viscous heating generates a profile for the third
derivatives of temperature of order $t^{-17/8}$ for the family below.

\subsubsection{Quadratic comparison with the linearized solution}

The nonlinear correction to the horizontal viscous heating is
controlled by the following quadratic estimate.

\Needspace{18\baselineskip}
\begin{lemma}[Quadratic error in the basic localization class]
\label{lem:new-quadratic-error}
There are $\delta_*,C>0$, depending only on the fixed coefficients
and equilibrium, such that the following holds. Let $U$ be the
global solution with $\mathcal J_2(U_0)\le\delta_*$, and let
$U^{\rm lin}=\mathscr S(t)U_0$ solve the full homogeneous linearized
system with the same initial datum. Set $r=u-u^{\rm lin}$ and
$R_W=(q-q^{\rm lin},\mathbb Qr,\tau-\tau^{\rm lin})$.
For every $t\ge0$,
\begin{align}\label{eq:new-quadratic-error}
 \sup_{s\ge0}\|U(s)-U^{\rm lin}(s)\|_{H^1}&\le C\mathcal J_2^2,\notag\\*
 \|r(t)\|_2+(1+t)^{1/2}\|\nabla_hr(t)\|_2
       &\le C\mathcal J_2^2(1+t)^{-1/2},\notag\\*
 \|\partial_3r(t)\|_2
       &\le C\mathcal J_2^2(1+t)^{-1/2}\sqrt{\log(e+t)},\notag\\*
 \|R_W(t)\|_2+(1+t)^{1/2}\|\nabla R_W(t)\|_2
       &\le C\mathcal J_2^2(1+t)^{-3/4},\notag\\*
 \|r_3(t)\|_{L^\infty_3L^2_h}+\|\partial_3r_3(t)\|_2
       &\le C\mathcal J_2^2(1+t)^{-1}.
\end{align}
\end{lemma}
\begin{proof}
\emph{Step 1. Source integrability and bounded-time estimates.}
Write $F=(F_1,F_2,F_3)$ for the unprojected source in \eqref{eq:perturbed}.
On $0\le s\le1$, \eqref{eq:H2-NLsource} and
\eqref{eq:H2-thermal-source} give
\[
 \|F(s)\|_{H^1}
 \le C\mathcal J_2\cD_2(U(s))^{1/2}.
\]
Its time integral is $O(\mathcal J_2^2)$ by Cauchy--Schwarz and
\eqref{eq:new-base-energy}. On $[1,\infty)$, sum the time-integrated $H^1$ estimates
\eqref{eq:r8-F2-H1-window} for $F_2$, and integrate
Lemma~\ref{lem:r10-scalar-sources} for $F_1,F_3$.
Consequently,
\begin{equation}\label{eq:r11-integrated-physical-source}
                 \int_0^\infty\|F(s)\|_{H^1}\dd s
                       \le C\mathcal J_2^2.
\end{equation}
The multiplier bound \eqref{eq:blockmultiplier} and the $H^1$
contraction of $S_h(t)\mathbb P$ give boundedness of the full linear
semigroup on $H^1$.
Duhamel's formula proves the first line of
\eqref{eq:new-quadratic-error}. It also handles all bounded
times in the remaining lines, since
$H^1(\R^3)\subset L^\infty_3L^2_h$.
We henceforth take $t\ge2$.

\smallskip\noindent
\emph{Step 2. Compressible variables.}
The exact split \eqref{eq:r9-early-flux} gives
$F_2=N_0+\partial_3B_0$, with
\[
 \int_0^\infty
   \bigl(\|(F_1,F_3)(s)\|_1+\|N_0(s)\|_1+\|B_0(s)\|_1\bigr)\dd s
 \le C\mathcal J_2^2 .
\]
For $s\ge1$, the respective powers are $7/4$, $3/2$, and
$5/4$. On $[0,1]$, the displayed split involves at most
second derivatives; $L^2\cdot L^2\subset L^1$ and the
original energy bound suffice. At low spatial frequencies,
one vertical derivative on the momentum column gives an
additional heat factor $(t-s)^{-1/2}$.
For $j=0,1$, the early part is therefore bounded by
\[
 C\mathcal J_2^2
       \bigl(t^{-3/4-j/2}+t^{-5/4-j/2}\bigr).
\]
This estimate uses the integrable flux pieces, not
$F_2\in L^1(0,\infty;L^1)$.

On $s\in[t/2,t]$ use the physical $L^2$ bound
$\|F(s)\|_2\le C\mathcal J_2^2(1+s)^{-7/4}$, from
Lemmas~\ref{lem:r10-momentum-source} and
\ref{lem:r10-scalar-sources}. The low-frequency contribution is
\[
 C\mathcal J_2^2t^{-7/4}
       \int_0^{t/2}(1+a)^{-j/2}\dd a
       \le C\mathcal J_2^2t^{-3/4-j/2},
                  \qquad j=0,1.
\]
At high spatial frequencies, \eqref{eq:component-Green-matrix}
gives, for both zero and one derivative, the kernel
$Ce^{-c_*a}(1+a^{-1/2})$ on the physical $L^2$ momentum source and the $H^1$ scalar sources. The early part is
$O(\mathcal J_2^2e^{-c_*t/4})$ by
\eqref{eq:r11-integrated-physical-source}; the late part is
$O(\mathcal J_2^2t^{-7/4})$ because this kernel is integrable.
This proves the fourth line of \eqref{eq:new-quadratic-error},
without a third density derivative.

\smallskip\noindent
\emph{Step 3. Horizontal and first vertical derivatives of the solenoidal velocity.}
Write
\[
 r_v=v-S_h(t)\mathbb Pu_0
     =v-S_h(t)\mathbb Pu_\infty
        +S_h(t)\mathbb P(u_\infty-u_0).
\]
The effective error \eqref{eq:new-effective-errors} and
the quadratic $L^2_3L^1_h$ bound
\eqref{eq:new-effective-datum} yield, for $j=0,1$,
\[
                \|\nabla_h^jr_v(t)\|_2
                    \le C\mathcal J_2^2t^{-(1+j)/2}.
\]
Add the compressible error to obtain the second line.
Proposition~\ref{prop:r8-vertical-correction} gives the third
line for $r_v$; the potential derivative is faster by
the estimate just proved.

\smallskip\noindent
\emph{Step 4. Mixed-norm estimates for the third velocity component.}
The early $N_0$ contribution uses the $L^1$ to
$L^\infty_3L^2_h$ kernel of power $1$ in
\eqref{eq:new-third-kernels-all}.
For the early vertical flux use
\[
 \partial_3(\mathbb PB_0)_3=-\Div_h(\mathbb PB_0)_h,
 \qquad
 \|B_0(s)\|_{L^\infty_3L^1_h}
 \le \|u_3(s)\|_{L^\infty_3L^2_h}
         \|u(s)\|_{L^\infty_3L^2_h}
 \le C\mathcal J_2^2(1+s)^{-11/8}.
\]
Here $\|u\|_{L^\infty_3L^2_h}
\le C\|u\|_2^{1/2}\|\partial_3u\|_2^{1/2}$.
For $\xi_h\in\R^2\setminus\{0\}$ and $i,j\in\{1,2\}$,
the vertical inverse Fourier transforms of the horizontal
projection rows are explicitly
\[
\mathcal F_3^{-1}[\mathbb P_{ij}(\xi_h,\cdot)](z)
    =\delta_{ij}\delta_0(z)
        -\frac{\xi_i\xi_j}{2|\xi_h|}e^{-|\xi_h||z|},
 \qquad \mathcal F_3^{-1}[\mathbb P_{i3}(\xi_h,\cdot)](z)
    =-\frac{\mathrm i\xi_i}{2}
                          \operatorname{sgn}(z)e^{-|\xi_h||z|}.
\]
Indeed, the inverse transform of
$(|\xi_h|^2+\xi_3^2)^{-1}$ is
$(2|\xi_h|)^{-1}e^{-|\xi_h||z|}$; differentiating it gives
the formula for $\mathbb P_{i3}$. The total variations of these
measures are bounded by an absolute constant, since
$|\xi_i\xi_j|/|\xi_h|^2\le1$ and
$|\xi_i|/|\xi_h|\le1$. In particular the Dirac term is retained.
Also
$\sup_{\xi_h,z}|\mathcal F_hG(\xi_h,z)|
 \le\|G\|_{L^\infty_3L^1_h}$.
Apply the vertical convolution at each horizontal frequency
and then horizontal Plancherel:
\[
\|\Div_hS_h(a)(\mathbb PG)_h\|_{L^\infty_3L^2_h}^2
 \le C\|G\|_{L^\infty_3L^1_h}^2
       \int_{\R^2}|\xi_h|^2e^{-2\nu_ha|\xi_h|^2}\dd\xi_h
 \le Ca^{-2}\|G\|_{L^\infty_3L^1_h}^2.
\]
Thus
\[
 \|\Div_hS_h(a)(\mathbb PG)_h\|_{L^\infty_3L^2_h}
       \le Ca^{-1}\|G\|_{L^\infty_3L^1_h},\qquad a>0.
\]
The zero horizontal frequency is a null set and does not
affect Plancherel. This is a bound for the differentiated heat
operator, without an endpoint mixed-space boundedness claim
for the projection alone.

The early flux is consequently $O(\mathcal J_2^2t^{-1})$.
This use of a mixed norm of the source also covers the Dirac part
of the projection kernel.
For the late part, the other third-row bound in
\eqref{eq:new-third-kernels-all} gives
\[
 C\mathcal J_2^2t^{-7/4}\int_0^{t/2}a^{-1/4}\dd a
                         \le C\mathcal J_2^2t^{-1}.
\]
For the potential error, vertical interpolation between its
$L^2$ and gradient bounds gives the same trace rate.
Finally $\partial_3(r_v)_3=-\Div_h(r_v)_h$ and the
compressible gradient bound prove the last derivative
estimate. All constants depend only on the fixed coefficients
and the small-data threshold, with no vertical mixed derivative
of the initial velocity required.
\end{proof}

\subsubsection{The thermal profile and lower bounds}

Let $\phi:\R^2\to\R^2$ be a nonzero real Schwartz field
with $\Div_h\phi=0$, and let $\psi$ be a nonzero real
Schwartz function on $\R$.
One may take $\phi=(-\partial_2,\partial_1)e^{-|X|^2}$
and $\psi(Z)=e^{-Z^2}$.
For $\delta>0$ and $\ell\ge2$, set
\begin{equation}\label{eq:r9-thermal-packet}
 \eta_\ell=\frac{\delta}{\log(e+\ell)},\qquad
 q_{0,\ell}=\tau_{0,\ell}=0,\qquad
 u_{0,\ell}(x_h,z)=\eta_\ell \ell^{-7/4}(\phi(x_h/\ell),0)\psi(\ell^{1/2}z).
\end{equation}
Direct scaling gives
\[
 \|U_{0,\ell}\|_{H^2}\simeq\eta_\ell,\quad
 \|u_{0,\ell}\|_{L^2_3L^1_h}\simeq\eta_\ell,\quad
 \|u_{0,\ell}\|_1=O(\eta_\ell \ell^{-1/4}).
\]
The first vertical mixed norm is comparable to
$\eta_\ell \ell^{1/2}$, so these data are not uniformly small
in $\mathcal J_{2,1}$.
The linearized solution with these initial data is the horizontal heat flow of $u_{0,\ell}$, with both scalars identically zero.

\begin{proposition}[Strong thermal profile for a concentrating family]

Let $\phi\in\mathscr S(\R^2;\R^2)$ and
$\psi\in\mathscr S(\R)$ be nonzero real-valued functions with
$\Div_h\phi=0$. There is $\delta_*>0$ such that, for every fixed
$0<\delta\le\delta_*$, the global solutions $U_\ell$ with initial
data \eqref{eq:r9-thermal-packet} satisfy
\begin{equation}\label{eq:r9-thermal-strong-profile}
 \eta_\ell^{-2}\ell^5\partial_3^3\tau_\ell(\ell^2,\ell X,\ell^{-1/2}Z)
 \longrightarrow
 -\frac{\mu}{\bar\rho\bar\theta c_v\chi}
       |\nabla_Xe^{\nu_h\Delta_X}\phi(X)|^2\partial_Z(\psi^2)(Z)
\end{equation}
strongly in $L^2_{X,Z}$.
For $T_3=\nabla^3$ or $T_3=\partial_3^3$, the pointwise and
time-integrated norms have the common positive limit
\begin{equation}\label{eq:r9-thermal-sharp-constant}
 \begin{aligned}
 \lim_{\ell\to\infty}\eta_\ell^{-2}\ell^{17/4}
                  \|T_3\tau_\ell(\ell^2)\|_2
 &=\lim_{\ell\to\infty}\eta_\ell^{-2}\ell^{17/4}
                  \|T_3\tau_\ell\|_{L^2(\ell^2,\ell^2+1;L^2)}\\
 &=\frac{\mu}{\bar\rho\bar\theta c_v\chi}
       \big\||\nabla_Xe^{\nu_h\Delta_X}\phi|^2\big\|_{2,X}
       \|\partial_Z(\psi^2)\|_{2,Z}>0.
 \end{aligned}
\end{equation}
The threshold $\delta_*$ depends only on $\phi,\psi$ and the fixed
physical coefficients.
\end{proposition}
\begin{proof}
For this family $\mathcal J_2\le C\eta_\ell$, so
Lemma~\ref{lem:new-quadratic-error} directly gives
\[
 \|\nabla_h(u_\ell-u_\ell^{\rm lin})(t)\|_2
       \le C\eta_\ell^2(1+t)^{-1}.
\]
For the mixed error, use
Proposition~\ref{prop:r8-vertical-correction} with $j=1$ on the
shear component. The potential part is bounded by the
$k=2$ case of \eqref{eq:new-H2-main-rates}
and the corresponding linear estimate \eqref{eq:Scestimate}.
Together they give
\[
 \|\nabla_h\partial_3(u_\ell-u_\ell^{\rm lin})(t)\|_2
 \le C\eta_\ell^2(1+t)^{-1}\sqrt{\log(e+t)}
                      +C\eta_\ell(1+t)^{-7/4}.
\]
These constants do not involve the large first mixed norm.
At $t=\ell^2$ the horizontal gradient of the nonlinear and linear velocities has $L^\infty$ norm at most $C\eta_\ell \ell^{-11/4}$,
and the linear mixed gradient has norm at most
$C\eta_\ell \ell^{-9/4}$.
Expanding the differentiated heating square gives
\begin{equation}\label{eq:r9-thermal-heating-comparison}
 \begin{aligned}
 &\|\nabla(|\nabla_hu_\ell|^2-|\nabla_hu_\ell^{\rm lin}|^2)(\ell^2)\|_2\\
 &\quad\le C\eta_\ell^3\ell^{-17/4}
       +C\eta_\ell^3\ell^{-19/4}\sqrt{\log(e+\ell)}
       +C\eta_\ell^2\ell^{-25/4}
       =o(\eta_\ell^2\ell^{-17/4}).
 \end{aligned}
\end{equation}
For the vertical derivative use the total horizontal
$L^\infty$ gradient with the mixed error, and the linear
mixed $L^\infty$ gradient with the horizontal error.
The horizontal derivative is faster by the heating estimates
of Section~\ref{sec:r9-thermal-upper}.

At this time the linear heating is exactly
\[
 \eta_\ell^2\ell^{-11/2}
       |\nabla_Xe^{\nu_h\Delta_X}\phi(x_h/\ell)|^2\psi(\ell^{1/2}z)^2.
\]
In the rescaled Fourier variables the pure vertical
third-derivative Newton multiplier is
$-i\zeta^3/(\zeta^2+\ell^{-3}|\xi_h|^2)$.
It tends to $-i\zeta$ and is bounded by $|\zeta|$.
Schwartz decay gives $L^2$ dominated convergence.
Equation \eqref{eq:r9-thermal-heating-comparison} permits
replacement by the nonlinear heating. The remainder in
\eqref{eq:r9-thermal-profile}, after normalization, is at
most $C\eta_\ell^{-1}\ell^{-1/4}\to0$.
This proves \eqref{eq:r9-thermal-strong-profile}.
For the ordered full tensor, the Newton norm equals the
heating gradient norm; its horizontal part has the extra
factor $\ell^{-3/2}$. Only the vertical part remains in the limit.
The limiting factors in \eqref{eq:r9-thermal-sharp-constant}
are nonzero. On $t=\ell^2+s$, $0\le s\le1$, the horizontal
heat time is $1+s/\ell^2$, and all estimates and convergence
are uniform. This also proves the time-integrated assertion.
\end{proof}

\paragraph{Fixed-data optimality.}
Consider a small $\mathcal J_2$ ball in the closed subspace
$\{(0,u_0,0):(u_0)_3=0,\ \Div_h(u_0)_h=0\}$.
For integers $k,m,N\ge1$, consider the set of initial data for which, at some integer $\ell$ with $\ell^2\ge N$, both the pointwise and time-integrated norms, multiplied by
$(1+\ell^2)^{17/8+1/k}$, exceed $m$.
The initial-data norm tends to zero, whereas
its normalized lower bound is
$c_*\delta^2\ell^{2/k}/[\log(e+\ell)]^2\to\infty$.
The countable-intersection argument following
Corollary~\ref{cor:r10-translation} makes these sets open and dense and gives
one datum satisfying the failure for every $\beta>17/8$,
at arbitrarily late times.
Its linearized temperature is zero, and its nonlinear temperature
satisfies the fixed-data $o(t^{-17/8})$ bound in
Theorem~\ref{thm:decay-H2}(iii).
This proves Theorem~\ref{thm:decay-H2}(iv)(a).

\subsection{Sharpness of the potential-velocity exponent}
\label{sec:r9-potential-counter}

The interaction of a broad thermodynamic background with a thin
horizontal shear layer produces a potential-velocity derivative of
order $t^{-2}$. The reference solution retains the density, velocity,
and temperature responses of the layer.

Choose the fixed functions
\[
 g_0(X,Y)=e^{-|X-(1,0)|^2-Y^2},\qquad
 \phi(X)=(-\partial_2,\partial_1)e^{-|X|^2},\qquad
 \psi(Z)=\partial_Z^3e^{-Z^2}.
\]
The first three moments of $\psi$ vanish and
$\partial_Z^{-2}\psi=\partial_Ze^{-Z^2}$.
For $\delta>0$, $\ell\ge2$, and $0<a_\ell\le\ell^{-2}$, put
$\eta_\ell=\delta/\log(e+\ell)$ and define
\begin{equation}\label{eq:r9-potential-data}
 \begin{aligned}
 q_{0,\ell}(x)&=\eta_\ell \ell^{-3}g_0(x/\ell),
 \qquad u_{0,\ell}(x)=\eta_\ell a_\ell^{3/2}\ell^{-1}
                       (\phi(x_h/\ell),0)\psi(z/a_\ell),\\
 \tau_{0,\ell}(x)&=-q_{0,\ell}(x)
       -\frac{\eta_\ell^2a_\ell^{7/2}}{\chi \ell^5}
           e^{\eta_\ell \ell^{-3}g_0(x_h/\ell,0)}
          [\phi(x_h/\ell)\cdot\nabla_Xg_0(x_h/\ell,0)]
                          \partial_Z^{-2}\psi(z/a_\ell).
\end{aligned}
\end{equation}
The small thermal correction matches the initial heat balance.

\begin{proposition}[Weak potential-velocity profile for a thin layer]

There is $\delta_*>0$ such that, for every fixed
$0<\delta\le\delta_*$, there exists a sequence
$0<a_\ell\le\ell^{-2}$ for which the data
\eqref{eq:r9-potential-data} satisfy
\[
 c_*\eta_\ell\le\mathcal J_{2,1}(U_{0,\ell})\le C\eta_\ell
\]
and the corresponding global solutions satisfy
\begin{equation}\label{eq:r9-potential-weak}
 \eta_\ell^{-2}\ell^5a_\ell^{1/2}
          \partial_3^3(w_\ell)_3(\ell^2+s,\ell X,a_\ell Z)
       \rightharpoonup G(X,Z)\quad\hbox{in }L^2_{X,Z},
                    \qquad 0\le s\le1,
\end{equation}
uniformly in $s$ against each fixed $L^2$ test function, where
\begin{equation}\label{eq:r9-potential-limit}
 G(X,Z)=-(e^{\nu_h\Delta_X}\phi)(X)\cdot
       \nabla_X(e^{\frac{\chi}{1+\alpha}\Delta}g_0)(X,0)
                    \,\psi''(Z)\ne0 .
\end{equation}
The Laplacian on $g_0$ acts in all three variables $(X_1,X_2,Y)$
before restriction to $Y=0$. In particular,
\[
 \liminf_{\ell\to\infty}\eta_\ell^{-2}\ell^4
       \inf_{0\le s\le1}\|\partial_3^3(w_\ell)_3(\ell^2+s)\|_2
       \ge\|G\|_2>0 .
\]
The same lower bound holds for
$\eta_\ell^{-2}\ell^4
\|\partial_3^3(w_\ell)_3\|_{L^2(\ell^2,\ell^2+1;L^2)}$.
The constants $c_*,C,\delta_*$ depend only on the fixed profiles
and physical coefficients.
\end{proposition}

\begin{proof}
\emph{Step 1. The broad background solution.}

Fix $\ell\ge2$ and $0<\eta\le\delta_*$. Let $U_b=(q_b,u_b,\tau_b)$ be
the original nonlinear solution with data
$(\eta \ell^{-3}g_0(x/\ell),0,-\eta \ell^{-3}g_0(x/\ell))$.
Reflection symmetry makes $q_b,\tau_b,u_{b,h}$ even in $z$
and $(u_b)_3$ odd.
\emph{Uniform estimates for the background solution.} Fix an integer $M\ge6$ large enough
for the derivatives below. For $0\le j\le M-2$ and $0\le t\le2\ell^2$,
\begin{equation}\label{eq:r9-broad-background}
 \begin{aligned}
 \|\nabla^j(q_b,\tau_b)(t)\|_\infty&\le C_j\eta \ell^{-3-j},\\
 \|\nabla^ju_b(t)\|_\infty&\le C_j\eta \ell^{-4-j},
 \qquad \|\nabla^j\partial_tU_b(t)\|_\infty\le C_j\eta \ell^{-5-j}.
\end{aligned}
\end{equation}
Moreover, with a fixed full high-frequency cutoff,
\begin{equation}\label{eq:r9-background-high}
 \|P_Hq_b(t)\|_{H^2}\le C\eta \ell^{-11/2},\qquad
 \eta^{-1}\ell^3q_b(\ell^2s,\ell\,\cdot)
       \longrightarrow e^{\frac{\chi}{1+\alpha}s\Delta}g_0
\end{equation}
in $H^M$ as $\ell\to\infty$, uniformly for $0\le s\le2$ and
$0<\eta\le\delta_*$. Constants may depend on the fixed $M$ and
$g_0$, but not on $\ell$ or $\eta$.

To prove \eqref{eq:r9-broad-background}--\eqref{eq:r9-background-high},
divide the unknowns by $\eta \ell^{-3}$ and put $s=t/\ell^2$, $y=x/\ell$.
The linear acoustic terms have coefficient $\ell$ and cancel
in the symmetrized integer-order energy with weights
$c^2e^{\tau_b},1,c^2e^{\tau_b}/\alpha$.
The transport and nonlinear pressure coefficients are
$O(\eta \ell^{-2})$; variable diffusion and heating have
coefficients $O(\eta \ell^{-3})$.
For $m=M,M+2$, the normalized energy satisfies
\[
 E_m'+c_*(\|\nabla_hu\|_{H^m_y}^2+
             \|\Div_yu\|_{H^m_y}^2+\|\nabla_y\tau\|_{H^m_y}^2)
       \le C_m\eta \ell^{-2}(E_m+E_m^3).
\]
Diffusion commutators are absorbed by the positive
dissipation. Thus $E_m$ is bounded on $0\le s\le2$,
uniformly in $\ell$ and $\eta$. The $H^{M+2}_y$ bound controls the
source in the difference equation between the nonlinear and linearized solutions with the same initial data in $H^M_y$;
the symmetrized difference energy gives error $C_M\eta\ell^{-2}$
on $0\le s\le2$.
At physical frequency $\xi/\ell$ the entropy root is
$-\chi|\xi|^2/((1+\alpha)\ell^2)+O(|\xi|^4/\ell^4)$.
Its velocity component, and the difference between the
initial vector $(1,0,-1)$ and the entropy projection,
are $O(|\xi|/\ell)$. Uniform semigroup bounds, a low-frequency
expansion, and the Schwartz tail give the linear velocity
bound $C_m/\ell$ and the entropy-heat limit with error $C_m/\ell$.
The scaled equations bound the time derivatives.
Restoring variables proves \eqref{eq:r9-broad-background}.
The physical fourth-derivative $L^2$ norm is
$O(\eta \ell^{-11/2})$, which controls the high-frequency
$H^2$ norm in \eqref{eq:r9-background-high}.
\Needspace{3\baselineskip}
\smallskip\noindent\emph{Step 2. The linearized layer equations.}

In the response equations all background coefficients
are evaluated at $z=0$, after any indicated physical derivative.
In particular $\Delta\tau_b$ retains its vertical second derivative.
Set $k=(c^2/\nu_c)e^{q_b+\tau_b}$ and
$a_b=\partial_3(u_b)_3(t,x_h,0)$.
Let $r$ denote the response of the logarithmic density and let
$b_h$ denote the horizontal velocity response. They solve
\begin{equation}\label{eq:r9-potential-response}
 \left\{\begin{aligned}
 &r_t+u_{b,h}\cdot\nabla_hr+a_bZ\partial_Zr+kr
       =-b_h\cdot\nabla_hq_b,\\
 &(b_h)_t+u_{b,h}\cdot\nabla_hb_h+a_bZ\partial_Zb_h
       -\nu_he^{-q_b}\Delta_hb_h\\
 &\qquad+b_h\cdot\nabla_hu_{b,h}
       +r(\partial_tu_b+u_b\cdot\nabla u_b)_h|_{z=0}=0,
 \end{aligned}\right.
\end{equation}
with $r(0)=0$ and $b_h(0,x_h,Z)=\phi(x_h/\ell)\psi(Z)$.
For each fixed pair of nonnegative integers $m,n$, the linear
parabolic--transport system has a smooth solution in the weighted
Sobolev norms
\[
 \|(1+|X|^2)^{n/2}(1+Z^2)^{(n+m-|\gamma|)/2}
                   \partial_X^\gamma\partial_Z^j(r,b_h)\|_2,
                          \qquad |\gamma|+j\le m,
                          \quad X=x_h/\ell.
\]
Indeed, commuting a horizontal derivative with $a_bZ\partial_Z$
lowers the number of horizontal derivatives on the unknown by at
least one. The corresponding norm has one additional $Z$ weight,
which controls the factor $Z$. Commuting $\partial_Z^j$ contributes
only $j a_b\partial_Z^j$. Weighted energy estimates, ordered by the
number of horizontal derivatives, therefore close the finite family.
The zero, first, and second $Z$ moments remain zero:
after integration, $Z\partial_Z$ multiplies the $j$th
moment by $-(j+1)$, so the moments solve homogeneous systems.
Therefore the decaying primitives
\begin{equation}
 b_3=\partial_Z^{-1}(kr-\Div_hb_h),\qquad
                   \Div_hb_h+\partial_Zb_3=kr
\end{equation}
are well-defined. The temperature response is determined by
\begin{equation}\label{eq:r9-potential-sigma}
 \begin{aligned}
 \chi e^{-q_b}\partial_Z^2\sigma
 ={}&b_h\cdot\nabla_h\tau_b+\alpha kr
       +\chi e^{-q_b}r(\Delta\tau_b+|\nabla\tau_b|^2)\\
 &-\frac{e^{-q_b-\tau_b}}{\bar\rho\bar\theta c_v}
 \left[2\mu\nabla_hu_{b,h}:\nabla_hb_h
       +2(\lambda+\mu)(\Div u_b)kr-rQ(u_b)\right]\big|_{z=0}.
 \end{aligned}
\end{equation}
Its right side has zero zeroth and first moments.
The formula
$\partial_Z^{-2}f(Z)=\int_{-\infty}^Z(Z-Y)f(Y)\dd Y$
and these zero moments bound the primitive by two more
polynomial weights, and give decay at both ends.

Weighted energy and the strictly positive $k$ give
$b_h(t,\ell X,Z)=O(1)$ and $r(t,\ell X,Z)=O(\eta \ell^{-4})$
in each fixed weighted Sobolev norm above.
The rescaled coefficients other than damping differ from
horizontal heat flow by $O(\eta \ell^{-3})$, so
\begin{equation}
 b_h(\ell^2s,\ell X,Z)=e^{\nu_hs\Delta_X}\phi(X)\psi(Z)+O(\eta \ell^{-3}).
\end{equation}
For the density response put $r_*=-k^{-1}b_h\cdot\nabla_hq_b$.
Equations \eqref{eq:r9-broad-background} and
\eqref{eq:r9-potential-response} give
$\|(\partial_t+u_{b,h}\cdot\nabla_h+a_bZ\partial_Z)r_*\|
\le C\eta \ell^{-6}$ in the same fixed norms.
Damping therefore gives
$\|r-r_*\|\le C\eta \ell^{-4}e^{-c_*t}+C\eta \ell^{-6}$.
Derivatives and weights follow in order, since derivatives
of $k$ are $O(\eta \ell^{-3})$.
At $t=\ell^2+s$ this proves, uniformly for $0\le s\le1$,
\begin{equation}\label{eq:r9-potential-r-profile}
 \eta^{-1}\ell^4r(\ell^2+s,\ell X,Z)
 \longrightarrow
 -\frac{\nu_c}{c^2}(e^{\nu_h\Delta_X}\phi)(X)\cdot
       \nabla_X(e^{\frac{\chi}{1+\alpha}\Delta}g_0)(X,0)\psi(Z).
\end{equation}
The trace at $Y=0$ is legitimate after taking a fixed
additional background derivative.

\smallskip\noindent\emph{Step 3. Residual estimates and comparison.}

For the remainder of this step, keep $\ell\ge2$ and
$0<\eta\le\delta_*$ fixed, and take $0<a\le\ell^{-2}$. Set
\[
 \varepsilon=\frac{\eta a^{3/2}}{\ell},\qquad Z=z/a.
\]
The reference field is
\begin{equation}
 \widetilde q=q_b+\varepsilon r,\qquad
 \widetilde u=u_b+\varepsilon(b_h,ab_3),\qquad
 \widetilde\tau=\tau_b+\varepsilon a^2\sigma.
\end{equation}
Here $r,b_h,b_3,\sigma$ are evaluated at $(t,x_h,Z)$, whereas
the background is evaluated at $(t,x_h,z)$.
Constants denoted by $C_\ell$ in this step may depend on $\ell$,
the fixed profiles, and the physical coefficients, but are uniform
for $0<a\le\ell^{-2}$ and $0<\eta\le\delta_*$. Their value may
increase from one estimate to the next.
At time zero $b_3=0$, and
\eqref{eq:r9-potential-sigma} gives exactly the thermal
correction in \eqref{eq:r9-potential-data}. Thus
$\widetilde U(0)=U_0$.

\paragraph{Cancellation of the leading linear terms.}
The perturbation divergence is $\varepsilon kr$. The largest
vertical linear terms in momentum are
\[
 -\nu_ce^{-q_b(t,x_h,aZ)}\frac{\varepsilon}{a}k\partial_Zr
 +c^2e^{\tau_b(t,x_h,aZ)}\frac{\varepsilon}{a}\partial_Zr.
\]
They cancel at $z=0$. Evenness of the background scalars makes
their sum $\varepsilon a$ times a uniformly bounded weighted
function of $(x_h,Z)$.
The horizontal terms cancel by the exact background identity
\[
 -\nu_ce^{-q_b}\nabla_h(kr)+c^2e^{\tau_b}\nabla_hr
                   +re^{-q_b}(\mathcal Lu_b)_h
       =r(\partial_tu_b+u_b\cdot\nabla u_b)_h.
\]
This explains the last term in \eqref{eq:r9-potential-response}.
The linear thermal terms of order $\varepsilon$ cancel by
\eqref{eq:r9-potential-sigma}.

\paragraph{Taylor remainders and layer interactions.}
For an even background coefficient $F$, Taylor's formula gives
\[
 F(t,x_h,aZ)-F(t,x_h,0)
 =a^2Z^2\int_0^1(1-s)\partial_3^2F(t,x_h,saZ)\dd s.
\]
Also,
$a^{-1}(u_b)_3(t,x_h,aZ)-a_bZ=O(a^2Z^3)$.
After multiplication by the layer responses, the powers of $Z$
are controlled by the weighted estimates above. To display the
mass residual explicitly, suppress $(t,x_h)$
in background arguments and use the density response equation:
\[
 \begin{aligned}
 R_q={}&\varepsilon[u_{b,h}(aZ)-u_{b,h}(0)]\cdot\nabla_hr
       +\varepsilon[a^{-1}(u_b)_3(aZ)-a_bZ]\partial_Zr\\
 &+\varepsilon b_h\cdot[\nabla_hq_b(aZ)-\nabla_hq_b(0)]
       +\varepsilon a b_3\partial_3q_b(aZ)
 +\varepsilon^2(b_h\cdot\nabla_hr+b_3\partial_Zr).
\end{aligned}
\]
The first four terms have size $O_\ell(\varepsilon a^2)$ in
weighted $L^2_{x_h,Z}$, and the last two have size
$O_\ell(\varepsilon^2)$.

In momentum, the unbalanced vertical acceleration, horizontal
diffusion of $ab_3$, and vertical gradient of $a^2\sigma$ have
size $O_\ell(\varepsilon a)$. The horizontal linear remainders have
size $O_\ell(\varepsilon a^2)$. The largest nonlinear viscosity
term is a bounded coefficient times
\[
 \frac{\varepsilon^2}{a}\,r\partial_Z(kr).
\]
It comes from the density coefficient multiplying the vertical
derivative of the layer divergence. The remaining quadratic
momentum terms have size at most $O_\ell(\varepsilon^2/a)$.

For temperature, the layer Laplacian is
\[
 \Delta(\varepsilon a^2\sigma)
   =\varepsilon\partial_Z^2\sigma
      +\varepsilon a^2\Delta_h\sigma.
\]
The first term is included in \eqref{eq:r9-potential-sigma}; the
second and the linear material derivative have size
$O_\ell(\varepsilon a^2)$. The quadratic part of
$Q(\widetilde u)$ contributed by the layer is exactly
\[
 \varepsilon^2\left\{
 \mu\bigl(|\nabla_hb_h|^2+a^2|\nabla_hb_3|^2\bigr)
       +(\lambda+\mu)|kr|^2\right\},
\]
which has size $O_\ell(\varepsilon^2)$. Density coefficient
differences multiplying the layer Laplacian have the same size.
The exponential Taylor remainders satisfy
$|e^h-1-h|\le C|h|^2$ on the fixed range of the perturbations;
their products with the displayed linear operators have no larger
order. This accounts for every term in the three logarithmic
residuals.

The conversion from weighted layer coordinates to physical $L^2$
uses
$\|f(x_h,z/a)\|_{L^2_x}=a^{1/2}\|f\|_{L^2_{x_h,Z}}$.
The principal orders are therefore
\begin{center}\small
\begin{tabularx}{\linewidth}{@{}>{\RaggedRight\arraybackslash}p{.39\linewidth}>{\centering\arraybackslash}p{.17\linewidth}>{\RaggedRight\arraybackslash}X@{}}
\toprule
Residual terms & Layer amplitude & Physical $L^2$ bound\\ \midrule
Vertical linear remainders & $\varepsilon a$ & $C_\ell\eta a^3$\\
Largest nonlinear viscosity term & $\varepsilon^2/a$ & $C_\ell\eta^2a^{5/2}$\\
Horizontal and mass linear remainders & $\varepsilon a^2$ & $C_\ell\eta a^4$\\
Quadratic transport and layer heating & $\varepsilon^2$ & $C_\ell\eta^2a^{7/2}$\\
\bottomrule
\end{tabularx}
\end{center}
Consequently, uniformly for $0\le t\le\ell^2+1$,
\begin{equation}\label{eq:r9-potential-residual}
 \begin{aligned}
 \|(R_q,(R_u)_h,R_\tau)(t)\|_2
       &\le C_\ell a^{1/2}
             (\varepsilon a^2+\varepsilon^2+\varepsilon^2/a),\\
 \|(R_u)_3(t)\|_2&\le C_\ell a^{1/2}
             (\varepsilon a+\varepsilon^2/a),\\
 \int_0^{\ell^2+1}\|(R_q,R_u,R_\tau)\|_2\dd t
       &\le C_\ell(\eta a^3+\eta^2a^{5/2}).
 \end{aligned}
\end{equation}

\paragraph{Relative energy and the choice of thickness.}
The differences of the reference velocity gradient and stress
divergence from the background values are bounded by
$C_\ell\varepsilon/a$. The corresponding difference of the
physical temperature second derivatives is bounded by
$C_\ell\varepsilon$. In particular,
$\int_0^{\ell^2+1}A_{\rm ref}(t)\dd t\le C_\ell$
uniformly as $a\downarrow0$.
Reference density and temperature remain in a fixed
positive interval; conversion to physical residuals
only multiplies \eqref{eq:r9-potential-residual} by
bounded reference coefficients.
Estimate~\eqref{eq:r8-relative-stability} yields
\begin{equation}\label{eq:r9-potential-true-error}
 \sup_{0\le t\le \ell^2+1}\|U(t)-\widetilde U(t)\|_2
       \le C_\ell(\eta a^3+\eta^2a^{5/2})=o(\eta a^2).
\end{equation}

The velocity initial norms have sizes
\[
 \|u_0\|_{H^2}\simeq\eta,\quad
 \|u_0\|_1=O(\eta a^{5/2}\ell),\quad
 \|u_0\|_{L^2_3L^1_h}=O(\eta a^2\ell),\quad
 \|\partial_3u_0\|_{L^2_3L^1_h}=O(\eta a\ell).
\]
The background has $\mathcal J_2$ size $O(\eta)$,
and the thermal correction has $H^2$ size
$O(\eta^2a^2\ell^{-4})$ and smaller $L^1$ size.
Thus $a\le \ell^{-2}$ gives $\mathcal J_{2,1}\simeq\eta$.
The second vertical mixed norm is comparable to $\eta\ell$.
With $\eta=\eta_\ell$, the normalization of the density error in
$L^2_{X,Z}$ gives
\[
 \sup_{0\le s\le1}
 \big\|\eta_\ell^{-2}\ell^5a^{-3/2}
        (q-\widetilde q)(\ell^2+s,\ell X,aZ)\big\|_{2,X,Z}
 \le C_\ell\ell^4
                   \left(\frac{a}{\eta_\ell}+\sqrt a\right).
\]
Here the full change of variables has Jacobian
$\dd x=\ell^2a\,\dd X\dd Z$.
For each fixed $\ell$ choose $a_\ell>0$ so that
\begin{equation}\label{eq:r9-potential-thickness}
 a_\ell\le \ell^{-2},\qquad
       C_\ell(a_\ell+\eta_\ell a_\ell^{1/2})\le\eta_\ell \ell^{-5}.
\end{equation}
The choice is possible because $\ell$, $\eta_\ell$, and $C_\ell$
are fixed while $a_\ell\downarrow0$. It makes the normalized error
at most $\ell^{-1}$; no uniform bound on $C_\ell$ as
$\ell\to\infty$ is required.
Equations \eqref{eq:r9-potential-r-profile} and
\eqref{eq:r9-potential-true-error} then show that
$\eta_\ell^{-2}\ell^5a_\ell^{-3/2}(q_\ell-q_b)(\ell^2+s,\ell X,a_\ell Z)$
converges strongly in $L^2_{X,Z}$ to the right side of
\eqref{eq:r9-potential-r-profile}, uniformly in $s$.

\smallskip\noindent\emph{Step 4. The weak derivative limit.}

For the data \eqref{eq:r9-potential-data},
\eqref{eq:new-data-sizes} and \eqref{eq:r9-potential-thickness} give
$\mathcal J_2\le\mathcal J_{2,1}\le C\eta_\ell$.
The exact elliptic identity is
\[
 (\nu_c+\nu_h\Delta_h\Delta^{-1})\Div u
       =c^2(q+\tau)+(-\Delta)^{-1}(\Div F_2-\Div u_t).
\]
By \eqref{eq:r9-high-time-final} and \eqref{eq:r9-improved-sources},
the high-frequency part of $\partial_3^3w_3$ differs from
\[
 c^2\partial_3^4\Delta^{-1}
       (\nu_c+\nu_h\Delta_h\Delta^{-1})^{-1}P_Hq
\]
by at most $C\mathcal J_2e^{-c_*t}
+C\mathcal J_2^2(1+t)^{-9/4}$ in $L^2$.
The low-frequency bound \eqref{eq:r9-low-third} gives
$\|\partial_3^3(w_L)_3(t)\|_2
 \le C\mathcal J_2(1+t)^{-9/4}$.
On the rescaled variables the multiplier acting on
$\partial_Z^2q$ is
\[
 \frac{c^2\zeta^2}{\zeta^2+(a_\ell/\ell)^2|\xi_h|^2}
 \left(\nu_c+
       \frac{\nu_h(a_\ell/\ell)^2|\xi_h|^2}
              {\zeta^2+(a_\ell/\ell)^2|\xi_h|^2}\right)^{-1}
       \longrightarrow c^2/\nu_c.
\]
It is uniformly bounded, and the fixed high cutoff tends
to one almost everywhere after rescaling.
Move this multiplier and two $Z$ derivatives onto a
Schwartz test function and use the strong density limit.
The normalized background contribution is at most
$C\eta_\ell^{-1}\ell^{-3/2}$ by \eqref{eq:r9-background-high}.
At~$t=\ell^2+s$, $0\le s\le1$, multiplication by
$\eta_\ell^{-2}\ell^4$ bounds the exponential and quadratic
high-frequency errors by
$C\eta_\ell^{-1}\ell^4e^{-c_*\ell^2}$ and $C\ell^{-1/2}$, respectively.
The low-frequency estimate \eqref{eq:r9-low-third} contributes
$C\eta_\ell^{-1}\ell^{-1/2}$.
All three remainders tend to zero for
$\eta_\ell=\delta/\log(e+\ell)$.
The complete normalized third derivative is bounded in
$L^2$ by \eqref{eq:r9-w-third}.
Distributional convergence therefore extends by density
to weak $L^2$ convergence against every fixed test function,
uniformly for $s\in[0,1]$.
The factors $\nu_c/c^2$ cancel and give
\eqref{eq:r9-potential-weak}--\eqref{eq:r9-potential-limit}.
The horizontal rotational Gaussian and the shifted
background Gaussian have a nonzero dot product proportional
to $X_2$ times a positive Gaussian, and $\psi''\ne0$.
Testing with $G/\|G\|_2$ gives the positive lower limit
and the time-integrated bound.
\end{proof}

\paragraph{Fixed-data optimality.}
Consider a small $\mathcal J_{2,1}$ ball in the
closed subspace
$\{(q_0,u_0,\tau_0):(u_0)_3=0,\ \Div_h(u_0)_h=0\}$.
For integers $k,m,N\ge1$, require both derivative norms at an integer $\ell^2\ge N$, multiplied by $(1+\ell^2)^{2+1/k}$,
to exceed $m$.
The initial-data norms tend to zero, and their normalized
lower bound is
$c_*\delta^2\ell^{2/k}/[\log(e+\ell)]^2\to\infty$.
The countable-intersection argument following
Corollary~\ref{cor:r10-translation} gives one fixed datum satisfying
all these lower-bound conditions, hence the asserted failure for every $\beta>2$.
The linearized solution satisfies $\|\nabla^3w^{\rm lin}(t)\|_2=O(t^{-9/4})$
by \eqref{eq:blockmultiplier} and Appendix~\ref{app:r13-green-matrix}.
At high frequencies, the slow mode is bounded by
$Ce^{-c_*t}(\varrho^2|\widehat q_0|+\varrho|\widehat d_0|+|\widehat\tau_0|)$,
while parabolic smoothing supplies the remaining derivative of the
fast modes; hence only $U_0\in H^2$ is required.
The upper bound \eqref{eq:r9-w-third} and its
fixed-data $o(t^{-2})$ refinement remain valid.
This proves Theorem~\ref{thm:decay-H2}(iv)(b).

\section{The exact relative-energy identity}\label{app:r10-relative}

Let the solution and reference fields satisfy the assumptions of
Lemma~\ref{lem:r10-relative}, and retain its definitions of the
residuals, ratios, and relative potentials. For smooth fields, the
three identities below sum to \eqref{eq:r8-relative-identity}; the
regularization argument in that lemma extends the calculation to
the stated energy class. Write
$D_u=\partial_t+u\cdot\nabla$ and
$D_V=\partial_t+V\cdot\nabla$ in this appendix.
The continuity equation for the solution implies
$\frac{\dd}{\dd t}\int\rho F=\int\rho D_uF$
whenever the relative integrals are finite.
The reference equations give
\[
 D_VV=r^{-1}(\Div\mathbb S^{\rm ed}(V)-\nabla(Rr\Theta)+f_m),
 \quad
 D_V\Theta=\frac{Q(V)+\kappa\Delta\Theta+f_\theta}{rc_v}
                         -\frac{R}{c_v}\Theta\Div V .
\]

\subsection{Kinetic relative energy}

Test the difference of the two material accelerations with
$\rho z$. The viscous part is
$-\int Q(z)-\int(b-1)z\cdot\Div\mathbb S^{\rm ed}(V)$.
For the pressure, integration by parts gives the exact relation
\[
 \int \bigl(R\rho\theta\Div z+bz\cdot\nabla(Rr\Theta)\bigr)\dd x
 =\int\bigl(R\rho(\theta-\Theta)\Div z
                         -R\rho\Theta z\cdot\nabla\log b\bigr)\dd x.
\]
Consequently,
\begin{equation}\label{eq:r10-relative-kinetic}
 \begin{aligned}
 \frac{\dd}{\dd t}\int\tfrac12\rho|z|^2\dd x+\int Q(z)\dd x
 =\int\bigl[&
 -\rho z\otimes z:\nabla V
 -(b-1)z\cdot\Div\mathbb S^{\rm ed}(V)\\
 &+R\rho(\theta-\Theta)\Div z
 -R\rho\Theta z\cdot\nabla\log b-bz\cdot f_m\bigr]\dd x .
 \end{aligned}
\end{equation}

\subsection{Density relative potential}

Since $\partial_\rho H_\rho=\log b$ and
$\partial_rH_\rho=1-b$, use both mass equations to obtain
\[
 \frac{\dd}{\dd t}\int R\Theta H_\rho\dd x
 =R\int\bigl[
 H_\rho D_V\Theta+\rho\log b\,z\cdot\nabla\Theta
       +\rho\Theta z\cdot\nabla\log b
       +\Theta(1-b)f_\rho\bigr]\dd x .
\]
Substituting the reference heat equation gives
\begin{equation}\label{eq:r10-relative-density}
 \begin{aligned}
 \frac{\dd}{\dd t}\int R\Theta H_\rho\dd x
 =\int\biggl[&
 \frac{RH_\rho}{rc_v}(Q(V)+\kappa\Delta\Theta+f_\theta)
       -\frac{R^2\Theta}{c_v}H_\rho\Div V\\
 &+R\rho\log b\,z\cdot\nabla\Theta
       +R\rho\Theta z\cdot\nabla\log b
       +R\Theta(1-b)f_\rho\biggr]\dd x .
 \end{aligned}
\end{equation}
The term with $\nabla\log b$ cancels the kinetic pressure
term. In particular no reference-density-gradient norm
remains in the final stability coefficient.

\subsection{Temperature relative potential}

Here $\partial_\theta K_\theta=1-d^{-1}$ and
$\partial_\Theta K_\theta=-\log d$.
Use $D_u\Theta=D_V\Theta+z\cdot\nabla\Theta$ and both heat
equations. Separate $\Div u=\Div V+\Div z$ to obtain
\begin{equation}\label{eq:r10-relative-temperature}
 \begin{aligned}
 \frac{\dd}{\dd t}\int c_v\rho K_\theta\dd x
 =\int\bigl[&
 (1-d^{-1})Q(u)-b\log d\,Q(V)\\
 &+\kappa(1-d^{-1})\Delta\theta-\kappa b\log d\,\Delta\Theta
 -R\rho(\theta-\Theta)\Div z-R\rho K_\theta\Div V\\
 &-c_v\rho\log d\,z\cdot\nabla\Theta-b\log d\,f_\theta\bigr]\dd x .
\end{aligned}
\end{equation}
The first divergence term cancels its counterpart in
\eqref{eq:r10-relative-kinetic}.
For the heat terms, write $\theta=d\Theta$ and integrate by parts:
\[
 \int(1-d^{-1})\Delta\theta\dd x
 =-\int\Theta|\nabla\log d|^2\dd x
                    -\int\nabla\log d\cdot\nabla\Theta\dd x .
\]
The part $-\kappa\int\log d\,\Delta\Theta$ cancels the last
cross term. The remaining coefficient is therefore exactly
$-\kappa(b-1)\log d\,\Delta\Theta$.

\subsection{Relative-energy identity for the eddy stress}

The constitutive tensor gives
$\mathbb S^{\rm ed}(V):\nabla z
=\mu\nabla_hV:\nabla_hz+(\lambda+\mu)\Div V\Div z$.
Its quadratic expansion is valid even though the stress
matrix is not symmetric.
Adding \eqref{eq:r10-relative-kinetic},
\eqref{eq:r10-relative-density}, and
\eqref{eq:r10-relative-temperature} gives
\begin{align}\label{eq:r8-relative-identity}
 E_{\rm rel}'
 &+\int[d^{-1}Q(z)+\kappa\Theta|\nabla\log d|^2]\dd x\notag\\
 =\int\big[&
 -\rho z\otimes z:\nabla V
 -(b-1)z\cdot\Div\mathbb S^{\rm ed}(V)
 -\rho(c_v\log d-R\log b)z\cdot\nabla\Theta\notag\\
 &-\left(R\rho K_\theta+\frac{R^2\Theta}{c_v}H_\rho\right)\Div V
 +\frac{RH_\rho}{rc_v}(Q(V)+\kappa\Delta\Theta)\notag\\
 &-(d^{-1}-1+b\log d)Q(V)
 -2(d^{-1}-1)\mathbb S^{\rm ed}(V):\nabla z
 -\kappa(b-1)\log d\,\Delta\Theta\notag\\
 &+R\Theta(1-b)f_\rho-bz\cdot f_m
       +\frac{RH_\rho-c_v\rho\log d}{rc_v}f_\theta\big]\dd x .
\end{align}
The residual coefficient in the last term is the sum of
$RH_\rho/(rc_v)$ from the density potential and
$-b\log d$ from the temperature potential.
All terms on the right are either quadratic relative
remainders or residuals paired with a relative difference.
This is the identity used in Lemma~\ref{lem:r10-relative}.

\section{Nonlinear source estimates}\label{app:r10-sources}

\subsection{A time-integrated estimate for the projected momentum source}

The following tail bound for the projected momentum source follows
from the basic decay estimates. It supplies the source term for
enhanced mixed third-order decay and the corresponding Gaussian
profile.

\Needspace{7\baselineskip}
\begin{lemma}[Tail bound for the projected momentum source]
\label{lem:r9-source-tail}
Let $U$ be the global solution with $\mathcal J_2(U_0)$ sufficiently
small, and set
\[
 \mathcal R_{\rm sh}=(e^q-1)u_t+e^q(u\cdot\nabla u).
\]
There is a constant $C$, depending only on the fixed coefficients
and equilibrium, such that
\begin{equation}\label{eq:r9-source-tail}
 \|\nabla_h\partial_3\mathcal R_{\rm sh}\|_{L^2(t,\infty;L^2)}
       \le C\mathcal J_2^2(1+t)^{-3/2},\qquad t\ge1.
\end{equation}
For each fixed datum the left side is $o(t^{-3/2})$.
\end{lemma}
\begin{proof}
Write $Z=\partial_3u$, $V=u_t$.
For $(e^q-1)V$, expand one horizontal and one vertical derivative.
The highest $V$ term is bounded by
$\|e^q-1\|_\infty\|\nabla_h\partial_3V\|_2$,
and the highest coefficient term by
$\|\nabla^2e^q\|_2\|V\|_\infty$.
The two cross terms use \eqref{eq:r9-trace}.
Equation~\eqref{eq:r9-time-L2} gives
$\|V(t)\|_2\le C\mathcal J_2(1+t)^{-5/4}$.
Together with \eqref{eq:r9-V-point}, \eqref{eq:r9-Y-point}, and
$\|f\|_\infty\le C(\|f\|_{H^1}+\|\nabla_hf\|_{H^1})$, this yields
$\|V\|_{L^2(b,b+1;L^\infty)}
\le C\mathcal J_2(1+b)^{-5/4}$ for $b\ge1$.
The decay of $q$ and its first two derivatives therefore gives
\[
 \|\nabla_h\partial_3[(e^q-1)V]\|_{L^2(b,b+1;L^2)}
       \le C\mathcal J_2^2(1+b)^{-11/4}.
\]
Summing the squared estimates over the intervals
$[t+n,t+n+1]$, $n\ge0$, gives
\[
 \|\nabla_h\partial_3[(e^q-1)V]\|_{L^2(t,\infty;L^2)}
       \le C\mathcal J_2^2(1+t)^{-9/4}.
\]
No $\partial_3^2V$ is used.

For $g=u\cdot\nabla u$, the two critical products in
$\nabla_h\partial_3g$ are
\[
       (\nabla_hu_3)\partial_3^2u,\qquad
       u_3\nabla_h\partial_3^2u.
\]
By \eqref{eq:r9-compressible-lips}, the first has pointwise
power $2$ times the bounded energy norm of $\partial_3^2u$.
The second is bounded in the time-integral tail by
$C\mathcal J_2(1+t)^{-3/2}
\|\nabla_h\partial_3^2u\|_{L^2(t,\infty;L^2)}$.
The remaining pairings all have pointwise exponent $17/8$
or better:
\[
 \begin{array}{c|c}
 \|\nabla_hZ\|_2\|\nabla_hu\|_\infty&3/4+11/8\\
 \|Z\|_3\|\nabla_h^2u\|_6&3/8+7/4\\
 \|u_h\|_\infty\|\nabla_h^2Z\|_2&7/8+5/4\\
 \|\nabla_h\partial_3u_3\|_6\|Z\|_3&7/4+3/8\\
 \|\partial_3u_3\|_\infty\|\nabla_hZ\|_2&11/8+3/4 .
 \end{array}
\]
The two $L^6$ estimates use directional Sobolev and
\[
 \nabla_h^2\partial_3u_3=\nabla_h^2\Div w-\Div_h\nabla_h^2u_h,
 \qquad
 \nabla_h\partial_3^2u_3
       =\nabla_h\partial_3\Div w-\Div_h\nabla_hZ .
\]
All derivatives on the right are already controlled.
The two critical products and the remaining pointwise bounds give
\[
 \begin{aligned}
 \|\nabla_h\partial_3g\|_{L^2(t,\infty;L^2)}
 &\le C\mathcal J_2(1+t)^{-3/2}
       \sup_{s\ge t}\|\partial_3^2u(s)\|_2\\
 &\quad+C\mathcal J_2(1+t)^{-3/2}
       \|\nabla_h\partial_3^2u\|_{L^2(t,\infty;L^2)}\\
 &\quad+C\mathcal J_2^2(1+t)^{-13/8}\\
 &\le C\mathcal J_2^2(1+t)^{-3/2}.
 \end{aligned}
\]

For the coefficient $e^q$, direct expansion also gives
\[
 \begin{aligned}
 \|g(t)\|_2&\le C\mathcal J_2^2(1+t)^{-7/4},
 \qquad \|\partial_3g(t)\|_2\le C\mathcal J_2^2(1+t)^{-3/2},\\
 \|\nabla_hg(t)\|_2+\|\nabla_h^2g(t)\|_2
       &\le C\mathcal J_2^2(1+t)^{-2}.
\end{aligned}
\]
Apply the same anisotropic $L^\infty$ inequality and integrate
these bounds and the preceding mixed-derivative estimate:
\[
 \|g\|_{L^2(t,\infty;L^\infty)}
       \le C\mathcal J_2^2(1+t)^{-1}.
\]
The highest coefficient term consequently satisfies
\[
 \|[\nabla_h\partial_3(e^q)]g\|_{L^2(t,\infty;L^2)}
       \le C\mathcal J_2^3(1+t)^{-11/4},
\]
using $\|\nabla^2e^q(t)\|_2\le
C\mathcal J_2(1+t)^{-7/4}$.
For the first cross term, \eqref{eq:r9-trace} and time H\"older give
\[
 \begin{aligned}
 \|(\nabla_he^q)\partial_3g\|_2
 &\le C\|\nabla^2e^q\|_2\|\partial_3g\|_2^{1/2}
                                  \|\nabla_h\partial_3g\|_2^{1/2},\\
 \|(\nabla_he^q)\partial_3g\|_{L^2(t,\infty;L^2)}
 &\le C\mathcal J_2^2\left(\int_t^\infty(1+s)^{-10}\dd s\right)^{1/4}
               \|\nabla_h\partial_3g\|_{L^2(t,\infty;L^2)}^{1/2}
 \le C\mathcal J_2^3(1+t)^{-3}.
 \end{aligned}
\]
For the second, the same trace estimate gives
$\|(\partial_3e^q)\nabla_hg\|_2
\le C\|\nabla^2e^q\|_2\|\nabla_hg\|_2^{1/2}
\|\nabla_h^2g\|_2^{1/2}\le C\mathcal J_2^3(1+t)^{-15/4}$;
its $L^2(t,\infty;L^2)$ norm is at most
$C\mathcal J_2^3(1+t)^{-13/4}$.
Together with the bound for $\nabla_h\partial_3g$ and the
boundedness of $e^q$, these estimates prove
\eqref{eq:r9-source-tail}.
For fixed data, $\sup_{s\ge t}\|\partial_3^2u(s)\|_2\to0$
and the total-dissipation tail tends to zero.
The two critical products therefore vanish after normalization
by the stated decay rate; all other terms have strictly faster algebraic decay.
\end{proof}

\subsection{An index of the source estimates}\label{sec:source-estimate-index}

The table gives the first proof and applications of each source bound.
All bounds have coefficient $C\mathcal J_2^2$ and hold for $t\ge1$;
powers refer to norms, and intervals indicate time $L^2$ norms.

\begingroup\small
\begin{longtable}{@{}>{\RaggedRight\arraybackslash}p{.34\linewidth}>{\centering\arraybackslash}p{.12\linewidth}>{\centering\arraybackslash}p{.09\linewidth}>{\RaggedRight\arraybackslash}p{\dimexpr.45\linewidth-6\tabcolsep\relax}@{}}
\toprule
Spatial norm & Time norm & Power & First proof; applications\\ \midrule
\endfirsthead
\toprule
Spatial norm & Time norm & Power & First proof; applications\\ \midrule
\endhead
\bottomrule
\endfoot
\bottomrule
\endlastfoot
$F_2$ in $L^2_3L^1_h,\ L^2$ &
Pointwise & $5/4,\ 7/4$ & Lemma~\ref{lem:r10-momentum-source};
(Proposition~\ref{prop:r8-H1-scattering}).\\
$F_2$ in $H^1$ &
$L^2(t,t+1)$ & $3/2$ & \eqref{eq:r8-F2-H1-window};
(Propositions~\ref{prop:r8-windows} and~\ref{prop:r8-vertical-correction}).\\
$F_2$ in $H^1$ &
Pointwise & $3/2$ & \eqref{eq:r9-F2-H1};
(Proposition~\ref{prop:r9-mixed-third}).\\
$\nabla_hF_2$ in $L^2$ &
$L^2(t,t+1)$ & $49/24$ & \eqref{eq:r8-horizontal-sources};
(Proposition~\ref{prop:r8-windows}).\\
$\nabla_hF_2,\Div F_2$ in $L^2$ &
Pointwise & $49/24$ & \eqref{eq:r9-point-sources-pre};
(Proposition~\ref{prop:r9-horizontal-third}).\\
$\nabla F_1,\nabla_hF_2,\Div F_2$ in $L^2$ &
Pointwise & $9/4$ & \eqref{eq:r9-improved-sources};
(\eqref{eq:r9-high-time-final} and~\eqref{eq:r9-thermal-profile}).\\
$F_1,F_3$ in $H^1$ &
Pointwise & $17/8$ & Lemma~\ref{lem:r10-scalar-sources};
(\eqref{eq:r11-integrated-physical-source}).\\
$F_3$ in $H^1$ &
$L^2(t,t+1)$ & $17/8$ & \eqref{eq:r9-source-thermal-old};
(Lemma~\ref{lem:r9-time-derivatives}).\\
$F_{2,t}$ in $L^2$ &
$L^2(t,t+1)$ & $9/4$ & \eqref{eq:r9-F2t-improved};
(\eqref{eq:r9-high-time-final}).\\
$N_0,B_0$ in $L^1$ &
Pointwise & $3/2,\ 5/4$ & \eqref{eq:r9-early-flux};
(\eqref{eq:r9-low-third} and Lemma~\ref{lem:new-quadratic-error}).\\
$N_1,B_1$ in $L^2$ &
$L^2(t,t+1)$ & $21/8,\ 5/2$ & \eqref{eq:r9-late-flux};
(\eqref{eq:r9-low-third}).\\
$\nabla_h\partial_3\mathcal R_{\rm sh}$ in $L^2$ &
$L^2(t,\infty)$ & $3/2$ & Lemma~\ref{lem:r9-source-tail};
(Lemma~\ref{lem:r13-enhanced-third} and
Section~\ref{sec:r13-mixed-profile}).\\
$F_2$ in $H^1$ &
$L^2(t,\infty)$ & $1$ & \eqref{eq:r13-effective-source-tails};
(\eqref{eq:r8-H1-scattering} and
Proposition~\ref{prop:r8-vertical-profile}).\\
\end{longtable}
\endgroup

The $49/24$ and $9/4$ bounds are used at successive stages.
Here $F_2=N_0+\partial_3B_0$ and
$F_{2,t}=N_1+\partial_3B_1$ are the decompositions in
\eqref{eq:r9-early-flux} and \eqref{eq:r9-late-flux}.
The projected source and horizontal viscous flux are estimated
in \eqref{eq:new-profile-source-ledger}.

\section{Refined compressible Green-matrix estimates}
\label{app:r13-green-matrix}

\begin{proof}[Proof of Lemma~\ref{lem:refined-compressible-green}]
At low frequencies, \eqref{eq:blockmultiplier} gives the required
decay. At high frequencies, we separate the slow density mode from
the two fast modes of \eqref{eq:fourierblock} to obtain the bound
in \eqref{eq:refined-second-Duhamel} with only one derivative on $F_2$.

\smallskip\noindent\emph{Step 1. High-frequency spectral decomposition.}
Write \(A(\xi)=\mathfrak a(\xi)\varrho^2\), where
\begin{equation*}
 \nu_c\leq\mathfrak a(\xi)\leq\nu_c+\nu_h,
 \qquad \varrho=|\xi|.
\end{equation*}
The matrix in \eqref{eq:fourierblock} is
\begin{equation*}
 B(\varrho,\mathfrak a)=
 \begin{pmatrix}
 0&\varrho&0\\
 -c^2\varrho&\mathfrak a\varrho^2&-c^2\varrho\\
 0&\alpha\varrho&\chi\varrho^2
 \end{pmatrix},
\end{equation*}
and
\begin{equation*}
p(z;\varrho,\mathfrak a):=\det(zI+B)
 =z^3+(\mathfrak a+\chi)\varrho^2z^2+
 \{\mathfrak a\chi\varrho^4+c^2(1+\alpha)\varrho^2\}z
 +c^2\chi\varrho^4.
\end{equation*}
Set \(y=\varrho^{-2}\).  Division by \(\varrho^4\) gives
\begin{equation*}
 y^2z^3+(\mathfrak a+\chi)yz^2+
 \{\mathfrak a\chi+c^2(1+\alpha)y\}z+c^2\chi=0.
\end{equation*}
The implicit-function theorem, uniformly for
\(\mathfrak a\in[\nu_c,\nu_c+\nu_h]\), gives constants
\(R_0,c_*,C_*>0\) and one root \(-\lambda_0\) such that, for
\(\varrho\geq R_0\),
\begin{equation}\label{eq:slow-root-bounds}
 c_*\leq\Re\lambda_0\leq|\lambda_0|\leq C_*,\qquad
 \left|\lambda_0-\frac{c^2}{\mathfrak a}\right|
 \leq C\varrho^{-2}.
\end{equation}
The other two roots \(-\lambda_\pm\) satisfy
\begin{equation*}
 \Re\lambda_\pm\geq c_*\varrho^2,\qquad
 |\lambda_\pm|\leq C_*\varrho^2.
\end{equation*}

\smallskip\noindent\emph{Step 2. Bounds for the Green matrix.}
The contour representation below remains uniform when the two
fast eigenvalues coincide. Choose
\(0<d_0<\frac14\min\{\nu_c,\chi,c^2/(\nu_c+\nu_h)\}\)
and a fixed positively oriented
contour \(\Gamma_f\subset\{\Re\zeta\leq-d_0\}\) enclosing
\[
 \{-\mathfrak a:\mathfrak a\in[\nu_c,\nu_c+\nu_h]\}
 \cup\{-\chi\}
\]
at positive distance, but not enclosing \(0\).  After
\(z=\varrho^2\zeta\),
\begin{equation*}
 \varrho^{-6}p(\varrho^2\zeta;\varrho,\mathfrak a)
 =\zeta(\zeta+\mathfrak a)(\zeta+\chi)+O(\varrho^{-2})
\end{equation*}
uniformly on \(\Gamma_f\).  Hence, after increasing \(R_0\),
\begin{equation*}
 |p(z;\varrho,\mathfrak a)|\geq c_*\varrho^4
 \quad\text{on }|z+c^2/\mathfrak a|=d_0,
 \qquad
 |p(\varrho^2\zeta;\varrho,\mathfrak a)|\geq c_*\varrho^6
 \quad\text{on }\Gamma_f.
\end{equation*}
The first contour encloses the slow root, and the scaled contour
encloses both fast roots. The resulting resolvent bounds remain
uniform at \(\mathfrak a=\chi\).

Let \(K(t,\varrho,\mathfrak a)=e^{-tB(\varrho,\mathfrak a)}\).  Direct
cofactor expansion gives
\begin{equation}\label{eq:resolvent-adjugate}
\begin{aligned}
 &(zI+B)^{-1}\\
 &=p(z;\varrho,\mathfrak a)^{-1}
 \begin{pmatrix}
 (z+\mathfrak a\varrho^2)(z+\chi\varrho^2)+\alpha c^2\varrho^2
 &-\varrho(z+\chi\varrho^2)&-c^2\varrho^2\\
 c^2\varrho(z+\chi\varrho^2)&z(z+\chi\varrho^2)&c^2\varrho z\\
 -\alpha c^2\varrho^2&-\alpha\varrho z
 &z(z+\mathfrak a\varrho^2)+c^2\varrho^2
 \end{pmatrix}.
\end{aligned}
\end{equation}
The slow root is simple, and
\[
 P_0=\frac{\operatorname{adj}(-\lambda_0I+B)}
              {\pa_zp(-\lambda_0;\varrho,\mathfrak a)},\qquad
 |\pa_zp(-\lambda_0;\varrho,\mathfrak a)|\geq c_*\varrho^4.
\]
The only additional cancellation needed in its adjugate is the
\((3,3)\) entry:
\[
 (-\lambda_0)(-\lambda_0+\mathfrak a\varrho^2)+c^2\varrho^2
 =\varrho^2(c^2-\mathfrak a\lambda_0)+\lambda_0^2=O(1)
\]
by \eqref{eq:slow-root-bounds}. The other cofactor bounds follow
entrywise from \eqref{eq:resolvent-adjugate}. Consequently,
\begin{equation*}
 |P_0(\varrho,\mathfrak a)|\leq C
 \begin{pmatrix}
 1&\varrho^{-1}&\varrho^{-2}\\
 \varrho^{-1}&\varrho^{-2}&\varrho^{-3}\\
 \varrho^{-2}&\varrho^{-3}&\varrho^{-4}
 \end{pmatrix}.
\end{equation*}
\Needspace{8\baselineskip}
\smallskip\noindent\emph{The fast modes.}
On the fast contour, \(dz=\varrho^2\,d\zeta\), and the entries of the
adjugate in \eqref{eq:resolvent-adjugate} have respective orders
\[
 \begin{pmatrix}
 O(\varrho^4)&O(\varrho^3)&O(\varrho^2)\\
 O(\varrho^3)&O(\varrho^4)&O(\varrho^3)\\
 O(\varrho^2)&O(\varrho^3)&O(\varrho^4)
 \end{pmatrix}.
\]
The Dunford formula and \(\Re\zeta\leq-d_0\) therefore imply
\begin{equation*}
 |K(t,\varrho,\mathfrak a)(I-P_0)|\leq Ce^{-c_*\varrho^2t}
 \begin{pmatrix}
 1&\varrho^{-1}&1\\
 \varrho^{-1}&1&\varrho^{-1}\\
 1&\varrho^{-1}&1
 \end{pmatrix}.
\end{equation*}
Consequently,
\begin{equation}\label{eq:component-Green-matrix}
 |K(t,\varrho,\mathfrak a)|\leq Ce^{-c_*t}
 \begin{pmatrix}
 1&\varrho^{-1}&\varrho^{-2}\\
 \varrho^{-1}&\varrho^{-2}&\varrho^{-3}\\
 \varrho^{-2}&\varrho^{-3}&\varrho^{-4}
 \end{pmatrix}
 +Ce^{-c_*\varrho^2t}
 \begin{pmatrix}
 1&\varrho^{-1}&1\\
 \varrho^{-1}&1&\varrho^{-1}\\
 1&\varrho^{-1}&1
 \end{pmatrix},\qquad \varrho\geq R_0.
\end{equation}
The constants in \eqref{eq:component-Green-matrix} are uniform in
\(\mathfrak a\in[\nu_c,\nu_c+\nu_h]\).

\smallskip\noindent\emph{Step 3. Source estimates and reconstruction.}
The Fourier source in the middle component equals
\begin{equation*}
 \widetilde F_2=|\xi|^{-1}i\xi\cdot\widehat F_2,\qquad
 |\xi|^m|\widetilde F_2|\leq|\xi|^m|\widehat F_2|,\quad m\geq0.
\end{equation*}
To estimate two derivatives of the physical solution, multiply each
row of \eqref{eq:component-Green-matrix} by \(\varrho^2\).
The second row controls \(\nabla\Div\mathbb Qu_F\), because
\(\widehat{\nabla\Div\mathbb Qu_F}=i\xi|\xi|\widehat d_F\).
Using \eqref{eq:parabolic-half-gain} for the fast-mode contribution
of $F_2$ yields
\begin{align}
 &\varrho^2\{|\widehat{P_Hq_F}|+|\widehat{P_H\tau_F}|+|\widehat{P_Hd_F}|\}
 \notag\\*
 &\quad\leq C\int_0^t e^{-c_*(t-r)}
 \{\langle\xi\rangle^2(|\widehat F_1|+|\widehat F_3|)
   +\langle\xi\rangle|\widehat F_2|\}(r)\,\dd r
 \notag\\*
 &\qquad+C\int_0^t(t-r)^{-1/2}e^{-c_*(t-r)}
 \langle\xi\rangle|\widehat F_2(r)|\,\dd r,
 \label{eq:high-frequency-source-bound}
\end{align}
where
\begin{equation}\label{eq:parabolic-half-gain}
 \sup_{\varrho\geq R_0}\varrho e^{-c_*\varrho^2s}
 \leq Cs^{-1/2}e^{-c_* s/2},\qquad s>0.
\end{equation}
Plancherel and Minkowski applied to
\eqref{eq:high-frequency-source-bound} prove the high-frequency part of
\eqref{eq:refined-second-Duhamel}.  On \(\operatorname{supp}\chi_\ell\),
\eqref{eq:blockmultiplier} gives
\begin{equation}\label{eq:low-frequency-second-source}
 \|\nabla^2P_L(q_F,\mathbb Qu_F,\tau_F)(t)\|_2
 \leq C\int_0^t(1+t-r)^{-7/4}\|F(r)\|_1\,\dd r.
\end{equation}
Since \(\widehat{\nabla\Div\mathbb Qu_F}=i\xi|\xi|\widehat d_F\),
\eqref{eq:low-frequency-second-source} also controls the low-frequency
part of \(\nabla\Div\mathbb Qu_F\). Combining
\eqref{eq:high-frequency-source-bound} and
\eqref{eq:low-frequency-second-source} proves
\eqref{eq:refined-second-Duhamel}. Finally, applying
\eqref{eq:ScDuhamel} with $k=j=0,1$ proves
\eqref{eq:refined-low-order-Duhamel}.
\end{proof}

\medskip
\noindent\textbf{Declaration of AI usage.}
ChatGPT (OpenAI) was used during the preparation of this manuscript
to assist with proof development, literature searches, and language editing. The authors are responsible for the mathematical
arguments, references, and final content of the manuscript.


\begin{thebibliography}{99}

\bibitem{BahouriCheminDanchin}
H. Bahouri, J.-Y. Chemin, and R. Danchin,
\emph{Fourier Analysis and Nonlinear Partial Differential Equations},
Grundlehren Math. Wiss., vol.~343, Springer, Heidelberg, 2011,
\href{https://doi.org/10.1007/978-3-642-16830-7}
{doi:10.1007/978-3-642-16830-7}.

\bibitem{BaiGuoXu}
H. Bai, Z. Guo, and L. Xu,
\emph{Global well-posedness of classical solutions to the isentropic
compressible Navier--Stokes equations with anisotropic viscous-stress tensor},
Acta Math. Sinica (Chinese Ser.) \textbf{69} (2026), 185--226,
\href{https://doi.org/10.12386/A20240059}{doi:10.12386/A20240059}.

\bibitem{BieFangWangZhou}
Q. Bie, H. Fang, S. Wang, and Y. Zhou,
\emph{Stability and sharp decay for the 3D incompressible anisotropic
Navier--Stokes equations with fractional horizontal dissipation},
J. Differential Equations \textbf{463} (2026), 114167,
\href{https://doi.org/10.1016/j.jde.2026.114167}
{doi:10.1016/j.jde.2026.114167}.

\bibitem{BreschJabin}
D. Bresch and P.-E. Jabin,
\emph{Global existence of weak solutions for compressible Navier--Stokes
equations: thermodynamically unstable pressure and anisotropic viscous stress
tensor},
Ann. of Math. (2) \textbf{188} (2018), 577--684.

\bibitem{CheminAnisotropic2000}
J.-Y. Chemin, B. Desjardins, I. Gallagher, and E. Grenier,
\emph{Fluids with anisotropic viscosity},
ESAIM Math. Model. Numer. Anal. \textbf{34} (2000), 315--335,
\href{https://doi.org/10.1051/m2an:2000143}{doi:10.1051/m2an:2000143}.

\bibitem{CheminZhang}
J.-Y. Chemin and P. Zhang,
\emph{On the global wellposedness to the 3-D incompressible anisotropic
Navier--Stokes equations},
Comm. Math. Phys. \textbf{272} (2007), 529--566,
\href{https://doi.org/10.1007/s00220-007-0236-0}
{doi:10.1007/s00220-007-0236-0}.

\bibitem{ChikamiDanchin}
N. Chikami and R. Danchin,
\emph{On the well-posedness of the full compressible Navier--Stokes system in
critical Besov spaces},
J. Differential Equations \textbf{258} (2015), 3435--3467,
\href{https://doi.org/10.1016/j.jde.2015.01.012}
{doi:10.1016/j.jde.2015.01.012}.

\bibitem{Danchin}
R. Danchin,
\emph{Global existence in critical spaces for compressible Navier--Stokes
equations},
Invent. Math. \textbf{141} (2000), 579--614,
\href{https://doi.org/10.1007/s002220000078}{doi:10.1007/s002220000078}.

\bibitem{DanchinFull}
R. Danchin,
\emph{Global existence in critical spaces for flows of compressible viscous
and heat-conductive gases},
Arch. Ration. Mech. Anal. \textbf{160} (2001), 1--39,
\href{https://doi.org/10.1007/s002050100155}{doi:10.1007/s002050100155}.

\bibitem{DanchinXuBarotropic}
R. Danchin and J. Xu,
\emph{Optimal time-decay estimates for the compressible Navier--Stokes
equations in the critical \(L^p\) framework},
Arch. Ration. Mech. Anal. \textbf{224} (2017), 53--90,
\href{https://doi.org/10.1007/s00205-016-1067-y}
{doi:10.1007/s00205-016-1067-y}.

\bibitem{DanchinXuFull}
R. Danchin and J. Xu,
\emph{Optimal decay estimates in the critical \(L^p\) framework for flows of
compressible viscous and heat-conductive gases},
J. Math. Fluid Mech. \textbf{20} (2018), 1641--1665,
\href{https://doi.org/10.1007/s00021-018-0381-6}
{doi:10.1007/s00021-018-0381-6}.

\bibitem{Deckelnick}
K. Deckelnick,
\emph{Decay estimates for the compressible Navier--Stokes equations in
unbounded domains},
Math. Z. \textbf{209} (1992), 115--130,
\href{https://doi.org/10.1007/BF02570825}{doi:10.1007/BF02570825}.

\bibitem{FeireislNovotny}
E. Feireisl and A. Novotn\'y,
\emph{Singular Limits in Thermodynamics of Viscous Fluids},
Advances in Mathematical Fluid Mechanics, Birkh\"auser, Basel, 2009.

\bibitem{FeireislNovotnyRelative}
E. Feireisl and A. Novotn\'y,
\emph{Weak-strong uniqueness property for the full
Navier--Stokes--Fourier system},
Arch. Ration. Mech. Anal. \textbf{204} (2012), 683--706,
\href{https://doi.org/10.1007/s00205-011-0490-3}
{doi:10.1007/s00205-011-0490-3}.

\bibitem{FengHongWangWu}
Z. Feng, G. Hong, Y. Wang, and J. Wu,
\emph{Global stability and anisotropic large-time behavior of the
three-dimensional compressible Navier--Stokes equations with eddy diffusion},
\href{https://arxiv.org/abs/2607.06264}{arXiv:2607.06264}, 2026.

\bibitem{GaoPengWuYao}
J. Gao, L. Peng, J. Wu, and Z.-A. Yao,
\emph{Global uniform regularity and vanishing vertical viscosity limit for the
compressible Navier--Stokes equations in the half-space},
\href{https://arxiv.org/abs/2605.23233}{arXiv:2605.23233}, 2026.

\bibitem{GuoWang}
Y. Guo and Y. Wang,
\emph{Decay of dissipative equations and negative Sobolev spaces},
Comm. Partial Differential Equations \textbf{37} (2012), 2165--2208.

\bibitem{Iftimie}
D. Iftimie,
\emph{The resolution of the Navier--Stokes equations in anisotropic spaces},
Rev. Mat. Iberoam. \textbf{15} (1999), 1--36.

\bibitem{JiLuoJiang}
R. Ji, W. Luo, and L. Jiang,
\emph{Stability of the 3D incompressible Navier--Stokes equations with
fractional horizontal dissipation},
Appl. Math. Comput. \textbf{448} (2023), 127934,
\href{https://doi.org/10.1016/j.amc.2023.127934}
{doi:10.1016/j.amc.2023.127934}.

\bibitem{JiWuYang}
R. Ji, J. Wu, and W. Yang,
\emph{Stability and optimal decay for the 3D Navier--Stokes equations with
horizontal dissipation},
J. Differential Equations \textbf{290} (2021), 57--77,
\href{https://doi.org/10.1016/j.jde.2021.04.026}
{doi:10.1016/j.jde.2021.04.026}.

\bibitem{KatoPonce}
T. Kato and G. Ponce,
\emph{Commutator estimates and the Euler and Navier--Stokes equations},
Comm. Pure Appl. Math. \textbf{41} (1988), 891--907,
\href{https://doi.org/10.1002/cpa.3160410704}
{doi:10.1002/cpa.3160410704}.

\bibitem{KawashimaThesis}
S. Kawashima,
\emph{Systems of a hyperbolic--parabolic composite type, with applications
to the equations of magnetohydrodynamics},
Doctoral dissertation, Kyoto University, 1984,
\href{https://doi.org/10.14989/doctor.k3193}{doi:10.14989/doctor.k3193}.

\bibitem{KobayashiShibata}
T. Kobayashi and Y. Shibata,
\emph{Decay estimates of solutions for the equations of motion of compressible
viscous and heat-conductive gases in an exterior domain in \(\mathbb R^3\)},
Comm. Math. Phys. \textbf{200} (1999), 621--659,
\href{https://doi.org/10.1007/s002200050543}{doi:10.1007/s002200050543}.

\bibitem{MatsumuraNishida}
A. Matsumura and T. Nishida,
\emph{The initial value problem for the equations of motion of viscous and
heat-conductive gases},
J. Math. Kyoto Univ. \textbf{20} (1980), 67--104.

\bibitem{ShizutaKawashima}
Y. Shizuta and S. Kawashima,
\emph{Systems of equations of hyperbolic--parabolic type with applications to
the discrete Boltzmann equation},
Hokkaido Math. J. \textbf{14} (1985), 249--275.

\bibitem{TemamZiane}
R. Temam and M. Ziane,
\emph{Some mathematical problems in geophysical fluid dynamics},
in: Handbook of Mathematical Fluid Dynamics, Vol.~III,
North-Holland, Amsterdam, 2004, pp.~535--657.

\bibitem{WangGuoAnisotropic}
Y. Wang and Z. Guo,
\emph{Global well-posedness to the three-dimensional compressible
Navier--Stokes equations with anisotropic viscous stress tensor},
Nonlinear Anal. \textbf{261} (2025), 113898,
\href{https://doi.org/10.1016/j.na.2025.113898}
{doi:10.1016/j.na.2025.113898}.

\bibitem{WuZhu}
J. Wu and Y. Zhu,
\emph{Global solutions of 3D incompressible MHD system with mixed partial
dissipation and magnetic diffusion near an equilibrium},
Adv. Math. \textbf{377} (2021), 107466,
\href{https://doi.org/10.1016/j.aim.2020.107466}
{doi:10.1016/j.aim.2020.107466}.

\bibitem{XuDu}
L. Xu and M. Du,
\emph{Global well-posedness of classical solutions to the isentropic
compressible Navier--Stokes equations with anisotropic viscous-stress tensor in
3D bounded domains},
Discrete Contin. Dyn. Syst. Ser. B \textbf{36} (2026), 463--502,
\href{https://doi.org/10.3934/dcdsb.2026041}{doi:10.3934/dcdsb.2026041}.

\bibitem{XuZhang}
L. Xu and P. Zhang,
\emph{Enhanced dissipation for the third component of 3D anisotropic
Navier--Stokes equations},
J. Differential Equations \textbf{335} (2022), 464--496,
\href{https://doi.org/10.1016/j.jde.2022.07.005}
{doi:10.1016/j.jde.2022.07.005}.

\end{thebibliography}
\end{document}